\documentclass[11pt]{amsart}
\usepackage{hyperref}
\usepackage[utf8]{inputenc}
\usepackage[nohug]{ptdiagrams}
\usepackage{xcolor}
\usepackage{amsthm}
\usepackage{booktabs}
\usepackage{geometry}
\usepackage[demo]{graphicx}
\usepackage{caption}
\usepackage{subcaption}
\usepackage{quiver}
\usepackage[mathscr]{eucal}
\usepackage{amsfonts, amsmath, amsthm, amssymb, latexsym}
\usepackage{enumitem}
\usepackage{mathtools}
\usepackage{tikz}
\usepackage{tikz-cd}
\tikzcdset{scale cd/.style={every label/.append style={scale=#1},
    cells={nodes={scale=#1}}}}
    
\usepackage{cleveref}
\usetikzlibrary{positioning,calc,decorations.markings,arrows.meta}
\usepackage{tikz-3dplot}
\usetikzlibrary{shapes, calc,shadings}

\newcommand{\edge}[2]{\{#1,\,#2\}}
\definecolor{thmblue}{RGB}{230,240,255}
\definecolor{thmframe}{RGB}{70,130,180}
\definecolor{defgreen}{RGB}{235,250,235}
\definecolor{defframe}{RGB}{60,140,60}
\definecolor{remyellow}{RGB}{255,250,230}
\definecolor{remframe}{RGB}{180,160,60}
\definecolor{excolor}{RGB}{250,235,245}
\definecolor{exframe}{RGB}{150,50,100}

\newtheorem{theorem}{Theorem}[section]
\newtheorem{lemma}[theorem]{Lemma}

\newtheorem{prop}[theorem]{Proposition}
\newtheorem{cor}[theorem]{Corollary}
\theoremstyle{definition}
\newtheorem{defi}[theorem]{Definition}
\newtheorem{example}[theorem]{Example}
\theoremstyle{remark}
\newtheorem{remark}[theorem]{Remark}

\theoremstyle{conjecture}
\newtheorem{conjecture}[theorem]{Conjecture}

\numberwithin{equation}{section}

\newcommand{\U}{\mathscr{U}}
\newcommand{\onto}{\,\,\twoheadrightarrow\,\,}
\newcommand{\la}{\label}
\newcommand{\colim}{\mathrm{colim}}

\newcommand{\Ker}{\mathrm{Ker}}

\newcommand{\Spec}{\mathrm{Spec}}

\newcommand{\Specc}{\mathrm{cSpec}}
\newcommand{\cSpec}{\mathrm{cSpec}}
\newcommand{\RSpec}{\mathrm{RSpec}}
\newcommand{\Map}{\mathrm{Map}}

\newcommand{\Ho}{\mathrm{Ho}}
\newcommand{\End}{\mathrm{End}}

\newcommand{\bcolim}{\boldsymbol{\mathrm{colim}}}
\newcommand{\blim}{\boldsymbol{\mathrm{lim}}}
\newcommand{\hocolim}{\mathrm{hocolim}}
\newcommand{\holim}{\mathrm{holim}}

\newcommand{\ux}{\underline{x}}
\newcommand{\Sym}{\mathrm{Sym}}

\def\deg{\mathrm{deg}\,}

\newcommand{\Pc}{\mathscr{P}}

\newcommand{\bQ}{\boldsymbol{Q}}
\newcommand{\id}{\mathrm{id}}
\newcommand{\into}{\hookrightarrow}

\newcommand{\Top}{\mathrm{Top}}
\newcommand{\Aff}{\mathrm{Aff}}
\newcommand{\dAff}{\mathrm{dAff}}
\newcommand{\cAff}{\mathrm{cAff}}

\newcommand{\nn}{\natural \natural}
\newcommand{\Ast}{\mathop{\scalebox{1.6}{\raisebox{-0.2ex}{$\ast$}}}}
\newcommand{\red}{\textcolor{red}}

\def\sd{\mathrm{sd}}
\def\I{\mathcal{I}}
\def \Dc {\mathscr{D}}
\def\sS{\mathscr{S}}
\def\Sp{\mathscr{S}\mathrm{p}}
\def\bpi{\boldsymbol{\pi}}
\def\bp{\boldsymbol{p}}
\def\bq{\boldsymbol{q}}
\def\bG{\mathbb G}
\def\bB{\mathbb B}
\def\bT{\mathbb T}
\def\bF{\mathbb F}
\def\bS{\mathbb S}
\def\O{\mathbb O}

\def\bbP{\boldsymbol{P}}
\def\bA{\mathbb A}
\def\bV{\boldsymbol{V}}
\def\bX{\boldsymbol{X}}
\def\bU{\boldsymbol{U}}

\def\bP{\mathbb P}
\def\bE{\mathbb E}
\def\bH{\boldsymbol{H}}
\def\NH{\mathfrak{N}}

\def\C{\mathbb C}
\def\F{\mathbb F}

\def\c{\mathbb{C}}
\def\e{\boldsymbol{e}}
\def\A{\mathcal{A}}

\def\M{\mathcal{M}}
\def\cQ{\mathscr{Q}}
\def\bcQ{\boldsymbol{\mathscr{Q}}}
\def\bcH{\boldsymbol{\mathscr{H}}}
\def\G{\mathcal{G}}

\def\bL{\boldsymbol{L}}
\def\bR{\boldsymbol{R}}

\def\sd{\mathrm{sd}}
\def\a{\mathbb{A}}
\def\k{\mathbb{K}}
\def\Z{\mathbb{Z}}
\def\Q{\mathbb{Q}}

\def\N{\mathbb{N}}
\def\P{\mathcal{P}}

\def\R{\mathcal{R}}
\def\D{\mathcal{D}}

\def\h{\mathfrak{h}}

\def\sl2{{\mathfrak{s}\mathfrak{l}}_2}
\def\vreg{V_{\rm{reg}}}

\def\k {\Bbbk}

\def\H{\mathcal{H}}

\def\Com{\mathrm{Com}}
\def\Hom{\mathrm{Hom}}
\def\mod{\mathrm{mod}}
\def\Mod{\mathrm{Mod}}

\def\im{\mathrm{Im}}
\def\Fun{\mathrm{Fun}}
\def\CAlg{\mathrm{CAlg}}
\def\Sch{\mathrm{Sch}}
\def\Comm{\mathrm{Comm}}
\def\sComm{\mathrm{sComm}}
\def\sset{\mathrm{sSet}}
\def\ccdga{\mathcal{CDGA}}
\def\cdga{\mathrm{CDGA}}

\def\Cat{\mathtt{Cat}}

\def \cF{\mathcal{F}}
\def \wcF{\widetilde{\mathcal{F}}}
\def\cO{\mathcal{O}}
\def\cA{\mathcal{A}}
\def\cR{\mathcal{R}}
\def\Cc{\mathscr{C}}
\def\M{\mathcal{M}}
\def\Mc{\mathscr{M}}
\def\Sc{\mathscr{S}}
\def\dDel{\mathbf{\Delta}}
\begin{document}

\title{
Quasi-flag manifolds and moment graphs}
%


\author{Yuri Berest}
\address{Department of Mathematics,
Cornell University, Ithaca, NY 14853-4201, USA}
\email{berest@math.cornell.edu}
\author{Yun Liu}
\address{Department of Mathematics,
University of Iowa,
Iowa city, IA 52240, USA}
\email{yun-liu@uiowa.edu}
\author{Ajay C. Ramadoss}
\address{Department of Mathematics,
Indiana University,
Bloomington, IN 47405, USA}
\email{ajcramad@iu.edu}
%

\begin{abstract}
In this paper, we introduce and study a new class of topological $G$-spaces generalizing the classical flag manifolds $  G/T $ of compact connected Lie groups. These spaces --- which we call the $m$-{\it quasi-flag manifolds} $ F_m = F_m(G,T) $ --- are topological realizations of the algebras $ Q_k(W) $ of $k$-quasi-invariant polynomials of the Weyl group $ W = W_G(T)$ in the sense that their (even-dimensional) $G$-equivariant cohomology $ H_G(F_m, \c) $ is naturally isomorphic to $ Q_k(W) $, where $ m $ is a $W$-invariant  multiplicity function on the system of roots of $W$ and $ k = \frac{m}{2}$ or $ \frac{m+1}{2}$ depending on whether $m$ is even or odd. Many topological properties and algebraic structures related to the flag manifolds can be extended to quasi-flag manifolds: in particular, we show that
\begin{itemize}
    \item the classical Borel (`double coinvariant') presentation and the so-called GKM presentation for the rational $T$-equivariant cohomology of $G/T$ generalize naturally to $ F_m(G,T) $;
    this implies that $ H_T(F_m, \Q) $ is a free module over $ H^*(BT, \Q) $, whenever $m$ is even or odd;
    \item there are natural actions of the nil-Hecke (aka Demazure) algebra and the 
    rational double-affine Hecke (aka Cherednik) algebra on equivariant cohomology of quasi-flag 
    manifolds; this gives a topological meaning to representation-theoretic results 
    of \cite{BEG03} and \cite{BC11};
    \item the quasi-flag manifolds carry a $G$-equivariant topological $W$-action 
    such that $ H^*_G(F_m, \Q)^W \cong H^*(BG, \Q)\,$  for all $ m$; this generalizes Borel's Isomorphism Theorem.   
\end{itemize}
We obtain these and other cohomological results by constructing rational algebraic models of quasi-flag manifolds in terms of coaffine stacks --- a certain kind of derived stacks introduced by B. To\"en \cite{To06} and J. Lurie  \cite{DAGVIII, DAGXIII} to provide an algebro-geometric framework for rational homotopy theory. The coaffine stacks of quasi-invariants that we introduce in this paper are natural `derived' extensions of the ordinary varieties of quasi-invariants: in fact, both are obtained from the geometric 
representation of $W$ by the same `quasi-quotient' construction  --- one in the $\infty$-category of coaffine stacks, while the other in the ordinary category of affine schemes.
Besides cohomology, we also compute the equivariant $K$-theory of quasi-flag manifolds and extend some of the above results to the multiplicative setting. Motivated by $K$-theory, we propose an algebraic definition of rings of exponential quasi-invariants and quasi-covariants of $W$.

On the topological side, our approach is motivated by the classical work on homotopy decompositions of classifying spaces of compact Lie groups (see, e.g., \cite{JM92, JMO92, JMO94}); however,
the diagrams that we use in our decompositions do not arise from collections of subgroups of $G$ but rather from {\it moment graphs} --- combinatorial objects introduced in a different area of topology called the GKM theory \cite{GKM98}. We obtain the quasi-flag manifolds from $G/T$ in two steps: first, we apply a homotopy-theoretic `$m$-thickening' construction to $G/T$, and then use a relative version of Quillen's plus-construction to make the resulting spaces simply connected.
Our approach  applies to more general spaces than $G/T$. For example, one can start with a partial flag manifold $G/P$ or an arbitrary $G$-manifold $M$ that satisfies the GKM geometric conditions for the restricted $T$-action. In this last case, the role of the Weyl group $W$ --- or rather, its root system --- is played by the moment graph $\Gamma$ of $M$, to which we can now associate a generalized ring $Q_k(\Gamma)$ of quasi-invariants. We will study these quasi-GKM spaces $M_m(G,T)$ `with multiplicities' and the associated rings $Q_k(\Gamma)$ in our subsequent paper.
\end{abstract}
\maketitle

\tableofcontents

\section{Introduction and main results}
\la{S1}

\subsection{Introduction}
Let $ G $ be a compact connected Lie group, $ T \subseteq G $ a maximal torus in $G $, and $ W = W_G(T)$ the associated Weyl group. Let $ F_0(G,T) = G/T $ denote the classical flag manifold of 
$G$.  It is well known that the $G$-equivariant, $T$-equivariant and ordinary 
cohomology rings of $G/T$ with coefficients in the rational numbers are given by
\begin{eqnarray}
H^*_G(G/T,\,\Q) & \cong & H^*(BT,\,\Q) \, \cong \, \Q[V]\,, \nonumber\la{In1}\\*[1ex]  
H^*_T(G/T,\,\Q) &\cong & \{(p_w)\in\prod_{w \in W}\!\Q[V]\,:\, p_{s_\alpha w} \equiv p_w\ \mod\, \langle \alpha \rangle \,\,,\,\forall \,\alpha \in \A\,\} \,,\la{In2}\\
H^*(G/T,\,\Q) & \cong & \Q[V]/\langle \Q[V]^W_+\rangle\,, \la{In3}\nonumber
\end{eqnarray}
where $ V $ is the $\Q$-vector space $ \pi_1(T) \otimes \Q \cong H_2(BT, \,\Q) $ carrying a reflection representation of $W$, $ \Q[V] = \Sym_{\Q}(V^*) $ is the polynomial algebra on $V$,
$ \A  $ is a set of linear forms (roots) $ \alpha \in V^* $ defining the
reflection hyperplanes $ H_{\alpha} = \Ker(\alpha)$ of $W$, and $ s_{\alpha} \in W $ are the corresponding reflections generating $W$. Formula \eqref{In2} is usually referred to as the GKM presentation of  the $T$-equivariant cohomology --- after the influential work of Goresky, Kottwitz, MacPherson \cite{GKM98} --- although, in the case of flag manifolds, it was known earlier (see, e.g., \cite{KK87}, \cite{Ara89} and historic notes in \cite{AF23}). 

The $T$-equivariant cohomology of $G/T$ has another classical description that goes back to the work of Borel \cite{Bo53}. By Borel's Theorem, the trivial map $ G/T \to {\rm pt} $ induces an injective algebra homomorphism on the $G$-equivariant cohomology:
$\,H^*(BG, \Q) \into H^*(BT,\Q)\,$, whose image is precisely the subring of $W$-invariant polynomials in $\Q[V]$:
\begin{equation}
\la{In4}
 H^*(BG, \Q)\,\cong\, H^*(BT, \Q)^W\,\cong\, \Q[V]^W.
\end{equation}
With this identification, there is a canonical isomorphism
\begin{equation}
\la{In5}
 H_T^*(G/T, \Q)\,\cong\, \Q[V] \otimes_{\Q[V]^W} \Q[V]\ ,
\end{equation}
which is called the  Borel presentation of $ H_T^*(G/T, \Q)$.

The goal of the present paper is to introduce a new class of topological $G$-spaces generalizing the classical flag manifolds of compact connected Lie groups. These spaces --- which we call the $m$-{\it quasi-flag manifolds} $ F_m = F_m(G,T) $ --- are topological realizations of algebras of quasi-invariant polynomials of $W$ in the sense that their (even-dimensional) rational cohomology rings are 
\begin{eqnarray}
H^{\rm ev}_G(F_m,\,\Q) & \cong & Q_k(W)\,, \nonumber \la{In6}\\*[1ex]  
H^{\rm ev}_T(F_m,\,\Q) &\cong & \bQ_{2k+1}(W)\,,\la{In7}\\*[1ex]
H^{\rm ev}(F_m,\,\Q) & \cong & Q_k(W)/\langle \Q[V]^W_+\rangle\,,\nonumber \la{In8}
\end{eqnarray}
where $ k = \frac{m}{2}$ or $ \frac{m+1}{2}$ depending on whether $m$ is even or odd.
The commutative algebras $ Q_k(W) $ and $ \bQ_m(W) $ that feature in the above isomorphisms are defined (over a field $ \k $ of characteristic zero for an arbitrary finite, not necessarily crystallographic, Coxeter group $W$) by
\begin{eqnarray}
Q_k(W) &:= & \{p \in \k[V]\,:\, s_{\alpha}(p) \equiv p\ \mod\, \langle \alpha\rangle^{2 k_{\alpha}}\ ,\ \forall \, \alpha \in \A\}\ , \la{In9} \\*[1ex]
\bQ_m(W) &:=& \{\,(p_w)\in\prod_{w \in W}\!\k[V]\,:\, p_{s_\alpha w} \equiv p_w\ \mod\, \langle \alpha \rangle^{m_{\alpha}} \,\,,\,\forall \,\alpha \in \A\,\} \,,\la{In10}
\end{eqnarray}
where $ k, m: \A \to \Z_+ $ are $W$-invariant functions assigning integral multiplicities ($k_{\alpha} $ and $m_{\alpha}$, respectively) to the reflection hyperplanes $ H_\alpha $ of $W$. 
We denote by $ \M(W) $ the set of all such multiplicity functions and regard it as a poset with the partial order: $\,m \le m'\,$ $\,\stackrel{}{\Leftrightarrow}\,$ $\,m_\alpha \le m_{\alpha}'\,$ for all $ \alpha\ \in \A $. The quasi-flag manifolds are then indexed by $\M(W)$,
forming a diagram of $G$-spaces 
\begin{equation}
\la{Inm4}
F_*(G,T):\ \M(W) \to \Top^G \ , \quad m \mapsto F_m(G,T)\,,
\end{equation}
all the above cohomology isomorphisms being natural in $ m $. This implies, in particular, that for any pair of multiplicities $\ m \le m' \,$ in $ \M(W)$, the natural maps of $G$-spaces 
$$ 
G/T \to F_m(G,T) \to F_{m'}(G,T) \to {\rm pt} 
$$ 
induces {\it injective} algebra homomorphisms:
$$
H^*(BG, \Q) \,\into\, H^{\rm ev}_G(F_{m'},\,\Q) \,\into\, H^{\rm ev}_G(F_m,\,\Q) \,\into\, H^*(BT,\,\Q)\,.
$$
In addition, each space $ F_m(G,T) $ carries a natural $W$-action, so that \eqref{Inm4} is actually
a diagram of $(G \times W)$-spaces and the above cohomology homomorphisms are $W$-equivariant. Restricting to $W$-invariants yields
$$
H^*(BG, \Q)\,\cong\, H^\ast_G(F_m,\,\Q)^W\,,\ \forall\, m \in \M(W)\,,
$$
which can be viewed as a generalization of Borel's classical isomorphism \eqref{In4}.

When $m=0$, the isomorphisms \eqref{In7} formally reduce to the classical ones \eqref{In2}. Formula \eqref{In7} for the $T$-equivariant cohomology is thus a generalization of the classical GKM presentation to arbitrary multiplicities. There is also a natural generalization of the Borel presentation \eqref{In5}, which takes the form
\begin{equation}
\la{Inm5}
 H_T^{\rm ev}(F_m, \Q)\,\cong\, \Q[V] \otimes_{\Q[V]^W} Q_k(W)\ ,
\end{equation}
where the relation between $m$ and $k$ is the same as in \eqref{In7}.

The algebras $Q_k(W) $ defined by \eqref{In9} are called the {\it $k$-quasi-invariants} of $W$. They first appeared in mathematical physics in the early 1990s  (see \cite{CV90, CV93}), and since then have received a good deal of attention and found applications in other areas:  most notably,  representation theory, algebraic geometry and combinatorics (see \cite{FV02}, \cite{EG02b}, \cite{Ch02}, \cite{BEG03}, \cite{FeV03}, \cite{GW06}, \cite{BM08}, \cite{T10}, \cite{BC11}, \cite{BEF20}, \cite{Gri21}).  

In contrast, the algebras $ \bQ_m(W) $ defined by \eqref{In10} seem to be much less studied. We call these algebras the {\it $m$-quasi-covariants} of $W$. To the best of our knowledge, in  representation theory, they were originally introduced in the paper \cite{BC11}, where they were defined for an arbitrary complex reflection group $W$, but in a different (`non-GKM') form and only as modules, not as commutative algebras (see Section~\ref{S2.2} for the comparison of our definition of $ \bQ_m(W) $ with that of \cite{BC11}).
The key observation of \cite{BC11} is that
$ \bQ_{m}(W) $ for $\,m = 2k\,$ is naturally a module 
over the
rational double-affine Hecke algebra $ \bH_{k}(W) $ associated to $W$. This result does not extend to
odd $m$: instead, we prove that $\bQ_{2k+1}(W) $ is a module over the classical 
nil-Hecke algebra $ \NH(W) $ of $W$ (see Theorem~\ref{Qfatodd}). In view of \eqref{In7}, this means that the $T$-equivariant cohomology of $ F_m(G,T)$ carries a natural (left) action of $ \NH(W)$ for all $ m \in \M(W) $. Finally, the connection to representation theory allows us to prove that 
$ \bQ_m(W) $ is a {\it free} module over $ \k[V] $ whenever $ m = 2k $ or $ m = 2k+1 $ for all $ k \in \M(W)$ (see Theorem~\ref{ThFree}). 


The problem of topological realization of algebras of quasi-invariants was posed in  \cite{BR1}.
%
%
%
We refer the reader to 
that paper 
for motivation and 
a general axiomatic formulation of the problem.
In terms of results,  \cite{BR1} focused entirely on the rank one case, solving the realization problem for $ G = SU(2) $ and its homotopy-theoretic generalizations --- the Rector spaces (aka the `fake' Lie groups of type $SU(2)$). The present paper can thus be viewed as a sequel to that work; however, it can be read independently as the methods we use are different
from those of \cite{BR1}. We should mention that, in the case of $SU(2)$,
an explicit geometric construction of spaces of quasi-invariants was also given earlier in the unpublished preprint \cite{FF}, which will appear as an Appendix in the updated version of \cite{BR1}.

The spaces $ F_m(G,T) $ are obtained from $G/T$ in a natural way --- by a homotopy-theoretic `gluing' construction that we call the {\it $m$-simplicial thickening} (see Section~\ref{S5}). 
This construction applies naturally to other spaces than $G/T$: for example, one can start with a partial flag manifold $G/P$, or more generally, with an arbitrary $G$-manifold $M$ that satisfies the so-called GKM conditions as a $T$-manifold (see \cite{GKM98}). In this last case, the role of the root system of $W$ is played by the {\it GKM graph} $\Gamma $ of $M$, to which we can now associate the rings of generalized quasi-invariants, $Q_k(\Gamma)$, and quasi-covariants, $\bQ_m(\Gamma)$ motivated by formulas \eqref{In7}. We briefly sketch this generalization in Section~\ref{S1.8}, deferring details to our forthcoming paper. 

Finally, we would like to mention another  generalization that deserves a further investigation. Our algebras $ \bQ_m(W) $ and $ \bQ_m(\Gamma) $ can be viewed as a special case of the so-called {\it generalized spline algebras} introduced in \cite{GTV16}  (see Remark~\ref{Rspline}). Motivated by various problems in applied mathematics, discrete geometry and combinatorics, the classical spline algebras and their generalizations have been extensively studied in recent years, mostly by methods of commutative algebra (see \cite{LST24} for a nice recent survey). The algebraic results about $\bQ_m(W)$ proved in this paper rely crucially on the connection to representation theory of 
rational Cherednik (aka double affine Hecke algebras) developed in \cite{BEG03} and \cite{BC11}. It would be interesting to see if this connection extends to other, more general classes of spline algebras.

\vspace*{0.8ex}

We now proceed with a detailed summary of the contents of the paper,
focusing on main ideas and results. We begin with the first two sections that provide algebraic prerequisites for our topological constructions.
\vspace*{-1ex}
\subsection{Schemes of quasi-invariants}
In Section~\ref{S2}, we re-examine the  algebraic definitions of quasi-invariants
and quasi-covariants from a geometric point of view. We describe a universal  (`quasi-quotient') construction in the category of $\k$-schemes that takes as an input a $W$-scheme $X$ --- or rather its categorical quotient 
$\, X \to X/\!/W\,$ --- and yields as an output a diagram (tower) of schemes
$\,X_*(W): \M(W) \to \Aff_{\k}\,$ defined over $ X/\!/W $.
The construction proceeds in two steps: first, we use a relative (algebro-geometric) 
version of the classical Ganea construction \cite{Ga65} (see also \cite{BR1}) to produce schemes with `cuspidal' singularities along the 
subschemes $ X^{W_\alpha} \subset X $ of multiplicity $k_{\alpha}$, and then, we `glue' these singular schemes by taking a colimit over the family of all elementary reflection subgroups $ \{W_\alpha\} $ of $W$.
We apply this `quasi-quotient' construction in two cases: to the geometric representation $V$ itself, viewed as a $W$-scheme, and its canonical $W$-free resolution $ \bV = V \times  W $. In the first case, we obtain the 
tower of classical varieties of quasi-invariants $\, V_k(W) = \Spec\,Q_k(W) $ 
(see Theorem~\ref{ThVQI}), while in the second case --- 
the varieties of quasi-covariants $ \bV_{\! m}(W) = \Spec\,\bQ_m(W) $ (see Theorem~\ref{ThVQ2}).
The relation between these varieties, which is a simple consequence of our formalism, is given by the canonical isomorphism (Lemma~\ref{bVV}):
\begin{equation}
\la{In11}
\bV_{\! m}(W)/\!/W \cong V_{[\frac{m}{2}]}(W)
\end{equation}
that holds for all multiplicities $ m \in \M(W) $. This is equivalent to an isomorphism of algebras %
\begin{equation}
\la{In111}
\bQ_{m}(W)^W \cong Q_{[\frac{m}{2}]}(W)\,,
\end{equation}
which was originally established for $ m = 2k $ (as an isomorphism of $ \k[V]^W$-modules)  in \cite{BC11}.

More subtle --- and indeed more surprising --- is the existence of another natural isomorphism:
\begin{equation}
\la{In12}
\bV_{2k+1}(W) \,\cong\,  V \times_{V/\!/W} V_k(W)
\end{equation}
that holds only for {\it odd} multiplicities $m =2k+1$ (see Proposition~\ref{bQodd}).
Algebraically, \eqref{In12} can be restated in the form
\begin{equation}
\la{In13}
\bQ_{2k+1}(W)\, \cong\,  \k[V] \otimes_{\k[V]^W} Q_k(W)
\end{equation}
underlying the Borel presentation \eqref{Inm5} of the $T$-equivariant cohomology of quasi-flag manifolds.

Now, the main result of Section~\ref{S2} reads
\newtheorem*{ThFree}{Theorem~\ref{ThFree}}
\begin{ThFree}
Assume that $ m \in \M(W) $ is either even or odd $($i.e. $m_{\alpha} = 2 k_{\alpha} $ or $ m_{\alpha} = 2 k_{\alpha} + 1 $ for all $ \alpha \in \A$$)$. Then $ \bQ_m(W) $ is a free module over $ \k[V] $ of rank $|W|$.  In particular, $\bV_m $ are Cohen-Macaulay varieties.
\end{ThFree}
For $\, m=2k$, Theorem~\ref{ThFree} follows directly (modulo the simple Lemma~\ref{fQGKM}) from the main results of \cite{BC11}, while, for $ m =2k+1$, it is a consequence of the isomorphism \eqref{In13} and the fundamental fact that the algebras $Q_k(W) $ are free over the invariants $ \k[V]^W $ for all $k \in \M(W)$ (see \cite{EG02b, BEG03}). It is natural to expect that Theorem~\ref{ThFree} holds true for {\it all} $W$-invariant multiplicities (the remaining case to prove is the one of `mixed' multiplicities: when some of the $m_{\alpha}$'s are even and some are odd). At the moment, we do not know whether the result of Theorem~\ref{ThFree}
extends to mixed multiplicities for an arbitrary group $W$, but we have verified this fact for all $W$ of rank 2 (see Proposition~\ref{prop:QC-rank2-free}).

In Section~\ref{S2.4}, we refine another important observation of \cite{BC11}: namely, the fact that the algebras $ \bQ_{m}(W) $ for $ m = 2k$ carry a natural module structure over the rational Cherednik (aka double affine Hecke) algebra $ \bH_k(W) $ (see Theorem~\ref{Qfat}). Proving this fact amounts to showing that  $ \bQ_{2k}(W) $ is stable under the action of the Dunkl differential-difference operators $T_{\xi,k}$. Unfortunately, when $ m \not= 2k $, this is no longer true: the 
action of Cherednik algebras on $Q_{2k}(W) $ constructed in \cite{BC11} does not extend to  odd multiplicities. Instead, we find that, for $m=2k+1$, the algebras $\bQ_{m}(W)$ are stable under the action of the classical divided difference (aka Demazure) operators $\nabla_{\alpha} $, which 
implies
\newtheorem*{Qfatodd}{Theorem~\ref{Qfatodd}}
\begin{Qfatodd} 
For all $ k \in \M(W) $, the algebras $\bQ_{2k+1}(W) $ carry a  left module structure over the nil-Hecke algebra $\NH(W)$ of $W$, extending the natural module structure over $\k[V]$.
\end{Qfatodd}
Just as Theorem~\ref{ThFree}, this result is motivated by topology: it shows that the even-degree $T$-equivariant cohomology of quasi-flag manifolds $F_m(G,T) $ is a module over the nil-Hecke algebra $ \NH(W)$ of the Weyl group $W$. We note that Theorem~\ref{Qfatodd} does {\it not} hold for $ m = 2k \,$: the algebras $ Q_{2k}(W)$
are not stable under the action of Demazure operators. We address this apparent dichotomy between
even and odd multiplicities in Remark~\ref{HvsNH}, where we show that the two cases can be actually unified by realizing the nil-Hecke algebra $\NH(W)$ as a `quasi-classical limit' of $ \bH_k(W)$.

We close this section with two explicit examples, illustrating 
Theorem~\ref{ThFree} and its extension to mixed multiplicities
(Proposition~~\ref{prop:QC-rank2-free}) in the case of dihedral groups.

\begin{example}
Let $W = S_3$ act on $V = \{(x_1, x_2, x_3) : x_1 + x_2 + x_3 = 0\} \cong \k^2$ by permuting coordinates. This is the standard 2-dimensional representation of $S_3$. 
There are 3 positive roots $\alpha_{ij} = x_i - x_j$ and 3 reflections $s_{ij}=(ij)$ transposes $x_i$ and $x_j$ for $1 \le i < j \le 3$.
We abbreviate
\[
s_1=(12),\; s_2 = (23),\; s_3 =s_1s_2s_1=s_2s_1s_2 = (13),
\qquad
a = \alpha_{12},\; b = \alpha_{23},\; c  = a + b = \alpha_{13}\,.
\]
Then $W = \{e,\, s_1,\, s_2,\, s_1s_2,\, s_2s_1,\, s_3\}$. 
The Bruhat moment graph of $S_3$ is shown on Figure~\ref{GammaS3}.
The $W$-action on the reflection hyperplanes is transitive, so multiplicity $m$ (resp. $k$) is give by a single number. 
The algebra of quasi-invariants $Q_k(S_3)$ is given by 
\[
Q_k(S_3) \,=\, \bigl\{p \in \k[V] \,:\, s_{ij}(p) \equiv p \pmod{\langle x_i-x_j \rangle^{2k}},\;\forall\, 1\leq i<j\leq 3 \;\bigr\}.
\]
and the algebra of quasi-covariants $\bQ_m(S_3)$ is given by 
\begin{equation}
\la{eq:QmS3}
\bQ_m(S_3) \,=\, \Bigl\{(p_w)_{w \in S_3} \in \prod_{w \in S_3} \k[V] \,:\,
p_w \equiv p_{s_{ij} w} \pmod{\langle x_i - x_j \rangle^{m}},\;\forall\, 1\leq i<j\leq 3 \;\Bigr\}.
\end{equation}
If we set $\deg(V^\ast)=1$, the Hibert series of $\bQ_m(W)$ for $W=S_3$ is given by 
\[
\text{Hilb}\bigl(\bQ_m(S_3)\bigr) \,=\, \frac{q_m(t)}{(1-t)^2}\ ,
\]
where the numerator is 
\begin{equation}
\label{eq:S3numerators}
p_{2k}(t) \,=\, 1 + 4\,t^{3k} + t^{6k},
\qquad
p_{2k+1}(t) \,=\, 1 + 2\,t^{3k+1} + 2\,t^{3k+2} + t^{6k+3}.
\end{equation}
Solving the divisibility conditions \eqref{eq:QmS3} degree by degree, 
we can compute a free basis $\{\sigma_w, w\in S_3\}$ for $\bQ_m(S_3)$ explicitly. 
For $m=1,2,3,4$, these bases are presented in the tables below. Here, the components
are ordered as $\,(e,\, s_1,\, s_2,\, s_1s_2,\, s_2s_1,\, s_3)\,$, and the 
$a,\,b,\,c$ denote the polynomials $x_1 - x_2 $, $\,x_2 -x_3 $ and $\,x_1 - x_3$, respectively. Note it is always true that $\sigma_1 = (1,1,1,1,1,1)$ and $\sigma_{s_3} = (0,0,0,0,0,a^mb^mc^m)$, so we skip them.
In particular, we recover the basis of equivariant Schubert classes for $m=1$. 
\begin{table}[ht]
\centering
\footnotesize
\begin{tabular}{cc|cccccc}
\toprule
Basis & deg & $1$ & $s_1$ & $s_2$ & $s_1s_2$ & $s_2s_1$ & $s_3$ \\
\midrule
$\sigma_{s_1}$ & 1 & $0$ & $a$ & $0$ & $a$ & $c$ & $c$ \\
$\sigma_{s_2}$ & 1 & $0$ & $0$ & $b$ & $c$ & $b$ & $c$ \\
$\sigma_{s_1s_2}$ & 2 & $0$ & $0$ & $0$ & $ac$ & $0$ & $ac$ \\
$\sigma_{s_2s_1}$ & 2 & $0$ & $0$ & $0$ & $0$ & $bc$ & $bc$ \\
\bottomrule
\end{tabular}
\caption{Basis for $\bQ_1(S_3)$: equivariant Schubert classes. }
\label{tab:basis-m1}
\end{table}
%

%
\begin{table}[ht]
\centering
\footnotesize
\begin{tabular}{c|cccccc}
\toprule
Degree & $1$ & $s_1$ & $s_2$ & $s_1 s_2$ & $s_2 s_1$ & $s_3$ \\
\midrule
3 & $0$ & $a^2(a{+}2b)$ & $0$ & $a^2(a{+}2b)$ & $ac^2$ & $ac^2$ \\
3 & $0$ & $0$ & $b^2(2a{+}b)$ & $bc^2$ & $b^2(2a{+}b)$ & $bc^2$ \\
3 & $0$ & $a^2 b$ & $0$ & $a^2 b$ & $bc^2$ & $bc^2$ \\
3 & $0$ & $0$ & $ab^2$ & $ac^2$ & $ab^2$ & $ac^2$ \\
\bottomrule
\end{tabular}
\caption{Basis for $\bQ_2(S_3)$.}
\label{tab:basis-m2}
\end{table}
%

%
\begin{table}[ht]
\centering
\footnotesize
\begin{tabular}{c|cccccc}
\toprule
Degree & $1$ & $s_1$ & $s_2$ & $s_1 s_2$ & $s_2 s_1$ & $s_3$ \\
\midrule
4 & $0$ & $a^3(a{+}2b)$ & $0$ & $a^3(a{+}2b)$ & $c^3(a{-}b)$ & $c^3(a{-}b)$ \\
4 & $0$ & $0$ & $b^3(2a{+}b)$ & $c^3(b{-}a)$ & $b^3(2a{+}b)$ & $c^3(b{-}a)$ \\
5 & $0$ & $a^2(c^3{-}b^3)$ & $0$ & $a^2(c^3{-}b^3)$ & $a^2 c^3$ & $a^2 c^3$ \\
5 & $0$ & $0$ & $b^2(c^3{-}a^3)$ & $b^2 c^3$ & $b^2(c^3{-}a^3)$ & $b^2 c^3$ \\
\bottomrule
\end{tabular}
\caption{Basis for $\bQ_3(S_3)$.}
\label{tab:basis-m3}
\end{table}
%
%
\begin{table}[ht]
\centering
\footnotesize
\begin{tabular}{c|cccccc}
\toprule
Degree & $e$ & $s_1$ & $s_2$ & $s_1 s_2$ & $s_2 s_1$ & $s_3$ \\
\midrule
6 & $0$ & $a^5(a{+}2b)$ & $0$ & $a^5(a{+}2b)$
      & $c^4(a^2{-}2ab{+}2b^2)$ & $c^4(a^2{-}2ab{+}2b^2)$ \\
6 & $0$ & $0$ & $b^5(2a{+}b)$ & $c^4(2a^2{-}2ab{+}b^2)$
      & $b^5(2a{+}b)$ & $c^4(2a^2{-}2ab{+}b^2)$ \\
6 & $0$ & $a^4 b(2a{+}5b)$ & $0$ & $a^4 b(2a{+}5b)$
      & $b(2a{-}3b)\,c^4$ & $b(2a{-}3b)\,c^4$ \\
6 & $0$ & $0$ & $ab^4(5a{+}2b)$ & $-a(3a{-}2b)\,c^4$
      & $ab^4(5a{+}2b)$ & $-a(3a{-}2b)\,c^4$ \\
\bottomrule
\end{tabular}
\caption{Basis for $\bQ_4(S_3)$.}
\label{tab:basis-m4}
\end{table}
\end{example}

\newpage
\begin{example}
Consider the Weyl group of type $B_2$ acting on $V = \k^2$ in the natural way (as a dihedral group of order~8). It is generated by two simple reflections $t$ and $s_y$, which are given (in standard coordinates) by $t\colon (x,y) \mapsto (y,x)$ and $s_y\colon (x,y) \mapsto (x, -y)$. There are 4 positive roots: two short roots $x, y$ and two long roots $x \pm y$. The other two reflections are $s_x=ts_yt$ and $s_{xy}=s_x s_y$.
The $W$-action on the reflection hyperplanes has two orbits so multiplicities are given by two numbers $m=(m_s,m_l), k=(k_s,k_l)$ associated to short and long roots. We write the elements of $W$ as $\{1,\, \,t, \,s_y, \, r=s_x t, \, r^3=s_y t, s_x, \, s_{xy}, \,r^2\}$

For $k=(k_s, k_l) \in \Z_{\ge 0}^2$, the algebra of quasi-invariants of multiplicity $k$ is 
\[
Q_{(k_s, k_l)}(B_2) = \bigl\{f \in \k[x,y] :
  \alpha^{2k_\alpha} \mid \bigl(f - s_\alpha(f)\bigr)
  \;\text{for every reflection $s_\alpha$}\bigr\},
\]
where $k_\alpha = k_s$ for short roots 
and $k_\alpha = k_l$ for long roots.
%
Explicitly, the conditions are:
\begin{alignat*}{2}
x^{2k_s} &\mid \bigl(f(x,y) - f(-x,y)\bigr), &\qquad
y^{2k_s} &\mid \bigl(f(x,y) - f(x,-y)\bigr), \\
(x{-}y)^{2k_l} &\mid \bigl(f(x,y) - f(y,x)\bigr), &\qquad
(x{+}y)^{2k_l} &\mid \bigl(f(x,y) - f(-y,-x)\bigr).
\end{alignat*}

For $m=(m_s, m_l) \in \Z_{\ge 0}^2$, the algebra of quasi-covariants of multiplicity $m$ is 
\begin{equation}
\la{QmB2}
  \bQ_{(m_s,m_l)}(B_2) = \bigl\{(p_w)_{w \in B_2} \in \k[x,y]^8 :
  \alpha^{m_c} \mid (p_w - p_{s_\alpha w})
  \;\forall \alpha \in \cA, \forall w\in B_2\bigr\},
\end{equation}
where $m_c = m_s$ for short roots and $m_c = m_l$ for long roots.

The Hilbert series numerator of $\bQ_{(m_s,m_l)}(B_2)$ for even or odd multiplicities are 
\[
p_{2k}(t) = 1 + 6t^{4k} + t^{8k}, 
\qquad
p_{2k+1}(t) = 1 + 2t^{4k+1} + 2t^{4k+2} + 2t^{4k+3} + t^{8k+4}.
\]

Solving the divisibility condition in \eqref{QmB2}, we can obtain a basis for $\bQ_{(m_s,m_l)}(B_2)$. 
%
%
\begin{table}[ht]
\centering
\footnotesize
\begin{tabular}{c|cccccccc}
\toprule
Degree & $1$ & $t$ & $s_y$ & $s_xt$ & $s_yt$ & $s_x$ & $s_{xy}$ & $r^2$\\
\midrule
1 & $x{+}y$ & $x{+}y$ & $x$ & $y$ & $x$ & $y$ & $0$ & $0$\\
1 & $x{-}y$ & $0$ & $x$ & $x$ & $-y$ & $-y$ & $x{-}y$ & $0$\\
2 & $xy$ & $xy$ & $0$ & $0$ & $0$ & $0$ & $xy$ & $xy$\\
2 & $\lambda$ & $0$ & $\lambda$ & $0$ & $0$ & $\lambda$ & $0$ & $\lambda$\\
3 & $y\lambda$ & $0$ & $0$ & $0$ & $0$ & $y\lambda$ & $0$ & $0$\\
3 & $x\lambda$ & $0$ & $x\lambda$ & $0$ & $0$ & $0$ & $0$ & $0$\\
\bottomrule
\end{tabular}
\caption{Basis for $Q_{(1,1)}(B_2)$, the equivariant Schubert basis; $\lambda = (x-y)(x+y)$. 
}
\label{tab:B211}
\end{table}
\begin{table}[ht]
\centering
\footnotesize
\begin{tabular}{c|cccccccc}
\toprule
Degree & $1$ & $t$ & $s_y$ & $s_xt$ & $s_yt$ & $s_x$ & $s_{xy}$ & $r^2$\\
\midrule
2 & $xy$ & $xy$ & $0$ & $0$ & $0$ & $0$ & $xy$ & $xy$\\
3 & $yq_+$ & $yq_+$ & $0$ & $yq_-$ & $0$ & $yq_-$ & $0$ & $0$\\
3 & $yq_-$ & $0$ & $0$ & $0$ & $yq_+$ & $yq_+$ & $yq_-$ & $0$\\
3 & $xq_+$ & $xq_+$ & $xq_-$ & $0$ & $xq_-$ & $0$ & $0$ & $0$\\
3 & $xq_-$ & $0$ & $xq_+$ & $xq_+$ & $0$ & $0$ & $xq_-$ & $0$\\
4 & $q_-q_+$ & $0$ & $q_-q_+$ & $0$ & $0$ & $q_-q_+$ & $0$ & $q_-q_+$\\
\bottomrule
\end{tabular}
\caption{Basis for $Q_{(1,2)}(B_2)$; $q_\pm = (x\pm y)^2$. 
}
\label{tab:B212}
\end{table}
\begin{table}[ht]
\centering
\footnotesize
\begin{tabular}{c|cccccccc}
\toprule
Degree & $1$ & $t$ & $s_y$ & $s_xt$ & $s_yt$ & $s_x$ & $s_{xy}$ & $r^2$\\
\midrule
2 & $\lambda$ & $0$ & $\lambda$ & $0$ & $0$ & $\lambda$ & $0$ & $\lambda$\\
3 & $y^2\mu_+$ & $x^2\mu_+$ & $0$ & $0$ & $x^2\mu_+$ & $y^2\mu_+$ & $0$ & $0$\\
3 & $y^2\mu_-$ & $0$ & $0$ & $x^2\mu_-$ & $0$ & $y^2\mu_-$ & $x^2\mu_-$ & $0$\\
3 & $x^2\mu_+$ & $y^2\mu_+$ & $x^2\mu_+$ & $y^2\mu_+$ & $0$ & $0$ & $0$ & $0$\\
3 & $x^2\mu_-$ & $0$ & $x^2\mu_-$ & $0$ & $y^2\mu_-$ & $0$ & $y^2\mu_-$ & $0$\\
4 & $x^2y^2$ & $x^2y^2$ & $0$ & $0$ & $0$ & $0$ & $x^2y^2$ & $x^2y^2$\\
\bottomrule
\end{tabular}
\caption{Basis for $Q_{(2,1)}(B_2)$; $\mu_\pm = x \pm y$, $\lambda = \mu_-\mu_+$. 
}
\label{tab:B221}
\end{table}
\begin{table}[ht]
\centering
\footnotesize
\begin{tabular}{c|cccccccc}
\toprule
Degree & $1$ & $t$ & $s_y$ & $s_xt$ & $s_yt$ & $s_x$ & $s_{xy}$ & $r^2$\\
\midrule
4 & $x^2y^2$ & $x^2y^2$ & $0$ & $0$ & $0$ & $0$ & $x^2y^2$ & $x^2y^2$\\
4 & $q_-q_+$ & $0$ & $q_-q_+$ & $0$ & $0$ & $q_-q_+$ & $0$ & $q_-q_+$\\
4 & $yq_+(2x{-}y)$ & $yq_+(2x{-}y)$ & $x^2y(2x{-}3y)$ & $-y^4$ & $x^2y(2x{-}3y)$ & $-y^4$ & $0$ & $0$\\
4 & $xq_+(x{-}2y)$ & $xq_+(x{-}2y)$ & $x^2(x^2{-}3y^2)$ & $-2xy^3$ & $x^2(x^2{-}3y^2)$ & $-2xy^3$ & $0$ & $0$\\
4 & $yq_-(2x{+}y)$ & $0$ & $2x^3y$ & $2x^3y$ & $y^2(y^2{-}3x^2)$ & $y^2(y^2{-}3x^2)$ & $yq_-(2x{+}y)$ & $0$\\
4 & $xq_-(x{+}2y)$ & $0$ & $x^4$ & $x^4$ & $xy^2(3x{+}2y)$ & $xy^2(3x{+}2y)$ & $xq_-(x{+}2y)$ & $0$\\
\bottomrule
\end{tabular}
\caption{Basis for $Q_{(2,2)}(B_2)$; $q_\pm = (x\pm y)^2$. 
}
\label{tab:B222}
\end{table}
\begin{table}[ht]
\centering
\footnotesize
\begin{tabular}{c|cccccccc}
\toprule
Degree & $1$ & $t$ & $s_y$ & $s_xt$ & $s_yt$ & $s_x$ & $s_{xy}$ & $r^2$\\
\midrule
5 & $c_+^3h_-$ & $c_+^3h_-$ & $x^3g_1$ & $-y^3g_2$ & $x^3g_1$ & $-y^3g_2$ & $0$ & $0$\\
5 & $c_-^3h_+$ & $0$ & $x^3g_1$ & $x^3g_1$ & $y^3g_2$ & $y^3g_2$ & $c_-^3h_+$ & $0$\\
6 & $x^3y^3$ & $x^3y^3$ & $0$ & $0$ & $0$ & $0$ & $x^3y^3$ & $x^3y^3$\\
6 & $c_-^3c_+^3$ & $0$ & $c_-^3c_+^3$ & $0$ & $0$ & $c_-^3c_+^3$ & $0$ & $c_-^3c_+^3$\\
7 & $y^3c_+^3n_-$ & $y^3c_+^3n_-$ & $0$ & $y^3c_-^3n_+$ & $0$ & $y^3c_-^3n_+$ & $0$ & $0$\\
7 & $y^3c_-^3n_+$ & $0$ & $0$ & $0$ & $y^3c_+^3n_-$ & $y^3c_+^3n_-$ & $y^3c_-^3n_+$ & $0$\\
\bottomrule
\end{tabular}
\caption{Basis for $Q_{(3,3)}(B_2)$, with $c_\pm = x \pm y$, $h_\pm = x^2 \pm 3xy + y^2$, $g_1 = x^2 - 5y^2$, $g_2 = 5x^2 - y^2$ and $n_\pm = 3x \pm y$. 
}
\label{tab:B233}
\end{table}
\end{example}

\newpage
\subsection{Stacks of quasi-invariants}
The geometric construction of quasi-invariants developed in Section~\ref{S2} is universal in the sense that it involves only natural operations in the category of schemes expressed in terms of limits and colimits. It thus generalizes naturally to other categories of geometric nature, allowing one to define `varieties of quasi-invariants' in more sophisticated settings. 
The most relevant example for us (see Sections~\ref{S3}) is the $\infty$-category $ \cAff_{\k} $ of {\it coaffine stacks} introduced by B.To\"en \cite{To06} and studied further (in the context of spectral algebraic geometry) by J. Lurie \cite{DAGVIII, DAGXIII}. This $\infty$-category is a natural derived extension of the category of affine schemes $\Aff_{\k}$ that is, in a sense, dual to the more familiar $\infty$-category $\dAff_{\k}$ of derived affine schemes. Both $ \dAff_{\k}$ and $ \cAff_{\k} $ 
are subcategories of the $\infty$-category of derived stacks $ {\rm dSt}_{\k}$ containing $\Aff_{\k}$, but modelled on different (`opposite') kinds of commutative algebras: when $\k$ is a field of characteristic zero\footnote{which we tacitly assume in this paper unless stated otherwise.},  $\dAff_{\k}$ is determined by the {\it negatively} graded cochain dg algebras, $ \cdga_{\k}^{\le 0}$, while $\cAff_{\k}$ by the {\it positively} graded coconnective dg algebras,  $ \cdga_{\k}^{\ge 0}$. Naturally, we say that a coaffine stack is simply connected if the corresponding dg algebra $A$ satisfies $ H^{1}(A)= 0 $.
Now, the key observation of To\"en and Lurie is that, under natural finiteness assumptions, the $\infty$-category $ \cAff^1_{\Q} $ of simply connected coaffine stacks is equivalent to the classical {\it rational} homotopy category of simply connected spaces (see Appendix~\ref{AB}, Theorem~\ref{SuTh}). Intuitively speaking, the coaffine stacks are derived stacks whose $\Q$-points represent the simply connected topological spaces of finite rational homotopy type in much the same way as the $\C$-points of usual schemes of finite type over $\C$ represent the classical complex algebraic varieties. Thus, to every simply connected topological space $X$ of finite rational type we can assign a coaffine stack --- called the {\it rational model} of $X$ --- that determines
(and is determined by) $X$ uniquely, up to rational equivalence. 

In Sections~\ref{S3.2} and \ref{S3.3}, we apply a derived (i.e., homotopy-coherent) version of 
the `quasi-quotient' construction of Section~\ref{S2} to define coaffine stacks that will serve as {\it rational} models for our quasi-flag manifolds and their homotopy quotients (see Theorem~\ref{MainT}). We construct the coaffine stacks $\, V_k^c(W)  $ and  $\, \bV_{\! m}^c(W) $ starting with the {\it same} $W$-schemes $V$ and $ \bV $ as in Section~\ref{S2} but viewing them as objects in $ \cAff_{\k}$. This amounts simply to putting the standard cohomological grading on the polynomial algebras $\k[V]$ and $ \k[\bV] $ making them 
positively graded cdga's with trivial differentials (i.e. objects in $ \cdga_{\k}^{\ge 0}$). 
Working in $\cAff_\k$ then amounts to replacing basic scheme-theoretic
operations (functors) with their derived analogues: i.e., using homotopy limits and homotopy colimits in $\cAff_{\k}$ in place of ordinary limits and colimits.
We describe the commutative dg algebras $Q_k^*(W)$ and $\bQ_m^*(W)$ representing the coaffine stacks  $\, V_k^c(W)  $ and  $\, \bV_m^c(W) $ in explicit terms, using polynomial differential forms on algebraic simplices (see \eqref{QcW} and \eqref{tbWQ}, respectively). We also compute the cohomology of these dg algebras and show that their even-dimensional subalgebras coincide (up to doubling the grading) with the ordinary algebras of quasi-invariants and quasi-covariants defined by \eqref{In9} and \eqref{In10} (see Theorem~\ref{Hcomp}): 
$$
H_*^{\rm ev}[Q_k^*(W)]\,\cong\, Q_k(W)\ ,\quad
H_*^{\rm ev}[\bQ_m^*(W)]\,\cong\, \bQ_m(W)\,,\quad \forall\,k,m\in\M(W)\,.
$$
We emphasize that, in general (when $k,m\not=0$), the dg algebras
$Q_k^*(W)$ and $\bQ_m^*(W)$ also have non-vanishing {\it odd} cohomology, which we  
compute in an explicit form as well (see \eqref{HQodd}). Thus, the dg algebras $Q_k^*(W)$ and $\bQ^*_m(W)$ --- or rather their cohomology algebras --- provide a nontrivial generalization 
of the classical rings of quasi-invariants defined by \eqref{In9} and \eqref{In10}. 

For the purpose of comparison, in Section~\ref{S3.1}, we perform our construction
in the $\infty$-category of derived affine schemes, $\dAff_k$. Again, we apply it to the same $W$-schemes $V$ and $\bV $ as in Section~\ref{S2}  --- but now regarded as objects in $ \dAff_{\k}$. In contrast to $\cAff_\k$, we find that, in $\dAff_{\k}$, our construction does not yield new geometric objects: we simply obtain derived affine schemes that are equivalent to the ordinary schemes of quasi-invariants (see Theorem~\ref{Tdc}). This shows that, in the setting of derived algebraic geometry, the definition of `varieties of quasi-invariants' depends on the category of derived rings that we use to model our spaces.

\subsection{Homotopy decomposition of a flag manifold over the moment graph} 
We define quasi-flag manifolds in terms of canonical homotopy decompositions over certain index categories $ \Cc^{(m)}(\Gamma)_{hW}$ attached to the Weyl group $ W $ and a $W$-invariant multiplicity function $ m \in \M(W)$. Although our definition is homotopy theoretic in nature, the diagrams of spaces that we use are geometric: they are built from finite products of orbits of complex algebraic tori acting on $ G/T$. 

To set the stage for our general definition we begin with the case of flag manifolds ($ m=0 $). In this case, the category $ \Cc^{(0)}(\Gamma)_{hW} = \Cc(\Gamma)_{hW} $ is determined by the classical {\it Bruhat moment graph} $ \Gamma = \Gamma(\cR_W) $ of the root system $\cR_W $ of $W$. Recall (see, e.g., \cite{Car94, Fie12}) that $ \Gamma $
is a simple unoriented graph with vertex set $ V_\Gamma = W $ and the edge 
set $ E_\Gamma = \{e(s_\alpha, w)\ :\ \alpha \in \cR_+ \,,\, w \in W\} $, where the edge $\, e(s_{\alpha}, w) \,$ connects the vertices $w$ and $ s_\alpha w $ for each reflection $ s_\alpha \in W $ (see Figure~\ref{GammaS3}). 
Combinatorially, we can view 
$\Gamma$ as a one-dimensional simplicial complex on the vertex set $ V_{\Gamma} $
and the face set $ \Sigma_{\Gamma} = V_{\Gamma} \cup E_{\Gamma}\,$: this, in turn, defines a poset and hence a small category that we denote by $ \Cc(\Gamma)$. The objects of 
$ \Cc(\Gamma) $ are the faces of $\Gamma $ (i.e., the elements of $ \Sigma_\Gamma $), 
while the (non-identity) morphisms are of the form $\,w \leftarrow e(s_{\alpha}, w) \to  s_{\alpha} w\,$, where $w \in W $ and $ \alpha \in \cR_+ $. Thus, $ \Cc(\Gamma)$ is the opposite face category of $ \Gamma $. 
\begin{figure}[ht]
\begin{center}
\begin{tikzpicture}[
  vertex/.style={
    rectangle, rounded corners=3pt, draw=black!70, thick,
    fill=black!4, inner sep=4pt, minimum height=20pt,
    font=\small
  },
  aedge/.style={thmframe, semithick},
  bedge/.style={defframe, semithick},
  gedge/.style={exframe!80!black, semithick},
  elabel/.style={font=\scriptsize, fill=white, inner sep=1.5pt, rounded corners=1pt},
  scale=0.8
]
  \node[vertex] (e)   at (0,0)       {$1$};
  \node[vertex] (s1)  at (-2.4,1.6)  {$s_1$};
  \node[vertex] (s2)  at (2.4,1.6)   {$s_2$};
  \node[vertex] (s12) at (-2.4,3.8)  {$s_1 s_2$};
  \node[vertex] (s21) at (2.4,3.8)   {$s_2 s_1$};
  \node[vertex] (w0)  at (0,5.4)     {$s_3$};


  \draw[aedge] (e)   -- node[elabel]  {$e(s_1,1)$} (s1);
  \draw[aedge] (s2)  -- node[elabel, pos=0.25]  {$e(s_1,s_2)$} (s12);
  \draw[aedge] (s21) -- node[elabel]  {$e(s_1,s_2s_1)$} (w0);

  \draw[bedge] (e)   -- node[elabel] {$e(s_2,1)$} (s2);
  \draw[bedge] (s1)  -- node[elabel, pos=0.25] {$e(s_2,s_1)$} (s21);
  \draw[bedge] (s12) -- node[elabel] {$e(s_2,s_1s_2)$} (w0);

  \draw[gedge] (e)   -- node[elabel, pos=0.25] {$e(s_3,1)$} (w0);
  \draw[gedge] (s1)  -- node[elabel]  {$e(s_3,s_1)$} (s12);
  \draw[gedge] (s2)  -- node[elabel] {$e(s_3,s_2)$} (s21);

\end{tikzpicture}
\end{center}
\caption{The Bruhat moment graph of $S_3$}
\label{GammaS3}
\end{figure}
The graph $ \Gamma $ has a natural geometric interpretation in terms of the flag manifold $G/T$ (see \cite{Car94}). To describe it we equip $ G/T $ with the structure of a complex algebraic variety:  we fix a complex reductive group $ \bG $, containing $G$ as its maximal compact subgroup, and identify $\, G/T  = \bG/\bB \,$ via the natural inclusion $ G \into \bG $, where $ \bB $ is a Borel subgroup in $ \bG $. The Bruhat graph $ \Gamma $ describes then the structure of $0$- and $1$-dimensional $\bT$-orbits, where $ \bT \subset \bB $ is a maximal complex algebraic torus containing $T$ as a compact abelian subgroup. Specifically, the vertices of $ \Gamma$ are in natural bijection with the  $0$-dimensional orbits (i.e. the  $\bT$-fixed points) in $G/T$, while the edges $ e(s_{\alpha}, w) \in E_{\Gamma} $ correspond to the complex one-dimensional $\bT$-orbits, which we denote by $ \O_{e(s_{\alpha}, w)} $.
Two vertices, say $w$ and $w'$, are connected by an edge $e$  in $ \Gamma $ iff there is
a complex one-dimensional orbit $ \O_e $ in $G/T$ such that $ \overline{\O}_e \!\setminus\! \O_e = \{w,w'\}$, where $ \overline{\O}_e $ denotes the (Zariski) closure of $ \O_e $.
This interpretation of $\Gamma $ allows us to define the canonical (tautological) functor 
%
\begin{equation*}
 \la{InF0}
\cF_{0}\,:\ \Cc(\Gamma)\,\to \, \Top^T\,,\quad w \mapsto \{w\}\ ,\quad e(s_{\alpha}, w) \mapsto \O_{e(s_{\alpha}, w)}\,,
\end{equation*}
which we call the $T$-orbit functor on $ \Cc(\Gamma)$.

Now, crucial to our construction is the fact that the category $ \Cc(\Gamma) $ carries a (strict) $W$-action induced by the natural action of $W$ on the moment graph $\Gamma$. Explicitly, this action is defined by
\begin{equation}
\la{InW1} 
 w \mapsto gw\ , \quad e(s_{\alpha}, w) \mapsto e(g s_{\alpha} g^{-1},\, gw)\ , \quad \forall\,g \in W\,.
\end{equation}
Taking the homotopy quotient of $ \Cc(\Gamma) $ with respect to \eqref{InW1} (in the category $ \mathsf{Cat} $ of small categories), we define a larger category $ \Cc(\Gamma)_{hW} $ with the property $ B[\Cc(\Gamma)_{hW}] \simeq  B[\Cc(\Gamma)]_{hW} $. There is an explicit canonical model for this category given by the Grothendieck construction: the objects of $\Cc(\Gamma)_{hW}$ 
are the same as in $ \Cc(\Gamma)$, while the morphisms are given by pairs 
$(g,\,\varphi) $, where $ g \in W $ and $ \varphi \in \Hom_{\Cc(\Gamma)}(c, \,gc')$ for $ c,c' \in \Cc(\Gamma)$, see \eqref{ChW} and \eqref{morChW}.

Following \cite{DM16}, we say that a diagram of spaces indexed by $\Cc(\Gamma)$ 
(i.e, a functor $ \Cc(\Gamma) \to \Top $) admits a {\it $W$-structure} if it extends to 
the category $ \Cc(\Gamma)_{hW} $. 
The $T$-orbit functor \eqref{InF0} itself does {\it not} admit a natural $W$-structure; however, the corresponding (induced) functor $ G \times_T \cF_{0}\,:\ \Cc(\Gamma) \to \Top^G $ with values in
$G$-spaces does (see Lemma~\ref{lem:W-fun-heq}). The $W$-structure on $ G \times_T \cF_{0} $ comes from the natural action of the normalizer $N_{\bG}(\bT) $ of the algebraic maximal torus on the $\bT$-equivariant $1$-skeleton of $G/T$, 
first studied in \cite{Kn03} (see Remark~\ref{pplus0}). Thus, we get a well-defined functor 
\begin{equation}
 \la{InFG0}
G \times_T \cF_{0}\,:\ \Cc(\Gamma)_{hW}\,\to \, \Top^G\,.
\end{equation}
The next theorem is the main result of Section~\ref{S4}.
\newtheorem*{hDec}{Theorem~\ref{hDec}}
\begin{hDec}
For every prime $p \not= 2 $, there is a $G$-equivariant mod-$p$ cohomology isomorphism 
\begin{equation}
\la{InGcohiso}
\hocolim_{\Cc(\Gamma)_{hW}} [\,G \times_T \cF_{0}\,] \,\simeq \, G/T\,, 
\end{equation}
which, abusing terminology, we refer to as a {\rm homotopy decomposition} of $G/T$ over  $ \Gamma$.
\end{hDec}
The result of Theorem~\ref{hDec} is not quite immediate: in particular, there is no natural map $ G \times_T \cF_{0} \to  G/T$ inducing the equivalence \eqref{InGcohiso}. Instead,
we construct a `zig-zag' of natural maps of $W$-diagrams:
\begin{equation}
\la{Inpq}
G \times_T \cF_0 \,\xleftarrow{\bar{p}}\, G \times_T \wcF_0 \,\xrightarrow{\bar{q}}\, G/T
\end{equation}
where $ \,\bar{p} \,$ is a weak homotopy equivalence, while $\,\bar{q}\,$ fits in a canonical fibration sequence whose fiber is $p$-acylic at every prime except $ p \not=2$ (see Lemma~\ref{lem:act-fib}). The functor $\wcF_{0} $ that appears in the above correspondence plays an important role in this paper: from an abstract homotopy-theoretic point of view, it can be viewed as a  cofibrant replacement (`resolution') of the $T$-orbit functor $ \cF_0 $ in category of $\Cc(\Gamma)$-diagrams (see  Lemma~\ref{Lcofres}); somewhat surprisingly, it admits a natural geometric construction:
\begin{equation}
 \la{InFFP0}
\wcF_{0}: \,\Cc(\Gamma) \,\to\, \Top^{T}\ ,\quad w \mapsto \bigvee_{e(s_{\alpha}, w) \in E_{\Gamma}(w)}\!\! \!\overline{\O}_{e(s_{\alpha}, w)}\!\setminus\!\{s_{\alpha} w\}\ ,\quad e(s_{\alpha}, w) \mapsto \O_{e(s_{\alpha}, w)}\,,
\end{equation}
where $ E_{\Gamma}(w) $ denotes the subset of all edges of $\Gamma $ emanating from the 
vertex $w \in V_{\Gamma}$.

\subsection{The definition of quasi-flag manifolds}\la{S1.5}
The homotopy decomposition of $G/T$ given by Theorem~\ref{hDec} may not be very useful on its own\footnote{Notice that $G/T$ itself is `hidden' in the homotopy colimit \eqref{InGcohiso} as the space assigned  to the vertex objects in the diagram \eqref{InFG0}.}, but it serves us as a motivation for the general definition of quasi-flag manifolds. This proceeds in two 
steps.

\subsubsection{Simplicial thickening of moment graphs} 
The first step consists in replacing the moment graph $\Gamma $ --- or rather its face category $ \Cc(\Gamma) $ --- with the cell category $\Cc^{(m)}(\Gamma) $ of a higher dimensional CW-complex $\Gamma^{(m)} $ depending on a given multiplicity function $ m \in \M(W)$. To construct this complex we will regard $m$ as a function $\, e(s_{\alpha}, w) \mapsto m_{\alpha} \,$ on the  edges of 
$\Gamma $. Note that, topologically, the graph $\Gamma$ is  obtained by gluing (`concatenating' together) the geometric $1$-simplices: $\,\Gamma = \Delta^1 \amalg_{\Delta^0} \Delta^1 \amalg_{\Delta^0} \,\ldots \,\amalg_{\Delta^0} \Delta^1\,$, one $\Delta^1$ for each edge $ e(s_{\alpha}, w) \in E_\Gamma $. We can `thicken' this complex in a natural way by replacing each of its $1$-simplices $ \Delta^1 = S \Delta^0 $ by the unreduced suspension (`bicone') $ S \Delta^{m_{\alpha}} $ over the $ m_{\alpha}$-dimensional simplex $\Delta^{m_{\alpha}}$, where $m_\alpha$ corresponds to the chosen edge $ e(s_{\alpha}, w) \in E_{\Gamma}$. As a result, we get an $(1+\max\{m\})$-dimensional CW-complex: 
\begin{equation*}
\la{InGm}
\Gamma^{(m)} \,:= \, S\Delta^{m_{\alpha_1}}\,\amalg_{\Delta^0}\, S\Delta^{m_{\alpha_2}}\,\amalg_{\Delta^0}\, \ldots \, \amalg_{\Delta^0} S\Delta^{m_{\alpha_r}}   
\end{equation*}
in which two ‘fat edges' $ S \Delta^{m_{\alpha_i}}$ and $ S\Delta^{m_{\alpha_{i+1}}}$
are glued to each other at their apex points
according to the incidence relations of the original graph $\Gamma$ (see Figure~\ref{bicone-graph}).

\begin{figure}[htbp]
\centering
\begin{tikzpicture}[
  knot/.style={circle, fill=white, draw=black!65, semithick,
               inner sep=0pt, minimum size=4pt},
  lead/.style={thin, draw=black!45},
  scale=1
]
  \def\coneR{7.5pt}\def\waistR{2.6pt}
  \def\bicone#1#2#3{%
    \coordinate (bcA) at (#2);
    \coordinate (bcB) at (#3);
    \coordinate (bcM) at ($(bcA)!.5!(bcB)$);
    \pgfmathanglebetweenpoints{\pgfpointanchor{bcA}{center}}%
                              {\pgfpointanchor{bcB}{center}}%
    \edef\bcang{\pgfmathresult}%
    \pgfpointdiff{\pgfpointanchor{bcA}{center}}%
                 {\pgfpointanchor{bcB}{center}}%
    \pgfgetlastxy{\bcx}{\bcy}%
    \pgfmathsetlengthmacro{\bcL}{veclen(\bcx,\bcy)/2}%
    \begin{scope}[shift={(bcM)}, rotate=\bcang]
      \shade[left color=#1!40, right color=#1!14]
         (-\bcL,0pt) -- (0pt,\coneR) -- (0pt,-\coneR) -- cycle;
      \shade[left color=#1!14, right color=#1!40]
         ( \bcL,0pt) -- (0pt,\coneR) -- (0pt,-\coneR) -- cycle;
      \shade[left color=#1!26, right color=#1!14]
         (0pt,0pt) ellipse [x radius=\waistR, y radius=\coneR];
      \draw[#1!50!black, densely dotted, thin]
         (0pt,\coneR) arc[start angle=90, end angle=270,
                          x radius=\waistR, y radius=\coneR];
      \draw[#1!65!black, semithick]
         (0pt,-\coneR) arc[start angle=-90, end angle=90,
                           x radius=\waistR, y radius=\coneR];
      \draw[#1!65!black, semithick]
         (-\bcL,0pt) -- (0pt,\coneR) -- (\bcL,0pt) -- (0pt,-\coneR) -- cycle;
    \end{scope}}
 
  \coordinate (C) at (0,0);
  \foreach \i in {1,...,6}{\coordinate (V\i) at ({60*\i-30}:2.45);}
 
  \bicone{thmframe}{C}{V1}          \bicone{defframe}{C}{V2}
  \bicone{exframe!80!black}{C}{V3}  \bicone{thmframe}{C}{V4}
  \bicone{defframe}{C}{V5}          \bicone{exframe!80!black}{C}{V6}
  \bicone{exframe!80!black}{V1}{V2} \bicone{thmframe}{V2}{V3}
  \bicone{defframe}{V3}{V4}         \bicone{exframe!80!black}{V4}{V5}
  \bicone{thmframe}{V5}{V6}         \bicone{defframe}{V6}{V1}
 
  \foreach \p in {C,V1,V2,V3,V4,V5,V6}{\node[knot] at (\p) {};}
 
  \draw[lead] ($(V6)!.5!(V1)$) -- ++(0.90,-0.18)
     node[anchor=west, font=\scriptsize, text=black!60] {$S\Delta^{m_{\alpha}}$};
  \draw[lead] (V1) -- ++(0.62,0.60)
     node[anchor=south west, font=\scriptsize, text=black!60] {$w$};
 
  \node[font=\small, anchor=west] at (2.75,-1.05)
    {\textcolor{thmframe}{\rule{12pt}{1.5pt}}\;$\alpha$};
  \node[font=\small, anchor=west] at (2.75,-1.65)
    {\textcolor{defframe}{\rule{12pt}{1.5pt}}\;$\beta$};
  \node[font=\small, anchor=west] at (2.75,-2.25)
    {\textcolor{exframe!80!black}{\rule{12pt}{1.5pt}}\;$\gamma$};
\end{tikzpicture}
\caption{A `thickened' planar graph.}
\label{bicone-graph}
\end{figure}
In Section~\ref{subsec:regCW}, we will show that $\Gamma^{(m)}$ carries a regular CW-structure and define $ \Cc^{(m)}(\Gamma) $ as the cell category of $\Gamma^{(m)}$. Furthermore, we refine this structure by making $\Gamma^{(m)}$ a polyhedral complex, which we call a {\it spherical lune complex} $\Gamma_{\rm sph}^{(m)}$, and identify $ \Cc^{(m)}(\Gamma) \cong \Cc[\Gamma_{\rm sph}^{(m)}] $ (see Theorem~\ref{thm:lune}).

Alternatively, the category $ \Cc^{(m)}(\Gamma) $ can be defined in a  algebraic way by 
$$
\Cc^{(m)}(\Gamma) \, :=\, \Cc(\Gamma) \smallint \sd^{(m)}(\Gamma)
$$
where ``$\,\smallint\,$'' stands for the classical Grothendieck construction (see \ref{Grothcon}) applied to the functor 
\begin{equation}
\la{Insdm}
\sd^{(m)}(\Gamma):\ \Cc(\Gamma)\,\to\, {\sf Cat}\ , \quad e(s_{\alpha}, w) 
 \mapsto {\rm sd}(\Delta^{m_{\alpha}})^{\rm op}\,. 
\end{equation} 
We call \eqref{Insdm} the {\it $m$-simplicial thickening} functor on $\Cc(\Gamma)$: it assigns to every edge object $e(s_{\alpha}, w)$ of $\Cc(\Gamma)$ the opposite poset of the barycentric subdivision $ {\rm sd}(\Delta^m)$ of the $m$-simplex $\Delta^m $ (i.e., the poset of all faces of $ \Delta^m $ ordered by reverse inclusion), while mapping each vertex object to $\Delta^0$.
Explicitly,  $\Cc^{(m)}(\Gamma)$ comprises the objects of two kinds:
the elements of $W$ (vertices) and the collections of `fat edges' $ \{e_{\sigma}(s_{\alpha}, w)\}_{\sigma \subseteq [m_{\alpha}]} $ indexed by the finite subsets of $ [m_{\alpha}] = \{0,1,2,\ldots, m_{\alpha}\}\,$.
The (nontrivial) morphisms in $\Cc^{(m)}(\Gamma)$ are of the form: $\, w \leftarrow e_{\sigma}(s_{\alpha},w) \rightarrow s_{\alpha} w \,$ for all $\,\sigma \subseteq [m_{\alpha}]$, and $ e_{\sigma}(s,w)\to e_{\tau}(s', w') $ for all $\tau \subseteq \sigma $ whenever $ s' = s $ and $ w'=w$ or $w' =sw$ (see Figure~\ref{CGamma}). 

Now, the $T$-obit functor \eqref{InF0} extends naturally to $\Cc^{(m)}(\Gamma) $ by 
\begin{equation}
\la{InFm}
\cF_m\,:\ \Cc^{(m)}(\Gamma)\,\to\, \Top^T\ ,\quad w \mapsto \{w\}\ ,\quad 
e_{\sigma}(s_\alpha, w) \,\mapsto\,\O^{\times \sigma}_{e(s_\alpha,w)}\ ,
\end{equation}
where $\,\O^{\times \sigma}_{e(s_\alpha,w)}  $ denotes the product of $\,|\sigma|$ copies of the $\bT$-orbit $\O_{e(s_\alpha,w)}$ indexed by the elements of the subset $ \sigma \subseteq [m_{\alpha}]$, with  $T \subset \bT $ acting diagonally. The morphisms $ e_{\sigma}(s_{\alpha}, w) \to e_{\tau}(s_{\alpha}, w')$ are taken by \eqref{InFm} to the canonical projections $\,\O^{\times \sigma}_{e(s_\alpha,w)} \onto \O^{\times \tau}_{e(s_\alpha,w)}\,$ for $ \sigma \supseteq \tau $. 
Since $ m \in \M(W) $ is $W$-invariant, the $m$-thickening functor \eqref{Insdm} admits a natural $W$-structure with respect to the $W$-action \eqref{InW1} on $ \Cc(\Gamma)$: this, in turn, makes $\Cc^{(m)}(\Gamma)$ a $W$-category. Just as in the case $m=0$, the $T$-orbit functor \eqref{InFm} itself does not admit a $W$-structure, but the induced $G$-equivariant functor $ G \times_T \cF_m $ does (see Lemma~\ref{WstrT}). Thus, we have a well-defined functor
\begin{equation}
\la{InFGm}
G \times_T \cF_m\,:\ \Cc^{(m)}(\Gamma)_{hW}\,\to\,\Top^G\ ,\quad 
e_{\sigma}(s_\alpha, w) \,\mapsto\, G \times_T \O^{\times \sigma}_{e(s_\alpha,w)}\,.
\end{equation}
We emphasize that the functor \eqref{InFGm} depends inherently on the geometry of $G/T$ via the $W$-structure on $ \cF_m $ induced by the left action of $N_G(T)$ on $G/T$ that normalizes the $T$-action. 

This leads us to our first definition motivated by Theorem~\ref{hDec}:
\begin{equation}
\la{InFmdef1}
F_m(G,T)\,:=\, {\rm hocolim}_{\Cc^{(m)}(\Gamma)_{hW}}\,[\,G \times_T \cF_m ]\,
\end{equation}
Cohomologically and $K$-theoretically, the spaces  \eqref{InFmdef1} 
are legitimate models for $m$-quasi-flag manifolds: in particular, all results stated in the beginning of this Introduction (including the isomorphisms \eqref{In6}-\eqref{In8}) hold true 
for \eqref{InFmdef1}. However, these spaces turn out to be non-simply connected
(see Proposition~\ref{fgrProp}), and hence {\it cannot} represent (rationally) the coaffine stacks defined in Section~\ref{S3}. We need to refine the definition \eqref{InFmdef1} making the spaces $F_m(G,T)$ simply connected. This brings us to the second step of our construction.

\subsubsection{$p$-plus construction}
In homotopy theory, there is a well-known procedure --- called the Quillen plus construction \cite{Q70} --- that eliminates the fundamental group of a path-connected space without changing its cohomology. 
This construction applies only to spaces whose fundamental group is perfect, 
which is, unfortunately, not the case for  \eqref{InFmdef1}. 
In the recent paper \cite{BLO21}, C. Broto, R. Levi and B. Oliver proposed a natural (relative) generalization of the  plus construction that does apply in our situation. Given a commutative ring $R$ and a connected CW-complex $X$, the authors of \cite{BLO21} established necessary and sufficient conditions for the existence of a plus construction $ X \to X^{R+}$ on $X$ {\it relative to}  $R$ (see Theorem~\ref{plusthm}). When $ R = \F_p $ ($p$ prime) or $ R = \Q $ ($p=0$), we refer to this relative plus construction as {\it $p$-plus construction} and, abusing notation, will denote it simply by $X^{+}$.

It turns out that our spaces $ F_m(G,T) $ satisfy the Broto-Levi-Oliver conditions for all values of $p$, except $p=2$. In fact, somewhat surprisingly, the $p$-plus construction can be performed on $ F_m(G,T) $ simultaneously for all these values of $p$ and, moreover, naturally in $m$. To be precise, in Section~\ref{S5.3}, we construct a diagram of spaces 
$$ 
F_*^{+}(G,T): \ \M(W) \to \Top^G \ ,\quad m \mapsto F_m^{+}(G,T)\,,
$$ 
together with a canonical morphism of diagrams $q_*: F_\ast(G,T) \to  F_*^{+}(G,T) $ such that, for all $ m \in \M(W) $, the maps $ q_m:\, F_m(G,T) \to  F_m^{+}(G,T)\,$ provide (simultaneous) $p$-plus constructions for $p=0$ and all primes $p > 2$. The key observation behind this construction is
Proposition~\ref{fgrProp}, which shows that the fundamental groups of all spaces $ F_m(G,T) $ are naturally isomorphic to the fundamental group of the category $ \Cc(\Gamma)_{hW} $ (in particular, they are independent of $m$).  To make these spaces simply connected we can therefore modify the definition \eqref{InFmdef1} by modifying our construction of the $m$-thickened moment categories $ \Cc^{(m)}(\Gamma)_{hW} $. For each $ m \in \M(W) $, the category $ \Cc(\Gamma)_{hW} $ canonically
embedds (as a full subcategory) into $ \Cc^{(m)}(\Gamma)_{hW} $, and  we can define a new category $ \Cc^{(m)}(\Gamma)^{+}_{hW} $ by homotopically contracting the image of this 
embedding to the terminal (point) object $\ast$ in $\Cat $. Just as in the case of topological spaces, such a contraction amounts to taking the homotopy cofiber of the inclusion functor
 $ \Cc(\Gamma)_{hW} \into \Cc^{(m)}(\Gamma)_{hW} $, which, in turn, amounts to erecting a cone over  $ \Cc(\Gamma)_{hW} $ inside $ \Cc^{(m)}(\Gamma)_{hW} $. Explicitly, the category  $ \Cc^{(m)}(\Gamma)^{+}_{hW} $ is represented by the Grothendieck construction
 \begin{equation}
 \la{InCGG-}
 \Cc^{(m)}(\Gamma)^{+}_{hW} \,=\,\I \smallint \{\ast \leftarrow \Cc(\Gamma)_{hW} 
 \into  \Cc^{(m)}(\Gamma)_{hW}\} \,,
 \end{equation}
where $\, \I = \{1 \leftarrow 0 \rightarrow 2 \}\,$ is the standard category parametrizing pushouts. The functors on \eqref{InCGG-} can be naturally represented by diagrams consisting
of three functors (defined on each of the categories $ \Cc(\Gamma)_{hW}$, $\, \Cc^{(m)}(\Gamma)_{hW} $ and $\,\ast\,$) and two natural transformations (see Lemma~\ref{funIG}). 
It turns out that our basic functor $G \times_T \cF_m $ admits a canonical extension to the category
$ \Cc^{(m)}(\Gamma)^{+}_{hW} $:
\begin{equation}
 \la{InFun-}
(G\times_T \cF_m)^{+}: \ \Cc^{(m)}(\Gamma)^{+}_{hW} \to \Top^G 
\end{equation}
defined by the following diagram
\begin{equation}
\la{IndiagFm-}
%
\begin{tikzcd}[scale cd= 0.9]
	\ast && {{\mathscr C}(\Gamma)_{hW}} && {{\mathscr C}^{(m)}(\Gamma)_{hW}} \\
	\\
	&& {\Top^G}
	\arrow["{G/T}"', curve={height=6pt}, from=1-1, to=3-3]
	\arrow[from=1-3, to=1-1]
	\arrow["", from=1-3, to=1-5]
	\arrow["{G\times_T\widetilde{\mathcal F}_0}"{description}, from=1-3, to=3-3]
	\arrow["{G \times_T\mathcal{F}_m}", curve={height=-6pt}, from=1-5, to=3-3]
	\arrow["{\tilde{q}}"', shift right=5, shorten <=24pt, shorten >=24pt, Rightarrow, from=3-3, to=1-1]
	\arrow["{\tilde{p}}", shift left=5, shorten <=27pt, shorten >=27pt, Rightarrow, from=3-3, to=1-5]
\end{tikzcd}
\end{equation}
Note that besides \eqref{InFGm}, the diagram \eqref{IndiagFm-} naturally 
incorporates the $G$-equivariant functor $ G \times_T \wcF_0 $ induced by 
geometric resolution \eqref{InFFP0} of the $T$-orbit functor $ \cF_0 $
and also the two canonical natural transformations \eqref{Inpq} 
related to the homotopy decomposition of $G/T$ in Theorem~\ref{hDec}.
We can now finally state our main definition.
\newtheorem*{maindef}{Definition~\ref{maindef}}
\begin{maindef}
For $m\in \M(W) $, we define the $m$-{\rm quasi-flag manifold} associated to $ (G,T)$ by
\begin{equation}
\la{InFmdef2}
F^{+}_m(G,T) \,:=\, \hocolim_{ \Cc^{(m)}(\Gamma)^{+}_{hW}}(G \times_T{\cF_m)}^+\, 
\end{equation}
where $(G \times_T{\cF_m)}^+$ is the diagram on $ \Cc^{(m)}(\Gamma)^{+}_{hW} $ represented by \eqref{IndiagFm-}.
\end{maindef}
Assuming that the Lie group $G$ is simply connected, we prove (see Corollary~\ref{corpplus}) that all
quasi-flag manifolds $F_m^{+}(G,T)$ are simply connected and have the same mod-$p$ cohomology groups as $F_m(G,T)$ for all $p\not=2$. We also allow ourselves to abuse terminology and 
refer to the spaces $ F_m(G,T) $ defined by \eqref{InFmdef1} as {\it non-simply connected models} of  $m$-quasi-flag manifolds. We may (and will) use $F_m(G,T)$ and $ F_m^+(G,T)$ interchangeably in cohomological computations.

Despite their intimidating appearance, the homotopy colimits \eqref{InFmdef1} and \eqref{InFmdef2} can be actually calculated. In Section~\ref{S5.4}, we construct explicit
geometric models for these spaces by gluing together even-dimensional spheres --- the dimensions depending on the multiplicities $m_\alpha$ --- with gluing data determined by the combinatorics of the Bruhat graph (see \eqref{realFm} and \eqref{rmodel}). A geometrically minded reader may consider these models as definitions of quasi-flag manifolds.

\subsection{Cohomology of quasi-flag manifolds} In Section~\ref{S6}, we prove our main result --- Theorem~\ref{MainT} --- that provides rational algebraic models for the $m$-quasi-flag manifolds and their homotopy quotients in terms of coaffine stacks constructed in Section~\ref{S3}. This theorem can be viewed as a solution to the topological realization problem for rings of quasi-invariants posed in \cite{BR1}.

Assume that $G$ is a compact simply connected Lie group with maximal torus
$T \subseteq G $ and the Weyl group $ W = W_G(T) $. For $ m \in \M(W)$, write
$ F_m^+ = F_m^+(G,T) $ for the $m$-quasi-flag manifold \eqref{InFmdef2} associated to $(G,T)$, and denote by
$ X_m^+ := EG \times_G F_m^+ $ and $\bX_m^+ := ET \times_T F^+_m $ the Borel quotient spaces of homotopy $G$-oribts and $T$-orbits in $ F^{+}_m$, respectively. By Corollary~\ref{corpplus}, these spaces are simply connected, hence  there are associated coaffine stacks: $\cSpec \,[\,C_{\Q}^*(F^+_m)\,] $, $\,\cSpec \,[\,C_{\Q}^*(X^+_m)\,] $ and $\,\cSpec \,[\,C_{\Q}^*(\bX^+_m)\,]$, where $\,C^*_{\Q}(X) = C^*(X, \Q)\,$ stands for the rational cochain complex of a space $X$  equipped with the structure of an $\bE_{\infty}$-algebra over $\Q$ (see Appendix~\ref{AB}).
\newtheorem*{MainT}{Theorem~\ref{MainT}}
\begin{MainT}
For all $ m \in \M(W) $, there are natural equivalences in $\cAff_{\Q} $:\begin{eqnarray}
\cSpec \,[\,C_{\Q}^*(F^+_m)\,] & \simeq & F^c_{[\frac{m+1}{2}]}(W)\,, \nonumber
\\
\cSpec \,[\,C_{\Q}^*(X^+_m)\,] & \simeq & V^c_{[\frac{m+1}{2}]}(W)\,,\la{IncoaffX}\\
\cSpec \,[\,C_{\Q}^*(\bX^+_m)\,] & \simeq & \bU^c_{[\frac{m+1}{2}]}(W)\,, \nonumber
\end{eqnarray}
where  $V^c_{k}(W)$, $\,F^c_{k}(W)$ and 
$\, \bU^c_{\!k}(W) $ are the coaffine stacks of quasi-invariants defined in Section~\mbox{\rm \ref{S3.2}} $($see \eqref{defvmc}, \eqref{Fck} and \eqref{bvck}, respectively$)$.
\end{MainT}      

The proof of Theorem~\ref{MainT} --- together with preliminary results ---
occupies almost the whole of Section~\ref{S6}. The main problem is of  homotopy-theoretic nature, and its solution given in Section~\ref{S6} may be of independent interest. Recall that our quasi-flag manifolds are defined in terms of canonical homotopy decompositions over the simplicially thickened moment categories $ \Cc^{(m)}(\Gamma)_{hW}$ and $ \Cc^{(m)}(\Gamma)^{+}_{hW}$. On the other hand, the algebraic coaffine stacks that appear in the equivalences \eqref{IncoaffX} 
are constructed in Section~\ref{S3} by gluing diagrams of stacks over the poset $ \Sc(W) $ of elementary reflection subgroups of $W$. To identify these objects we thus need to compare the homotopy colimits over the categories  $ \Cc^{(m)}(\Gamma)_{hW} $ and $ \Sc(W) $. The problem is that these last two categories are not directly related: there is no natural functor in either direction. Instead, we construct a sequence (`zig-zag') of five functors: 
%
\begin{equation}
\la{Inzig}
\begin{tikzcd}[scale cd=0.8]
	& {\Sc(W)\smallint \Phi} && {\Cc_F(\Pc)} && {\Cc^{(m)}(\Gamma)_{hW}} \\
	{\Sc(W)} && {\sd(\Pc)^{\rm op}} && {\overline{\Cc(\Gamma)}_{hW}}
	\arrow[from=1-2, to=2-1]
	\arrow[from=1-2, to=2-3]
	\arrow[from=1-4, to=2-3]
	\arrow["{p_F}", from=1-4, to=2-5]
	\arrow[from=1-6, to=2-5]
\end{tikzcd}
\end{equation}
where the three functors pointing to the left are cocartesian fibrations, while the two pointing to the right are homotopy cofinal (i.e., preserve all homotopy colimits under restriction).

The most interesting in \eqref{Inzig} is the functor $p_F $. Its target $ \overline{\Cc(\Gamma)}_{hW} $ is a canonical skeletal subcategory  of our quotient moment category $ \Cc(\Gamma)_{hW} $ (see Figure~\ref{skC}), and $ \Pc = \Pc(\overline{\Cc(\Gamma)}_{hW}) $ is the {\it universal poset} under $ \overline{\Cc(\Gamma)}_{hW} $. By definition, $\Pc$ comes with a canonical functor $ F: \ \overline{\Cc(\Gamma)}_{hW} \to \Pc $, which is initial among all functors from $\overline{\Cc(\Gamma)}_{hW}$ to poset categories in $\Cat$. The functor $F$ determines naturally a simplicial presheaf $\, {\rm Sing}(NF): \,\sd(\Pc)^{\rm op} \to \sset \,$ on the barycentric subdivision of the poset $ \Pc $, which, in turn, corresponds --- under Lurie's Straightening Equivalence \cite[3.2.5]{HTT} --- to the cocartesian fibration $ \Cc_{F}(\Pc) \to \sd(\Pc)^{\rm op} $ that appears as the middle arrow in \eqref{Inzig}. The category $ \Cc_{F}(\Pc) $ thus defined comes with a natural functor $ p_F:\, \Cc_{F}(\Pc) \to \overline{\Cc(\Gamma)}_{hW} $ that has the important property of being homotopy cofinal. We deduce this property from a general Cofinality Theorem~\ref{SlThmInf}  that holds for an arbitrary $\infty$-functor $ F: \Cc \to N\Pc $ from an $\infty$-category $\Cc$ to the nerve of a poset. We defer the proof of Theorem~\ref{SlThmInf}
and its corollaries  to Appendix~\ref{AC}.

Now, using the  functors \eqref{Inzig}, we can transform the homotopy decomposition \eqref{InFmdef1} over the category $ \Cc^{(m)}(\Gamma)_{hW} $, defining the space $F_m(G,T)$, to a homotopy decomposition over the reflection category $\Sc(W) $ (see Proposition~\ref{natdec}):
\begin{equation}
\la{IndecFm}
F_m(G,T)\,\simeq\, \hocolim_{\Sc(W)} \, (G \times_T \cF_m)^\natural\,.
\end{equation}
The functor $ (G \times_T \cF_m)^\natural:\, \Sc(W) \to \Top^G $ in \eqref{IndecFm}  is defined by
\begin{equation*}
\la{InGNW}
W_0 \,\mapsto\, G/T\ ,\quad W_\alpha \,\mapsto\, 
G \times_{N_{\alpha}}\bigl[(G_{\alpha}/T) \ast (T/T_{\alpha})^{\ast\, m_{\alpha}}\bigr]\,,
\end{equation*}
where $ T_\alpha \subset T $ is the singular subtorus in $T$ corresponding to the root $\alpha $,
$G_\alpha := C_G(T_\alpha) $ is the centralizer of $T_{\alpha} $ in $G$, and $ N_\alpha = G_\alpha \cap N_G(T) $ is the normalizer of $T$ in $ G_{\alpha} $ (all of which are compact subgroups in $G$).
The star `$ \ast $' in \eqref{InGNW} stands for the usual (topological) join operation, with $ (T/T_{\alpha})^{\ast\,m_\alpha}$ denoting the $m_{\alpha}$-fold join equipped with the 
diagonal action of $ N_{\alpha} $.  Having the decomposition \eqref{IndecFm} in hand,  one can deduce
the result of Theorem~\ref{MainT} by standard rational homotopy-theoretic argument using the Bousfield-Gughenheim functor $A_{\rm PL}  $ (For reader's convenience, we review the definition and formal properties of the $A_{\rm PL}$ in Appendix~\ref{AB}).

Theorem~\ref{MainT} has a number of important consequences. First, in combination with results of Section~\ref{S3}, it implies
\newtheorem*{CorHev}{Corollary~\ref{CorHev}}
\begin{CorHev}
Let $\k$ be a field of characteristic $0$. Then, for all $ m \in \M(W) $, there are natural isomorphisms of graded algebras
    \begin{eqnarray*}
    \la{Inevencohom}
H^{\rm ev}(F^+_m,\,\k) \,\cong\, H^{\rm ev}(F_m,\,\k) & \cong & Q_{k}(W)/\langle\k[V]^W_+\rangle\ ,\nonumber\\*[1ex]
H_G^{\rm ev}(F^+_m,\,\k) \,\cong\, H_G^{\rm ev}(F_m,\,\k) & \cong & Q_{k}(W)\ ,\\*[1ex]
    H_T^{\rm ev}(F^+_m,\,\k) \,\cong\, H_T^{\rm ev}(F_m,\,\k) & \cong & \bQ_{2k+1}(W)\nonumber\ ,
\end{eqnarray*}
where $\,k = [\frac{m+1}{2}] \in \M(W)\,$.
\end{CorHev}
Note that Corollary~\ref{CorHev} includes the isomorphisms \eqref{In7} given in the beginning of the Introduction.
Furthermore, together with Theorem~\ref{ThFree} and Theorem~\ref{Qfatodd}, it implies
\newtheorem*{CorHev2}{Corollary~\ref{CorHev2}}
\begin{CorHev2}
For all $ m \in \M(W) $, $\, H_T^{\rm ev}(F^{+}_m,\,\k) \cong H_T^{\rm ev}(F_m,\,\k) \,$ is a free $ H^*(BT,\, \k)$-module of rank $ |W|$. With 
identification $ H^*(BT,\,\k) = \k[V]\,$,
this module structure extends to a left module structure over the nil-Hecke algebra $ \NH_\k(W) $.
\end{CorHev2}
Observe that the $T$-equivariant cohomology of the quasi-flag manifolds
provides a topological interpretation for the algebras \eqref{In10} only
for odd multiplicities $ m = 2k+1$. To extend this interpretation to all $ m \in \M(W) $,
we consider the $T$-space
\begin{equation*}
\tilde{F}_m(T)\,:=\,
\hocolim_{ \Cc^{(m)}(\Gamma)}(\cF_m)\,,
\end{equation*}
where $ \cF_m $ is the generalized $T$-orbit functor defined by \eqref{InFm}.
It is easy to check that the associated $G$-space $\,\tilde{F}_m(G,T) := G \times_T  \tilde{F}_m(T)\,$ carries a {\it free} $W$-action, with $\,\tilde{F}_m(G, T)/W \cong F_m(G,T)\,$. Hence, we have the canonical $W$-bundle over $ F_m(G,T)$:
\begin{equation}
\la{InWbundle}
W \to \tilde{F}_m(G,T) \to F_m(G,T)
\end{equation}
In the end of Section~\ref{S6.3}, as yet another consequence of Theorem~\ref{MainT}, we prove
\newtheorem*{tilProp}{Proposition~\ref{tilProp}}
\begin{tilProp} 
For all $ m \in \M(W) $, there are natural isomorphisms of algebras
\begin{equation}
\la{InbQtop}
H_G^{\rm ev}(\tilde{F}_m(G,T),\,\k) \,\cong\, H_T^{\rm ev}(\tilde{F}_m(T),\,\k)\,\cong\, \bQ_{m+1}(W)\,,
\end{equation}
which fit in the commutative diagram
\begin{equation*}
    \begin{diagram}[small]
H_G^{\rm ev}(F_m(G,T),\,\k) & \rTo & H_G^{\rm ev}(\tilde{F}_m(G,T),\,\k)^W\\
\dTo & & \dTo\\
Q_{[\frac{m+1}{2}]} & \rTo & \bQ_{m+1}(W)^W 
\end{diagram}
\end{equation*}
where all arrows are isomorphisms, the top one being induced topologically by the bundle map 
\eqref{InWbundle}, while the bottom is the isomorphism \eqref{In111}.
\end{tilProp}
%

In Section~\ref{S6.4}, we construct and study a natural action of the Weyl group $W$ on the spaces $ F_m(G,T) $, and extend another well-known result of A. Borel  on cohomology of $G/T$ to quasi-flag manifolds. In its standard formulation (see, e.g., \cite[Theorem 1.1]{JM92}), Borel's Theorem asserts that, for any compact connected Lie group $G$, the natural map $\, p_N: BN \to BG \,$ induced by the inclusion of the normalizer $ N = N_G(T) $ of $T$ in $G$ is a $ \mod\,p $ cohomology isomorphism, provided $ p $ does not divide
the order of $W$. We can reformulate this result in terms of $G/T$ as follows. Recall that 
$G/T$ carries the following $W$-action commuting with natural $G$-action:
\begin{equation}
\la{WG/T}
W \times G/T \to G/T \ ,\quad (w, gT) \mapsto g\dot{w}^{-1} T \,,
\end{equation}
where $ \dot{w} \in N$ such that $ w = \dot{w}\, T \in W $. Although the action \eqref{WG/T} is known to be non-algebraic, it is topologically well-defined and, in fact, free. Hence, the space of its homotopy orbits can be easily calculated:
$$
(G/T)_{hW} \,\simeq\, (G/T)/W = (G/T)/(N/T) \,\cong \,G/N\,,
$$
which is precisely the homotopy fiber of the Borel map $p_N$. Thus, Borel's
Theorem can be restated by saying that the homotopy quotient $\,(G/T)_{hW}\,$ is $p$-acyclic for all primes $\, p \not|\ |W|\,$. We  prove

\newtheorem*{WFmTh}{Theorem~\ref{WFmTh}}
\begin{WFmTh}
For all $ m \in \M(W) $, the $G$-spaces  $ F_m(G,T)$ and $ F^+_m(G,T) $ carry a natural $W$-action induced by \eqref{WG/T}. The associated homotopy quotients
$ F_m(G,T)_{hW} $ and $ F^+_m(G,T)_{hW} $ are $p$-acyclic for all primes $\, p \not|\ |W|\,$.
Consequently, there are natural mod $p$ cohomology isomorphisms
\begin{equation}
\la{XmW}
X_m(G,T)_{hW} \,\xrightarrow{\sim} \, X^+_m(G,T)_{hW}\, \xrightarrow{\sim}\, BG
\end{equation}
where $X_m $ and $X^+_m $ are the spaces of homotopy $G$-orbits in $ F_m$ and $ F^+_m $, respectively.
\end{WFmTh}      

\noindent
Note that, for $ \k = \Q $ or $ \k = \F_p $ with $\, p $ prime to $ |W|\,$, \eqref{XmW} 
implies  (see Corollary~\ref{Coroddcoh})
$$
H^*_G(F^+_m,\,\k)^W \cong H^*_G(F_m,\,\k)^W \cong H^*(BG, \,\k) \cong \k[V]^W\,,
$$
which shows, in particular, that the $G$-equivariant cohomology of $F_m $ and $ F_m^+ $ over 
$\k$ contains no $W$-invariants in {\it odd}\, dimensions:
$\, H^{\rm odd}_G(F^+_m,\,\k)^W \cong H^{\rm odd}_G(F_m,\,\k)^W = 0\, $.

\subsection{Equivariant $K$-theory and exponential quasi-invariants}
In Section~\ref{S7}, we compute the $G$- and $T$-equivariant topological $K$-theory
of quasi-flag manifolds and --- motivated by these computations --- define
exponential analogues of the algebras $Q_k(W)$ and $ \bQ_m(W)$. We also extend a number of 
algebraic results from Section~\ref{S2} to the exponential setting. 

The main result of Section~\ref{S7} --- Theorem~\ref{KGFm} --- yields the ring isomorphisms
\begin{eqnarray}
K_G[F^{+}_m(G,T)] & \cong & K_G[F_m(G,T)]\,\cong\,  \cQ_{k}(W)\,, \la{Inexp}\\*[1ex]
 K_T[F_m^{+}(G,T)] &\cong & K_T[F_m(G,T)]\, \cong \,\bcQ'_{2k+1}(W)\,, \la{Inexpb}
\end{eqnarray}
which hold (upon tensoring with $\Z[\frac{1}{2}]$) for all $m \in \M(W)$ with $ k = [\frac{m+1}{2}]\,$. The commutative rings $ \cQ_k(W) $ and $ \bcQ'_m(W) $ that appear in the above isomorphisms are defined as follows. 

For any $ k \in \M(W) $, we set ({\it cf.} Definition~\ref{Expqi})
\begin{equation}
\la{IncQ_k}
\cQ_k(W)\, := \,\{f \in R(T)\ :\ s_{\alpha}(f) \equiv f\ \mod\ \langle (1 - e^{\alpha})(1-e^{\bar{\alpha}})^{2 k_{\alpha}}\rangle \ ,\ \forall\ \alpha \in \R_+\}\,,
\end{equation}
where $ R(T) $ denotes the representation ring of $T$ (i.e., the group ring $\Z[\hat{T}]$ of the character lattice $ \hat{T} := \Hom(T, U(1)) $ of the maximal torus) and, for $ \alpha \in \R_+$, the  $ \,\bar{\alpha} $ equals $ \frac{1}{2}\alpha $ or $ \alpha $ depending on whether 
$ \frac{1}{2}\alpha $ belongs to $ \hat{T} $ or not. 
We call \eqref{IncQ_k} the ring of {\it exponential quasi-invariants of $W$ of multiplicity $ k \in \M(W)$}. We should mention that the exponential analogues of rings of quasi-invariants (in fact, their $q$-deformations) appeared earlier in work of O. Chalykh \cite{Ch00, Ch02} (see Remark~\ref{R202}).  We also remark that the above definition \eqref{IncQ_k}, although stated in a different form, is equivalent to the one given in \cite{BR1} (see Remark~\ref{R22}).

Next, we introduce exponential analogues of the ring $\bQ_m(W)$ of $m$-quasi-covariants. Recall that $ \bQ_m(W) $ behaves differently depending on whether the multiplicity is even ($m=2k$) or odd $(m=2k+1)$. In the exponential
case, it is natural to define two {\it different} rings: $\, \bcQ_m(W) \,$ and $\, \bcQ_m'(W) $, the first being the analogue of $\bQ_{2k}(W)$ and the second of $\bQ_{2k+1}(W)$. The explicit definitions of these rings are ({\it cf.} Definition~\ref{expqmring}):
\begin{eqnarray*} 
\bcQ_m(W) & := & \bigg\{ \sum_{w \in W} p_w \otimes w \in R(T) \otimes \Z[W]\ :\ p_{s_{\alpha}w} \equiv p_w \ \mod\, \langle(1 - e^{\bar{\alpha}})^{m_\alpha}\rangle \,,\ \forall\,  \alpha \in \R_+\,\bigg\}\,,\la{Inbcqm2}\\*[1ex]
\bcQ_m'(W) & := & \bigg\{ \sum_{w \in W} p_w \otimes w \in R(T) \otimes \Z[W]\ :\ p_{s_{\alpha}w} \equiv p_w \ \mod\, \langle(1-e^{{\alpha}})(1 - e^{\bar{\alpha}})^{{m_\alpha}-1}\rangle \,,\ \forall\,  \alpha \in \R_+\,\bigg\}\,\la{Inbcqm1}
\end{eqnarray*}
It is easy to check that $ \bcQ_m'(W) $ is actually a subring of $\bcQ_m(W)$ such that $\,\bcQ_m'(W) = \bcQ_m(W) \cap \bcQ_1'(W)\, $ for all $ m \in \M(W)$. Moreover, both rings contain $ R(T) $ as a subring (under the obvious
map $\,f \mapsto \sum f \otimes w \,$), and both are stable under the (diagonal) action of $W$. Thus, for all $ m\in\M(W) $, we have natural $W$-equivariant inclusions
\begin{equation} 
\la{Inrtmodstr}
R(T) \,\subseteq \, \bcQ'_{m}(W) \,\subseteq\, \bcQ_m(W)\,\subseteq\,R(T) \otimes \Z[W]\,.
\end{equation}
When restricted to the (diagonal) $W$-invariants these  become
$$
R(T)^W \,\subseteq \, \bcQ'_{m}(W)^W \,=\, \bcQ_m(W)^W\,\subseteq\,R(T)\,.
$$
Now, in Lemma~\ref{QmulHalf}, we prove that
\begin{equation}
\la{InInvcQ}  
 \bcQ'_m(W)^W \,=\, \bcQ_m(W)^W \,\cong \, \cQ_{[\frac{m}{2}]}(W) \,,
\end{equation}
and further, in Proposition~\ref{bcQodd}, we establish the ring isomorphisms 
\begin{equation} \la{Inphimulk} 
R(T) \otimes_{R(T)^W} \cQ_k(W) \ \stackrel{\sim}{\to}\ \bcQ'_{2k+1}(W)\ ,\quad 
f \otimes g \mapsto \sum_{w \in W} f w(g) \otimes w\ ,
\end{equation}
for all $\,m,k \in \M(W)\,$. Formulas \eqref{InInvcQ} and \eqref{Inphimulk} are exponential analogues of the isomorphisms \eqref{In111} and \eqref{In12}, respectively.

In Section~\ref{S7.1}, we also prove some representation-theoretic results on nil-Hecke and double affine Hecke actions on  $ \bcQ_m'(W) $ and 
$\bcQ_m(W)$. Specifically, we show that, for $ m = 2k+1 $, the abelian group $ \bcQ_{m}'(W) $ carries a natural action of the Demazure nil-Hecke ring $\D(W) $ (see Theorem~\ref{cQfatodd}), while for $ m = 2k $, the $\k$-vector space
$\,\bcQ_{m}(W) \otimes \k\,$ is naturally a module over the {\it trigonometric} Cherednik algebra $ \bH^{\rm trig}_{k}(W)\,$
defined over $\k$ (see Theorem~\ref{cQfat}).
These results are exponential analogues of Theorem~\ref{Qfatodd} and Theorem \ref{Qfat}; their proofs, albeit parallel, are algebraically a bit more delicate.

Finally, there is a natural question regarding the generalization of Theorem \ref{ThFree}. Although we do not address this question in the present paper, we feel that it is safe to make the following
\newtheorem*{freeQmul}{Conjecture~\ref{freeQmul}}
\begin{freeQmul}
For all $k \in \M(W)$, after $($possibly$)$ inverting $ |W|$, the following is true:
\begin{enumerate}

\item[$(i)$] $\bcQ_{2k}(W)$ is a free module  of rank $|W|$ over $R(T)$,
\item[$(ii)$] $\cQ_{k}(W)$ is a free module of rank $|W|\,$ over $R(T)^W$,
\item[$(iii)$] $\bcQ'_{2k+1}(W)$ is a free module  of rank $|W|$ over $R(T)$,
\end{enumerate}
where the $R(T)$-module structure on $\bcQ_{2k}(W)$ and $\bcQ_{2k+1}(W)$ is defined via the inclusion \eqref{Inrtmodstr}.
\end{freeQmul}
We refer the reader to the end of Section~\ref{S7.1} for further discussion of this conjecture.

\subsection{Vista: GKM spaces with multiplicities}\la{S1.8}
The construction of $m$-quasi-flag manifolds given in this paper 
generalizes naturally to a broader class of spaces called GKM manifolds. 
A {\it GKM manifold} is a compact $T$-manifold $M$ equipped with a $T$-invariant complex (or almost complex) structure that satisfies the following conditions: (1) $M$ is $T$-equivariantly formal, (2) the set $M^T$ of $T$-fixed points in $M$ is finite, and (3) the isotropy representation of $T$ on the tangent space $ T_x(M)$ at each fixed point $x \in M^T $ is nondegenerate (in the sense that the weights of this representation are pairwise linearly independent). The cohomological properties of such a manifold --- in particular, its $T$-equivariant cohomology --- are determined combinatorially by a GKM moment graph $ \Gamma = (V_{\Gamma}, E_{\Gamma}) $ together with a labeling (called 
{\it axial}\,) function $ \alpha: E_{\Gamma} \to \Lambda(T)\setminus\{0\} $, where $ \Lambda(T)  \subset V^* $ is viewed as a lattice in the dual Lie algebra $\,V^* = {\rm Lie}(T)^*$ of $T$. The vertex set $ V_{\Gamma} $ of $\Gamma$
is simply $ M^T$, while the edge set $ E_{\Gamma} $ consists of (complex) closures of the one-dimensional $T$-orbits in $M$; the axial function $ \alpha: E_{\Gamma} \to \Lambda(T) $ is defined by taking the weights of the tangential isotropy representations of $T$ at fixed points of $M$ (see \cite{GZ01}). Introduced in \cite{GKM98} the class of GKM manifolds includes many important spaces that play a role both in geometry and representation theory (see, e.g., the recent surveys \cite{Fie12}, \cite{Fr24}, \cite{Tym24}, \cite{GKZ24} and references therein).
%
%
Of particular relevance to us are the so-called {\it GKM manifolds with extended $G$-action}: these are  manifolds equipped with a (smooth effective) action of a compact connected Lie group $G$, whose restriction to a maximal torus $ T \subseteq G$ satisfies the above geometric conditions. A prototypical example is $ G/T$ but there are many other interesting examples (see \cite{GHZ06}, \cite{Ku10, Ku11}, \cite{Wie12}, \cite{Ka15}, \cite{GM17}, \cite{KuM17}). If $ M $ is a GKM manifold with extended $G$-action, its moment graph $ \Gamma $ --- and hence the associated moment category $ \Cc(\Gamma)$ --- carry a natural $W$-action, the corresponding axial function $ \alpha: E_{\Gamma} \to \Lambda(T)  $ being $W$-equivariant. Motivated by our basic example of the Bruhat graph, we define a {\it multiplicity function} for $\Gamma $ to be a $W$-invariant function 
$ m: E_{\Gamma} \to \Z_{+} $, $\,e \mapsto m_e$, satisfying, in addition, the property that
$ m_{e'} = m_e $ whenever $ \alpha(e') = \pm\,\alpha(e)$ (i.e., $ m $ is required to be constant on the fibers of the composite map $\,E_{\Gamma} \xrightarrow{\alpha} \Lambda(T) \to \Lambda(T)/\{\pm 1\} \,$). We denote the set of all such multiplicity functions by $ \M(\Gamma) $. 

Now, starting with a GKM $G$-manifold $M$ and applying {\it mutatis mutandis} the procedure of Section~\ref{S5}, we can define a family of $G$-spaces $ M_m(G,T) $ and $ M^+_m(G,T) $ indexed by
the multiplicities $ m \in \M(\Gamma) $ that generalize our $m$-quasi-flag manifolds $F_m(G,T)$ and
$F_m^+(G,T)$. The $G$- and $T$-equivariant even-dimensional cohomology of these spaces can be computed explicitly in terms of the following commutative algebras 
\begin{eqnarray*}
Q_k(\Gamma)\!\!\!\! &:= & \!\!\!\{(p_{\bar{x}}) \in \prod_{\bar{x} \in M^T\!/W} \!\! \k[V]^{W_{x}}:\, \bar{w}(p_{\bar{x}}) \equiv p_{\bar{x}'}\ \mod\, \langle \alpha(e)\rangle^{2 k_e},\ \forall \, \bar{w} \in W/W_{x} : e(wx,x') \in E_{\Gamma}\}\ , \la{In90} \\*[1ex]
\bQ_m(\Gamma)\!\!\!\! &:=&\!\!\! \{\,(p_x)\in\prod_{x \in M^T}\!\k[V]\,:\, p_{x'} \equiv p_{x}\ \mod\, \langle \alpha(e) \rangle^{m_e},\ \forall \,e(x,x') \in E_{\Gamma}\,\} \,,\la{In100}
\end{eqnarray*}
where $ \bar{x} \in M^T\!/W $ stands for the $W$-orbit of a fixed point $x \in M^T $ and
$\,W_x := \{w \in W \ :\ wx=x\}\,$ denotes its stabilizer in $W$. In the case $ M = G/T $, the 
above formulas specialize to \eqref{In9} and \eqref{In10}, and, in general, one can easily check that 
$ Q_k(\Gamma) \cong \bQ_m(\Gamma)^W $ whenever $ m = 2k $. It is therefore natural to think of
(and refer to) the rings $ Q_k(\Gamma) $ and $ \bQ_m(\Gamma) $ as generalized quasi-invariants and quasi-covariants, respectively. The algebraic properties of these rings are similar to those of the classical ones; moreover, the rings $Q_k(\Gamma)$ seem to be closely related to other interesting generalizations of quasi-invariants that arise in mathematical physics ({\it cf.} \cite{BC23}).
We will study the spaces $ M_m(G,T) $ and the associated rings $ Q_k(\Gamma)$ and $\bQ_m(\Gamma)$ in our subsequent paper.

\subsection{Appendices}
The paper contains three appendices included for reader's convenience. 

In Appendix~\ref{AA}, we review the basic homotopy theory of $W$-diagrams (aka $W$-functors) that we use systematically in the present paper. The notion of a $W$-diagram was introduced in topology in \cite{JS01} and \cite{VF99}. In recent years, it has been developed into a general model-categorical framework for equivariant homotopy theory in the work of E. Dotto and K. Moi  \cite{DM16, Dotto16, Dotto17, Moi20}. The starting point is the following simple question.
Let $\,F: \,\Cc \to \Mc \,$ be a diagram of spaces (schemes, stacks, derived stacks, or 
more abstractly, objects in some homotopical category $\Mc$), and let $ X := \holim_{\Cc}(F) $ or $ X := \hocolim_{\Cc}(F) $, whichever exists in $\Mc$.
Assume that the indexing category $ \Cc $ is equipped with an action of a discrete group $W$. What extra structure do we need to put on $F$ to ensure that $X$ is a $W$-space, i.e. $X$ carries a natural $W$-action compatible with the action of $W$ on $\Cc$? If the $W$-action on $\Cc$ is trivial, then $X$ a $W$-space iff $F$ is a diagram of $W$-spaces (i.e., the functor $F$ takes values in the category $ \Mc^W$ of $W$-spaces and $W$-equivariant maps in $\Mc$). In general, a $W$-structure on
$F$ that makes it a $W$-diagram is encoded by a family of natural transformations (indexed by the elements of $W$) satisfying certain compatibility conditions  (see Definition~\ref{Wfun}). In our context, the usefulness of this notion is that, when we define a space $X$ by `gluing' a family of spaces $ \{X_c = F(c)\}_{c \in \Cc}$ via a $W$-diagram $F$, we get a $W$-space, even though the spaces $ X_c $ carry no $W$-action or are equipped only with an action of proper subgroups of $W$. In Appendix~\ref{AA}, we review main definitions and constructions related to $W$-functors, following mostly \cite{DM16}, and prove a few technical results needed for the present paper.  
We choose to work in the general setting of $W$-model categories (not just with topological spaces) as in the main body of the paper we use $W$-diagrams in different homotopical 
contexts (spaces, dg algebras, derived schemes, coaffine stacks, etc.)

In Appendix~\ref{AB}, we give a quick introduction to derived algebraic geometry (DAG), focusing on the To{\"e}n-Lurie theory of coaffine stacks and its application to rational homotopy theory.
Recall that rational homotopy theory is concerned with the study of homotopy types of topological spaces up to rational equivalence (that is, modulo torsion in cohomology and homotopy groups). There is a well-known classical approach to the subject, due to D. Sullivan \cite{Su77}, that uses
commutative cochain dg algebras and polynomial differential forms defined over $\Q$ as algebraic models for rational homotopy types.  Sullivan's theory has been recently re-interpreted and incorporated into the general DAG framework in the work B. Toën \cite{To06} and J. Lurie \cite{DAGVIII, DAGXIII}. In the To{\"e}n-Lurie approach,  the rational homotopy types are modeled geometrically by certain kinds of derived stacks called the coaffine stacks. Roughly speaking, the coaffine stacks are derived stacks defined over $\Q$, whose $\Q$-points represent the simply connected topological spaces of finite rational homotopy type in much the same way as the $\C$-points of usual schemes of finite type over $\C$ represent the classical complex varieties. We found this geometric approach to rational homotopy theory to be particularly well-suited for our constructions: our main results are  most naturally expressed in terms of coaffine stacks.
To the best of our knowledge --- besides the original references --- 
no expository account of the To\"en-Lurie theory is currently available in the literature: for readers' convenience, we decide to include one in Appendix~\ref{AB}. 

Finally, in Appendix~\ref{AC}, we prove a general result --- Cofinality Theorem~\ref{SlThmInf} --- that provides a canonical homotopy decomposition for diagrams of spaces over an arbitrary $\infty$-category $\Cc$ equipped with a functor to a poset. This theorem is one of the main technical tools of the paper: it plays a key role in the proof of Theorem~\ref{MainT} and several other computations in Section~\ref{S6}. To the best of our knowledge, Theorem~\ref{SlThmInf} is a new result that has not appeared in the literature in this form and generality; 
however, its important special case, where $\Cc $ is the nerve of an ordinary (small) category,  was established earlier in the work of J. S{\l}omi{\'n}ska \cite{Sl01} (see Theorem~\ref{SlThm} in Section~\ref{AC4}), and, in fact,  various consequences and specializations of S{\l}omi{\'n}ska's theorem, including Corollary~\ref{EIA}, had been known in homotopy theory  before \ref{SlThm} (see, for example, \cite{DK85, JM89, Sl91, JMO94, Gr02}).
Although Theorem~\ref{SlThmInf} is a natural generalization of S{\l}omi{\'n}ska's Theorem (from 1-categories to $\infty$-categories), our proof --- and, in fact, the very approach to this result based on J. Lurie's Straightening--Unstraightening Principle \cite{HTT} --- are quite different from  \cite{Sl01}. In particular, our arguments do not specialize to S{\l}omi{\'n}ska's  when $\Cc$ is an ordinary category\footnote{For comparison, the interested reader may see the proof of S{\l}omi{\'n}ska's theorem contained in the first {\tt arXiv} version of the present paper.}. 


\subsection*{Acknowledgements}{\footnotesize
We would like to thank O. Chalykh, P. Etingof, and M. Feigin for many stimulating conversations, interesting questions and comments that helped us to shape basic ideas of the current work.
We also thank O. Chalykh for clarifying to us the relation between exponential quasi-invariants
and their $q$-deformations introduced in \cite{Ch02} and S. Kuroki for his interesting comments and references related to our work. Part of this research was conducted while the first author was visiting the Simons Laufer Mathematical Sciences Institute (formerly known as MSRI) as Research Professor for the program ``Noncommutative Algebraic Geometry'' in Spring 2024. He is very grateful to the Institute for its hospitality and financial support under the NSF Grant DMS-1928930. The work of Yu.B. was also partially supported by the Simons Collaboration Grant 712995. }

\section{Varieties of quasi-invariants}\la{S2}
In this section, we give a geometric construction of the algebras of quasi-invariants $Q_k(W) $ and $ \bQ_m(W)$ --- or rather the associated varieties $V_k(W)$ and $ \bV_{\! m}(W)$ --- in terms of universal operations in the category of schemes. The advantage of this approach is that it extends naturally to other objects of geometric nature, such as derived schemes and derived stacks, allowing one to define `varieties of quasi-invariants' in more sophisticated categories. (Two such extensions will be treated in the next section.)  On the algebraic side, the main result of this section is Theorem~\ref{ThFree} which shows that $ \bQ_m(W)$ is a free module over $ \k[V] $
(or equivalently, $ \bV_m(W) $ is a Cohen-Macaulay variety), provided the multiplicities $m_{\alpha} $ are all either even or odd. In the last section, we define a natural action of the classical Demazure (aka nil-Hecke) algebra $\NH(W) $ on $ \bQ_m(W) $ in the case when $m = 2k+1$ is odd (see Theorem~\ref{Qfatodd}) and compare it to the action of the Cherednik
(aka the rational double affine Hecke) algebra $ \bH_k(W) $ on $\bQ_m(W)$ when $ m = 2k$.

\vspace{2ex}

\noindent
{\it Notation.}
Throughout this section, $\k$ is a field of characteristic zero.
$ \Comm_\k $ denotes the category of commutative $\k$-algebras, and  $ \Aff_{\k} $ the corresponding category of affine $\k$-schemes. Recall that the categories $\Comm_{\k}$ and $\Aff_{\k}$ are formally dual (i.e., anti-equivalent) to each other, with inverse equivalences given by $ \Spec: \Comm_{\k}^{\rm op} \leftrightarrows \Aff_{\k}: \k[\,\mbox{--}\,] $. Both these categories are
complete and cocomplete, the above equivalences interchanging limits and colimits.
If $V$ is a $\k$-vector space, we write $ V^* = \Hom_{\k}(V,\k)$ for the $\k$-linear dual of $V$ and $ \k[V] := \Sym_{\k}(V^*)$  for the $\k$-algebra of polynomial functions on $V$. Abusing notation, we denote the affine $\k$-scheme $ \Spec\, \k[V] $ by $ V $, and similarly $ H = \Spec\,\k[H] $ for any linear subspace $ H \subseteq V$.

\subsection{Classical varieties of quasi-invariants} \la{S2.1}  
Let $W$ be a finite Coxeter group, and let $ V $ denote a finite-dimensional
$\k$-vector space that carries a reflection representation of $W$. As in the Introduction, for each reflecting hyperplane $H_{\alpha} $ of $W$, we fix a linear form $ \alpha \in V^* $
 such that $ H_{\alpha} = \Ker(\alpha) $ and write $ \A = \{H_\alpha\} $ (or simply $ \A = \{\alpha\} $) for the collection of all such hyperplanes in $V$. The group $W$ acts naturally on $ \A $ by permuting the $H_\alpha$'s,
and we denote by $ W_{\alpha} $ the (pointwise) stabilizers of $H_\alpha $ in $W$. Further, we introduce the set $ \M(W)$ of $W$-invariant functions $ k: \A \to \Z_+ $, $ \alpha \mapsto k_{\alpha} $ with values in nonnegative integers, which we write as $ k = \{k_\alpha\}_{\alpha \in \A} $ and call the multiplicities of $ W $. Note that  $ \M(W)$ carries a natural partial order induced from $ \Z$, with $ k \leq k' $ defined by the condition that
$ k_{\alpha} \leq k_{\alpha}' $ for all $ \alpha \in \A $.

Now, fix a hyperplane  $ H_{\alpha} \in \A $ and consider two natural maps in $ \Aff_{\k} $: the categorical quotient map $ p_{\alpha}:\, V \to V/\!/W_{\alpha} $ and its restriction\footnote{As we will see in Section~\ref{S2.2},
the map $ p_{\alpha}' $ should be actually written as $ H_{\alpha}/\!/W_{\alpha} \to V/\!/W_{\alpha}$ to indicate that we consider $H_{\alpha}$ together with the action of $W_{\alpha}$. Since this last action is trivial, there is a canonical isomorphism $ H_{\alpha} \xrightarrow{\sim} H_{\alpha}/\!/W_{\alpha}$ in $\Aff_{\k}$ that identifies $ p_{\alpha}': H_{\alpha} \to V/\!/W_{\alpha} $ with the induced map given above.}
$\, p_{\alpha}':\, H_{\alpha} \into V 
\to V/\!/W_{\alpha} $  to $H_{\alpha} \,$, and form the (categorical) fibered product 
\begin{equation}
\la{FP1}
V \times_{V/\!/W_{\alpha}} H_{\alpha} := \lim \{V \to V/\!/W_{\alpha} \leftarrow H_{\alpha}\}
\end{equation}
By definition, the scheme \eqref{FP1} comes with two canonical maps --- the projections onto $V$ and $H_{\alpha}$ ---
that form the cartesian square  in $\Aff_{\k}$:
\begin{diagram}[small]
V \times_{V/\!/W_{\SEpbk\alpha}} H_{\alpha} & \rTo & H_{\alpha}\\
\dTo\ &  & \dTo\\
V    & \rTo & V/\!/W_{\alpha}
\end{diagram}
We can now use these two maps to define the fibered sum (coproduct) of $ V $ and $ H_{\alpha} $
which fits in the {\it cocartesian} square in $ \Aff_{\k} $: 
\begin{diagram}[small]
V \times_{V/\!/W_{\alpha}} H_{\alpha} & \rTo & H_{\alpha}\\
\dTo\ &  & \dTo\\
V    & \rTo^{\quad} & V_{1}{\NWpbk}(W_{\alpha})
\end{diagram}
By analogy with topology, we refer to $ V_1(W_{\alpha})$ as a {\it relative join} of the schemes $V$ and $H_\alpha$ over $V/\!/W_{\alpha}$ --- more precisely, the join of the  morphisms $ p_{\alpha}$ and $ p_{\alpha}'$  --- and denote it by 
\begin{equation}
\la{V1}
V_1(W_\alpha) =  V \ast_{V/\!/W_{\alpha}} H_{\alpha} := \colim \{V \leftarrow V \times_{V/\!/W_{\alpha}} H_{\alpha} \to H_{\alpha}\}
\end{equation}
Geometrically, the scheme \eqref{V1} should be thought of as a `gluing' of $ V $ and $ H_{\alpha} $ along their common (scheme-theoretic) intersection in $V/\!/W_{\alpha}$. Since one of the maps, $\, p_{\alpha}' $, factors through the inclusion $H_\alpha \into  V $, this gluing results effectively in creating a `cusp' along the subscheme $ H_\alpha $ in $V$.

By the universal property of pushouts, the scheme $ V_1(W_\alpha) $ comes equipped with two canonical maps $\, V \to V_1(W_\alpha) \,$ and $\,H_{\alpha} \to V_1(W_\alpha) \,$, and there is a third  map $ p_{\alpha, 1}:\,V_1(W_{\alpha}) \to V/\!/W_{\alpha} $ induced by $ p_{\alpha} $ and $ p_{\alpha}' = p_{\alpha}|_{H_\alpha} $ in \eqref{FP1}. In addition, $ V_1(W_\alpha) $ carries a natural $W_{\alpha}$-action given by the geometric action of $W_\alpha$ on $V$ that makes $\, V \to V_1(W_\alpha) \,$ a $ W_\alpha$-equivariant map. This last map induces an isomorphism
of quotient schemes 
$ V/\!/W_{\alpha} \xrightarrow{\sim} V_1(W_\alpha)/\!/W_{\alpha} $, its inverse being induced by the map $ p_{\alpha,1} $. 
We can now repeat the above construction to define the scheme $ V_2(W_\alpha) := V_1(W_\alpha) \ast_{V/\!/W_{\alpha}} H_{\alpha} $ by taking the  join of 
$ p_{\alpha, 1} $ with its restriction $ p_{\alpha,1}':= p_{\alpha,1}|_{H_\alpha}$.
%
%

Iterating this way, for each hyperplane $ H_{\alpha} \in \A $, we get by induction a sequence of $W_{\alpha}$-schemes
\begin{equation}
\la{Vk}
V_{k+1}(W_\alpha) := V_{k}(W_\alpha) \ast_{V/\!/W_{\alpha}} H_{\alpha}\ , \quad k\ge 0\,,
\end{equation}
together with $W_\alpha$-equivariant maps $ p_{\alpha, k}:  V_{k}(W_\alpha) \to V/\!/W_{\alpha} $ 
and  $ \pi_{\alpha, k}: V_{k}(W_\alpha)  \to  V_{k+1}(W_\alpha) $ that form a tower of fibrations in $\Aff_{\k}$ over $V/\!/W_{\alpha}$:
\begin{equation}
\la{towV}
V = V_0(W_\alpha) \xrightarrow{\pi_{\alpha, 0}} V_1(W_\alpha) \xrightarrow{\pi_{\alpha, 1}} \ldots
\xrightarrow{} V_{k}(W_\alpha)  \xrightarrow{\pi_{\alpha, k}} V_{k+1}(W_\alpha) \to\, \ldots
\end{equation}
As explained above, the maps in \eqref{towV} induce natural isomorphisms of the categorical quotients $ V_k(W_\alpha)/\!/W_{\alpha}$ for all $ k \ge 0$.
%
%

Our goal is to `glue' the $W_{\alpha}$-schemes $V_{k_{\alpha}}(W_{\alpha})$ attached to different hyperplanes $H_{\alpha}$ into a single scheme $V_k(W)$ carrying a natural action of the group $W$. To this end, we construct a small category $ \Sc(W) $ equipped with a (strict) action of $W$, and organize the collection of schemes $ \{V_{k_{\alpha}}(W_{\alpha})\}$ into a {\it $W$-diagram} on $ \Sc(W) $ (see Definition~\ref{Wfun}). The colimit of this $W$-diagram will then carry an action of the entire group $W$, even though each of the schemes $V_{k_{\alpha}}(W_{\alpha})$ 
is equipped only with an action of  a proper subgroup of $  W $
({\it cf.} Lemma~\ref{Wob}). The idea of this gluing construction --- which we use systematically in the present paper in various contexts --- comes from equivariant homotopy theory. We refer the reader to Appendix~\ref{AA}, where recall main definitions and review abstract results needed for our applications.

To define $\Sc(W)$ we simply take the set $  \{W_{\alpha}\}_{\alpha \in \A} $ of all elementary reflection subgroups of $W$ and make it into a {\it family}\footnote{Recall that a {\it family of subgroups} of a group $W$ is defined to be a (non-empty) set $ \Sc $ of subgroups closed under inclusion ($ H \in \Sc $ and $ H' \leq H $ $\Rightarrow$ $ H' \in \Sc$)
and conjugation ($ H \in \Sc $ $\Rightarrow$ $ wHw^{-1} \in \Sc $ for all $w \in W $). This notion is a natural extension of the notion of a group from a topological point of view: in particular, for any  $\Sc $, there is a  well-defined classifying space $ B_{\Sc}(W)$ and the universal $W$-fibration $ E_{\Sc}(W) \to B_{\Sc}(W)$ that specializes to the classical universal $W$-bundle $ EW \to BW $ when $\Sc$ consists of the single (trivial) subgroup of $W $ (see, e.g., \cite{Lu05}).} 
by adjoining the trivial subgroup $W_0:=\{e\}$: 
\begin{equation} 
\la{sw} 
\Sc(W)\,:=\,\{W_0\} \cup \{W_{\alpha}\}_{\alpha \in \A} \ . 
\end{equation}
A family of subgroups is naturally a $W$-poset, on which the partial order is given by inclusion and the $W$-action by conjugation of subgroups. Thus, we can view $\Sc(W)$ as a $W$-category, with objects being the subgroups of $\Sc(W)$, the morphisms being the inclusions, and the $W$-action $\,a(g): \Sc(W) \to \Sc(W)\,$ defined by conjugation: $ W_{\alpha} \mapsto g W_{\alpha} g^{-1}$. The multiplicity functions $ k \in \M(W) $ can be viewed as $W$-invariant maps of posets $\,\Sc(W) \to \Z_{\ge 0} $, $\, W_\alpha \mapsto k_{\alpha} $, which we extend to the trivial subgroup in the trivial way, setting $ k(W_0) := 0 $. 

Now, for a fixed $ k\in \mathcal{M}(W)$, we assign 
\begin{equation} 
\la{cvm} 
\mathcal{V}_k(W)\,:\ \Sc(W) \to \Aff_{\k}\ ,
\quad W_{\alpha} \mapsto V_{k_{\alpha}}(W_{\alpha})\ ,\quad
W_0 \mapsto V_0(W_\alpha)\,.
\end{equation}
Note that $ V_0(W_\alpha) = V $ for all $ \alpha \in \A $, hence we can
extend \eqref{cvm} to a functor on $ \Sc(W)$ by assigning to each inclusion $ W_0 \hookrightarrow W_{\alpha}$  the composition of the first $k_{\alpha}$ maps from the tower \eqref{towV}:
\begin{equation} 
\la{cvm1} 
\pi_{\alpha, k} := \pi_{\alpha,k_{\alpha}-1} \circ \pi_{\alpha,k_{\alpha}-2} \circ 
\ldots \circ \pi_{\alpha, 0}\,:\ V \,\to \, V_{k_{\alpha}}(W_{\alpha})\,.
\end{equation}
Thus, $\mathcal{V}_k(W)$ can be visualized as the following diagram in $\Aff_{\k}$:
%
%
\begin{equation}
\la{GPO}
\begin{tikzcd}[scale cd=0.85, column sep=small, row sep=small,
                 every arrow/.append style={shorten <=2pt, shorten >=2pt}]
	&& {V_{k_{\alpha}}(W_\alpha)} \\
	{V_{k_{\nu}}(W_\nu)} &&&& {V_{k_{\beta}}(W_\beta)} \\
	&& V \\
	{V_{k_{\lambda}}(W_\lambda)} &&&& {V_{k_{\gamma}}(W_\gamma)} \\
	&& {V_{k_{\delta}}(W_\delta)}
	\arrow["{\pi_{\alpha, k}}"', pos=0.4, tail reversed, no head, from=1-3, to=3-3]
	\arrow["{\pi_{\nu, k}}"', pos=0.4, tail reversed, no head, from=2-1, to=3-3]
	\arrow["{\pi_{\beta, k}}"', pos=0.4, tail reversed, no head, from=2-5, to=3-3]
	\arrow["{\pi_{\lambda, k}}"', pos=0.4, tail reversed, no head, from=4-1, to=3-3]
	\arrow["{\pi_{\gamma, k}}"', pos=0.4, tail reversed, no head, from=4-5, to=3-3]
	\arrow["{\pi_{\delta, k}}"', pos=0.4, tail reversed, no head, from=5-3, to=3-3]
\end{tikzcd}
\end{equation}
We claim
\begin{prop} \la{cvmWfun}
The diagram \eqref{GPO} defined by \eqref{cvm} and \eqref{cvm1} is a $W$-diagram on $\Sc(W)$.
\end{prop}
We will use the following formal lemma, the proof of which 
is a (trivial) exercise to the reader.
\begin{lemma} 
\la{jointhm}
Any commutative diagram in $\Aff_{\k}$ 
\begin{equation} \la{cdsch}
\begin{diagram}[small]
X & \rTo^p & S &\lTo^{q} & Y\\
\dTo^f & & \dTo^g & & \dTo_h\\
X' & \rTo^{p'} & S' &\lTo^{q'} & Y'\\
\end{diagram}
\end{equation}
induces a commutative diagram on relative joins
\begin{equation} \la{cdrelj}
\begin{diagram}[small]
X & \rTo & X \ast_S Y &\rTo & S\\
\dTo^f & & \dTo^{f_1} & & \dTo_g\\
X' & \rTo & X' \ast_{S'} Y' &\rTo & S'\\
\end{diagram}
\end{equation}
Moreover, if each of the two inner squares in \eqref{cdsch} is cartesian in $\Aff_{\k}$, then so are the two inner squares in \eqref{cdrelj}.
\end{lemma}
\begin{proof}[Proof of Proposition~\ref{cvmWfun}]
For any element $g \in W$, we have $g \, W_{\alpha}\, g^{-1}=W_{g(\alpha)}$ in $W$. Moreover, the linear automorphism $g: V \to V$ restricts to $g|_{H_{\alpha}}:H_{\alpha} \to H_{g(\alpha)}$ and induces the map $\bar{g}:V/\!/W_{\alpha} \to V/\!/W_{g(\alpha)}$. Hence, we have the following commutative diagram in $\Aff_{\k}$:
\begin{equation} \la{cdschh}
\begin{diagram}[small]
V & \rTo^{p_\alpha} & V/\!/W_{\alpha} & \lTo^{p_\alpha'} & H_{\alpha}\\
\dTo^g & & \dTo^{\bar{g}} & & \dTo_{g|_{H_\alpha}}\\
V & \rTo^{p_{g(\alpha)}} & V/\!/W_{g(\alpha)} & \lTo^{p_{g(\alpha)}'} & H_{g(\alpha)}\\
\end{diagram}
\end{equation}
Applying Lemma \ref{jointhm} to \eqref{cdschh} yields the commutative diagram
\begin{equation} \la{cdm1a}
\begin{diagram}[small]
V & \rTo^{\pi_{\alpha,0}} & V_1(W_\alpha) & \rTo & V/\!/W_\alpha\\
\dTo^g & & \dTo^{g_1} & & \dTo_{\bar{g}}\\
V & \rTo^{\pi_{g(\alpha),0}} & V_1(W_{g(\alpha)}) & \rTo & V/\!/W_{g(\alpha)}\\
\end{diagram}
\end{equation}
Now, replacing the first commutative square in diagram \eqref{cdsch} by the second one in diagram \eqref{cdm1a}, we can apply Lemma \ref{jointhm} again to get the commutative diagram
\begin{equation} \la{cdm2a}
\begin{diagram}[small]
V_1(W_\alpha) & \rTo^{\pi_{\alpha,1}} & V_2(W_\alpha) & \rTo & V/\!/W_\alpha\\
\dTo^{g_1} & & \dTo^{g_2} & & \dTo_{\bar{g}}\\
V_1(W_{g(\alpha)}) & \rTo^{\pi_{g(\alpha),1}} & V_2(W_{g(\alpha)}) & \rTo & V/\!/W_{g(\alpha)}\\
\end{diagram}
\end{equation}
Continuing this way, we construct a sequence of maps of schemes $g_{k}:V_{k}(W_{\alpha}) \to V_{k}(W_{g(\alpha)})$, $k \geqslant 0$ which actually give a morphism of towers \eqref{towV} from $\alpha$ to $g(\alpha)$. 

Now, to construct a $W$-action of the functor $\mathcal{V}_m(W): \Sc(W) \to \Aff_{\k}$, we need to define natural transformations $\theta_g: \mathcal{V}_k(W) \to \mathcal{V}_k(W) \circ a(g) $, one for each $g \in W$. We set $(\theta_g)_{W_0}:=g: V \to V$ and $(\theta_g)_{W_{\alpha}}:=g_{k_{\alpha}}:V_{k_{\alpha}}(W_\alpha) \to V_{k_{\alpha}}(W_{g(\alpha)}) = V_{k_{g(\alpha)}}(W_{g(\alpha)})$ where we use the fact that the multiplicities $ k=\{k_{\alpha}\}_{\alpha \in \A}$ are $W$-invariant. The second property of Definition \ref{Wfun} amounts to checking that the following diagram 
$$ 
\begin{diagram}[small]
V_{k_{\alpha}}(W_\alpha) & \rTo^{h_{k_{\alpha}}} & V_{k_{\alpha}}(W_{h(\alpha)})\\
  & \rdTo_{\!\!(gh)_{k_{\alpha}}\ } & \dTo_{g_{k_{\alpha}}}\\
  & & V_{k_{\alpha}}\!(W_{gh(\alpha)})\\
\end{diagram}
$$
commutes for all $g,h \in W$. This follows formally from the above construction of the maps $g_i$ and the fact that the linear maps $g:V \to V$ come from the action of $W$ on $V$.
\end{proof}

\begin{remark}
Lemma \ref{jointhm} admits a useful homotopy-theoretic generalization (known as the Join Theorem) in the context of Quillen model categories (see \cite{Do93}). In that context, commutative (resp., cartesian) squares are replaced with homotopy commutative (resp., homotopy cartesian) ones, and Lemma \ref{jointhm} holds under the additional assumption that the base change map $g:S \to S'$ satisfies the so-called Cube Axiom. The Cube Axiom holds in the usual (model) category of topological spaces, hence  Lemma \ref{jointhm} can be used to define maps between topological joins.
\end{remark}

Now, having the functor \eqref{cvm} in hand, we define
\begin{equation} 
\la{defvm} 
V_k(W)\,:=\, {\rm colim}_{\Sc(W)}\,[\mathcal{V}_k(W)] \ .
\end{equation}
By Lemma \ref{Wob} and Proposition \ref{cvmWfun}, $V_k(W)$ is a $W$-scheme, even though the functor $\mathcal{V}_k(W)$ does not take values in the category of $W$-schemes. The next theorem is the main result of this section.
\begin{theorem} \la{ThVQI}
For all $k \in \mathcal{M}(W)$, there are natural isomorphisms of $W$-schemes
$$ 
V_k(W)\,\cong\, \Spec\,Q_k(W)\ .
$$
Thus, $V_k(W)$ are the varieties of $W$-quasi-invariants of multiplicity $k$ $($as defined in \cite{BEG03}$)$.
\end{theorem}
\begin{proof}
By definition, the fibered product \eqref{FP1} in $\Aff_{\k} $ is given by
\begin{equation}
\la{prodf}
V \times_{V/\!/W_{\alpha}} H_{\alpha} \,=\,\Spec \bigl(\k[V] \otimes_{\k[V]^{W_\alpha}} \k[H_\alpha]\bigr)
\end{equation}
%
%
As a closed subscheme of $V$, the hyperplane $ H_{\alpha} $ is described by an isomorphism
$ \k[H_{\alpha}] \cong \k[V]/(\alpha) $. On the other hand, the canonical algebra maps 
induce an isomorphism
\begin{equation}\la{iso1}
\begin{diagram}[small]
\k[V]^{W_\alpha} & \rInto & \k[V]        & \rOnto & \k[V]/(\alpha^2)\\
\dOnto            &        & \ruInto(2,4) &       & \dOnto \\
\k[V]^{W_{\alpha}}\!/(\alpha^2)  &       &  \rTo(2,4)^{\cong} &       & \k[V]/(\alpha)
\end{diagram}
\end{equation}
that describes $ H_\alpha $ as a closed subscheme in $ V/\!/W_{\alpha}$.
We can use this last isomorphism to identify
\begin{equation}
\la{tensp}
\k[V] \otimes_{\k[V]^{W_\alpha}} \k[H_{\alpha}] \,\cong\,
\k[V] \otimes_{\k[V]^{W_\alpha}}\, \k[V]^{W_\alpha}\!/(\alpha^2)\,\cong\, \k[V]/(\alpha^2)\,,
\end{equation}
describing the fibered product   $ V \times_{V/\!/W_{\alpha}} H_{\alpha}$ (see \eqref{prodf}) as a closed subscheme in $V$.
With this description, the canonical map $ V \times_{V/\!/W_{\alpha}} H_{\alpha} \to V $ corresponds to the algebra projection $ \k[V] \onto \k[V]/(\alpha^2) $.
Hence the scheme $V_1(W_\alpha)$ defined by \eqref{V1} can be described by
\begin{eqnarray*}\la{coprodf}
V_1(W_\alpha) 
&\cong& \Spec \left(\, \lim \{\,\k[V] \to \k[V] \otimes_{\k[V]^{W_\alpha}} \k[H_\alpha]  \leftarrow \k[H_{\alpha}]\}\,\right) \\*[1ex]
&\cong& \Spec \left(\k[V] \times_{\k[V]/(\alpha^2)} \k[H_\alpha]\right)\nonumber
\end{eqnarray*}
where in the first step we use the isomorphism \eqref{prodf} and in the 
last step the isomorphism \eqref{tensp}. Note that under  \eqref{tensp},
the natural algebra map $ \k[H_\alpha] \to \k[V] \otimes_{\k[V]^{W_\alpha}} \k[H_\alpha] $ corresponds to the `unnatural' map
\begin{equation} 
\la{iota}
\iota:\, \k[V]/(\alpha) \, \xrightarrow{\sim} \, \k[V]^{W_{\alpha}}\!/(\alpha^2) \,\into\, \k[V]/(\alpha^2) \end{equation}
where the first arrow is the inverse of the induced isomorphism in  \eqref{iso1}. The map \eqref{iota} is explicitly given by
the formula
\begin{equation} \la{iota1}
 p \ \mbox{mod}\,(\alpha) \ \mapsto \  p_+\ \mbox{mod}\,(\alpha^2)
\end{equation}
where $\,p_+ := \frac{1}{2}(p + s_\alpha(p)) \,$ denotes its symmetrization (averaging) of $ p \in \k[V] $  over the reflection subgroup $W_\alpha$. Using the fact that $\, s_{\alpha}(p) \equiv p\,\mbox{mod}\,(\alpha)\,$ for all $ p \in \k[V]$, one can easily check that formula \eqref{iota1} indeed gives a well-defined algebra homomorphism, which is a section of the canonical projection $ \k[V]/(\alpha^2) \onto \k[V]/(\alpha)$. It follows that
\begin{eqnarray}
\la{fibiso}
\k[V] \times_{\k[V]/(\alpha^2)} \k[H_\alpha] 
&\cong&
\{(p, \bar{q}) \in \k[V] \times \k[V]/(\alpha)\ :\
p \equiv q_+\,\mbox{mod}\,(\alpha^2)\} \nonumber\\
&=&
\{(p, \bar{p}) \in \k[V] \times \k[V]/(\alpha)\ :\
p \equiv p_+\,\mbox{mod}\,(\alpha^2)\} \nonumber\\
&=&
\{(p, \bar{p}) \in \k[V] \times \k[V]/(\alpha)\ :\
s_\alpha(p) \equiv p\,\mbox{mod}\,(\alpha^2)\}\nonumber
\end{eqnarray}
The last identification shows that the natural algebra map $ \k[V] \times_{\k[V]/(\alpha^2)} \k[H_\alpha] \to \k[V] $ is injective, and its image is precisely the subring of $W_\alpha$-invariants $ Q_1(W_\alpha) $. Thus 
$\,V_1(W_\alpha) \cong \Spec\,Q_1(W_\alpha)\,$.

Similarly, arguing by induction,  we  identify ({\it cf.} \eqref{prodf} and \eqref{tensp})
$$
V_{k}(W_\alpha) \times_{V/\!/W_{\alpha}} H_{\alpha} \, \cong\,  
\Spec \left(Q_{k}(W_\alpha) \times_{\k[V]^{W_{\alpha}}} \k[H_\alpha]\right)
\,\cong\, \Spec\,[Q_{k}(W_\alpha)/(\alpha^2)]
$$
Hence
$$
V_{k+1}(W_\alpha) = V_{k}(W_\alpha) \ast_{V/\!/W_{\alpha}} H_{\alpha} \cong 
\Spec \left(Q_{k}(W_\alpha) \times_{Q_{k}(W_\alpha)/(\alpha^2)} \k[H_\alpha]\right)
$$
The same calculation as above then shows that the natural algebra map 
$$ 
Q_{k}(W_\alpha) \times_{Q_{k}(W_\alpha)/(\alpha^2)} \k[H_\alpha] \into Q_{k}(W_\alpha) 
$$ 
is injective, and its image is given by
$$
Q_{k}(W_\alpha) \times_{Q_{k}(W_\alpha)/(\alpha^2)} \k[H_\alpha] \cong Q_{k}(W_\alpha) \times_{Q_{k}(W_\alpha)/(\alpha^2)} \k[V]^{W_\alpha}\!/(\alpha^2)\,\cong\,
Q_{k+1}(W_\alpha)
$$
for all $ k \ge 0 $.
Finally, with above isomorphisms we can formally compute \eqref{defvm}:
\begin{eqnarray}\la{VmW}
V_k(W) &=& {\rm colim}_{\Sc(W)}\,\{\Spec\,[Q_{k_{\alpha}}(W_\alpha)]\}\nonumber\\
&\cong& \Spec\,\{{\rm lim}_{\Sc(W)^{\rm op}}\,[Q_{k_{\alpha}}(W_\alpha)]\}\nonumber\\
&\cong& 
\Spec \,\bigl\{\bigcap_{Q_0(W)} Q_{k_{\alpha}}(W_\alpha)\bigr\}\\
&=& \Spec\, Q_k(W)\,.\nonumber
\end{eqnarray}
This finishes the proof of the theorem.
\end{proof}
\begin{remark} We mention three possible generalizations of the definition \eqref{defvm}.
\vspace*{1ex}

$(1)$
To build the `global' scheme of quasi-invariants, $V_k(W)$,
one can try to glue the {\it diagrams} of `local' schemes \eqref{towV} instead of their individual members. To do this 
one should replace the indexing category $ \Sc(W) $ with a larger  category $ \Sc_k(W) $, depending on the given multiplicity $k \in \M(W) $. The latter can be formally defined via the Grothendieck construction (see Definition~\ref{Grothcon}): $\, \Sc_k(W) := \Sc(W) \int \,[k]\,$, where the functor $ [k]: \Sc(W) \to \Cat $ assigns to the subgroup $W_\alpha$ the poset category $[k_{\alpha}] = \{0<1<2< \ldots < k_{\alpha}\}\,$. The objects of $ \Sc_k(W) $ are  the pairs $ 
(W_0,\, 0) $ and $ \{(W_{\alpha},\, i)\}_{i \in [k_{\alpha}]} $ for all $ \alpha \in \A$, and morphisms are the obvious maps: $(W_0,\, 0) \to (W_\alpha,\, i) $ for all $i \in [k_{\alpha}]$, and $\,(W_\alpha,\, i) \to (W_\alpha,\, i') \,$ for $ i \leq i'$. Then, for a fixed $k \in \M(W)$, there is a natural  $\Sc_k(W)$-diagram:
\begin{equation*} 
\la{cvmk} 
\tilde{\mathcal{V}}_k(W)\,:\ \Sc_k(W) \to \Aff_{\k}\ ,\quad (W_0,0) \mapsto V\ ,\quad
(W_{\alpha}, \,i) \mapsto V_{i}(W_{\alpha})\, ,\quad 0 \leq i \leq k_{\alpha}\,,
\end{equation*}
refining the $\Sc(W) $-diagram \eqref{cvm}.
We set $\, \tilde{V}_k(W) := \colim_{\Sc_k(W)}\,[\tilde{\mathcal{V}}_k(W)] \,$.
Although it looks more natural, this definition is actually equivalent to \eqref{defvm}:  since each poset $[k_{\alpha}]$ has the largest element (terminal object), we have canonical isomorphisms
$$
\tilde{V}_k(W)\,\cong\,\colim_{W_{\alpha} \in \Sc(W)}\, 
\colim_{i \in [k_{\alpha}]}\,[V_i(W_{\alpha})]\,\cong\, 
\colim_{W_{\alpha} \in \Sc(W)}\,[V_{k_{\alpha}}(W_{\alpha})]\, = \, V_k(W)
$$
for all $ k \in \M(W) $. By Thomason's Theorem \cite{ChS02}, the above argument extends to homotopical contexts (e.g., to diagrams of derived schemes or topological spaces), when the usual colimits are replaced by the homotopy ones.

$(2)$ A more interesting generalization arises if --- in place of  $ \Sc(W) $ ---  we consider the poset $ \mathcal{P}(W) $ of all {\it parabolic} subgroups of $W$. Then, we should replace the colimit in \eqref{defvm} with the left Kan extension functor $\,i_!: \,\Aff^{\Sc(W)}_{\k} \to \Aff^{\mathcal{P}(W)}_{\k} \,$ taken along the inclusion $i: \Sc(W) \into {\mathcal P}(W) $. When applied to the $\Sc(W)$-diagram \eqref{cvm} this functor yields the ${\mathcal P}(W)$-diagram  $\, i_! \mathcal{V}_k:\,{\mathcal P}(W) \to \Aff_{\k} $ that assigns to each parabolic subgroup $ W' $ of $W$ the scheme $ V_k(W') $ of $k$-quasi-invariants of $W'$, and to each inclusion of parabolic subgroups $ W' \subseteq W'' $ the morphism of schemes $  V_k(W') \to V_k(W'') $ corresponding to the natural inclusion of commutative algebras $ Q_k(W'') \into Q_k(W')$. Note that the definition \eqref{defvm} can be recovered from this diagram by
$\, i_! \mathcal{V}_k(W) = V_k(W) $, while $ i_! \mathcal{V}_k(W_0) = V $. This shows that our geometric
construction behaves well with respect to parabolic induction.

$(3)$ One can also apply the above construction in the category $ \Sch_{\k} $ of all (i.e., not necessarily affine) $\k$-schemes. The category $ \Sch_{\k} $ admits all finite limits but it is not well behaved with respect to colimits: in particular, arbitrary pushouts do not exist in $ \Sch_{\k} $ (see \cite[Tag 0ECH]{StackPr}). The proof of Theorem~\ref{ThVQI} shows that the schemes $V_{k}(W_\alpha) $ are defined by pushouts of closed embeddings of affine schemes, and hence they exist in $ \Sch_{\k} $ (see, e.g., \cite{Sch05}). On the other hand,  whether the diagrams \eqref{GPO} admit colimits in $\Sch_{\k} $ and whether such colimits agree with those in $ \Aff_{\k} $ seem to be subtle geometric questions that deserve a separate investigation.

\end{remark}

\subsection{Varieties of quasi-covariants} \la{S2.2} In classical invariant theory, given an action of a finite group $W$ on a vector space $V$ and a linear representation $ \tau $ of $W$, a {\it $W$-covariant of type $\tau $} on $V$ is defined to be a $W$-equivariant polynomial map $f: V \to \tau $ (see \cite{KP96}, Section 1.4). The space of  such maps is naturally identified with $ (\k[V] \otimes \tau)^W$. When $ V$ is a geometric representation of a finite (complex) reflection group $W$, there is a generalization of this notion proposed in \cite{BC11} under the name `$\tau$-valued quasi-invariants'.  Of special interest is the case when $ \tau = \k[W]$, the regular representation of $W$. In that case, the space of $W$-quasi-covariants --- or $\k[W]$-valued quasi-invariants in terminology of \cite{BC11} --- is a linear subspace $\bQ_{k}^{\rm BC}(W) $ of $ \k[V] \otimes \k[W] $,  closed under the (multiplication) action of $ \k[V] $ and the (diagonal) action of $W$, with key property that $ \boldsymbol{e} \bQ_{k}^{\rm BC}(W) = \boldsymbol{e} (Q_k(W) \otimes 1) $,
where $ \boldsymbol{e} \in \k[W] $ is the symmetrizing idempotent of $W$. 
In this section, we construct the modules of $W$-quasi-covariants, $\bQ_m(W)$, 
in a slightly greater generality than in \cite{BC11} (so that our $ \bQ_m(W) $ specialize to $ \bQ_k^{\rm BC}(W) $ defined in \cite{BC11} when $ m = 2 k $, see Corollary~\ref{BCcor}). We show that, for all $ m \in \M(W)$, these modules carry natural commutative algebra structures compatible with the action of $W$, and there are canonical isomorphisms of algebras (see Lemma~\ref{bVV})
$$
\bQ_m(W)^W \cong Q_{[\frac{m}{2}]}(W) 
$$ 
We define the algebras $\bQ_m(W)$ geometrically by constructing the corresponding affine schemes as iterated relative joins. The key idea is to apply the construction of Section~\ref{S2.1} not  to the $W$-scheme $V$ itself,  but to its canonical {\it $W$-free resolution} $\bV$. 

Recall that a $W$-scheme $ \bX $ is called {\it free} if $W$ acts freely on $ \bX  $ and the quotient map $ \bX \to \bX\!/W $ is a locally trivial fibration (in \'etale topology) with fiber $W$. Given an arbitrary $W$-scheme $X$, its free resolution is then defined to be a $W$-equivariant map
$ p: \bX \to X $, with $ \bX $ being a free $W$-scheme. For any $X$, there is a canonical (trivial) 
free resolution given by the projection $ X \times W \to X $, where $\, X \times W \,$ is equipped with the diagonal $W$-action.
 
Thus, to `resolve' the  $W$-action on $V$, we simply replace $V$ and $ H_\alpha $ with their trivial $W$-free resolutions: 
$ \bV := V \times W \onto V $ and $ \bH_{\alpha} := H_{\alpha} \times W \onto H_{\alpha}
\,$.
As schemes, $ \bV $ and $ \bH_{\alpha} $ are isomorphic to the disjoint unions of copies of $ V $  and $ H_{\alpha} $ indexed by the elements of $W$: hence, we have canonical isomorphisms
$$ 
\k[\bV] \,=\,  \prod_{w \in W }\k[V] \ ,\quad 
\k[\bH_{\alpha}] \,=\, \prod_{w \in W }\k[H_{\alpha}] \,, 
$$
with above projections corresponding to the diagonal inclusions.


Now,  as in Section~\ref{S2.1}, we fix $ \alpha \in \A $ and consider two natural maps $ \bp_{\alpha}:  \bV \to \bV\!/W_{\alpha} $ and $ \bp_{\alpha}':\,\bH_{\alpha}/W_{\alpha} \to \bV\!/W_{\alpha} $. Since $ \bV $ and $ \bH_{\alpha} $ are free $W$-schemes, the above (geometric) quotients coincide with the categorical ones: that is, $\,\bV\!/W_{\alpha} \cong \bV/\!/W_{\alpha} $ and $\bH_{\alpha}/W_{\alpha} \cong \bH_{\alpha}/\!/W_{\alpha}$.
Explicitly, the maps $\bp_{\alpha}$ and $ \bp_{\alpha}' $ can be written as
\begin{equation}\la{bps}
 V \times W \xrightarrow{\bp_{\alpha}} V \times_{W_\alpha} W \xleftarrow{\bp_{\alpha}'} H_{\alpha} \times_{W_{\alpha}} W
\end{equation}
Taking the relative join of  these maps in $\Aff_{\k}\,$, we define
\begin{equation} \la{bV1}
\bV_1(W_\alpha) \,:= \,\bV \ast_{\bV\!/W_{\alpha}}(\bH_{\alpha}/W_{\alpha}) \,=\, 
(V \times W) \ast_{(V \times_{W_\alpha} W)} (H_{\alpha} \times_{W_{\alpha}} W)
\end{equation}
By the universal property of joins, this scheme comes with three canonical maps $\, \bV \to \bV_1(W_{\alpha}) $, $\,\bH_{\alpha}/W_{\alpha} \to \bV_1(W_{\alpha}) $ and $\,\bV_1(W_{\alpha}) \to \bV\!/W_{\alpha} $. We take $\,\bp_{\alpha,1}\,$ to be the last map and $\, \bp_{\alpha,1}' \, $ to be its (pre)composition with the second. Forming the join of these two maps, we define 
 $ \bV_2(W_{\alpha}) := \bV_1(W_{\alpha}) \ast_{\bV\!/W_{\alpha}} \bH_{\alpha}/W_{\alpha} $, and then iterating further as in Section \ref{S2.1}, we set by induction
\begin{equation}
\la{bVm+1}
\bV_{m+1}(W_\alpha) \,:= 
 \bV_{m}(W_{\alpha}) \ast_{\bV\!/W_{\alpha}}(\bH_{\alpha}/W_{\alpha})\ ,\quad m\ge 0 \,.
\end{equation}
The next proposition gives an explicit description of (the coordinate rings of) these affine schemes.
\begin{prop} \la{LQwa} 
For each $ \alpha \in \A $ and for all $m \ge 0$, there are natural isomorphisms
$$
\bV_{m}(W_\alpha) \cong \Spec\,\bQ_{m}(W_\alpha)\,,
$$
where $\bQ_{m}(W_\alpha)$ are the subalgebras of 
$\, \k[V \times W] \,=\, \prod_{w \in W }\k[V] \,$ defined by 
\begin{equation} 
\la{ebqkwa} 
\bQ_{m}(W_\alpha)\,=\,\{ (p_w)_{w \in W} \in \k[V \times W]\,:\, p_w \equiv p_{s_\alpha w}\,{\rm mod}\,\langle\alpha\rangle^m \}\ .
\end{equation}
The natural maps $\pi_{\alpha,m}:\bV \to \bV_{m}(W_\alpha)$ correspond to the inclusions $\bQ_{m}(W_\alpha) \hookrightarrow \bQ_0=\k[V \times W]$.
\end{prop}
\begin{proof} 
The map $\,[v,w]_{W_{\alpha}} \mapsto (w^{-1} v,\, W_{\alpha}w)\,$  induces an isomorphism of schemes  
$$
V \times_{W_{\alpha}} W \,\xrightarrow{\sim}\, V \times (W_\alpha\backslash W)\,,
$$
where $W_\alpha\backslash W$ is the set of right cosets of $W_\alpha $ in $W$ and $  V \times (W_\alpha\backslash W) $ stands for the disjoint union of copies of $V$ indexed by  $W_\alpha\backslash W$. Similarly, we have an isomorphism $\, H_{\alpha} \times_{W_{\alpha}} W \xrightarrow{\sim} H_{\alpha} \times (W_\alpha\backslash W)  \,$ induced by the map $\, [v,w]_{W_{\alpha}} \mapsto (v, W_{\alpha} w) $ (which makes sense since 
$W_{\alpha}$ acts trivially on $H_{\alpha}$). With these isomorphisms, we can identify the maps \eqref{bps} as
\begin{equation} 
\la{bps1}
 V \times W \xrightarrow{\bp_{\alpha}} V \times (W_\alpha\backslash W)  \xleftarrow{\bp_{\alpha}'} H_{\alpha} \times (W_\alpha\backslash W)  
\end{equation}
where $\, \bp_{\alpha}(v,w) = (w^{-1} v, \,W_{\alpha}w) \,$ and 
$\, \bp_{\alpha}'(v, W_{\alpha}w) = (w^{-1} v,\, W_{\alpha}w) $.
It follows that $\,\bV \times_{\bV\!/W_{\alpha}}(\bH_{\alpha}/W_{\alpha})
\,\cong\,  H_{\alpha}  \times W \,$, with canonical maps
$\,\bV \leftarrow \bV \times_{\bV\!/W_{\alpha}} (\bH_{\alpha}/W_{\alpha}) \to \bH_{\alpha}/W_{\alpha}\,$ corresponding to the natural inclusion and  projection, respectively:
\begin{equation} 
\la{dig1}
V \times W \hookleftarrow  H_{\alpha} \times W  \onto H_{\alpha} \times (W_\alpha\backslash W) 
\end{equation}
The scheme $ \bV_1(W_\alpha) $ is thus the colimit of the diagram \eqref{dig1}.
To describe its coordinate ring we identify the vector spaces
\begin{equation}
\la{idenkwb}
\k[V \times W] \,=\, \prod_{w \in W} \k[V] \,\cong\, \k[V]\otimes \k[W]\ ,\quad 
(p_w)_{w \in W}\, \leftrightarrow \, \sum_{w \in W} p_w \otimes w\,,
\end{equation}
and think of their elements as functions $\,W \to \k[V]$, $\,w \mapsto p_w\,$. The algebra structure on $\k[V]\otimes \k[W]$ is then defined by the formula 
\begin{equation}
\la{prodVW}
 (p_w \otimes w) \cdot (q_{w'} \otimes w') = \delta_{w,w'} \,(p_w q_w \otimes w) 
\end{equation}
which corresponds to the usual (pointwise) multiplication of functions. With this identification, 
\begin{equation} \la{e18} 
\k[\bV_1(W_\alpha)] \,\cong\,  (\k[V]\otimes \k[W]) \times_{(\k[H_{\alpha}] \otimes \k[W])} (\k[H_\alpha] \otimes \k[W_{\alpha} \backslash W])\ .\end{equation}
To describe this fiber product  explicitly we identify $\k[H_\alpha] \cong \k[V]/\langle \alpha\rangle$ and write the elements of $\,\k[V] \otimes \k[W]\,$ and $\,  \k[H_\alpha] \otimes \k[W_\alpha \backslash W]\, $ respectively as 
\begin{equation*}
p = \sum_{w \in W} p_w \otimes w =\!\! \sum_{[w] \in W_\alpha\! \backslash W}(p_w \otimes w +  p_{s_\alpha w} \otimes s_{\alpha}w) \ ,\quad
\overline{q} = \sum_{[w] \in W_\alpha\! \backslash W}  \overline{q}_{[w]} \otimes [w]  \,,
\end{equation*}
where $\overline{q}_{[w]} \in \k[V]/\langle \alpha\rangle\,$. Then, a pair $\,(p,\overline{q})\,$ represents an element in $\k[\bV_1(W_\alpha)]$ if and only if
\begin{equation} \la{cond} 
\bar{p}_w=\bar{p}_{s_\alpha w}=q_{[w]} \,\in\, \k[V]/\langle \alpha\rangle\ ,\quad 
\forall w \in W\,,
\end{equation}
where $\bar{p}_w$ and $\bar{p}_{s_\alpha w}$ denote the images of $p_w,\,p_{s_\alpha w} \in \k[V]$ under the  projection $\k[V] \twoheadrightarrow \k[V]/\langle \alpha\rangle$. Equations \eqref{cond} show that $\overline{q}$ is uniquely determined by $p \in \k[V] \otimes \k[W]$ whenever $ (p,\overline{q}) \in \k[\bV_1(W_\alpha)]$; hence, the canonical map $\k[\bV_1(W_\alpha)] \to \k[V] \otimes \k[W]$ is injective,  and its image is 
\begin{equation} \la{Q1W} 
\bQ_1(W_\alpha) = \left\{\,\sum\, p_w \otimes w   \in \k[V] \otimes \k[W]\,:\, p_w \equiv p_{s_\alpha w}  \ {\rm mod}\,\langle\alpha \rangle\,\right\}
\end{equation}
This proves the claim for $m=1$. 

Now, the same argument, starting with $ \bQ_m(W_{\alpha})$ instead of $ \k[V \times W]$, shows that
$$ 
\bQ_{m+1}(W_\alpha) \,\cong\, \bQ_{m}(W_\alpha) \times_{\bQ_{m}(W_\alpha)/\langle \alpha\rangle} (\k[H_\alpha] \otimes \k[W_\alpha \backslash W]).
$$
To verify  \eqref{ebqkwa} we need to describe the above fiber product explicitly.
Observe that the natural map 
$\, 
\k[V]\bigl/\langle \alpha\rangle \otimes \k[W_\alpha\backslash W]\,\to\,
\bQ_m(W_\alpha)\bigl/\langle \alpha\rangle
\,$
is given by
%
%
$$
\sum_{[w] \in W_\alpha\! \backslash W}  \overline{q}_{[w]} \otimes [w]  \ \mapsto \ 
\sum_{w\in W} q_{[w]} \otimes w  \ (\mod\,\langle\alpha \rangle) \ ,
$$
where $\,q_{[w]}\,$ is a(ny) representative of $\,\overline{q}_{[w]} 
\in \k[V]/\langle \alpha\rangle $ in $ \k[V]$. Hence 
$$
(p, \,\overline{q}) \,= \,\biggl(\sum_{w\in W} p_w \,\otimes\, w,\
\sum_{[w]\in W_\alpha\backslash W}  \overline{q}_{[w]} \,\otimes\, [w] \biggr)
\;\in\; \bQ_m(W_\alpha)\,\times\,\bigl((\k[V]/\langle \alpha\rangle) \otimes \k[W_\alpha\backslash W]\bigr)
$$ 
represents an element in $\bQ_{m+1}(W_{\alpha})$ if there exists $\,\displaystyle{r=\sum_{w\in W}  r_w \otimes w \in \bQ_m(W_{\alpha})} \,$ such that
\begin{equation} \la{cond-m}
p_w - q_{[w]} = \alpha \cdot r_w, \quad p_{s_\alpha w}-q_{[w]} = \alpha \cdot r_{s_{\alpha} w} \ ,\quad 
\forall\, w \in W\,.
\end{equation}
Again, this implies that $\overline{q}$ is uniquely determined by $p\in \bQ_m(W_{\alpha})$: i.e., the canonical map $\bQ_{m+1}(W_{\alpha})\to \bQ_m(W_{\alpha})$ is injective. The composite map  $\bQ_{m+1}(W_{\alpha}) \to \bQ_m(W_{\alpha}) \to \k[V]\otimes\k[W]$ is then injective as well, and the image of $\bQ_{m+1}(W_{\alpha}) $ 
under this map is  
$$
\bQ_{m+1}(W_{\alpha}) = \left\{\sum\, w \otimes p_w   \in \k[V]\otimes\k[W]\,:\, p_w - p_{s_\alpha w} = \alpha \cdot (r_w - r_{s_\alpha w}) \text{ for some } r\in \bQ_m(W_{\alpha}) \right\}$$
By induction, we have $r_w \equiv r_{s_\alpha w} \ \mod \,\langle\alpha\rangle^m $. Hence, the condition on the components $p_w$ in $\bQ_{m+1}(W_{\alpha})$ can be rewritten as
$p_w \equiv p_{s_\alpha w} \ \mod \,\langle\alpha\rangle^{m+1}$, which proves formula \eqref{ebqkwa}.
\end{proof}
Next, observe that the schemes $ \bV_{m}(W_\alpha) $ are defined as iterated joins of $W_\alpha$-equivariant maps
\begin{equation}
\la{bvmWa}
\bV_{m}(W_\alpha) \,= 
\bV\ast_{\bV\!/W_{\alpha}}(\bH_{\alpha}/W_{\alpha}) \ast_{\bV\!/W_{\alpha}}\,
\stackrel{m}{\ldots}
\, \ast_{\bV\!/W_{\alpha}} (\bH_{\alpha}/W_{\alpha})
\end{equation}
and hence these are $W_\alpha$-schemes, with $W_{\alpha}$-action induced from $ \bV $. As in Section~\ref{S2.1}, we thus obtain a tower of $W_{\alpha}$-schemes
\begin{equation}
\la{towbV}
\bV = \bV_0(W_\alpha) \xrightarrow{\bpi_{\alpha, 0}} \bV_1(W_\alpha) \xrightarrow{\bpi_{\alpha, 1}} \ldots
\xrightarrow{} \bV_{m}(W_\alpha)  \xrightarrow{\bpi_{\alpha, m}} \bV_{m+1}(W_\alpha) \to \ldots
\to \bV\!/W_{\alpha} 
\end{equation}
interpolating between $ \bV$ and $\bV/\!/W_{\alpha} $. 
Fixing  $ m = (m_{\alpha})_{\alpha \in \A}\in \M(W)$, we can now define a functor 
$\,\boldsymbol{\mathcal{V}}_{m}(W):\ \Sc(W) \to \Aff_{\k} \,$ on the reflection category \eqref{sw}
 that assigns to the inclusion $ W_0 \leq W_{\alpha} $ the composition of the first $m_{\alpha}$ maps in the tower \eqref{towbV}: 
\begin{equation}
\la{bGPO}
\begin{tikzcd}[scale cd=0.85, column sep=small, row sep=small,
                 every arrow/.append style={shorten <=2pt, shorten >=2pt}]
	&& {\bV_{m_{\alpha}}(W_\alpha)} \\
	{\bV_{m_{\nu}}(W_\nu)} &&&& {\bV_{m_{\beta}}(W_\beta)} \\
	&& \bV \\
	{\bV_{m_{\lambda}}(W_\lambda)} &&&& {\bV_{m_{\gamma}}(W_\gamma)} \\
	&& {\bV_{m_{\delta}}(W_\delta)}
	\arrow["{\bpi_{\alpha, m}}"', pos=0.4, tail reversed, no head, from=1-3, to=3-3]
	\arrow["{\bpi_{\nu, m}}"', pos=0.4, tail reversed, no head, from=2-1, to=3-3]
	\arrow["{\bpi_{\beta, m}}"', pos=0.4, tail reversed, no head, from=2-5, to=3-3]
	\arrow["{\bpi_{\lambda, m}}"', pos=0.4, tail reversed, no head, from=4-1, to=3-3]
	\arrow["{\bpi_{\gamma, m}}"', pos=0.4, tail reversed, no head, from=4-5, to=3-3]
	\arrow["{\bpi_{\delta, m}}"', pos=0.4, tail reversed, no head, from=5-3, to=3-3]
\end{tikzcd}
\end{equation}
The same construction as in Proposition~\ref{cvmWfun} shows that \eqref{bGPO} is a $W$-diagram. Hence, by Lemma~\ref{Wob}, its colimit
\begin{equation} 
\la{defbvm} 
\bV_{\!m}(W)\,:=\, {\rm colim}_{\Sc(W)}[\boldsymbol{\mathcal{V}}_{m}(W)] 
\end{equation}
is a $W$-scheme. By Proposition~\ref{LQwa}, all maps $ \bpi_{\alpha, m}$ in  \eqref{bGPO} are represented by injective algebra homomorphisms $\bQ_{m_{\alpha}}(W_{\alpha}) \into \k[\bV] $. Hence, the above colimit in $ \Aff_\k$ is represented by 
the (set-theoretic) intersection of the commutative subalgebras $ \bQ_{m_{\alpha}}(W_{\alpha}) $ in $ \k[\bV] = \prod_{w \in W} \k[V] $ ({\it cf.} \eqref{VmW}). In parallel to Theorem~\ref{ThVQI}, we thus conclude
\begin{theorem} 
\la{ThVQ2}
For all $m \in \mathcal{M}(W)$, there are natural isomorphisms of $W$-schemes
\begin{equation}
\la{bVm}
\bV_{\!m}(W) \,\cong\, \Spec\, \bQ_m(W)\,,
\end{equation}
where 
\begin{equation} 
\la{ebqkwb} 
\bQ_{m}(W) = \{(p_w)_{w \in W} \in \k[V \times W]\,:\, 
p_w \equiv p_{s_\alpha w}\,{\rm mod}\,\langle\alpha\rangle^{m_{\alpha}},\ \forall\, \alpha \in \A  \} \,.
\end{equation}
%
\end{theorem}
It is easy to see that the algebras $\bQ_m(W)$ are finitely generated (and hence Noetherian) for all $ m $; we thus refer to the associated schemes $ \bV_m(W) $ as {\it varieties of quasi-covariants} of $W$. 

\subsection{Algebraic properties} \la{S2.3}
We will study the algebras $ \bQ_m(W) $ and the associated varieties $ \bV_m(W)$ in more detail. To this end we introduce the following three $W$-actions on $ \bV = V \times W\,$:
\begin{eqnarray}
W \times \bV \to \bV &, \quad & g \cdot (v,w) = (gv,\,gw)\,, \la{act1} \\
W \times \bV \to \bV &, \quad & g \cdot (v,w) = (v,\,gw)\,, \la{act2} \\
W \times \bV \to \bV &, \quad & g \cdot (v,w) = (v,\,wg^{-1})\,, \la{act3} 
\end{eqnarray}
which we refer to as diagonal, left and right, respectively. 
By duality, these define three  $W$-actions on the algebra $ \k[\bV] $, which, using the identification $\, \k[\bV] \cong \k[V] \otimes \k[W] \,$ (see \eqref{idenkwb}),  we denote as $g \otimes g $,$\, 1 \otimes g $ (acting on the left) and $ (1 \otimes g^{-1})$ (acting on the right), respectively. 

Note that the diagonal $W$-action \eqref{act1} is the one used in our construction in Section~\ref{S2.2}: by Theorem~\ref{ThVQ2}, it induces a $W$-action on the schemes 
$ \bV_m(W) $ and the algebras $\bQ_m(W) $. In the next lemma, we describe the categorical 
quotient of $ \bV_m(W)$ under this induced action, clarifying the relation to the varieties $V_k(W)$ introduced in Section~\ref{S2.1}. As mentioned in the Introduction, this result is a refinement of a key observation of \cite{BC11}.
\begin{lemma}
\la{bVV}
For all  $\, m \in \M(W) $, there are natural isomorphisms of $W$-schemes
\begin{equation}\la{isovv}
\bV_{\!m}(W)/\!/W \cong V_{[\frac{m}{2}]}(W)\,, 
\end{equation}
where the quotient is taken with respect to the diagonal $W$-action \eqref{act1} and 
the $W$-structure on $\bV_{\!m}(W)/\!/W$ is induced by the right action \eqref{act3}.
Equivalently, we have $W$-equivariant isomorphisms of commutative $\k$-algebras
\begin{equation}
\la{isoqq}
\bQ_{m}(W)^W \cong Q_{[\frac{m}{2}]}(W)\,.
\end{equation}
\end{lemma}
\begin{proof}
Note that the diagonal action \eqref{act1} is free, and its quotient 
$ \bV \to \bV\!/W $ can be represented by the action map
$\, V \times W \to V\,$, $\,(w,v) \mapsto w^{-1} v\,$ on $V$. Dually, this gives 
an algebra embedding
\begin{equation}
\la{twistemb}
\k[V] \into \k[\bV]\ ,\quad f \,\mapsto\, \sum_{w \in W} w(f) \otimes w    \ , 
\end{equation}
whose image is the subring of $W$-invariants in $ \k[\bV] \cong \k[V] \otimes \k[W]\,$:
\begin{equation}
\la{isoQm}
\k[V] \cong \k[\bV]^W \,\cong\, 
\bigl\{\,\sum_{w \in W} w(f) \otimes w\, \in \k[V] \otimes \k[W] :\ f \in \k[V] \,\bigr\}\,,
\end{equation}
where $ f \mapsto w(f)$ denotes the natural (contragradient) action of $W$ on $\k[V] = \Sym_{\k}(V^*)$. Using \eqref{ebqkwb}, the invariant subring $\,\bQ_m(W)^W \,=\,\bQ_m(W)\,\cap\, \k[\bV]^W \,$ of $ \bQ_m(W) $ can then be described as
$$
\bQ_m(W)^W \cong \{f \in \k[V]\,:\,w(f)\, -\,s_{\alpha}w(f)\,\equiv\,0 \ {\rm mod}\,\langle\alpha\rangle^{m_{\alpha}},\ \forall\, w \in W\,,\ \forall\,\alpha \in \A\}
$$
Since $m$ is $W$-invariant, the above conditions on $f$
can be simplified to $ s_{\alpha}(f) - f \equiv\, 0 \, {\rm mod}\,\langle\alpha\rangle^{m_{\alpha}}$ for all $ \alpha \in \A $, which, in turn, are equivalent to   $\, s_{\alpha}(f) - f \equiv\, 0 \, {\rm mod}\,\langle\alpha\rangle^{2 k_{\alpha}}\,$,
where $ k_{\alpha} = [m_{\alpha}/2] $ ({\it cf.} \cite[Example 5.1]{BC11}). Hence, 
for $\, k = [\frac{m}{2}] \,$, we have
$\, \bQ_{m}(W)^W \cong Q_{k}(W)\,$,
which proves \eqref{isoqq} and can be geometrically restated as \eqref{isovv}.

Note that the right action \eqref{act3} on $ \bV $ commutes with the diagonal one and, under the quotient map $ \bV \to \bV\!/W $, induces the natural (geometric) action of $W$ on $ V $. This action also restricts to all subschemes $ \bH_{\alpha} \subset \bV $, making both $ \bp_{\alpha}$ and $ \bp_{\alpha}'$  $W$-equivariant maps. Hence, all
$W_{\alpha}$-schemes $ V_{m_{\alpha}}(W_{\alpha}) $ inherit (by induction) the `right' $W$-action, so are, in fact, $(W_{\alpha}\times W)$-schemes. It follows that the schemes $ \bV_m(W)$ are $(W\times W)$-schemes, where the first factor of $ W \times W $ corresponds to the diagonal $W$-action and the second to the right action. The isomorphism \eqref{isovv} constructed above 
is then an isomorphism of $W$-schemes.
\end{proof}

Now, we state the main theorem of the this section. Recall that, for all $m \in \M(W)$, the algebra $\bQ_m(W) $ contains $ \k[V] $ as its subalgebra via the canonical inclusion
\begin{equation}
\la{canincl}
\k[V] \into \bQ_m(W)\ ,\quad f \mapsto \sum_{w \in W} f \otimes w \,.
\end{equation}
\begin{theorem}
\la{ThFree}
Assume that $ m \in \M(W) $ is either purely even or purely odd, i.e. $m_{\alpha} = 2 k_{\alpha} $
or $ m_{\alpha} = 2 k_{\alpha} + 1 $ for all $ \alpha \in \A$. Then 
$ \bQ_m(W) $ is a free module over $ \k[V] $ of rank $|W|$. 
In particular, $\bV_m $ are Cohen-Macaulay varieties.
\end{theorem}
We will deduce Theorem~\ref{ThFree} from the main results of \cite{BC11}. As a first step,
we give an alternative description of the module $\bQ_m$ (valid for all $ m \in \M(W)$), using
the left $W$-action defined by \eqref{act2}. Recall that this action defines 
an action of the group algebra $\k[W]$ on  $ \k[\bV] = \k[V] \otimes \k[W]$, which (following \cite{BC11}) we denote by `$ 1 \otimes g $' for $ g \in \k[W]$.

\begin{lemma} 
\la{fQGKM} 
For all $\,m \in \M(W)$, under the identification \eqref{idenkwb}, we have
$$ 
\bQ_{m}(W)\,=\, \,\{\,p \in \k[V] \otimes \k[W]\ :\ (1 \otimes e_{\alpha})\,p \,\equiv\,0\ \, {\rm mod}\,  (\langle\alpha\rangle^{m_{\alpha}} \otimes \k[W])\,,\ \forall\,\alpha \in \A\,\}\,, 
$$
where $\, e_{\alpha} := \frac{1}{2}(1 - s_{\alpha}) \in \k[W] $.
\end{lemma}
\begin{proof} For a fixed $ \alpha \in \A $, we can write the elements of $ \k[V] \otimes \k[W] $ as 
$$ 
p = \sum_{w \in W} p_w \otimes w = \sum_{[w] \in W_\alpha\! \backslash W}(p_w \otimes w +  p_{s_{\alpha}w} \otimes s_\alpha w) \,,
$$
where the last sum is taken over the set of all right cosets of $W_{\alpha}$ in $W$. 
It follows that
$$
(1 \otimes e_\alpha) p \,=\, \frac{1}{2}
\sum_{[w]\in W_{\alpha}\!\backslash W} [\,(p_w - p_{s_\alpha w}) \otimes w \,+\,
(p_{s_\alpha w}-p_w) \otimes s_{\alpha}w \,] \,.
$$
Hence $\, (1 \otimes e_{\alpha})\,p \in  \langle \alpha \rangle^{m_\alpha} \otimes \k[W] \,$
if and only if $p_w-p_{s_\alpha w} \in \langle \alpha\rangle^{m_\alpha}$ in $ \k[V] $ for all $w \in W$.
\end{proof}

Now, we recall the definition of the modules of `$\k W$-valued quasi-invariants' of complex reflection groups introduced in \cite{BC11}. When specialized to Coxeter groups, this definition reads (see {\it loc.cit.}, formula (3.8)):
\begin{equation}
\la{BCdef}
\bQ_k^{\rm BC}(W) := \,\{\, p \in  \k[V_{\rm reg}] \otimes \k[W] \ :\ (1 \otimes e_{\alpha})p \,\equiv\,0\ {\rm mod}\, (\langle\alpha\rangle^{2k_{\alpha}} \otimes \k[W]),\ \forall\,\alpha \in \A\,\}\,,   
\end{equation}
where $\,\k[V_{\rm reg}] := \k[V\!\setminus \!\cup_{\alpha \in\A} H_{\alpha}]\, $.
As already observed in \cite{BC11}, the conditions imposed on elements of $\, \bQ_k^{\rm BC}(W) \,$  in \eqref{BCdef} actually force $\, \bQ_k^{\rm BC}(W) \,$ to be a subspace of $\, \k[V] \otimes \k[W] \,$. Hence, Lemma~\ref{fQGKM} implies immediately
\begin{cor} \la{BCcor}
 $\, \bQ_k^{\rm BC}(W) = \bQ_{2k}(W) \,$    
for all $\,k \in \M(W)$.
\end{cor}
This shows that the algebras \eqref{ebqkwb} are natural generalizations of the modules of $\k W$-valued quasi-invariants introduced in \cite{BC11}.

Next, we examine the case when $ m = 2k+1 \,$ (i.e. $ m_{\alpha} = 2 k_{\alpha}+1 $ for all $ \alpha \in \A $). It is convenient to approach this case geometrically. To this end, we consider the categorical quotients $\bV/\!/W$ with respect to the first two $W$-actions \eqref{act1} and \eqref{act2} on $\bV = V \times W$ and observe that they can be both identified with $V$ in such a way that the quotient maps make the commutative diagram
\[
\begin{diagram}[small]
\bV & \rTo^{q_1} & V\\
\dTo^{q_2}\ &  & \dTo\\
V    & \rTo & V/\!/W
\end{diagram}
\]
where $ V \to V/\!/W $ is the quotient with respect to the geometric action of $W$ on $V$. By the universal mapping property of fiber products, we then get the canonical map in $\Aff_{\k} $
\begin{equation}
\la{dualphi}
q:\ \bV \to V \times_{V/\!/W} V\ ,\quad (v,\,w)\ \mapsto\ (v,\,w^{-1}v)\,.
\end{equation}
Dually, with identification \eqref{idenkwb}, this defines a homomorphism of commutative $\k$-algebras
\begin{equation}
\la{defphi}
\varphi: \ 
\k[V] \otimes_{\k[V]^W} \k[V] \,\to \, \k[V] \otimes \k[W]\ , \quad
f \otimes  g\,\mapsto\,  \sum_{w \in W}  f\, w(g) \otimes w
\end{equation}
Now, recall that, by construction, the varieties $V_k$ and $ \bV_m$ come with natural maps
$\pi_k: V \to V_k $ and $ \bpi_m: \bV \to \bV_{m} $ for all $k$ and $m$ in $ \M(W)$. We claim
\begin{prop}
\la{bQodd}
For all $ k\in \M(W) $, the composite map $\,(1 \times \pi_k)\,q:\, \bV \to V \times_{V/\!/W} V \to V \times_{V/\!/W} V_k \,$ factors though $ \bpi_{2k+1}: \bV \to \bV_{2k+1} $, inducing an {\rm isomorphism} of $\k$-schemes
\begin{equation}
\la{oddisosch} 
q_{k}:\ \bV_{\!2k+1} \xrightarrow{\sim}  V \times_{V/\!/W} V_k
\end{equation}
Equivalently, for all $ k\in \M(W) $, the map \eqref{defphi} restricts to an isomorphism of $\k$-algebras
\begin{equation}
\la{oddisoalg} 
\varphi_k:\ \k[V] \otimes_{\k[V]^W} Q_k(W)\,\xrightarrow{\sim}\,\bQ_{2k+1}(W)\,.
\end{equation}
\end{prop}
\begin{proof}
Clearly, $\, \varphi_k = q_k^* \,$, hence it suffices to prove \eqref{oddisoalg}.
First of all, we note that the map \eqref{defphi} being injective with image $ \bQ_1(W) $ 
is a well-known result that underlies two classical presentations of the $T$-equivariant cohomology $ H_T^*(G/T, \k) $ of the flag manifold  (see \cite{GKM98, KK87, AF23} or \cite[Proposition 3.5]{Kaji}). Thus, for $k=0$, we know that the map \eqref{oddisoalg} is an algebra isomorphism, and we denote its inverse by
$$
\psi := \varphi_0^{-1}:\ \bQ_1(W)\, \xrightarrow{\sim} \,\k[V] \otimes_{\k[V]^W} \k[V]\,. 
$$ 
Next, we observe that, for all $ k \in \M(W) $, the natural inclusions $ Q_k(W) \subseteq \k[V] $ 
induce {\it injective} algebra homomorphisms $\, \k[V] \otimes_{\k[V]^W} Q_k(W) \into \k[V] \otimes_{\k[V]^W} \k[V]\,$ (since $ \k[V]$ is free as a module over $ \k[V]^W$). Hence the maps $\varphi_k $ obtained by composing these homomorphisms 
with $\varphi$ are injective as well. We need only to check that these maps have correct images.

First, we show that 
\begin{equation}
\la{inj_k}
 \im(\varphi_k) \subseteq \bQ_{2k+1} \ \mbox{for all}\ k \in \M(W)\,.
\end{equation}
Given an element $ f \otimes g \in \k[V] \otimes_{\k[V]^W} Q_k $, we apply $ \,1 \otimes e_{\alpha} \,$ to its image $\, \varphi(f \otimes g) = \sum_{w \in W} f w(g) \otimes w \,$ 
in $ \k[V] \otimes \k[W] $ and find
$$
(1 \otimes e_{\alpha}) \cdot \varphi(f \otimes g) = \frac{1}{2} \sum_{w \in W }
f\,(w(g) - s_{\alpha}w(g)) \otimes w
$$
Since $ Q_k \subseteq \k[V] $ is stable under the natural action of $W$, for any $ g \in Q_k $, we have $\, w(g) - s_{\alpha} w(g) \equiv 0 \ \mod \,\langle\alpha \rangle^{2 k_{\alpha} +1} \,$, which shows that  $\,(1 \otimes e_{\alpha}) \cdot \varphi(f \otimes g) \,\in\, \langle\alpha \rangle^{2 k_{\alpha} +1} \otimes \k[W] \,$ for all $ f \in \k[V]$. Then, by Lemma~\ref{fQGKM} , we conclude that $\,\varphi(\k[V] \otimes_{\k[V]^W} Q_k) \subseteq \bQ_{2k+1}\,$, which is equivalent to \eqref{inj_k}.

To prove the opposite inclusion we note that the map \eqref{dualphi} is  equivariant with respect to the right $W$-action \eqref{act3} on $\bV$ and the natural action of $W$ on the second factor in $V \times_{V/\!/W} V$. Hence, for the dual map \eqref{defphi} (and all its restrictions $\varphi_k$), 
we have
\begin{equation} 
\la{equiv2} 
\varphi[(1 \otimes s) \cdot (f \otimes g)] \,=\, \varphi(f \otimes g) \cdot (1 \otimes s^{-1}) 
\,,\quad \forall\,s\in W\,.
\end{equation}
For the inverse map $ \psi = \varphi_0^{-1} $, this reads $\, (1 \otimes s) \cdot \psi(p) = \psi[p\cdot (1 \otimes s^{-1})]\,$, where $ p \in \bQ_1(W) $, which, in particular, implies
\begin{equation}
 \la{psiQ_k}
 (1 \otimes e_{\alpha})\cdot \psi(\bQ_{2k+1}) = \psi[\bQ_{2k+1}\cdot (1 \otimes e_{\alpha})] \,,\quad \forall\,\alpha \in \A\,.
\end{equation}
We now look closer at the subspace $\bQ_{2k+1}\cdot (1 \otimes e_{\alpha}) \subset \k[V] \otimes \k[W] $ for a fixed $ \alpha \in \A $. For an element $ p = \sum_{w\in W} p_w \otimes w \in \k[V] \otimes \k[W]$, we compute
$$
p\cdot (1 \otimes e_{\alpha}) = 
 \frac{1}{2} \sum_{w \in W} \,(p_w - p_{w s_{\alpha}}) \otimes w \,= \, 
\frac{1}{2} \sum_{w \in W} \,(p_w - p_{s_{w(\alpha)}w}) \otimes w
$$
Hence, if $\, p \in \bQ_{2k+1}\,$, the difference $\,p_w - p_{s_{w(\alpha)}w}$ is divisible by $ w(\alpha)^{2 k_{\alpha} +1 }\,$ for all $ w \in W $; then, setting  $ q_w := \frac{1}{2} (p_w - p_{s_{w(\alpha)}w})\, w(\alpha)^{- 2 k_{\alpha} - 1 } \in \k[V]$, we can write
\begin{equation}
\la{pQ_1}
p\cdot (1 \otimes e_{\alpha}) = 
\sum_{w \in W} \, w(\alpha)^{2k_{\alpha}+1} q_w\otimes w \,=\, \varphi(1 \otimes \alpha^{2k_{\alpha}+1}) \cdot \sum_{w \in W} \,(q_w \otimes w)
\end{equation}
where the `dot' product on the right denotes the (commutative) multiplication  in $ \k[V] \otimes \k[W]$, see \eqref{prodVW}. We claim that  
\begin{equation}
\la{wqwQ1}
\sum_{w \in W}\, q_w \otimes w\ \in \ \bQ_1(W)
\end{equation}
Indeed, for any root $\beta$, we have the equation
\begin{equation} 
\la{genddq}
2\,(w\alpha)^{2k_{\alpha}+1}\,(s_\beta w \alpha)^{2k_\alpha+1}\, (q_w-q_{s_{\beta}w})\,= \, (s_\beta w\alpha)^{2k_\alpha+1}(p_w-p_{ws_{\alpha}})-(w\alpha)^{2k_\alpha+1}(p_{s_{\beta}w}-p_{s_{\beta}ws_{\alpha}})
\end{equation}
If $\beta= c \,w(\alpha)$ for some $ c \in \k^*$, then $s_\beta = s_{w(\alpha)} $ and  $s_{\beta}w(\alpha)=-w(\alpha)$. Hence, in this case, the right-hand side of \eqref{genddq} --- and thus the difference $(q_w-q_{s_{\beta}w})$ --- vanish identically in $\k[V]$. 

On the other hand, since $\,(s_\beta w \alpha)^{2k_\alpha+1} \equiv (w\alpha)^{2k_{\alpha}+1}\,\mod\,\langle \beta\rangle\,$, the right-hand side of \eqref{genddq} is congruent (modulo $\beta $) to
$$
(w\alpha)^{2k_{\alpha}+1}\left[(p_w - p_{s_\beta w})\,-\,(p_{ws_{\alpha}} - p_{s_{\beta}(w s_{\alpha})}) \right]\,,
$$
which is, in turn, congruent to $0$, since $ p = \sum p_w \otimes w \in \bQ_{2k+1}$.
Thus, in general, \eqref{genddq} implies 
\begin{equation} 
\la{genddq1}
(w\alpha)^{2k_{\alpha}+1}\,(s_\beta w \alpha)^{2k_\alpha+1}\, (q_w-q_{s_{\beta}w})\,\equiv\, 0\ \mod\,\langle \beta\rangle
\end{equation}
Now, if $ \beta \not= c \,w(\alpha) $ for $ c \in \k^*$, then the linear forms $w(\alpha)$ and $s_\beta w(\alpha)$ are not divisible by $ \beta$ in $ \k[V] $. Hence, \eqref{genddq1} may hold only if $\,q_w-q_{s_\beta w}\,\equiv\, 0\ \mod\,\langle \beta\rangle\,$. This proves \eqref{wqwQ1}. Then, by \eqref{pQ_1} and \eqref{wqwQ1}, 
$$
\bQ_{2k+1} \cdot (1 \otimes e_{\alpha}) \,\subseteq  \,\varphi(1 \otimes \alpha^{2k_{\alpha}+1}) \cdot \bQ_1
$$
which, in combination with  \eqref{psiQ_k}, yields
\begin{equation} 
\la{bqtoqps2}
(1 \otimes e_{\alpha})\cdot \psi(\bQ_{2k+1}) \subseteq  (1 \otimes \alpha^{2k_{\alpha}+1}) \cdot \psi(\bQ_1) = \k[V]\otimes_{\k[V]^W} \langle \alpha \rangle^{2k_{\alpha}+1}\,.
\end{equation}
Choosing a basis $\{f_w\}_{w \in W}$ in $ \k[V]$ as a (free) module over $ \k[V]^W$, we can write
the elements of $ \k[V] \otimes_{\k[V]^W} \k[V] $ uniquely in the form $ \sum_{w\in W} f_w \otimes g_w $. Then, the above inclusion says that $\,\sum_{w \in W} f_w \otimes (g_w - s_{\alpha} g_w) \in \sum_{w \in W} f_w \otimes \langle \alpha \rangle^{2k_{\alpha}+1}\,$ whenever $ \sum_{w \in W} f_w \otimes g_w \in \psi (\bQ_{2k+1})$, which implies, by uniqueness, 
$ s_{\alpha}g_w \equiv g_w \,\mod\,\langle \alpha \rangle^{2k_{\alpha}+1} $ for all $ w \in W$. 
Since this holds for each $ \alpha \in \A $, we conclude that $ g_w \in Q_k(W) $. Hence 
\eqref{bqtoqps2} implies
$$
\psi(\bQ_{2k+1}) \subseteq \k[V] \otimes_{\k[V]^W} Q_k \,.
$$ 
Now, recall that $ \psi $ is an algebra isomorphism inverse to $ \varphi_0 $, hence
applying $\varphi_0 $, we get
\begin{equation*} 
\bQ_{2k+1} \subseteq \varphi_0(\k[V] \otimes_{\k[V]^W} Q_k) =  \varphi(\k[V] \otimes_{\k[V]^W} Q_k)\, =\, \im(\varphi_k)\ .
\end{equation*}
This proves the opposite inclusion to \eqref{inj_k} and thus completes the proof of the proposition.
\end{proof}
\begin{proof}[Proof of Theorem~\ref{ThFree}]
For $ m = 2k $, the result follows immediately from Corollary~\ref{BCcor} and 
\cite[Proposition 8.1]{BC11}. For $ m = 2k+1$, we observe that,  under the isomorphism \eqref{oddisoalg}, the natural inclusion $ \k[V] \into \bQ_m(W) $, see \eqref{canincl}, 
corresponds to the inclusion $ \k[V] \into \k[V] \otimes_{\k[V]^W} Q_k(W) $ as 
the first factor. Since $ Q_k(W) $ is a free module over $ \k[V]^W $ of rank $|W|$ (\cite[Theorem~1.1]{BC11}), the result follows from  Proposition~\ref{bQodd}.
\end{proof}
\begin{remark}
As shown in (the proof of) Lemma~\ref{bVV}, the algebra map \eqref{twistemb} restricts to the embedding $ Q_{[\frac{m}{2}]}(W)  
\into \bQ_m(W) $, identifying $Q_{[\frac{m}{2}]}(W)$ with the invariant subring $ \bQ_m(W)^W $ of 
$ \bQ_m(W)$ for the diagonal $W$-action. When $ m = 2k+1$, this embedding corresponds, under the isomorphism \eqref{oddisoalg}, to the natural inclusion 
$ Q_k(W) \into \k[V] \otimes_{\k[V]^W} Q_k(W) $ as the second factor. As $\k[V]$ is free over
$\k[V]^W$, we conclude that $ \bQ_{2k+1}(W) $ is also free over $ \bQ_{2k+1}(W)^W \cong Q_k(W) $ of rank $|W|$. We do not know whether this result holds true in general when $ m \not= 2k+1$.
\end{remark}
%
A natural question is whether the result of Theorem~\ref{ThFree} holds for `mixed' multiplicities: namely, whether $\bQ_m(W)$ remains a free module over $\k[V] $ when some of the $m_{\alpha}$'s are even and some are odd. The next proposition shows that this is indeed the case for all Coxeter groups of rank two.
\begin{prop}
\label{prop:QC-rank2-free}
Let $W$ be a rank two finite Coxeter group acting in its reflection representation $V = \k^2$.
Then, for any multiplicity function $m$, $\bQ_m(W)$ is a free $\k[V]$-module of rank~$|W|$.
\end{prop}

\begin{proof}
To simplify the notation we set $R :=\k[V]$, $K :=\operatorname{Frac}(R)$, $M :=\bQ_m(W)$, 
and $M_\alpha := \bQ_m(W_\alpha)$ for each  $\alpha \in \A$ (as defined in~\eqref{ebqkwa}). 
First, we show that the subalgebra $M_\alpha \subseteq \k[\bV]=\prod_{w \in W} \k[V]$ 
is a reflexive $\k[V]$-module. 
It fits into the short exact sequence 
\[ 
  0 \to M_\alpha \to \k[\bV] \to \k[\bV]/M_\alpha \to 0 
\]
where each $\k[\bV]/M_\alpha \cong \bigoplus_{W/W_\alpha} R/(\alpha^{m_\alpha})$, whose component $R/(\alpha^{m_\alpha})$ has only height-one associated prime $(\alpha)$. By the depth lemma (see \cite[Lemma 10.72.6]{StackPr}), for any prime $\mathfrak{p}\subset R$, 
\[
  \operatorname{depth}((M_\alpha)_{\mathfrak{p}})
  \;\ge\; \min \bigl(\operatorname{depth}(\k[\bV]_{\mathfrak{p}}),\;
  \operatorname{depth}\bigl((\k[\bV]/M_\alpha)_{\mathfrak{p}}\bigr) + 1\bigr) \,.
\]
As $\operatorname{depth}(\k[\bV]_{\mathfrak{p}}) = \operatorname{ht}(\mathfrak{p})$ and $\operatorname{depth}\bigl((\k[\bV]/M_\alpha)_{\mathfrak{p}} = \operatorname{ht}(\mathfrak{p}) - 1$,
we have
\(
  \operatorname{depth}((M_\alpha)_{\mathfrak{p}})
  \;\ge\; \operatorname{ht}(\mathfrak{p}),
\)
showing that $M_\alpha$ satisfies Serre's criterion $(S_2)$.
As $M_\alpha$ is also torsion free, it is reflexive. 

Next we show $M$ is reflexive. As both $M$ and $M_\alpha$ are finitely generated torsion-free module over the normal Noetherian domain $R$, the reflexive hull $M^{**}=\Hom_R(\Hom_R(M,R),R)$ satisfies $M^{**}= \bigcap_{\operatorname{ht}(\mathfrak{p})=1} M_{\mathfrak{p}} \subset M\otimes_R K$ and likewise for $M_\alpha$ (see \cite[Lemma 15.24.19]{StackPr}). As $M=\bigcap_{\alpha}M_\alpha$, 
\[
M^{**} = \bigcap_{\operatorname{ht}(\mathfrak{p})=1} M_{\mathfrak{p}} \subseteq \bigcap_{\operatorname{ht}(\mathfrak{p})=1} (M_\alpha)_{\mathfrak{p}} = M_\alpha^{**} = M_\alpha, \quad \forall \alpha\in\A.
\]
Hence $M^{**}\subseteq \bigcap_\alpha M_\alpha = M$, so $M=M^{**}$ is reflexive.

Since $\dim R = 2$, reflexivity gives $\operatorname{depth}(M_{\mathfrak{m}}) \geq 2$ at every maximal ideal (see \cite[Lemma 15.24.16]{StackPr}), whence $M$ is projective by the Auslander--Buchsbaum formula
\[
\operatorname{pd}(M_\mathfrak{m}) + \operatorname{depth}(M_\mathfrak{m}) = \dim R_\mathfrak{m} = 2.
\]
By Quillen-Suslin theorem, $M$ is a free module over $R$. The rank is $|W|$ as $M\otimes_R K \cong K^{|W|}$.

\end{proof}

\begin{remark}
The argument above shows that $\bQ_m(W)$ is still a reflexive $\k[V]$-module for rank~$\ge 3$, but in the geneal case, it only gives $\operatorname{pd}(M_\mathfrak{m}) \le \dim \k[V] - 2$, so we lack proof of freeness.
At this moment, we do not know if $\bQ_m(W)$ is free for mixed multiplicity when rank~$\ge 3$.
\end{remark}

\begin{remark}
\la{Rspline}
As we mentioned in the Introduction, our algebras $ \bQ_m(W)$ can be viewed as generalized spline algebras in the sense of \cite{GTV16}. Given a (finite simple) graph $ \Gamma = (V_{\Gamma}, E_{\Gamma}) $ together with a map $ e \mapsto I(e) $ assigning to each edge $ e \in E_{\Gamma} $
an ideal $ I(e) \subset R $ in a commutative ring $R$, the ring of generalized 
splines of $\Gamma$ is defined by
\begin{equation}
\la{splalg}
R_{\Gamma, I} := \bigl\{\, (p_v) \in \prod_{v \in V_{\Gamma}} R\ :\ p_u \equiv p_v\ \mod\, I(e)\ \, \mbox{for every edge}\ \, e = (u  v) \in E_{\Gamma}\bigr\}
\end{equation}
Obviously, if $ R = \k[V] $ and $ \Gamma = \Gamma(\R_W) $ is the Bruhat graph of  $W$ (see Section~\ref{S4.1}) with edge-labeling function $\, I: e(s_{\alpha}, w) \mapsto \langle \alpha \rangle^{m_{\alpha}}\,$, then we recover from \eqref{splalg} our definition of $\bQ_m(W)$.
When $R$ is an integral domain, it is proved in \cite[Corollary 5.2]{GTV16} that, for any graph $ \Gamma $, there is an {\it inclusion} of $R$-modules: $\, R^n \, \into \,R_{\Gamma, I} \,$, 
where $ R^n $ is the free module of rank $ n = |V_{\Gamma}| $. 
In our case, by Theorem~\ref{ThFree}, such an inclusion is an isomorphism.
\end{remark}

\subsection{Double affine and nil-Hecke algebra actions} \la{S2.4}
Let $\, \vreg := V \setminus \bigcup_{\alpha \in \A} H_{\alpha}\,$. We regard $ \vreg $ as an affine variety and denote by $ \k[\vreg] $ and $ \D(\vreg) $ the rings of regular functions and regular (algebraic) differential operators on $ \vreg $, respectively. We will often identify
$ \k[\vreg] $ with the subalgebra of $ \D(\vreg) $ consisting of zero order (multiplication) operators.
The action of $ W $ on $ V $ restricts to a free action on $ \vreg $, which induces a natural $W$-action   on $\k[\vreg]$ and $ \D(\vreg) $ by algebra automorphisms. We write $ \k[\vreg] \rtimes W $ and
$ \D(\vreg) \rtimes W$ to be the corresponding crossed products. As an algebra, $ \D(\vreg) \rtimes W $ is generated by its two subalgebras $ \k[W] $ and $ \D $: that is, by the elements of $ W $, the regular functions in $\, \k[\vreg] \,$ and the $k$-linear derivations $\, \partial_\xi $ on $ \k[\vreg] $ defined for each $\, \xi \in V $. Clearly, $ \k[\vreg] \rtimes W  $ is the subalgebra of $\D(\vreg) \rtimes W $ generated $W$ and the $\k[\vreg]$.

To define Cherednik (a.k.a. rational double affine Hecke) algebra $ \bH_k(W) $ we fix  $ k \in \M(W)$ and introduce the family of {\it Dunkl $($diffrential-difference$)$ operators} $ T_{\xi} = T_{\xi, k} $ in $ \D(\vreg) \rtimes W $ depending on $\xi \in V$:
\begin{equation}
\label{du} 
T_\xi := \partial_\xi-\sum_{\alpha \in \A} k_{\alpha} \,\frac{\alpha(\xi)}{\alpha}\,(1 - s_{\alpha})\ .
\end{equation}
The key property of the operators \eqref{du}, first established in \cite{Du89}, is that they commute in different directions: $\, T_{\xi}\,T_{\eta} - T_{\eta}\,T_{\xi} = 0 \,$  for all $\xi, \eta \in V $.
Another useful property, which is easy to check directly, is the $W$-equivariance: $\,w\,T_\xi = T_{w(\xi)}\,w\,$ for all $w \in W$. These two properties show that the assignment $\,\xi \mapsto T_\xi\,$ extends to an injective algebra homomorphism
\begin{equation}
\la{hom}
\k[V^*] \,\hookrightarrow \,\D(\vreg) \rtimes W \ ,\quad  \xi_1 \xi_2\ldots \xi_n \,\mapsto \,T_{\xi_1} T_{\xi_2} \ldots T_{\xi_n} \ .
\end{equation}
Identifying $ \k[V^*] $ with its image in $ \D(\vreg) \rtimes W $ under
\eqref{hom}, we now define the {\it Cherednik algebra}
$\bH_k(W)$ of $W$ as the subalgebra of $\D(\vreg) \rtimes W$ generated by $ \k[V] $,
$\,\k[V^*] $ and $ \k[W]$. 

The Cherednik algebras can be also defined abstractly, 
in terms of generators and relation (see, e.g., \cite{EG02, BEG03, GGOR03}). 
The definition makes sense for arbitrary ($ \k$-valued) multiplicities
$ k $ (not necessarily integers). The family $ \{\bH_k(W)\} $ can be then viewed as a deformation (in fact, the universal deformation) of the crossed product algebra $\, \bH_0(W) = \D(V)\rtimes W \,$. The algebra inclusion $\, \bH_k(W) \into \D(\vreg) \rtimes W
\,$ defined by \eqref{hom} is called the {\it Dunkl embedding}. The existence of such an embedding implies the PBW property for $ \bH_k(W) $, which says that the multiplication map:
$\, \k [V] \otimes \k W \otimes \k[V^*] \stackrel{\sim}{\to} \bH_k(W)\,$
is an isomorphism of vector spaces.

To define the action of $ \bH_{k}(W) $ on $ \bQ_{m}(W)$ we will first define the action of 
$ \D(\vreg) \rtimes W $ on the space $\, \k[\vreg] \otimes \k[W] $ of $ \k[W]$-valued functions on $ \vreg$. To this end, we consider
the canonical inclusion $\, \k[\vreg] \otimes \k[W] \into \D(\vreg) \rtimes W $ and identify its image
with the quotient $(\D(\vreg) \rtimes W)$-module $ (\D(\vreg) \rtimes W)/J $, where $ J $ is the left ideal of $ \D(\vreg) \rtimes W $
generated by the derivations $\, \partial_\xi \in \D(\vreg) \rtimes W $ for all $\, \xi \in V \,$. 
Explicitly, the action of $ \D(\vreg) \rtimes W$ on $ \k[\vreg] \otimes \k[W] $ is then given by
\begin{eqnarray}
&& g(f\otimes w) =  gf\otimes w\,,\quad g\in\k[\vreg]\ , \nonumber\\
\la{diffact}
&& \partial_\xi(f\otimes w) = \partial_\xi f\otimes w\,,\quad \xi\in V\ ,\\
\nonumber
&& s(f\otimes w)=s(f)\otimes sw\,,\quad s\in W\ .
\end{eqnarray}

Now, restricting scalars via the Dunkl embedding
$\,\bH_k(W) \into \D(\vreg) \rtimes W $ makes $ \k[\vreg]\otimes \k[W] $ an $ \bH_k(W)
$-module for all $ k \in \M(W)$. Following \cite{BC11}, we refer to this as the {\it differential representation} of $ \bH_k(W) $.  In this representation, the Dunkl operators \eqref{du}
act as follows
\begin{equation}
\la{dudiff}
T_{\xi}\bigl[\sum_{w \in W} f_w \otimes w\bigr] \,=\,
\sum_{w \in W}\, \bigl(\,\partial_{\xi}f_w - \sum_{\alpha \in \A} k_{\alpha} \frac{\alpha(\xi)}{\alpha}(f_{w} - s_{\alpha} f_{s_{\alpha} w})\,\bigr)  \otimes w  
\end{equation}
Note that the natural action of the Dunkl operators \eqref{du} on $\k[\vreg]$ preserves the subspace
of polynomials $ \k[V] \subset \k[\vreg]$ for all $k$. In contrast, the action \eqref{dudiff} does {\it not} preserve the subspace $\, \k[V] \otimes \k[W] \subset \k[\vreg]\otimes \k[W]\,$ if $ k \not= 0 $. Nevertheless, we have the next result that follows (by Corollary~\ref{BCcor}) from \cite[Theorem 3.4]{BC11}.
\begin{theorem}[\cite{BC11}]
\la{Qfat} 
The subspace $\bQ_{2k}(W) $ defined by \eqref{ebqkwb} for $m=2k$ is stable under the action \eqref{dudiff} and hence carries a left module structure over the Cherednik algebra $\bH_k(W)$ extending the natural $\k[V]$-module structure on $\bQ_{2k}(W) $ defined
via \eqref{canincl}.
\end{theorem}
Next, following \cite{BEG03}, we introduce the {\it spherical Cherednik subalgebra}  
$\, \e\bH_{k}\e \,$ of $\bH_k$, where the product is inherited from $ \bH_k $ but 
the identity element is different: namely, the symmetrizing idempotent $ \e := \frac{1}{|W|} \sum_{w\in W}w \in \bH_k $. The Dunkl embedding of $ \bH_k $ restricts to its spherical subalgebra, giving an algebra monomorphism 
$\,\e \bH_k \e \into \D(\vreg)^W \,$ called the spherical Dunkl embedding.
Theorem~\ref{Qfat}  implies that the image of this embedding coincides with the subalgebra $ \D(V_k)^W \subset \D(\vreg)^W  $ of $W$-invariant differential operators  on the variety $ V_k \,$ (see  \cite[Proposition 4.3]{BC11}), so that
$$
\e\bH_{k}\e \,\xrightarrow{\sim}\,\D(V_k)^W,
$$
which defines a natural differential action of $\e\bH_{k}\e $ on $Q_k$. Further, using the fact that the algebra $ \e\bH_{k}\e $ is Morita equivalent to $\bH_k$ (see \cite[Theorem 3.1]{BEG03}), one can deduce from 
Theorem~\ref{Qfat} that,   for all $ k \in \M(W) $, there are $\bH_k$-module isomorphisms
\begin{equation}
\la{Hkmodiso}
\bH_k \e \otimes_{\e \bH_k \e} Q_{k}(W)\,\cong\, \bH_k \e \otimes_{\e \bH_k \e} \e\bQ_{2k}(W)\,\xrightarrow{\sim}\, \bQ_{2k}(W)\,
\end{equation}
where the second arrow is given by the canonical multiplication-action map.

The result of Theorem~\ref{Qfat} does not extend to arbitrary multiplicities: 
the subspaces $ \bQ_m(W) \subset \k[\vreg] \otimes \k[W]$ are not stable under the action  of Dunkl operators \eqref{dudiff} --- and hence, do not carry a (differential) module structure over Cherednik algebras --- unless $ m = 2k$. Instead, for $ m = 2k+1$, we can define on $ \bQ_m(W) $ a natural action of the {\it Demazure $($nil-Hecke$)$ algebra} $ \NH(W) $ of $W$. By analogy with $ \bH_k(W) $, we introduce $ \NH(W) $ as the subalgebra of $ \k[\vreg] \rtimes W  $ 
generated by the elements of $\k[V] $, $\, W$ and the family of {\it Demazure $($divided difference$)$ operators} 
 $\,\nabla_{\alpha} := \frac{1}{\alpha}(1 - s_{\alpha}) \in \k[\vreg] \rtimes W\,$ defined for  $ \alpha \in \A $. Choosing a basis of `simple roots' $S = \{\alpha_1, \ldots, \alpha_n\}$ in $\A$ we can write 
 the elements $w \in W $ as (reduced) products $ w = s_{\alpha_{i_1}} s_{\alpha_{i_2}} \ldots s_{\alpha_{i_k}}$  with $ \alpha_{i} \in S $. Then, we  define
$$
\nabla_w := \nabla_{\alpha_1} \nabla_{\alpha_2} \ldots \nabla_{\alpha_n} \,\in\, \NH(W) \,.
$$
It is well known that the operators $ \nabla_w $, while depending on the choice $S$ of simple roots, do not depend on the reduced decomposition of $w$; moreover, they satisfy the relations 
$$
\nabla_w \nabla_{w'} = \left\{\begin{array}{cc}
  \nabla_{w w'}   & \mbox{if}\quad l(ww') = l(w) + l(w')  \\
   0  & \mbox{\rm otherwise}
\end{array}
\right.
$$
where $l: W \to \Z_+ $ is the length function on $W$ associated to $S$. 
It follows from the above relations that $\NH(W) $ is a free (left) module over $ \k[V]$ with basis $ \{\nabla_w\}_{w \in W}$  (see, e.g., \cite[Section 4, Corollary 1]{D73}): in particular, every element of $ \NH(W) $  can be written uniquely in the form $ \sum_{w\in W}\, a_w \nabla_w\,$, where $ a_w \in \k[V] $. This is a PBW property for the algebra $ \NH(W)$.

Now, to state the analog of Theorem~\ref{Qfat} for $ m =2k+1$ we consider the action of $\NH(W)$
on $ \k[\vreg] \otimes \k[W] $ obtained by restricting the natural (left) 
multiplication-action of $ \k[\vreg]\rtimes W $. Explicitly, the Demazure operators act
on $ \k[\vreg] \otimes \k[W] $ by\footnote{This action corresponds topologically to the so-called
{\it left} action on $T$-equivariant cohomology of flag manifolds, see Remark~\ref{leftDem}.}
\begin{equation}
\la{demdiff}
\nabla_{\alpha}\bigl[\sum_{w \in W} f_w \otimes w\bigr] \,=\,
\sum_{w \in W}\, \frac{(f_{w} - s_{\alpha} f_{s_{\alpha} w})}{\alpha}  \otimes w  
\end{equation}

We claim
\begin{theorem}
\la{Qfatodd} 
The subspace $\bQ_{2k+1}(W) $ defined by \eqref{ebqkwb} for $m=2k+1$ is stable under the action \eqref{demdiff} and hence carries a left module structure over the Demazure algebra $\NH(W)$ extending the natural $\k[V]$-module structure on $\bQ_{2k+1}(W) $ defined
via \eqref{canincl}.
\end{theorem}
\begin{proof} Let $ \sum_{w \in W} f_w \otimes w \in \bQ_{2k+1}(W) $. First, observe that 
$\,f_w - s_{\alpha} f_{s_{\alpha} w} = (f_w - f_{s_{\alpha} w}) + (f_{s_{\alpha}w} - s_{\alpha} f_{s_{\alpha} w}) \equiv 0\ \mod \langle \alpha \rangle \,$. Hence, under \eqref{demdiff}, $\nabla_{\alpha}(\sum_{w \in W} f_w \otimes w) \in \k[W] \otimes \k[V]\,$ for all $ \alpha \in \A $. Next, for a pair 
$ \alpha, \beta \in \A $, consider the difference
\begin{equation}
\la{abdiff}
\frac{(f_{w} - s_{\alpha} f_{s_{\alpha} w})}{\alpha}\,-\,\frac{(f_{s_{\beta}w} - s_{\alpha} f_{s_{\alpha} s_{\beta} w})}{\alpha}\,=\,\frac{(f_w - f_{s_{\beta} w}) - s_{\alpha}(f_{s_{\alpha}w} - f_{s_{\alpha} s_{\beta}w})}{\alpha}\ .
\end{equation}
Since
$$ 
f_{s_{\alpha}w} - f_{s_{\alpha} s_{\beta}w} = (f_{s_{\alpha}w} - 
f_{s_{s_{\alpha}(\beta)} s_{\alpha} w})\, \equiv \,0 \ \mod \,\langle s_{\alpha}(\beta)^{2k_{\beta}+1} \rangle \,,
$$
we have $\,s_{\alpha}(f_{s_{\alpha}w} - f_{s_{\alpha} s_{\beta}w})\ \equiv \,0\
\mod \,\langle \beta^{2k_{\beta}+1} \rangle $. Hence, the numerator of \eqref{abdiff} 
is divisible by $ \beta^{2k_{\beta}+1} $ for all pairs $ \alpha, \beta \in \A $.
Now, if $ \beta \not= c \alpha $ for $ c \in \k^*$, this implies immediately that \eqref{abdiff} is 
also divisible by $ \beta^{2k_{\beta}+1} $, since $\alpha$ and $\beta$ are 
irreducible polynomials in $ \k[V] $. On the other hand, if $ \beta = c \alpha $ for some $ c \in \k^* $, then $ s_{\beta} = s_{\alpha} $ and $ k_{\alpha} = k_{\beta} $. 
It is easy to check that \eqref{abdiff} equals $ \frac{1}{\alpha}(1+s_{\alpha})(f_{w} - f_{s_{\alpha} w}) $ in $\k[\vreg]$ in that case. 
Since $\,f_w - f_{s_{\alpha} w} = \alpha^{2 k_{\alpha}+1} g_w\,$ for some $ g_w \in \k[V]$, we have
$ (1+s_{\alpha})(f_{w} - f_{s_{\alpha} w}) = \alpha^{2 k_{\alpha}+1} (g_w - s_{\alpha} g_w) $ and therefore $ \frac{1}{\alpha}(1+s_{\alpha})(f_{w} - f_{s_{\alpha} w})\,\equiv 0\ \mod \,\langle \alpha^{2k_{\alpha}+1} \rangle $. As $ \langle \alpha^{2k_{\alpha}+1} \rangle =  \langle \beta^{2k_{\beta}+1} \rangle $ for $ \beta = c \alpha $, we see that 
\eqref{abdiff} is divisible $ \beta^{2k_{\beta}+1} $ in this case as well. This proves the first claim of Theorem~\ref{Qfatodd}. The second follows from the proof of Lemma~\ref{bVV} (specialized to the case $ m = 2k+1$).
\end{proof}
Next, to compare the results of Theorem~\ref{Qfat} and Theorem~\ref{Qfatodd} we observe that the action of $ \NH(W) $ on $ \k[\vreg] \otimes \k[W] $ coincides with the differential action \eqref{diffact} if we regard $ \NH(W) $ as the sublagebra
of $ \D(\vreg) \rtimes W $ (via the natural inclusions $  \NH(W) \subset \k[\vreg]\rtimes W \subset 
\D(\vreg) \rtimes W $). Since
$ \nabla_{\alpha}\, \e = \frac{1}{\alpha} \e_{\alpha} \e = 0 $ in $ \D(\vreg) \rtimes W $ for all $ \alpha \in \A $, the PBW property of $\NH(W)$ implies 
\begin{equation}
\la{demspher}
 \NH(W)\e = \k[V]\e\quad \mbox{and} \quad   \e \,\NH(W) \e  \,=\,\e\, \k[V]^W\!\e \,\cong\, \k[V]^W .
\end{equation}
Then, combining the isomorphism \eqref{isoqq} (for $m=2k+1$) with \eqref{demspher}, we can rewrite the isomorphism \eqref{oddisoalg} of Proposition~\ref{bQodd} in the form
\begin{equation}
\la{NHmodiso}
\NH \e \otimes_{\e \NH \e} Q_k \,\cong\,\NH \e \otimes_{\e \NH \e} \e \bQ_{2k+1} \,\xrightarrow{\sim}\, \bQ_{2k+1}(W)\,
\end{equation}
where the second arrow is given by the canonical multiplication-action map. 
Now, it is easy to check that \eqref{NHmodiso} is actually an {\it isomorphism of
$\NH$-modules}, where the $\NH$-module structure is canonical on the left
and the one given by Theorem~\ref{Qfatodd} on the right. The isomorphism \eqref{NHmodiso} 
is obviously parallel to \eqref{Hkmodiso}: in fact,
there is a direct relation:
\begin{cor}
For all $ k \in \M(W) $, there is a commutative diagram of $\k[V]$-modules
\begin{equation}
\la{diaNH}
\begin{diagram}[small]
\NH \e \otimes_{\e \NH \e} Q_k(W) & \rTo^{\eqref{NHmodiso}\ } & \bQ_{2k+1}(W)\\
\dInto\ &  & \dInto\\
\bH_k \e \otimes_{e \bH_k \e} Q_k(W)  & \rTo^{\ \eqref{Hkmodiso}\ } & \bQ_{2k}(W)
\end{diagram}
\end{equation}
where the first map is induced by the inclusion 
$ \NH\e = \k[V] \e \subset \bH_k \e $ in $\D(\vreg) \rtimes W $  via  \eqref{demspher}.
\end{cor}
\begin{proof}
This follows from Proposition~\ref{bQodd} by a straightforward calculation.
\end{proof}

\begin{remark}
\la{HvsNH}
One can unify the two cases ($m=2k$ and $m=2k+1$) by changing the definition of the Dunkl embedding. To this end, in place of \eqref{du}, consider the family of operators
$$
T_{\xi,k}(t) \,:= \,t \cdot \partial_\xi-\sum_{\alpha \in \A} k_{\alpha} \,\frac{\alpha(\xi)}{\alpha}\,(1 - s_{\alpha})\ ,\quad \xi \in V\,,
$$
depending on the additional parameter $ t \in \k $. These operators still commute  in $ \D(\vreg) \rtimes W$ (for different $ \xi $'s) and are $W$-equivariant: together with $ \k[V]$ and $\k[W]$, they define --- for fixed $ (t,k) \in \k\times \M(W) $ --- a subalgebra $ \bH_{t,k}(W) $ in $ \D(\vreg) \rtimes W $, which is still called the Cherednik algebra.
For any $ t \not= 0 $, there is a natural isomorphism $ \bH_{t,k}(W) \cong \bH_{k'}(W) $, where 
$ k' = t^{-1}k$. On the other hand, for $ t=0 $, the structure of the algebra $\bH_{0,k}(W) $ and 
its representation theory are quite different from those of $ \bH_{k'}(W) $  (see \cite[Section 4]{EG02} or \cite[Section 5]{Rou05}). Now, a direct calculation (similar to that of \cite[Lemma 3.6]{BC11}) shows that the operators
$ T_{\xi,k}(t) $ preserve the subspace $ \bQ_{m}(W) \subset \k[\vreg] \otimes \k[W] $ under the differential action \eqref{diffact} if and only if the multiplicities $k$ and $m$ satisfy
\begin{equation}
\la{kmrel}
t \,m_{\alpha}\,=\,(1 + (-1)^{m_{\alpha}})\,k_{\alpha}\ ,\quad \forall\ \alpha \in \A\,.
\end{equation}
It follows from \eqref{kmrel} that if $ m_{\alpha} $ is even for {\it all} $\, \alpha \in \A $, then $\bQ_m(W) $ carries a differential action of $ \bH_{k'}(W) $ for $ k' = t^{-1}k = m/2$, which recovers precisely the result of Theorem~\ref{Qfat}. On the other hand, if there exists
an $\alpha' \in \A $ such that $ m_{\alpha'}$ is odd, then $ t=0 $. In this case, writing $\A^{\rm odd} \subseteq \A $ for the subset of all $\alpha$'s for which  $m_{\alpha} $ is odd, we see from \eqref{kmrel} that  $ k_{\alpha} $ may take any value when $ \alpha \in \A^{\rm odd}$ but must be $0$ otherwise. Thus, if $ \A^{\rm odd} \not= \varnothing $, then $\, \bQ_{m}(W) $ carries a differential action  $ \bH_{0, k}(W) \to \End[\bQ_m(W)] $, provided $ k_{\alpha} = 0 $ for all $ \alpha \in \A\setminus\A^{\rm odd}$. Under this last assumption on $k$, the assignment
$\, T_{\xi,k}(0) \,\mapsto\, - \sum_{\alpha \in A} \,k_{\alpha}\, \alpha(\xi)\, \nabla_{\alpha} \,$
extends to an algebra homomorphism $ \bH_{0,k}(W) \to \NH(W_{\rm odd}) $, where $\NH(W_{\rm odd}) $ is the nil-Hecke algebra associated to the parabolic subgroup $ W_{\rm odd} \subseteq W $ generated by the reflections $ s_{\alpha} $ with $ \alpha \in \A_{\rm odd} $. Moreover, it is easy to check that the above action  $\, \bH_{0,k}(W) \to \End[\bQ_m(W)] \,$ factors through this algebra homomorphism inducing
an action of $\NH(W_{\rm odd}) $ on $ \bQ_m(W) $. For $ m = 2k+1 $, we have $ \A_{\rm odd} = \A $ and $ W_{\rm odd} = W $, thus we get an action $ \NH(W) $ on $ \bQ_{m}(W) $ which 
coincides with the one described in Theorem~\ref{Qfatodd}.
\end{remark}

In Section~\ref{S7}, we will prove multiplicative analogs of these results for the Weyl groups.

\section{Stacks of quasi-invariants}\la{S3}
In this section, we extend the construction of  Section~\ref{S2} to the realm of Derived Algebraic Geometry (DAG).  To this end we replace the category $\Aff_{\k}$ of affine schemes by two ($\infty$-)subcategories of the ($\infty$-)category of derived stacks: the {\it derived affine schemes}, $ \dAff_\k $, 
and the {\it coaffine stacks}, $ \cAff_{\k} $,  which are, in a sense, dual to each other.  
While the derived affine schemes are an obvious derived extension of affine schemes, playing a role in many applications of DAG, the coaffine stacks seem to be less familiar objects. For reader's convenience, we include an Appendix~\ref{AB}, where we review the definition and basic properties of coaffine stacks and make precise their relation to rational homotopy theory.  The stacks of quasi-invariants constructed in Section~\ref{S3.2} 
will be particularly relevant for us as they will provide {\it rational}  models for our basic spaces --- the quasi-flag manifolds.\\

We shall keep the assumption that $\k $ is a field of characteristic $0$.

\subsection{Derived schemes of quasi-invariants} \la{S3.1} 
Over a field of characteristic zero,
the category $ \dAff_\k $  can be identified with  $ (\cdga_{\k}^{\leq 0})^{\rm op} $, 
the opposite  category of {\it negatively} graded commutative cochain dg $\k$-algebras  (see Appendix~\ref{AB}, Section~\ref{Sullivan}). The latter  is isomorphic to the category of positively graded commutative {\it chain} dg algebras and hence carries a natural (projective) model structure, where the weak equivalences are the quasi-isomorphisms of commutative dg algebras, and the fibrations are the dg algebra maps that are surjective in {\it strictly}\, negative cohomological degrees (see, e.g., \cite[Appendix B]{BKR13}).  
$ \dAff_{\k}$ is then equipped with the dual (injective) model structure. Note that the category $ \Comm_{\k} $ of ordinary commutative algebras
is naturally a subcategory of $\cdga_{\k}^{\leq 0}$, hence $\Aff_{\k}$ is naturally as a subcategory of $\dAff_{\k}$. 

Given an algebra $A \in  \cdga_{\k}^{\leq 0} $, we write
$ \RSpec(A) \in \dAff_\k $ for the derived affine scheme associted to $A$. Viewed as a derived  stack, $ \RSpec(A) $
represents the functor: $ \cdga_{\k}^{\leq 0} \to \sS $, $\, B \mapsto \Map(A, B) $, where `$ \Map $' stands for a simplicial mapping complex in $\cdga_{\k}^{\leq 0}$ ({\it cf}. \eqref{RSpec}). By definition,
the (homotopy) fibre products and coproducts  in $\dAff_{\k}$ are given by
\begin{equation}
\la{pullbh}
\RSpec(A)\times^h_{\RSpec(B)} \RSpec(C) \,=\, \RSpec(A \otimes^{\bL}_{B} C)\ ,
\end{equation}
\begin{equation}
\la{pushh}
\RSpec(A) \amalg^h_{\RSpec(B)} \RSpec(C) \,=\, \RSpec(A \times^{\bR}_{B} C)\ ,
\end{equation}
where $ \otimes_B^{\bL} $ and $ \times_{B}^{\bR} $ stand for the usual derived tensor and cartesian products in $ \cdga_{\k}^{\leq 0} $. 

As in Section~\ref{S2.1}, we fix a Coxeter group $W$ and its reflection representation $V$
over $\k$. Simplifying the notation, we will write\footnote{Note that this new notation is {\it not} in conflict with the one used in  Section~\ref{S2} (where we denoted  $ V = \Spec(\k[V])$, etc.), because all the commutative algebras $ \k[V]$, $\k[V]^W$, $\k[H_{\alpha}]$, $\ldots$ are free, and hence cofibrant objects in $ \cdga_{\k}^{\le 0}$.} 
$$ 
V = \RSpec(\k[V]) \ ,\quad
V/\!/W = \RSpec(\k[V]^W)\ , \quad
H_{\alpha} = \RSpec(\k[H_{\alpha}])\ ,\quad \ldots
$$
for the derived affine schemes associated to the algebras $ \k[V]$, $\,\k[V]^W$, $\,\k[H_{\alpha}]\,$, $\ldots$ viewed as objects in $\cdga_{\k}^{\leq 0}$. 

Now, constructing the derived schemes of quasi-invariants amounts to replacing basic categorical operations in $\Aff_k $ (used  in Section~\ref{S2.1} to construct the ordinary schemes of quasi-invariants) with their derived --- homotopy-invariant --- counterparts in $ \dAff_k$. Thus, for each $\alpha \in \A$, we first set
\begin{equation}
\la{FP1h}
V \times^h_{V/\!/W_{\alpha}} H_{\alpha} := \holim \{V \to V/\!/W_{\alpha} \leftarrow H_{\alpha}\}
\end{equation}
\begin{equation}
\la{V1h}
 V \ast^h_{V/\!/W_{\alpha}} H_{\alpha} := \hocolim \{V \leftarrow V \times^h_{V/\!/W_{\alpha}} H_{\alpha} \to H_{\alpha}\}
\end{equation}
\begin{equation}
\la{V2kh}
V^h_{k+1}(W_\alpha) := V^h_{k}(W_\alpha) \ast^h_{V/\!/W_{\alpha}} H_{\alpha}\ , \quad k\ge 0\,,
\end{equation}
which are derived extensions of formulas \eqref{FP1}, \eqref{V1} and \eqref{Vk}, respectively. Next, for a fixed multiplicity function $ k = (k_{\alpha})_{\alpha \in \A} \in \mathcal{M}(W)$, we form the diagram of derived schemes over the family \eqref{sw}  of reflection subgroups of $W$, extending 
\eqref{cvm}:
\begin{equation} 
\la{cvmh} 
\mathcal{V}^h_k(W)\,:\ \Sc(W) \to \dAff_{\k}\ ,
\quad W_{\alpha} \mapsto V^h_{k_{\alpha}}(W_{\alpha})\ ,\quad
W_0 \mapsto V^h_0(W_\alpha)\,,
\end{equation}
where $ V^h_0(W_\alpha) = V $ for all $ \alpha \in \A $. Finally, we define
\begin{equation} 
\la{defvmh} 
V^h_k(W)\,:=\, \hocolim_{\Sc(W)}[ \mathcal{V}^h_k(W) ] \ .
\end{equation}

Note that the homotopy pullbacks and  pushouts, defining the derived schemes
\eqref{FP1h}, \eqref{V1h} and  \eqref{cvmh}, are naturally equipped with $W_{\alpha}$-action, since 
their underlying diagrams take values in $W_\alpha$-schemes. Then, the same construction as 
in the proof of Proposition~\ref{cvmWfun} shows that \eqref{cvmh} is a $W$-functor, and hence, 
by Proposition~\ref{Whob}, the homotopy colimit
\eqref{defvmh} represents a derived affine $W$-scheme, i.e. an object of $\dAff_{\k}$ equipped with 
a $W$-action. We call $V_k^h(W) $ the {\it derived schemes of $W$-quasi-invariants  of multiplicity $k$}.\\

The next theorem, which is the main result of this section, asserts that the derived schemes of quasi-invariants defined above coincide with the classical (`underived') ones defined in Section~\ref{S2.1}.
\begin{theorem}
\la{Tdc}
For all $ k \in \M(W) $, there are $W$-equivariant $($weak$)$ equivalences in $ \dAff_{\k}$
$$
V^h_{k}(W)\,\simeq\, \RSpec[Q_k(W)] \,\simeq\, V_{k}(W)
$$
where $V_k(W) $ are the classical varieties of quasi-invariants viewed as objects in $\dAff_{\k}$.
\end{theorem}

This result can be reformulated by saying that the geometric 
construction  of schemes of quasi-invariants presented in Section~\ref{S2.1} is homotopy invariant (at least, in the context of affine DAG).

For the proof of Theorem~\ref{Tdc} we will need a simple lemma from abstract homotopical algebra.
Let $ \I $ denote the category with $ n+1 $ objects $\{0,1,2,\ldots, n\} $ and $n$ 
(non-identity) morphisms $\{i \to 0 \}$, one for each $ i =1,2,\ldots, n$.  A functor
on $ F: \I \to \M $ can thus be represented by a generalized pullback diagram of the form
\begin{equation}
\la{gpull}
\begin{tikzcd}[scale cd=0.8, column sep=small, row sep=small,
                 every arrow/.append style={shorten <=2pt, shorten >=2pt}]
	&& {F(1)} \\
	{F(n)} &&&& {F(2)} \\
	\vdots && {F(0)} && \vdots \\
	{F(i+1)} &&&& {F(i-1)} \\
	&& {F(i)}
	\arrow[from=1-3, to=3-3]
	\arrow[from=2-1, to=3-3]
	\arrow[from=2-5, to=3-3]
	\arrow[from=4-1, to=3-3]
	\arrow[from=4-5, to=3-3]
	\arrow[from=5-3, to=3-3]
\end{tikzcd}
%
\end{equation}
\begin{lemma}
\la{Lpull}
Let $\M$ be a model category, and let $\, F: \I \to \M $ be a diagram in $\M$ of shape \eqref{gpull}, where $F(0) \in \M $ is a fibrant object, and
all morphisms $ F(i) \to F(0) $ are fibrations in $\M$. Then the canonical map
$\,{\rm lim}_{\I}(F) \xrightarrow{\sim} \holim_{\I}(F) \,$ is a  weak equivalence in $ \M $.
If $F$ is a $W$-diagram, then both $ \lim_{\I}(F) $ and $ \holim_{\I}(F) $ are $W$-objects in $\M$, 
and the above equivalence is $W$-equivariant.
\end{lemma}
\begin{proof}
Note that $\I$ is a `very small' category in the sense of \cite{DS95} (see Definition~\ref{vscat}). Then, by  
\cite[Section 10.13]{DS95}, the category of $\I$-diagrams, $ \M^{\I} := \Fun(\I, \M) $, carries
two model structures: the projective, $\M^{\I}_{\rm proj}$, and the injective, $ \M^{\I}_{\rm inj}$, and the limits define a right Quillen functor $\,\lim_{\I}: \M^{\I}_{\rm inj} \to \M \,$
relative to the injective one. The $ \holim_{\I} $ is the right derived functor of 
$ \lim_{\I} $ and hence agrees with $ \lim_{\I}$ on fibrant objects in $ \M_{\rm inj}^{\I} $. It suffices to show the diagram \eqref{gpull} with $ F(0) $, being fibrant and all arrows $ F(i) \to F(0) $ being fibrations in $\M$ is a fibrant object in $\M_{\rm inj}^{\I}$.

For an arbitrary very small category $\I$,
the fibrant objects in $  \M^{\I}_{\rm inj} $ can be described as follows. For each object $ i \in \I $, let $ \I^{\circ}_{i/} $ denote the full subcategory of the coslice category $ \I_{i/}$ of $ \I $ under the object $i$ obtained by discarding the identity arrow $ \id_{i} $ from $\I_{i/}$. By restriction, any functor
$ F: \I \to \M $ defines a diagram  $ F|_{ \I^{\circ}_{i/}}: \I^{\circ}_{i/} \to \I \to \M $, 
$(i \to j) \mapsto F(j) $, whose limit we denote by $\delta_i(F) := \lim(F|_{\I^{\circ}_{i/}}) $. 
This comes with a natural map $ p_i: F(i) \to \delta_i(F) $ induced by the inclusion $  \I^{\circ}_{i/} \into \I_{i/} $.  Now, by \cite[10.13]{DS95}, a diagram 
$ F: \I \to \M $ is fibrant in $ \M^{\I}_{\rm inj} $ if the  maps 
$ p_i: F(i) \onto \delta_i(F) $ are fibrations in $ \M $ for all $ i \in \I $. 

For the category $\I$ of shape \eqref{gpull}, the category $  \I^{\circ}_{0/} = \varnothing $  is empty, while $\I^{\circ}_{i/} = \{\ast_i\} $  is a one-point category, the single object ``$ \ast_i $'' being the arrow $ \{i \to 0\} $ in $ \I $, for each $i=1,2,\ldots,n$. For a functor $F: \I \to \M $, we then have $ \delta_0(F) = \ast $ (the terminal object in $\M $)  and $ \delta_i(F) = F(0) $ for all $i= 1, 2,\ldots, n\,$. Hence, such $F$ is fibrant in $ \M^{\I}_{\rm inj}$ if $p_0: F(0) \to \ast $ and $ p_i: F(i) \to F(0) $ are fibrations in $ \M $. 

Finally, the last claim of the lemma follows from Lemma~\ref{vsdiagdual}
in Appendix~\ref{AA}.
\end{proof}
Now, we proceed with
\begin{proof}[Proof of Theorem~\ref{Tdc}]
We begin by computing the derived fibre product \eqref{FP1h}:
\begin{eqnarray*}
V \times^h_{V/\!/W_{\alpha}}\! H_{\alpha} &\simeq&  
\RSpec(\k[V] \otimes^{\bL}_{\k[V]^{W_\alpha}} \k[H_{\alpha}]) \\
&\simeq& \RSpec(\k[V] \otimes_{\k[V]^{W_\alpha}} \k[H_{\alpha}]) \\
&\simeq& \RSpec\,[\k[V]/(\alpha^2)]\,, 
\end{eqnarray*}
where the first equivalence is  formula \eqref{pullbh}, the second follows from the fact that
$\k[V]$ is a free module over $\k[V]^{W_{\alpha}}$, and the last one is given by \eqref{tensp}.
Next, using the above identification, we compute the derived relative join \eqref{V1h}:
\begin{eqnarray*}
V \ast^h_{V/\!/W_{\alpha}}\! H_{\alpha} &\simeq &  
\RSpec(\k[V] \times^{\bR}_{\k[V]/(\alpha^2)} \k[H_{\alpha}]) \\
&\simeq & \RSpec(\k[V] \times_{\k[V]/(\alpha^2)} \k[V]/(\alpha))\\
&\simeq & \RSpec \,[Q_{1}(W_{\alpha})]\,, 
\end{eqnarray*}
where the first equivalence is formula \eqref{pushh}, the second is due to the fact that
the canonical maps $\,\k[V] \to \k[V]/(\alpha^2) \leftarrow \k[V]/(\alpha)\,$ are fibrations in the model category $ \cdga_{\k}^{\leq 0} $ (since this last model category 
is fibrant, the homotopy fibre product is then equivalent to the ordinary one, by \cite[13.3.4]{Hir03}). The last isomorphism follows from \eqref{fibiso}.  The same argument works by induction for \eqref{V2kh},
giving the equivalences
\begin{equation}
\la{Vhka}
V^{h}_{k_{\alpha}}(W_{\alpha}) \simeq \RSpec\,[ Q_{k_{\alpha}}(W_{\alpha})]\ ,\quad \forall\, k_{\alpha} \ge 0 \,,
\end{equation}
where $ Q_{k_{\alpha}}(W_\alpha)$ are the ordinary rings of quasi-invariants of $ W_\alpha $ viewed as objects of $ \cdga_{\k}^{\le 0}$.

Now, in the final step, we take $\I = \Sc(W)^{\rm op} $ and apply Lemma~\ref{Lpull} to a diagram of dg algebras $ F: \Sc(W)^{\rm op} \to \cdga_{\k}^{\le 0} $ of the form \eqref{gpull}, with 
$  F(W_0) =  \k[V] $ and $ F(W_{\alpha}) = Q_{k_{\alpha}}(W_{\alpha}) $ and the morphisms
$ F(W_{\alpha}) \to F(W_0) $ in $ \cdga_{\k}^{\le 0} $ being the natural inclusions. Precisely, by \eqref{defvmh} and \eqref{Vhka}, we have
\begin{eqnarray*}
V^h_k(W) &\simeq & \hocolim_{\Sc(W)}\,[\,\RSpec(F)] \\ 
& \simeq & \RSpec\,[\,\holim_{\Sc(W)^{\rm op}}(F)]  \\
& \stackrel{(\ast)}{\simeq} &  \RSpec\,[\,\lim{}_{\Sc(W)^{\rm op}}(F)] \\
&\simeq & \RSpec \bigl(\bigcap_{Q_0(W)} Q_{k_{\alpha}}(W_{\alpha})\bigr)\\
& = & \RSpec\,[Q_k(W)]\,,
\end{eqnarray*}
where the key equivalence $(*)$ identifying $\, \holim_{\Sc(W)^{\rm op}}(F) \simeq \lim_{\Sc(W)^{\rm op}}(F)\,$ in $\cdga_{\k}^{\le 0} $ follows from Lemma~\ref{Lpull}. The conditions of this lemma 
are satisfied for the functor $F$ because every dg algebra in $\cdga_{\k}^{\le 0} $ is fibrant and any morphism $f: A \to B $, where $A$ and $B$ are {\it ordinary} algebras viewed as dg algebras with single component, is automatically a fibration in  $\cdga_{\k}^{\le 0} $. 

To complete the proof of Theorem~\ref{Tdc}, it remains to note that the last claim of Lemma \ref{Lpull} ensures that all the above equivalences are $W$-equivariant.
\end{proof}

We conclude this section with two remarks.

\begin{remark}
The assumption that $ \k $ has characteristic zero is not essential for Theorem~\ref{Tdc}. It is well known that, over fields of positive characteristic, the commutative dg algebras do not provide a correct model for the derived affine schemes: instead, one should work with {\it simplicial} commutative $\k$-algebras. Thus, if ${\rm char}(\k) > 0 $, we can identify  $\dAff_{\k} = (\sComm_{\k})^{\rm op} $, where $ \sComm_{\k} $ is equipped with standard Quillen model structure. In that model structure, the fibrations are represented by the morphisms
$ A_* \to B_* $ in $ \sComm_{\k} $ with the property that the induced maps $ A_* \to \pi_0(A_*) \times_{\pi_0(B_*)} B_* $ are degreewise surjective (see, e.g., \cite[Theorem~4.17]{GS07}). This implies that  any morphism $ A \to B $ of {\it discrete} (i.e., ordinary) commutative algebras in $ \sComm_{\k} $ is automatically a fibration. This allows us to verify the conditions of Lemma~\ref{Lpull} to conclude that the homotopy colimit  \eqref{defvmh} taken in  $ \sComm_{\k} $ is still equivalent to the ordinary colimit.
\end{remark}

\begin{remark}
The result of Theorem~\ref{Tdc} extends {\it mutatis mutandis} to the varieties of quasi-covariants $ \bV_{\! m}(W) $ defined in Section~\ref{S2.2}. 
\end{remark}

\subsection{Coaffine stacks of quasi-invariants}
\la{S3.2}
As explained in Appendix~\ref{AB}, the coaffine stacks are represented by {\it coconnective} dg algebras: that is, the positively graded cochain dg algebras $A$ with the property that $ H^0(A) \cong \k $. Such algebras do not form a model category (even if $\k$ is a field of characteristic zero); instead, we treat them as objects of an {\it $\infty$-category} which we denote by $ \cdga_\k^{>0} $. The $\infty$-category of coaffine stacks can then be identified with the opposite $\infty$-category: $ \, \cAff_{\k} := (\cdga_\k^{>0})^{\rm op} $ (see \eqref{DGcc}). Using the notation of \cite{DAGVIII}, we will write $ \cSpec(A) \in \cAff_{\k}$ for the coaffine stack represented by an algebra $ A \in \cdga_{\k}^{>0}$.

To perform our construction in $\cAff_{\k} $ we fix, as in Section~\ref{S2.1}, a reflection representation $V$ of a finite Coxeter group $W$ and consider its polynomial algebra $ \k[V] = \Sym_{\k}(V^*) $, but now we equip $ \k[V] $  with the {\it cohomological}  grading, setting
$$
\deg(x) = 2\ ,\quad \forall \,x \in V^*\,.
$$
Then we can regard $ \k[V] $ as an object of $\, \cdga_{\k}^{>0} $ (with trivial differential).
The algebras $\, \k[V]^{W}$, $\,\k[H_\alpha], \ldots $ inherit the grading from $\k[V]$ and hence
can also be viewed as objects in $\, \cdga_{\k}^{>0} $. This way we define the coaffine stacks:
$$ 
V^c := \cSpec\,\k[V] \ ,\quad
V^c/\!/W := \cSpec\,\k[V]^W\ , \quad
H^c_{\alpha} := \cSpec\,\k[H_{\alpha}]\ ,\quad \ldots
$$
Next, for each $\alpha \in \A$, we set
\begin{equation}
\la{FP1c}
V^c \times_{V^c/\!/W_{\alpha}} H^c_{\alpha} := \blim \{V^c \to V^c/\!/W_{\alpha} \leftarrow H^c_{\alpha}\}
\end{equation}
\begin{equation}
\la{V1c}
 V^c \ast_{V^c/\!/W_{\alpha}} H^c_{\alpha} := \bcolim \{V^c \leftarrow V^c \times_{V^c/\!/W_{\alpha}} H^c_{\alpha} \to H^c_{\alpha}\}
\end{equation}
\begin{equation}
\la{V2kc}
V^c_{k+1}(W_\alpha) := V^c_{k}(W_\alpha) \ast_{V^c/\!/W_{\alpha}} H^c_{\alpha}\ , \quad k\ge 0\,,
\end{equation}
where $\blim $ and $ \bcolim $ denote the pullback and the pushout in the $\infty$-category $ \cAff_{\k}$, 
i.e. the limit and colimit taken over the nerves of the categories indexing the pullback and pushout diagrams, respectively (see Section~\ref{B5}).

For a fixed multiplicity $ k = (k_{\alpha})_{\alpha \in \A} \in \mathcal{M}(W)$, we then consider the diagram of coaffine stacks over the family \eqref{sw}  of elementary reflection subgroups of $W$:
\begin{equation} 
\la{cvmc} 
\mathcal{V}^c_k(W)\,:\ \Sc(W) \to \cAff_{\k}\ ,
\quad W_{\alpha} \mapsto V^c_{k_{\alpha}}(W_{\alpha})\ ,\quad
W_0 \mapsto V^c\,,
\end{equation}
%
and define
\begin{equation} 
\la{defvmc} 
V^c_k(W)\,:=\, \bcolim_{{\mathcal N}\Sc(W)}[ \mathcal{V}^c_k(W) ] \ ,
\end{equation}
where $ \bcolim_{{\mathcal N}\Sc(W)} $ denotes the colimit taken in $ \cAff_\k$ over the nerve of the poset $\Sc(W)$. We call \eqref{defvmc} the {\it coaffine stack of $W$-quasi-invariants of multiplicity $k$.}

\begin{lemma}
\la{LcW}
The coaffine stacks $ V_k^c(W) $ are equipped with natural $W$-action induced from $ V^c $, and there are canonical $W$-equivariant maps
\begin{equation}\la{towerc}
V^c \to V_k^c(W) \to V_{k'}^c(W) \to  V^c/\!/W \ ,\quad \forall\, k' \ge k\,, 
\end{equation}
factoring the quotient map $\,V^c \to V^c/\!/W $ in $ \cAff_{\k}$.
\end{lemma}
\begin{proof}
Since $\cAff_{\k}$ is realized as a full $\infty$-subcategory of the $\infty$-category $ (\cdga_{\k}^{\geq 0})^{\rm op} $ with underlying simplicial model category $(\ccdga_{\k}^{\geq 0})^{\rm op}$ (see Section \ref{Sullivan}), we may identify (co)limits in $\cAff_{\k}$ with homotopy (co)limits in $\ccdga_{\k}^{\geq 0,{\rm op}}$ provided that the latter lie in $\cAff_{\k}$ (see Proposition \ref{PropA}). 
This is indeed the case by Theorem \ref{Tcc} and Theorem \ref{Hcomp} below. 

Note that for any $g \in W$, there is a map of pullback diagrams
$$
\begin{diagram}[small]
V^c & \rTo^{} & V^c/\!/W_{\alpha} & \lTo &  H^c_{\alpha}\\
 \dTo^{g} & & \dTo^{g} & & \dTo^{g}\\
V^c & \rTo^{} & V^c/\!/W_{g\cdot \alpha} & \lTo & H^c_{g \cdot \alpha}\\ 
\end{diagram}
$$
This yields maps $\vartheta_{g}:V_{k_\alpha}^c(W_\alpha) \to V^c_{k_\alpha}(W_{g\cdot \alpha})$ for all $k \geq 0$. The natural transformations $\vartheta_{g}: \mathcal{V}_k^c(W) \to \mathcal{V}_k^c(W) \circ \hat{g}$ make the functor $\mathcal{V}_k^c(W):\Sc(W) \to \ccdga_{\k}^{\geq0, {\rm op}}$ a $W$-functor. By Proposition \ref{Whob}, $V_k^c(W) \cong \hocolim_{\Sc(W)} \mathcal{V}^c_k(W)$ carries a $W$-action for all $k$.

Next, since $V^c_{k_\alpha+1}(W_\alpha) = V^c_{k_\alpha}(W_\alpha) \ast_{V^c/\!/W_\alpha} H^c_{\alpha}$, there is a (natural) map $V^c_{k_\alpha}(W_\alpha) \to V^c_{k_\alpha+1}(W_\alpha)$. Hence, for $k'>k$, there are natural transformations of $W$-functors $\mathcal{V}^c_{k}(W) \to \mathcal{V}^c_{k'}(W)$. Moreover, by \eqref{FP1c}, \eqref{V1c} and \eqref{V2kc}, for each $\alpha \in \A$ there are maps $V^c \to V^c_{k_{\alpha}} \to V^c/\!/W_{\alpha} \to V^c/\!/W$ for each multiplicity $k$. This yields a natural transformation of $V^c \to \mathcal{V}^c_k(W) \to V^c/\!/W$ of $W$-functors $\Sc(W) \to \ccdga^{\geq 0,{\rm op}}_{\k}$, where $V^c$ and $V^c/\!/W$ are constant $W$-functors with trivial $W$-action. For all $k'>k$ there are commutative diagrams of natural transformations (of $W$-functors)
$$\begin{diagram}[small]
  V^c & \rTo & \mathcal{V}^c_k(W) & \rTo & V^c/\!/W\\
  & \rdTo & \dTo & \ruTo &\\
  & & \mathcal{V}^c_{k'}(W) & &\\
\end{diagram}\ .$$
Hence, by Proposition \ref{Whob}, there are $W$-equivariant maps 
$$ V^c \to V^c_k(W) \to V^c_{k'}(W) \to V^c/\!/W$$
for all $k'>k$. This proves the desired lemma.
\end{proof}

Our next goal is to construct an explicit model for $ V_k^c(W) $. To this end, we will use  polynomial differential forms on algebraic simplices  (see Section~\ref{DGModel}). Recall that $\A $ is the set of reflection hyperplanes $H_{\alpha} $ of $W$ (equivalently, linear forms $ \alpha \in V^* $ such that $ H_{\alpha} = \Ker\,\alpha $), on which $W$ acts by permutations. Let $ \k\A$ denote the $\k$-linear permutation representation of $W$ with (unordered) basis $\A$, and let $\k[\A] := \Sym_{\k}[(\k\A)^*] $  be  the algebra of polynomial functions on $ \k\A $. Write $ \{x_{\alpha}\}_{\alpha \in \A } $ for the  generators of $ \k[\A]$ corresponding to the linear basis of $ (\k\A)^*$ dual to $\A  $, and let $ \dDel_{\k}(\A)$ denote the algebraic simplex on $\A $ with `barycentric coordinates' $ \{ x_{\alpha}\}$.
Thus, $ \dDel_{\k}(\A)$ is the hyperplane $\, \sum_{\alpha \in \A} x_{\alpha} = 1 \,$ in 
the $|\A|$-dimensional affine $\k$-space $ \bA^{|\A|}_{\k} = \Spec(\k[\A]) $. 
Now, define $\,\Omega_{\k}^*(\A) := \Omega^{\ast}_{\rm dR}[\dDel_{\k}(\A)] $ to be the 
algebraic de Rham complex on $ \dDel_{\k}(\A) $. Explicitly,
\begin{equation}
\la{OmegA} 
\Omega_{\k}^*(\A)\,\cong\, \frac{ \k[\{x_{\alpha}\}] \otimes \Lambda_{\k}(\{dx_{\alpha}\})}{(\sum_{\alpha \in \A} x_{\alpha} - 1,\,\sum_{\alpha \in \A} dx_{\alpha})} \,,
\end{equation}
where $\deg(x_{\alpha})=0$ and $\deg(dx_{\alpha})=1$, and the differential $d: \Omega_{\k}^*(\A) \to \Omega_{\k}^{*+1}(\A) $ is defined on generators by $ x_{\alpha} \mapsto  d x_{\alpha}$ and $d x_{\alpha} \mapsto 0$ for  $ \alpha \in \A $. Note that $\Omega_{\k}^*(\A)$ comes equipped with a natural $W$-action (induced by the permutation action of $W$ on $\A$) and admits
 $|\A|$ augmentation maps $ {\rm ev}_{\alpha}: \Omega_{\k}^*(\A) \to k \,$, one for each $ \alpha \in \A $, obtained by evaluation of differential forms at vertices $ x_{\alpha} = 1$ of the algebraic simplex $ \dDel_\k(\A) $. Given a dg algebra $B \in \cdga_{\k}^{\ge 0}$, we set
\begin{equation}\la{Bdel}
B^{\dDel_{\k}(\A)} := B \otimes \Omega_{\k}^*(\A)
\end{equation}
and put on \eqref{Bdel} the diagonal $W$-action whenever $B$ is a $W$-algebra. Extending the augmentation maps  ${\rm ev}_{\alpha}$ to $B^{\dDel_{\k}(\A)}$ and assembling them together, we define a $W$-equivariant dg algebra map:
\begin{equation}
\la{evmaps}
p:  \  B^{\dDel_{\k}(\A)} \onto B^{\A} 
\end{equation}
where
$ B^{\A}$ stands for the direct product of copies of $B$ indexed by the set $\A$ and equipped with the diagonal $W$-action  whenever $B$ is a $W$-algebra.

Next, we consider the subalgebras $ Q_{k_{\alpha}} = Q_{k_{\alpha}}(W_{\alpha}) $ of ordinary $W_{\alpha}$-quasi-invariants in $ \k[V]$ related to the hyperplanes in $\A$ but equipped with cohomological grading and the trivial differential (thus the $ Q_{k_{\alpha}}$ are viewed as objects in $\cdga_{\k}^{>0}$). Then, we define the dg algebra $Q_{k}^*(W)$ by pulling back the map \eqref{evmaps} (for $ B = \k[V]$) along the product of the natural inclusions $ Q_{k_{\alpha}} \into\, \k[V] $ taken over all $ \alpha \in \A$: 
thus, $Q_{k}^*(W)$ is defined by the fibre product
\begin{equation}\la{QcW}
Q_k^\ast(W)\,:=\, \biggl(\prod_{\alpha \in \A} Q_{k_{\alpha}} \biggr)
\bigtimes_{\k[V]^{\A}}  \k[V]^{\dDel_{\k}(\A)}
\end{equation}
or equivalently, by the cartesian square in $\cdga_{\k}^{>0}$:
\[
\begin{diagram}[small]
Q_k^\ast(W \SEpbk) & \rTo_{} &  \k[V]\otimes \Omega_{\k}^*(\A) \\
\dTo\ & & \dTo_{p}\\
\prod_{\alpha \in \A} Q_{k_{\alpha}}    & \rInto^{\quad {\rm can}} & \prod_{\alpha \in \A} \k[V]
\end{diagram}
\]
Note, since both maps in the above fibre product are $W$-equivariant, the dg algebra $Q_{k}^\ast(W)$ carries a natural $W$-action induced from the geometric action on $V$ and the permutation action on $\A$. 

\begin{theorem}
\la{Tcc}
For all $ k \in \M(W) $, there are $W$-equivariant equivalences in $\cAff_{\k}:$
\begin{equation}
\la{Weq}
V^c_k(W)\,\simeq\, \cSpec\,[\,Q_k^*(W)\,] 
\end{equation}
compatible with the $W$-equivariant maps \eqref{towerc}.
\end{theorem}
\begin{proof} 
We refer the reader to Appendix~\ref{AB}, Section~\ref{B5}, where we review definitions and explain how to compute limits and colimits in the $\infty$-category $\cAff_{\k} $ of coaffine stacks. 

The proof of Theorem~\ref{Tcc} is parallel to that of the proof of Theorem~\ref{Tdc}, except two things.
First, one should use Proposition~\ref{PropA} to convert  $\infty$-categorical (co)limits in $ \cAff_{\k} $ into classical {\it homotopy} (co)limits in a simplicial model structure on the category $ \ccdga_{\k}^{\ge 0} $ of non-negatively graded cochain dg algebras. Second, one should take into account the fact that, in this last model structure (originally constructed in \cite{BG76}) the homotopy limits behave differently from those in $ \cdga_{\k}^{\le 0} $ ({\it cf.} Example~\ref{EA2}).

For the first three steps of the proof --- computing \eqref{FP1c}, \eqref{V1c} and \eqref{V2kc} --- 
the different behavior of homotopy limits in $ \ccdga_{\k}^{\ge 0} $ and $  \cdga_{\k}^{\le 0} $
does not effect the result ({\it cf.} \eqref{Vhka}):
\begin{equation}
\la{Vcka}
V^{c}_{k_{\alpha}}(W_{\alpha}) \simeq \cSpec\,[ Q_{k_{\alpha}}(W_{\alpha})] \,,
\end{equation}
where $ Q_{k_{\alpha}}(W_\alpha)$ are equipped with cohomological grading (as subalgebras of $\k[V]$) and thus, are viewed as objects of the $\infty$-category $ \cdga_{\k}^{> 0}$.

Now, for the final step, we take $\I = \Sc(W)^{\rm op} $ and define the diagram of algebras $ F:\,\Sc(W)^{\rm op} \to \ccdga_{\k}^{\ge 0}\,$ 
of the form \eqref{gpull}:
\begin{equation}
\la{Fgpull}
\begin{tikzcd}[scale cd=0.8, column sep=tiny, row sep=small,
                 every arrow/.append style={shorten <=2pt, shorten >=2pt}]
	&& {Q_{k_{\alpha_1}}(W_{\alpha_1})} \\
	{Q_{k_{\alpha_r}}(W_{\alpha_r})} &&&& {Q_{k_{\alpha_2}}(W_{\alpha_2})} \\
	\vdots && {\k[V]} && \vdots \\
	{Q_{k_{\alpha_{i+1}}}(W_{\alpha_{i+1}})} &&&& {Q_{k_{\alpha_{i-1}}}(W_{\alpha_{i-1}})} \\
	&& {Q_{k_{\alpha_i}}(W_{\alpha_i})}
	\arrow[from=1-3, to=3-3]
	\arrow[from=2-1, to=3-3]
	\arrow[from=2-5, to=3-3]
	\arrow[from=4-1, to=3-3]
	\arrow[from=4-5, to=3-3]
	\arrow[from=5-3, to=3-3]
\end{tikzcd}
\end{equation}
with all maps being the natural inclusions. Then, in view of \eqref{Vcka}, by Proposition~\ref{PropA}, we can rewrite \eqref{defvmc} in the form
\begin{equation}
\la{hocvmc}
V^c_k(W) \,\simeq\,\cSpec\,[\,\holim_{\Sc(W)^{\rm op}}(F)\,]
\end{equation}
where the homotopy limit is taken in the Bousfield-Gugenheim model category $ \ccdga_{\k}^{\ge 0} $.
Unlike in $\cdga_{\k}^{\le 0} $, the conditions of Lemma~\ref{Lpull} do not hold for the diagram $F$ in  
$ \ccdga_{\k}^{\ge 0} $: hence the homotopy limit $\holim_{\I}(F)$ is not equivalent to $ \lim_{\I}(F) $. However, as shown in the proof of Lemma~\ref{Lpull}, $\holim_{\I}(F)$  can be identified with the ordinary limit, provided we replace $F$ by an equivalent fibrant diagram $ \tilde{F}$ over $\I$, in which each arrow is represented by a fibration in $ \ccdga_{\k}^{\ge 0} $. Thus, for every $ \alpha \in \A $, we choose a factorization $\, Q_{k_{\alpha}}(W_{\alpha})\stackrel{\sim}{\to} R_{\alpha} \onto \k[V] \,$, where the first arrow is an equivalence in $ \ccdga_{\k}^{\ge 0} $ and the second is a fibration, and define $ \tilde{F}: \I \to \ccdga_{\k}^{\ge 0} $ by $ \tilde{F}(W_{\alpha} \to W_0) = (R_{\alpha} \onto \k[V]) $. Then, we have
\begin{eqnarray*}
\holim_{\I}(F) 
&\simeq& \lim{\!}_{\I}(\tilde{F})\, \simeq \, \lim\,\bigl(\,\prod_{\alpha} R_{\alpha} \onto \k[V]^{\A}  \leftarrow \k[V]\bigr) \\
&\simeq& 
\holim\,\bigl(\prod_{\alpha} Q_{k_{\alpha}} \to \k[V]^{\A}  \leftarrow \k[V]\bigr)\\
&\simeq& 
\holim\,\bigl(\prod_{\alpha} Q_{k_{\alpha}} \to \k[V]^{\A} \twoheadleftarrow \k[V]^{\dDel_{\k}(\A)}\bigr)\\
&\simeq& 
\lim\,\bigl(\prod_{\alpha} Q_{k_{\alpha}} \,\onto \,\k[V]^{\A} 
\,\twoheadleftarrow \,\k[V]^{\dDel_{\k}(\A)}\bigr) \,,
\end{eqnarray*}
where $ \k[V]^{\dDel_{\k}(\A)} $ is defined by \eqref{Bdel}. It follows from \eqref{QcW} that $V_k^c(W) \simeq \cSpec\,[Q_k^*(W)]\,$. 

Finally, note that the fibrant replacement $\tilde{F}$ may be chosen to be a $W$-diagram (for example, if we take $\,R_{\alpha}=\{p(t)+q(t)dt \,\in\, \k[V] \otimes \k[t,dt]\,:\, p(0) \in Q_{k_\alpha}(W_\alpha)\}$). Then, applying Lemma~\ref{Lpull} to $\tilde{F}$ shows that the equivalence $V_k^c(W) \simeq \cSpec\,[Q_k^*(W)]\,$ is indeed $W$-equivariant.
\end{proof}

Note that, if $k=0$, then, by \eqref{QcW}, $\, Q^*_0(W) \cong  \k[V]^{\dDel_{\k}(\A)} = \k[V] \otimes \Omega^*_{\k}(\A) \,$, which shows that $ Q^*_0(W) \simeq \k[V]$ for any  $W$ (as the de Rham algebra $\Omega^*_{\k}(\A)$ is acyclic, the unit map $ \k \to \Omega^*_{\k}(\A)$ being a quasi-isomorphism). Furthermore, if $ W $ is of rank one, then $ Q_k^*(W) \cong Q_k(W) $ for all $ k \in \Z_{\ge 0}$. In general, however, $ Q_k^*(W) $ is {\it not} quasi-isomorphic to $Q_k(W)$.  Before proving this for arbitrary $W$, we examine a simplest example.
\begin{example}
\la{Er2}
Let $ W = W_1 \times W_2 $ be the product of cyclic groups $ W_1 \cong W_2 \cong \Z/2 $, 
acting by reflections in a vector space $V$. Write $H_{1}$ and $H_2$ for the corresponding
reflection hyperplanes in $V$. Fix a multiplicity $ k = (k_1, k_2 ) \in \Z_{\ge 0}^2 $ and 
consider the subalgebras $Q_{k_1}(W_1)$ and $ Q_{k_2}(W_2)$ of quasi-invariants of $W_1$ and 
$W_2$ in $ \k[V]$ of multiplicities $ k_1 $ and $ k_2 $, respectively. In this case, the 
de Rham algebra \eqref{OmegA} can be identified as $\, \Omega_{\k}^*(\A) \cong  
\Omega^*_{\rm dR}(\Delta_{\k}^1) \cong \Lambda_{\k}(t, dt) \cong \k[t] \oplus \k[t] dt \,$, 
with $W$ acting trivially on $ \Lambda_{\k}(t, dt)$. The map  \eqref{evmaps} 
is given explicitly by
$$
p:\ \k[V] \otimes \Lambda_{\k}(t,dt) \to \k[V] \times \k[V]\ , \quad 
f \otimes (h(t) + g(t) dt) \mapsto (f \otimes h(0), \, f \otimes h(1))
$$
where $ \deg(t) = 0 $ and $ \deg(dt) = 1 $. The dg algebra \eqref{QcW} can then be written as 
\begin{eqnarray*}
Q_k^*(W) &\cong & 
Q_{k_1}\times_{\k[V]} (\k[V] \otimes \Lambda_{\k}(t, dt)) \times_{\k[V]} Q_{k_2}\\*[1ex]
& \cong & 
( Q_{k_1}  \times_{\k[V]} {\k[V, t]} \times_{\k[V]} Q_{k_2})\,\oplus\, k[V,t] dt
\end{eqnarray*}
where $ \k[V,t] := \k[V] \otimes \k[t] $ and the double fiber product in the first summand 
is defined via the evaluation maps $ {\rm ev}_0, \,{\rm ev}_1: \k[V,t] \to \k[V] $ sending
$t \mapsto 0 $ and $ t \mapsto 1$, respectively. Thus, the cochains  in $ Q_k^*(W) $ of
{\it even} cohomological degree  are represented  by 
\begin{equation*}\la{Qev}
Q_k^{\rm ev}(W) = \{(f_1, h(t), f_2) \in Q_{k_1} \times \k[V,t] \times Q_{k_2}\,|\  
h(0) = f_1  \ \mbox{and}\  h(1) = f_2 \ \mbox{in}\ \k[V]\}
\end{equation*}
while $ Q_k^{\rm odd}(W) \cong \k[V,t] dt $. The differential on $ Q_k^c(W)$ is determined  
by the map 
\begin{equation}\la{dQc}
d: Q_k^{\rm ev}(W) \to Q_k^{\rm odd}(W)\ ,\quad
(f_1, h(t), f_2) \mapsto h'(t) dt\,,
\end{equation}
where $h'(t) $ denotes the  derivative of a polynomial function with respect to $t$. From this identification,
we see immediately that 
\begin{eqnarray}\la{Hev}
H^{\rm ev}[\,Q_k^*(W)\,] &\cong & \{(f_1, h(t), f_2)\,|\ h'(t) = 0,  
\,h(0) = f_1 \,,\,  h(1) = f_2 \ \mbox{in}\ \k[V]\}\\
&\cong& 
\{(f_1, h, f_2)\,|\ f_1 = h = f_2 \ \mbox{in}\ \k[V]\} \nonumber \\
&\cong& Q_{k_1} \,\cap\, Q_{k_2}\nonumber\\
& \cong & Q_k(W)\nonumber
\end{eqnarray}
 To compute the odd-dimensional cohomology we consider the integration map
\begin{equation}
\la{intmap}
Q^{\rm odd}_k(W) = \k[V,t] dt \,\to\, \k[V] \ , \quad g(t)dt \mapsto \int_0^1 g(t) dt
\end{equation}
The composition of \eqref{intmap} with differential map \eqref{dQc} is given  by
\begin{equation*}
(f_1, p(t), f_2) \mapsto \int_0^1 p'(t) dt = p(1) - p(0) = f_1 - f_2 \,\in \, Q_{k_1} + Q_{k_2}
\end{equation*}
Hence, up to shift in degree, \eqref{intmap} induces a well-defined map
\begin{equation}\la{Hodd}
H^{\rm odd}[\,Q_k^*(W)\,]\,\xrightarrow{\sim}\, \k[V]/(Q_{k_1} + Q_{k_2})
\end{equation}
which is easily seen to be an isomorphism. Thus, the even-dimensional part of the cohomology of $ Q^*_k(W) $ is isomorphic to  $Q_k(W)$, the algebra of ordinary polynomial $k$-quasi-invariants of $W$, while the odd-dimensional part is given by \eqref{Hodd} and does not vanish unless $ k_1 =0 $ or $k_2=0$.
\end{example}
We now show that these properties hold in general. 
\begin{theorem}
\la{Hcomp}
For all $ k \in \M(W)$, there is a natural isomorphism of graded algebras 
$$ 
H^{\rm ev}[\,Q_k^*(W)\,] \cong Q_k(W) 
$$
where $ Q_k(W) $ is equipped with cohomological grading. The `odd' cohomology $ H^{\rm odd}[\,Q_k^*(W)\,] $
is given explicitly by formula \eqref{HQodd}, which shows that $ H^{\rm odd}[\,Q_k^*(W)\,] \not= 0 $ in general.
\end{theorem}
\begin{proof}
First, we note that 
the forgetful functor $ \cdga_{\k}^{>0} \to \Com_{\k} $ from dg $\k$-algebras to cochain complexes of $\k$-vector spaces preserves cartesian squares. Hence, as a cochain complex, $ Q^*_k(W) $ can be represented by the cone of a moprhism  in $ \Com_{\k} $:
\begin{equation}\la{dimap}
\bigl(\prod_{\alpha \in \A} Q_{k_{\alpha}}\bigr) \oplus \, \k[V]^{\dDel_{\k}(\A)} \ \xrightarrow{i\,-\,p}\ \k[V]^{\A} 
\end{equation}
given by the difference of the two maps $i$ and $p$ in \eqref{QcW}. 
Next, observe that the canonical (diagonal) map $\, \nabla:\, \k[V] \to \k[V]^{\A}\,$ factors  in $ \Com_{\k} $ 
as
\begin{equation}\la{nabla}
\nabla:\,\k[V] \xrightarrow{\id \otimes 1} \, \k[V]^{\dDel_{\k}(\A)} \, \xrightarrow{p} \, \k[V]^{\A} \,,
\end{equation}
the first arrow being a quasi-isomorphism. Hence, in the derived category of cochain complexes, we can 
replace \eqref{dimap} with the simpler map
\begin{equation}\la{dinmap}
\bigl(\prod_{\alpha \in \A} Q_{k_{\alpha}}\bigr) \oplus \k[V]  \ \xrightarrow{i\,-\,\nabla}\ \k[V]^{\A} 
\ ,\quad ((f_{\alpha})_{\alpha \in \A}\ , \, f) \,\mapsto\, (f_{\alpha} - f)_{\alpha \in \A} \,,
\end{equation}
which we denote by $\delta := i-\nabla$. Thus, we have
$\,H^*[Q^*_k(W)] \cong H^*[{\rm cone}(\delta)] \,$,
and there is a long exact cohomology sequence associated to (the cone of) the map \eqref{dinmap}:
\begin{equation*}\la{longex}
\ldots \,\to\, H^n[Q_k^*(W)] \,\to\, (\prod_{\alpha \in \A} Q^n_{k_{\alpha}})  
\oplus \k[V]^n \,\xrightarrow{\delta^n}\, (\k[V]^{\A})^{n}  \,\to\, H^{n+1}[Q_k^*(W)] \,\to\, \ldots
\end{equation*}
From the above exact sequence, it follows immediately that 
\begin{eqnarray}
H^{\rm ev}[Q^*_k(W)] & \cong &
\Ker(\delta)\,\cong\, \bigcap_{\alpha \in \A} Q_{k_{\alpha}} \,\cong\, Q_k(W) \la{HQev}   \\
H^{\rm odd}[Q^*_k(W)] & \cong & {\rm Coker}(\delta)[-1]\,\cong\, (\k[V]^{\A}\!/S_k(\A))[-1] \la{HQodd}
\end{eqnarray}
where $S_k(\A) $ is the linear subspace of $\, \k[V]^{\A} = \bigoplus_{\alpha \in \A} \k[V] \,$ spanned by the elements ($ \A$-tuples) of the form $\,(f_{\alpha} - f)_{\alpha \in \A}\, $ with 
$\, f_{\alpha} \in Q_{k_{\alpha}} \,$ and $\, f \in \k[V]\,$.

Note that \eqref{HQev} is actually an isomorphism of graded {\it algebras}
induced by a canonical dg algebra map $\, Q_k(W) \to Q^*_k(W) \,$. To construct this last map
we realize $ Q_k(W) $ as the fibre product similar to \eqref{QcW}:
\begin{equation}\la{QW}
Q_k(W)\,=\, \biggl(\prod_{\alpha \in \A} Q_{k_{\alpha}} \biggr)
\bigtimes_{\k[V]^{\A}}  \k[V]
%
%
\end{equation}
where the map $p$ is replaced by the diagonal map $\nabla $. Since $\nabla$ factors through $p$ as a morphism in $ \cdga_{\k}$ (see \eqref{nabla}), we get a map of pullback diagrams, and hence, the canonical map $\, Q_k(W) \to Q^*_k(W) \,$.
The isomorphism \eqref{HQodd} is then an isomorphism of 
graded {\it  $ Q_k(W) $-modules}, where the action of $Q_{k}(W) $ is induced by the action on $ \k[A]^{\A}$ via the diagonal embedding 
$ Q_{k}(W) \subseteq \k[V] \into \k[V]^{\A} $. The subspace $ S_k(\A) \subseteq \k[V]^{\A}$
is stable under this action of $Q_k(W) $, since $\, Q_k(W) \subseteq Q_{k_{\alpha}} \,$ for all $ \alpha \in \A$; hence the  quotient $ \k[V]^{\A}\!/S_k(\A) $ is indeed an $ Q_k(W)$-module. With a little more effort, one can prove that the multiplication vanishes on $ H^{\rm odd}[Q^*_k(W)] $: the full cohomology algebra of $ Q_K^*(W) $ is thus given by the semidirect product:
$\,
H^*[Q^*_k(W)]\,\cong\, H^{\rm ev}[Q^*_k(W)] \ltimes H^{\rm odd}[Q^*_k(W)] 
\,$.
\end{proof}
Note that the canonical dg algebra map $\, Q_k(W) \to Q_k^*(W) \,$ constructed in the proof
of Theorem~\ref{Hcomp} induces a map of coaffine stacks $\,V^c_k(W) \to \cSpec[Q_k(W)]\,$.
This  map is similar to the map of derived affine schemes $\, V_k^h(W) \to \RSpec[Q_k(W)] \,$ from Section~\ref{S3.1}, that is, according to Theorem~\ref{Tdc}, an equivalence in $\dAff_\k $ for all $ k $. In contrast,  Theorem~\ref{Hcomp} shows that $\,V^c_k(W) \to \cSpec[Q_k(W)]\,$ is {\it not} an equivalence in $ \cAff_{\k}$ in general (when $ k \not=0 $).

\subsection{Coaffine stacks of quasi-covariants}
\la{S3.3}
Parallel to \eqref{defvmc}, for any $ m \in \M(W)$, we can define the {\it coaffine stacks $ \bV_m^c(W)$ of quasi-covariants}:
\begin{equation} 
\la{bdefvmc} 
\bV^c_{\!m}(W)\,:=\, \bcolim_{{\mathcal N}\Sc(W)}[ \boldsymbol{\mathscr{V}}^c_m(W) ] \,,
\end{equation}
where the functor $\,\boldsymbol{\mathscr{V}}^c_m(W)\,:\ \Sc(W) \to \cAff_{\k}\,$ is given by
\begin{eqnarray*}  
W_0 &\mapsto& \bV^c := V^c \times W\ ,\\*[1ex]
 W_{\alpha} &\mapsto& \bV^c_{m_\alpha}(W_\alpha) \,:= 
\bV^c\ast_{\bV^c\!/W_{\alpha}}(\bH^c_{\alpha}/W_{\alpha}) \ast_{\bV^c\!/W_{\alpha}}\,
\stackrel{m_\alpha}{\ldots}\, \ast_{\bV^c\!/W_{\alpha}} (\bH^c_{\alpha}/W_{\alpha})\,.
\end{eqnarray*}
with relative joins being defined in $\cAff_{\k}$. As an exercise, we leave to the reader to prove that the analogs of Theorem~\ref{Tcc} and Theorem~\ref{Hcomp} hold for $\bV^c_m(W) $:
in particular, we have
\begin{equation} 
\la{tbWck} 
\bV^c_{\!m}(W)\,\simeq\,{\rm cSpec}[\bQ^\ast_m(W)]  
\end{equation}
and
\begin{equation} 
\la{tbWQ} 
H^{\rm ev}[\bQ^\ast_m(W)] \cong \bQ_m(W) \,, 
\end{equation}
where 
\begin{equation}
\la{bbQmW}
\bQ_m^\ast(W)\,:=\, \biggl(\prod_{\alpha \in \A} \bQ_{m_{\alpha}} \biggr)
\bigtimes_{\k[\bV]^{\A}}  \k[\bV]^{\dDel_{\k}(\A)}\,.
\end{equation}
Again, the dg algebra $ \bQ_m^\ast(W) $ has nonvanishing odd cohomology: 
in fact, $\, H^{\rm odd}[\,\bQ_m^*(W)\,] \not= 0 \,$ even for $ m = 0$.

Our final goal is to define coaffine stacks that will serve as rational models for our basic spaces: the quasi-flag manifolds and their homotopy $T$-quotients.
Recall that, by Lemma~\ref{LcW}, the stacks $ V_k^c(W) $ come with natural $W$-equivariant maps $ p_k: V^c_k(W) \to V^c/\!/W $. Using these maps, for each $ k \in \M(W)$, we define 
\begin{eqnarray}
 F_k^c(W) &:=&  V_k^c(W) \times_{V^c/\!/W} \cSpec(\k)\ , \la{Fck}\\*[1ex]
\bU^c_k(W) &:=& V_k^c(W) \times_{V^c/\!/W} V^c\ ,\la{bvck}
\end{eqnarray}
where $\cSpec(\k) \to V^c/\!/W $ represents the canonical basepoint in $V^c/\!/W$ (the $W$-orbit of the origin), and the fibre products are taken in $ \cAff_{\k}$.

Equivalently, the stack \eqref{Fck} can be described as follows. Recall the `local' coaffine stacks $ V^c_{k_{\alpha}} = V^c_{k_{\alpha}}(W_\alpha) $ defined, for each $\alpha \in \A $, by \eqref{Vcka}. 
By Lemma~\ref{LcW}, the canonical projection $ p: V^c \to V^c/\!/W $ factors through  $ p_{k_{\alpha}}: V^c_{k_{\alpha}} \to V^c/\!/W $. We let  $F^c $ and $ F^c_{k_{\alpha}}$
the fibres of these maps:
$$ 
F^c := V^c \times_{V^c/\!/W} \cSpec(\k)
\ ,\quad 
F^c_{k_{\alpha}} := V^c_{k_{\alpha}} \times_{V^c/\!/W} \cSpec(\k) \,, 
$$
and define a diagram $ \mathcal{F}^c_k(W): \Sc(W) \to \cAff_{\k} $  by $\, W_0  \mapsto F^c $ and $ W_{\alpha} \mapsto F^c_{k_{\alpha}}$. Then, by a classical theorem of Puppe \cite{Pu74}, we can  deduce from \eqref{defvmc} that
\begin{equation} 
\la{deffmc} 
F^c_k(W)\,\simeq\, \bcolim_{{\mathcal N}\Sc(W)}[ \mathcal{F}^c_k(W) ] \,.
\end{equation}
Similarly, we can decompose \eqref{bvck} as
\begin{equation} 
\la{defbumc} 
\bU^c_k(W)\,\simeq\, \bcolim_{{\mathcal N}\Sc(W)}[ \boldsymbol{\mathscr{U}}^c_k(W) ] \ ,
\end{equation}
where the diagram $ \boldsymbol{\mathscr{U}}^c_k(W):\, \Sc(W) \to \cAff_{\k} $ is defined  
by 
$$
W_0  \mapsto V^c \times_{V^c/\!/W} V^c \ ,\quad 
W_{\alpha} \mapsto V^c_{k_{\alpha}} \times_{V^c/\!/W} V^c \,.
$$

The decompositions \eqref{deffmc} and \eqref{defbumc} allow us to give explicit presentations for the coaffine stacks $F_k^c(W)$ and $\bV^c_k(W)$ parallel to \eqref{QcW}. 
Let $\, \H_W(V) := \k[V]/(\k[V]_{+}^{W})\,$ denote the graded algebra of
polynomial $W$-coinvariants over $\k$, and let $\,\H_{k_{\alpha}}(W_\alpha) := 
Q_{k_{\alpha}}(W_{\alpha})/(\k[V]_{+}^{W})\,$ be the graded algebra of $k_{\alpha}$-quasi-coinvariants for each $ \alpha \in \A $. Note that there are canonical embeddings $ \H_{k_{\alpha}}(W_{\alpha}) \into \H_W(V)$ in $ \cdga^{>0}_{\k} \,$ induced by $ Q_{k_{\alpha}}(W_{\alpha}) \into \k[V] $. 
We define the cochain dg algebra $ \H^*_k(W) $ by taking the fibre product of $W$-equivariant maps in $ \cdga^{>0}_{\k} $:
\begin{equation}\la{HcW}
\H_k^\ast(W)\,:=\, \biggl(\prod_{\alpha \in \A} \H_{k_{\alpha}}\biggr)
\bigtimes_{\H_W(V)^{\A}}  \H_W(V)^{\dDel_{\k}(\A)}\,,
%
%
\end{equation}
where the dg algebra $\H_W(V)^{\dDel_{\k}(\A)}$ and the projection $p$ are defined by
\eqref{Bdel} and \eqref{evmaps}. 
Similarly, we set
\begin{equation}
\la{bPkW}
\bbP_k^\ast(W)\,:=\, \biggl(\prod_{\alpha \in \A} \bbP_{k_{\alpha}} \biggr)
\bigtimes_{\bbP_W(V)^{\A}}  \bbP_W(V)^{\dDel_{\k}(\A)}\,,
\end{equation}
where $\, \bbP_W(V) := \k[V] \otimes_{\k[V]^W} \k[V]\,$ and $\,{\bbP}_{k_{\alpha}} := Q_{k_{\alpha}}(W_{\alpha}) \otimes_{\k[V]^W} \k[V]$ in $ \cdga^{>0}_{\k} $. 

\vspace{1ex}

Now, as a consequence of Theorem~\ref{Tcc}, we get
\begin{cor}
\la{CHcomp}
For all $ k \in \M(W) $, there are $W$-equivariant equivalences in $\cAff_{\k}:$
\begin{eqnarray}
F^c_k(W) &\simeq& \cSpec\,[\,\H_k^*(W)\,]\,, \la{FckW}\\*[1ex]
{\bU}^c_k(W) &\simeq& {\cSpec}[\,\bbP^\ast_k(W)\,]\la{bWck}
\end{eqnarray}
compatible with \eqref{Weq} 
via the canonical maps $ F_k^c(W) \to V^c_k(W) $ and ${\bU}^c_k(W) \to V^c_k(W)$.
\end{cor}
Note that, by definition, $\, \bbP_k^\ast(W) \cong Q_k^*(W) \otimes_{\k[V]^W} \k[V] \,$,
hence
\begin{equation}
\la{evcohP_k}
H^{\rm ev}[\bbP_k^\ast(W)]\,\cong\, Q_k(W) \otimes_{\k[W]^W} \k[V]\,\cong\,\bQ_{2k+1}(W)\,,
\end{equation}
where the last isomorphism is the result of Proposition~\ref{bQodd}, see \eqref{oddisoalg}. Similarly, we have
\begin{equation}
\la{evcohH_k}
H^{\rm ev}[\H_k^\ast(W)]\,\cong\, Q_k(W) \otimes_{\k[W]^W} \k\,\cong\,
Q_{k}(W)/\langle \k[V]^W_+\rangle \,=\, \H_k(W)\,.
\end{equation}

\begin{remark} \la{tildebv}
Although  $\,H^{\rm ev}[\bbP_k^\ast(W)] \cong H^{\rm ev}[\bQ_{2k+1}^\ast(W)] $ for all $ k\in \M(W) $, the dg algebras $\bbP_k^\ast(W) $ and $ \bQ^*_{2k+1}(W) $ are not quasi-isomorphic:
in fact, $ H^{\rm odd}[\bbP_0^\ast(W)] \cong H^{\rm odd}(\k[V]\otimes_{\k[V]^W} \k[V]) = 0 $, while $ H^{\rm odd}[\bQ_{1}(W)] \not=0$. Thus, the coaffine stacks $ \bU_k^c(W) $ and $ \bV^c_{2k+1}(W) $ are not equivalent in $ \cAff_{\k}$ in general.
\end{remark}

\section{Homology decomposition of a flag manifold over the moment graph}
\la{S4}
In this section, we construct a natural (co)homology decomposition of the flag manifold $G/T$
(see Theorem~\ref{hDec}) that will serve as a motivation for our general definition of quasi-flag manifolds.

\subsection{Moment graphs and categories}\la{S4.1}
We begin by recalling the combinatorial definition of a moment graph (see, e.g., \cite{BM01, Fie12}). 
Let $ \Lambda $ be a lattice (i.e., a free abelian group of finite rank). An (unordered) {\it moment graph} over $\Lambda$ is a  simple unoriented graph $ \Gamma = (V, E)$ with vertex set $ V = V_{\Gamma}$ and edge set $ E = E_{\Gamma} $ given together with a function $\, \lambda: E \to 
\Lambda\!\setminus\!\{0\} \,$ called a labeling of $\Gamma$.
By forgetting the labeling, we may (and often will) think of $ \Gamma $ combinatorially as a one-dimensional simplicial complex with vertex set $V_{\Gamma}$ and the simplex set $ \Sigma_{\Gamma} = V_{\Gamma} \cup E_{\Gamma}\,$, and we will write $ \Cc(\Gamma) $ for the {\it opposite} face category
of $ \Gamma$.  Thus, by definition, the objects of $ \Cc(\Gamma) $ are the simplices $ \Sigma_{\Gamma} $, i.e.
the vertices and the edges of $\Gamma$, while the morphisms correspond to the reverse inclusions of faces, i.e. $ \{v \leftarrow e \to v'\} $, where $ e \in E_{\Gamma} $ and $ v,v' \in V_{\Gamma} $ are the adjacent vertices to $e$.

Our main example of a moment graph will be the {\it Bruhat  graph} $ \Gamma = \Gamma(\cR_W) $  of a Weyl group $W$ (see \cite{Car94, Fie12}).  Let $G$ be a compact connected Lie group with a maximal torus $T \subset G $
and associated Weyl group $ W = N_G(T)/T $. Let $ \R = \R_W$ be a root system of $W$ with  a choice of positive roots $ \R_+ \subset \R $. The Bruhat graph $ \Gamma = \Gamma(\cR_W) $ is a moment graph defined over the weight lattice $ \Lambda(T) = \Hom(T, U(1))$; its vertex set $ V_\Gamma = W $  and the edge set is $ E_\Gamma = \{e(s_\alpha, w)\, :\, \alpha \in \cR_+ \,,\, w \in W\} $, where the edge $\, e(s_{\alpha}, w) \,$ connects the vertices $w$ and $ s_\alpha w $ for each reflection $ s_\alpha \in W $  (see Figure~\ref{GammaS3}). Note that, since  $ \Gamma $ is an unoriented graph, we have the identifications $\,e(s_\alpha, w) = e(s_\alpha, s_\alpha w)\,$ for all $w \in W$ and all $\alpha \in \R_+$.
The graph $ \Gamma $ is labeled by the roots $ \alpha$ which we identify with characters of $ T $: thus, we define the function $ \lambda: E \to \Lambda(T) $ by $ e(s_\alpha, w) \mapsto \alpha $.


The Bruhat moment graph $ \Gamma $ has a simple geometric interpretation in terms of
the flag manifold $G/T$ (see \cite[Theorem F]{Car94}, and also \cite{GKM98, BM01}).
Let $ \bG = G_{\c} $ be  a complex connected reductive algebraic group containing $G$ as its maximal compact subgroup, and let $ \bB \subset \bG $ be its Borel subgroup containing a maximal complex torus $ \bT = (\c^*)^r $, with $T = G \cap \bT = G \cap \bB$ being the compact abelian subgroup in $\bT$. 
The quotient $ \bF := \bG/\bB $ carries the 
structure of a smooth projective algebraic variety over $\c$ referred to as a 
{\it flag variety} of $\bG$. The natural inclusion $ G \into \bG $ induces a homeomorphism:
\begin{equation}
\la{GTB}
G/T \,\cong\, \mathbb{F}
\end{equation}
where $\bF $ is viewed as a topological manifold with classical (metric) topology,
and we have 
\begin{equation}
\la{GTW}
W = N_G(T)/T = (N_G(\bT) \cap G)/(\bT \cap G) \cong N_{\bG}(\bT)/\bT\,,
\end{equation}
where $N_{\bG}(\bT) $ denotes the normalizer of $ \bT $ in $ \bG $.
We will often abuse the notation writing simply $ \bF = G/T $  and $ W = N_{\bG}(\bT)/\bT $
instead of  \eqref{GTB} and \eqref{GTW}. 
The Bruhat graph $ \Gamma $ describes the structure of the $0$- and $1$-dimensional complex
$\bT$-orbits in $ \bF $. Specifically, the vertices of $ \Gamma$ are in natural bijection with the $0$-dimensional $\bT$-orbits (i.e., the $\bT$-fixed points) in $\bF$. Indeed, 
under identification \eqref{GTB}, we have $ \bF^{\bT} \cong (G/T)^T = N_G(T)/T = W $,  and  we write $ \{w\} \in \bF^{\bT} $ for the fixed point corresponding to $w \in W$.
The edges $ E_{\Gamma} $ can be identified with the $\bT$-orbits of (complex) dimension one  in $ \bF $: we write $ \O_{e} $ for the  $\bT$-orbit corresponding to $ e = e(s_{\alpha}, w)\in E_\Gamma $.
The closure $ \overline{\O}_{e(s_{\alpha}, w)} $ of each  orbit  $ \O_{e(s_{\alpha}, w)} $ in $\bF$ is isomorphic to $\C\bP^1 $ and often referred to as a {\it $\bT$-invariant curve}.
Topologically, $ \overline{\O}_{e(s_{\alpha}, w)} $ is the $2$-sphere $S^2$ that contains two $\bT$-fixed points: $ \overline{\O}_{e(s_{\alpha}, w)} \cap {\bF}^{\bT} = \overline{\O}_{e(s_{\alpha}, w)}\setminus \O_{e(s_{\alpha}, w)} = \{w, s_{\alpha}w\}$, corresponding to the vertices adjacent to the edge $ e(s_{\alpha}, w)\in E_\Gamma$. Thus, in geometric terms, two
vertices $w$ and $w'$ are connected by an edge $e $  in $ \Gamma $ iff there is
a one-dimensional orbit $ \O_e $ in $\bF$ such that $ \overline{\O}_e = 
\O_e \cup \{w,w'\}$.

Next, we give a geometric interpretation of the moment category $ \Cc(\Gamma)$. 
To this end, we consider the $\bT$-equivariant $1$-skeleton $ \bF^{(1)} $ of the flag manifold
$\bF$. Recall that $ \bF^{(1)} $  is the union of all $ \bT$-orbits of complex dimension $0$ and $1$ in $\bF $: i.e., $ \bF^{(1)} = \cup_{e \in E_{\Gamma}} \overline{\O}_{e(s_{\alpha}, w)} $. As a $T$-space, it can be also defined topologically as 
 \begin{equation}
 \la{bF1}
 \bF^{(1)} \,=\,\{x\in\bF \, : \, \dim\,(T \cdot x) \leq 1\}
 \end{equation}
where $ \dim $ stands for the (real) dimension of $T$-orbits in $ \bF $  (see, e.g., \cite{AP93}). Thus, $\bF^{(1)} $ is a two-dimensional $T$-CW subcomplex of $\bF$ obtained by `gluing together' finitely many $2$-spheres $\bS^2$. We define on $ \bF^{(1)} $ the following `orbit relation':
\begin{equation}
\la{preorder}
x \,\lesssim\, y\quad \xLeftrightarrow{\rm def} \quad \overline{\bT\cdot x} \,\supseteq\,\bT\cdot y\ ,
\end{equation}
where the `bar' stands for the (Zariski) closure of a $\bT$-orbit in $\bF = \bG/\bB$. 
\begin{lemma}
\la{Lpreorder}
$ (\bF^{(1)},\,\lesssim) $ is a topological $\bT$-preorder, whose space of 
$\bT$-orbits is a $($discrete$)$ poset isomorphic to the moment category of $\Gamma:$
\begin{equation}
\la{quotientF}
\bF^{(1)}\!/\bT\,\cong\,\Cc(\Gamma)
\end{equation}
\end{lemma}
\begin{proof}
Recall that a preorder is a binary relation, which is reflexive and transitive (but not necessarily antisymmetric). These two properties obviously hold for \eqref{preorder}, and
the bijection \eqref{quotientF} follows from the above
geometric interpretation of $ \Gamma $. To see that this bijection is an isomorphism of posets we check that the $\bT$-orbits in $\bF^{(1)}$ are precisely the equivalence classes modulo the equivalence relation on $\bF^{(1)}$ associated to the preorder \eqref{preorder}: $\, x \sim y \,$ iff $\, x \lesssim y \,$ and $\,y \lesssim x \,$, and the induced partial order on $  \bF^{(1)}/\sim $ coincides under \eqref{quotientF} with the face partial order on $ \Cc(\Gamma)$. Thus, $ \Cc(\Gamma) $ is isomorphic to the universal poset under the preorder $ (\bF^{(1)},\,\lesssim) $.
\end{proof}
\begin{remark}
In (the proof of) Lemma~\ref{LemmaL4}, we will give another topological interpretation of the category $ \Cc(\Gamma) $ as the nerve diagram of an open covering of the space $ \F^{(1)}$, see \eqref{coverposet}.
\end{remark}
Note that the natural action of $\bG$ on $\bF$ restricts to an action of $ N_{\bG}(\bT) $ on $\bF^{(1)}$ that preserves 0- and 1-dimensional $\bT$-orbits respectively. Therefore this last action preserves the preorder \eqref{preorder}, and hence, under the isomorphism \eqref{quotientF} of Lemma~\ref{Lpreorder}, induces an action of $W = N_{\bG}(\bT)/\bT $ on the poset
$ \Cc(\Gamma)$. Explicitly, the action $\, W \times \Cc(\Gamma) \to \Cc(\Gamma)\,$ is given by
\begin{equation}
\la{Wact1} 
 w \mapsto gw\ , \quad e(s_{\alpha}, w) \mapsto e(g s_{\alpha} g^{-1},\, gw)\ , \quad g \in W\,.
\end{equation}

Now, if we regard $ \Cc(\Gamma) $ as a (small) category, then \eqref{Wact1} defines a strict action on it: that is, a functor $ W \to \Cat $, where $W$ is viewed as a category with the single object mapped to $\Cc(\Gamma)$ (see Appendix~\ref{AA}).
This allows us to enlarge $ \Cc(\Gamma)$ by taking its homotopy quotient in $ \Cat $. Specifically, we define the category $ \Cc(\Gamma)_{hW} $ by applying the classical Grothendieck construction to the functor $ W \to \Cat $ associated to the $W$-action \eqref{Wact1}:
\begin{equation}
\la{ChW}
\Cc(\Gamma)_{hW} := W \!\smallint \Cc(\Gamma).
\end{equation}
The category \eqref{ChW} will play a crucial role in the present paper. 
By definition, this category has the same objects as $ \Cc(\Gamma)$, while its morphisms are given by all pairs $(g,\,\varphi) $, where $ g \in W $ and $ \varphi \in \Hom_{\Cc(\Gamma)}(gc, \,c')$ for $ c,c' \in \Cc(\Gamma)$: thus, the non-identity morphisms in $\Cc(\Gamma)_{hW}$ are of the form:
\begin{equation}
\la{morChW}
w \to gw\ ,\quad  e(s_\alpha, w) \to e(gs_\alpha g^{-1}, gw)\ , \quad e(s_\alpha, w) \to gw\,,
\end{equation}
where $\,g,w \in W $ and $\alpha \in \cR_+$. Note that there is a canonical (inclusion) functor $ \Cc(\Gamma) \into \Cc(\Gamma)_{hW} $ that acts as the identity on objects while taking a morphism $ \varphi $ in $ \Cc(\Gamma) $ to $ (e, \varphi) $ in $ \Cc(\Gamma)_{hW}$. 
Using this functor, we will often regard $\Cc(\Gamma) $ as a (wide) subcategory of  
$ \Cc(\Gamma)_{hW}$. The category $\,  \Cc(\Gamma)_{hW}\,$ has a natural skeleton (i.e., a full subcategory $\,\overline{\Cc(\Gamma)}_{hW} $ that is equivalent to $ \Cc(\Gamma)_{hW}$ but contains no non-identity isomorphisms), which is described in Lemma~\ref{skCG} (see Fig.~\ref{skC}). Finally, we remark that our notation  \eqref{ChW} is justified by the well-known fact (see \eqref{ThomForm}) that 
$$ 
B[\Cc(\Gamma)_{hW}] \simeq B[\Cc(\Gamma)]_{hW}\,,
$$
where $B[\Cc(\Gamma)]$ is the classifying space of $ \Cc(\Gamma)$ and $B[\Cc(\Gamma)]_{hW}$ is 
the usual (topological) Borel construction on $B[\Cc(\Gamma)]$ with respect to the continuous 
action of $W$ induced by \eqref{Wact1}. 

\subsection{Homology decomposition theorem} \la{S4.2}
The above geometric interpretation of the graph $\Gamma $ allows us to define the canonical functor
\begin{equation}
 \la{F0}
\cF_{0}\,:\ \Cc(\Gamma)\,\to \, \Top^{T}\,,\quad w \mapsto \{w\}\ ,
\quad e(s_{\alpha}, w) \mapsto \O_{e(s_{\alpha}, w)}\,,
\end{equation}
which we call $ \eqref{F0} $ the {\it $T$-orbit functor} on $ \Cc(\Gamma) $.

Recall that a functor $ \cF $ on $ \Cc(\Gamma) $ is a {\it $W$-functor} if it extends to $ \Cc(\Gamma)_{hW} $. Such an extension is given by additional data: a family of natural transformation $ \{\vartheta_w: \cF \to \cF \circ \hat{w}\}_{w \in W} $ called a $W$-structure on $ \cF $ (see Definition~\ref{Wfun} and Lemma \ref{LA4}).  The $T$-orbit functor   \eqref{F0} does {\it not} admit a $W$-structure: it does not `lift' to $ \Cc(\Gamma)_{hW} $ in a natural way. To remedy this problem we consider the induced functor with values in $G$-spaces:
\begin{equation}
\la{GF0}
G \times_T \cF_0:\ \Cc(\Gamma) \to \Top^G \,,\quad w \mapsto G\times_T \{w\}\,,\quad e(s_{\alpha}, w) \mapsto G \times_{T} \O_{e(s_\alpha, w)}\,,
\end{equation}
obtained by composing $\cF_0 $ with the canonical functor $ G \times_T (-):\,\Top^T \to \Top^G $. 

\vspace*{1ex}

The next theorem is the main result of this section.
\begin{theorem} \la{hDec}
$(1)$ The functor \eqref{GF0} has a natural $W$-structure with respect to the $W$-action \eqref{Wact1} on $ \Cc(\Gamma)$.

$(2)$ For every prime\footnote{including `$p=0$', in which case `mod-$p$' means `with coefficients in $\Q$'.} $p \not= 2 $, there is a $G$-equivariant mod-$p$ cohomology isomorphism 
\begin{equation}
\la{Gcohiso}
\mathrm{hocolim}_{\Cc(\Gamma)_{hW}} [\,G \times_T \cF_{0}\,] \,\simeq \, G/T\,, 
\end{equation}
where $ G \times_T \cF_{0} $ is viewed as a functor $ \Cc(\Gamma)_{hW} \to \Top^G $
with the $W$-structure from $(1)$.
\end{theorem}
%

We will refer to \eqref{Gcohiso} as a {\it homology decomposition of $G/T$ over $\Gamma$}. 

\subsection{Proof of \Cref{hDec}}
\la{S4.3}

\vspace{1ex}

In preparation for the proof of Theorem~\ref{hDec}, we introduce another functor:
\begin{equation}
 \la{FFP0}
\wcF_{0}\,:\ \Cc(\Gamma)\,\to \, \Top^{T}\,,
\quad w \mapsto \bigvee_{e(s_{\alpha}, w) \in E_{\Gamma}(w)}\!\! \!\overline{\O}_{e(s_{\alpha}, w)}\!\setminus\!\{s_{\alpha} w\}\ ,
\quad e(s_{\alpha}, w) \mapsto \O_{e(s_{\alpha}, w)}\,,
\end{equation}
where each space $\overline{\O}_{e(s_{\alpha}, w)}\!\setminus\!\{s_{\alpha} w\}$
is pointed by $\{w\} \in \bF^{\bT} $ and the wedge sum runs over the set  $ E_{\Gamma}(w) := \{e \in E_{\Gamma}\,:\, w \in \overline{\O}_e\}$ of all edges emanating from the corresponding vertex $w \in  V_{\Gamma}$. On morphisms, the functor $ \wcF_0$ is defined by sending the arrow $ e(s_{\alpha}, w) \to w $ in $ \Cc(\Gamma)$ to the natural inclusion of the orbit $ \O_{e(s_{\alpha}, w)} $ into the corresponding summand  in the wedge sum over $ E_{\Gamma}(w)$. 

Note that there is a morphism of functors $\,\wcF_0 \to \cF_0 $ defined on edges $ e \in E_{\Gamma}$ as the identity map on $ \O_{e} $, while on vertices $ w \in V_{\Gamma}$, as the trivial map collapsing the wedge sum over $ E_{\Gamma}(w) $  to its basepoint $\{w\}$.  
Since each wedge sum over $ E_{\Gamma}(w) $ in \eqref{FFP0} is contractible, this 
morphism is  an (objectwise) homotopy equivalence: $ \wcF_0 \xrightarrow{\sim} \cF_0 $. We consider the induced 
morphism of functors to $G$-spaces, $\Cc(\Gamma) \to \Top^G$, 
which we denote by
\begin{equation}
\la{GTF}
\tilde{p}:\, G \times_{T} {\wcF}_0 \,\to\, G \times_{T} \cF_0. 
\end{equation}
\begin{lemma}
\la{lem:W-fun-heq}
The functors $ G \times_T \cF_0 $ and $\, G \times_T \wcF_0 $  carry canonical $W$-structures, the morphism \eqref{GTF} being a $G$-equivariant weak equivalence of $W$-functors.
\end{lemma}
\begin{proof}
For $w \in W $, we choose a representative $ \dot{w} \in N_G(T)$ such that
$ w = \dot{w}\, T$,  and define the natural transformation $\,\vartheta_w:\, G \times_T \cF_0 \to (G \times_T \cF_0) \circ \hat{w} \,$ by 
\begin{equation}
\la{eq:tildeF0-ver}
\vartheta_{w,v}:\ G \times_T \{v\} \to G \times_T \{wv\}\,,\quad
[(g,v)]_T \mapsto [(g \dot{w}^{-1}, wv)]_T
\end{equation}
where 
$ \{v\} \in \bF^{\bT} $ corresponds to $ v \in V_{\Gamma} = W $, and
\begin{equation}
\la{eq:tildeF0-edge}
\vartheta_{w,e}:\ G \times_T \O_e \to G \times_T \O_{we}\,,\quad
[(g,x)]_T \mapsto [(g \dot{w}^{-1},\, \dot{w} \cdot x)]_T
\end{equation}
where $ e \in E_{\Gamma} $, $ x \in \O_e $, and $ w\in W $ acts on $e$ as in \eqref{Wact1}. 
The `dot' action $ x \mapsto \dot{w} \cdot x $ on $\bT$-orbits in \eqref{eq:tildeF0-edge} is obtained by restricting the natural (left) action of $N$ on $ \bF^{(1)} $.
It is straightforward to check that the above formulas for $ \vartheta_{w,v} $ and $ \vartheta_{w,e} $ are independent of the choice of representatives $ \dot{w} $ in $N_G(T)$. Moreover, for each $ w \in W $, they indeed define a natural transformation of functors: for this, it suffices to note that the following diagrams commute in $\Top^G$
%
\[\begin{tikzcd}
	{G \times_T \O_{e(s_\alpha, v)}} & {G\times_T\{v\}} \\
	{G\times_T \O_{w \cdot e(s_{\alpha}, v)}} & {G \times_T\{wv\}}
	\arrow[from=1-1, to=1-2]
	\arrow["{\vartheta_{w,e}}"', from=1-1, to=2-1]
	\arrow["{\vartheta_{w,v}}", from=1-2, to=2-2]
	\arrow[from=2-1, to=2-2]
\end{tikzcd}\]
for all non-identity maps $ e(s_{\alpha}, v) \to v $ in $ \Cc(\Gamma)$. Finally,
verifying that the family  $ \{\vartheta_w\}_{w \in W} $ defines
a $W$-structure on the functor $ G \times_T \cF_0 $ amounts to verifying
the identities ({\it cf.} \eqref{diagact})
$$
\vartheta_{w_1 w_2, \,v} = \vartheta_{w_1, \,w_2 \cdot v} \,\circ\, \vartheta_{w_2, v} \ ,\quad
\vartheta_{w_1 w_2, \,e} = \vartheta_{w_1, \,w_2 \cdot e} \,\circ\, \vartheta_{w_2, e}
$$
for all $ v \in V_{\Gamma}$ and $ e \in E_{\Gamma}$, which is a trivial exercise in view of formulas \eqref{eq:tildeF0-ver} and \eqref{eq:tildeF0-edge}.

Similarly, we define the morphisms $ \tilde{\vartheta}_w: G \times_T \wcF_0 \to (G \times_T \wcF_0) \circ \hat{w} $: for each $ v \in V_{\Gamma}$,
$$
\tilde{\vartheta}_{w,v}:\ G \times_T \left(\bigvee_{e \in E_{\Gamma}(v)}\!\! \!\overline{\O}_{e(s_{\alpha}, v)}\!\setminus\!\{s_{\alpha} v\}\right)\ \to \ G \times_T \left(\bigvee_{w\cdot e \in E_{\Gamma}(v)}\!\! \!\overline{\O}_{w\cdot e(s_{\alpha}, v)}\!\setminus\!\{ws_{\alpha} v\}\right)\ 
$$
is given by $\,[(g,x)]_T \mapsto [(g \dot{w}^{-1}, \dot{w}x)]_T \,$,
while for each  $ e = e(s_{\alpha}, w) \in E_\Gamma $, $\,\tilde{\vartheta}_{w,e}\,$ by the same formula as \eqref{eq:tildeF0-edge}. 
Again, it is routine to check that $\tilde{\vartheta}_{w}$ is a natural transformation, and the family $\{\tilde{\vartheta}_{w}\}_{w \in W }$ defines a $W$-structure on the functor $ G \times_T \wcF_0 $.  
As remarked before Lemma~\ref{lem:W-fun-heq}, the morphism 
$ \wcF_0 \xrightarrow{\sim} \cF_0 $ is a natural weak equivalence of functors, and hence so is the induced morphism \eqref{GTF}. It suffices to note that \eqref{GTF} intertwines the families $ \{\vartheta_{w}\}_{w \in W } $ and $ \{\tilde{\vartheta}_{w}\}_{w \in W } $ and thus is a morphism of $W$-functors.
\end{proof}

Next, we observe that, for all $w\in V_{\Gamma} $ and $ e \in E_{\Gamma} $, there are canonical $G$-equivariant maps 
\[
\tilde{q}_w:\ G\times_T \left(\bigvee_{e \in E_{\Gamma}(w)}\!\! \!\overline{\O}_{e(s_{\alpha}, w)}\!\setminus\!\{s_{\alpha} w\}\right)\ \to\ 
\bF \quad ,\quad 
\tilde{q}_e:\ G\times_T \O_{e} \to\ 
\bF
\]
induced by the natural ($T$-equivariant) inclusions $\overline{\O}_{e(s_{\alpha}, w)}\!\setminus\!\{s_{\alpha} w\} \into \bF $ 
and $\,\O_e  \into \bF $. It is easy to check that these maps assemble to a morphism of $W$-functors 
\begin{equation}
\la{GTFq}
\tilde{q}:\,G\times_T \wcF_0 \to \bF\,, 
\end{equation}
where $G/T$ is regarded as a constant functor on $\Cc(\Gamma)$ equipped with the {\it trivial} $W$-structure. 

The morphism \eqref{GTFq} induces a $G$-equivariant map of spaces:
\begin{equation}
\la{eq:p-hocolim}
   q:\, \hocolim_{\Cc(\Gamma)_{hW}}(G\times_T \wcF_0) \to \bF
  \end{equation}
Our next goal is to show that the map \eqref{eq:p-hocolim} is a $p$-local weak homotopy equivalence at all primes except $p=2$. We begin by identifying the domain of this map, i.e. computing the homotopy colimit of $ G \times_T \wcF_0 $, in geometric terms.
 \begin{lemma}
 \la{LemmaL4}
     There is a natural weak equivalence of $G$-spaces 
     \begin{equation}
     \la{eq:we-1-skel}
         \hocolim_{\Cc(\Gamma)_{hW}}(G\times_T \wcF_0) \,\simeq\, G\times_N \bF^{(1)}
     \end{equation}
     where  $\bF^{(1)} $ is the $ T$-equivariant $1$-skeleton of $\bF$.
 \end{lemma}
 %
 \begin{proof}
 Recall that $\bF^{(1)} $ is defined topologically by \eqref{bF1}.
 The action of $N = N_G(T) $ on $\bF^{(1)}$ used in \eqref{eq:we-1-skel} is obtained by restricting the natural action of $G$ on $\bF = G/T$.  
     It is easy to see that this gives a well-defined continuous $N$-action on $\bF^{(1)}$, and moreover, the corresponding (diagonal) action on $ G \times \bF^{(1)}$,$\,(g,x) \mapsto (gn^{-1},nx)\,$, induces a  free action of $W $ on $\,G\times_T\bF^{(1)}\,$. 

    We will calculate the homotopy colimit \eqref{eq:we-1-skel} geometrically, using a simple covering argument. For each fixed point $ w \in \bF^{\bT} = W $, we let
     \begin{equation}
     \la{eq:cover}
    U_w := \bigvee_{e \in E_{\Gamma}(w)}\!\! \!\overline{\O}_{e(s_{\alpha}, w)}\!\setminus \{s_{\alpha} w\} \subset \bF^{(1)}
    \, 
    \end{equation}
    and denote by $ \U := \{U_w\}_{w \in W}$ the set of all such subsets of $ \bF^{(1)}$. Clearly, $\U $ is an open covering of $\bF^{(1)} $ stable under the $T$-action. Associated to $ \U $ is an abstract simplicial complex $ {\rm Nrv}(\U)$ (called
    the nerve of $\U$), with vertex set $W$ and the simplices 
    $$
    \Sigma_{\U} =\{I \subseteq W\ :\ U_I \not= \varnothing\}  
    $$
    indexing all non-empty intersections $ U_I :=\bigcap_{w\in I}U_w $ of the covering sets \eqref{eq:cover}. It is easy to see that  $U_I\neq \varnothing $  for $I\subseteq W$ iff either $I=\{w\}$ or $I=\{w,s_{\alpha}w\}$ for some $ w \in W $ and $ \alpha \in \R_+ $. (In the first case, we have $U_{\{w\}}=U_w$, while in the second, $\,U_I = U_w \cap U_{s_{\alpha}w}=\O_{e(s_{\alpha}, w)}$.) Hence, there is a natural isomorphism of 
    simplicial complexes: $\,{\rm Nrv}(\U) \cong \Gamma \,$, and hence an isomorphism
    of posets: 
    \begin{equation}
     \la{coverposet}
     \Sigma_{\U}^{\rm op}\, \cong \,\Cc(\Gamma)\,,
    \end{equation}
    which shows that $ \Cc(\Gamma) $ is the {\it nerve diagram} of the covering $\U $ ({\it cf.} \cite[Definition~2.3]{BKRR23}).
    
     Now, observe that, under the isomorphism \eqref{coverposet},  $\,\wcF_0\,$ becomes the standard functor associated to the covering \eqref{eq:cover}:
     \[ 
     \cF_{\U}\,:\, \Sigma_{\U}^{\rm op} \,\to\, \Top^T, \quad I \mapsto U_I
     \]
Hence, by \cite[Proposition 3.2]{DI04},
     \begin{equation}
     \la{eq:F1-hocolim}
         \hocolim_{\Cc(\Gamma)}(\wcF_0) \,\cong\, \hocolim_{\Sigma_{\U}^{\rm op}}\left(\cF_{\U}\right) \,\simeq\, \bF^{(1)}\,,
     \end{equation}
    and therefore
     \[\hocolim_{\Cc(\Gamma)}(G\times_T\wcF_0) \,\simeq\, G\times_T[\hocolim_{\Cc(\Gamma)}(\wcF_0)] \,\simeq\, G\times_T\bF^{(1)} \]\,
     By \Cref{lem:Thom}, we then conclude 
     \[\hocolim_{\Cc(\Gamma)_{hW}}(G\times_T\wcF_0) \,\simeq\, (G\times_T\bF^{(1)})_{hW} \,\simeq\, (G\times_T\bF^{(1)})\!/W \,\cong\, G\times_N\bF^{(1)}\,,
     \]
     where on the second step we use the (obvious) fact that $W$ acts freely on  $G\times_T\bF^{(1)}$.
 \end{proof}
\begin{remark}
\la{remcofres} The result of
Lemma~\ref{LemmaL4} --- specifically, the key equivalence~\eqref{eq:F1-hocolim} --- can be established `abstractly', using  model categorical arguments. In fact, one can show (see Lemma~\ref{Lcofres}) that
$ \wcF_0 $ is a cofibrant object (hence, a cofibrant replacement of $ \cF_0$) in the category $ \Top^{\Cc(\Gamma)} $ of $\Cc(\Gamma)$-diagrams equipped with the projective model structure. The homotopy colimit of $ \wcF_0$ is therefore equivalent to the ordinary colimit, which is easily seen to be  homeomorphic to $ \bF^{(1)}$. Thus,
\begin{equation}
\la{F0F0}
\hocolim_{\Cc(\Gamma)}(\cF_0) \,\simeq\,
\hocolim_{\Cc(\Gamma)}(\wcF_0)\, \simeq\,\colim_{\Cc(\Gamma)}(\wcF_0) \,\cong\, \bF^{(1)}\,.
\end{equation}
We mention that the last homeomorphism holds in general
for the nerve diagram of an arbitrary open covering 
(see \cite[Remark 2.4]{BKRR23}).
\end{remark}
As an application,
we will use the homotopy decomposition \eqref{F0F0} to clarify the question
about (the odd-dimensional part of) $T$-equivariant cohomology of $\F^{(1)}$ raised by A.~Knutson and I.~Rosu in the Appendix of \cite{KR03}. 
\begin{cor} \la{tgkm1sk}
Let $\k$ be any field of characteristic $\neq 2$. 
\begin{eqnarray*}
H^{\rm ev}_T(\F^{(1)},\,\k) & \cong & \big\{\sum_{w \in W} p_w \otimes w\in \k[V \times W]\,:\, p_w-p_{s_\alpha w} \in \langle \alpha \rangle\,,\ \forall \alpha \in \A\ ,\ \forall w \in W \big\}\,,\\
H^{\rm odd}_T(\F^{(1)},\,\k) & \cong & 
\frac{ \bigoplus_{\alpha\in \A}\, \k[V \times W]\,[-1]}{\langle (f_{\alpha}-f)_{\alpha \in \A}\ : \ f\in \k[V \times W]\,,\ f_{\alpha} \in \bQ_1(W_{\alpha})\rangle}\,, 
\end{eqnarray*}
where $ \bQ_1(W_{\alpha}) \subseteq \k[V \times W] $ is defined in \eqref{Q1W}.
\end{cor}
\begin{proof}
We will study the Bousfield-Kan spectral sequence (see, e.g., \cite[Theorem 1.1]{JMO94}) associated with the homotopy decomposition \eqref{F0F0}:
\begin{equation} \la{BKS} 
E_2^{pq}\,=\, {\varprojlim}_{\,\Cc(\Gamma)^{\rm op}}^p \big[\, H^q_T(\cF_0(\,\mbox{--}\,),\,\k)\,\big] \ \Rightarrow\ H^{p+q}_T(\F^{(1)},\,\k)\,.\end{equation}
First,  note that, for all $w \in W$ and all $e(s_\alpha,w) \in E_\Gamma$,  there are natural isomorphisms
\begin{equation*} \la{hobj}
H^\ast_T(\cF_0(w),\k)\, \cong\, H^*(BT,\,\k)\,\cong \,\k[V]\ ,\qquad H^\ast_T(\cF_0(e(s_\alpha,w)),\k)\, \cong\, H^*(BT_\alpha,\,\k)\,\cong \,\k[V]/\langle \alpha \rangle 
\end{equation*}
Since $\Cc(\Gamma)^{\rm op}$ is a poset of height $1$, the higher inverse limits $\,{\varprojlim}_{ \Cc(\Gamma)^{\rm op}}^p $  vanish for all $p>1$ (see \cite[Theorem 1.8]{Ch76}),
 hence the spectral sequence \eqref{BKS} degenerates at $ E_2$-term.
Furthermore, since $H^q_T(\cF_0(\mbox{--}),\k) = 0 $ for odd $q$, we get from \eqref{BKS} the isomorphisms
\begin{equation} \la{bk0}
H^{2i}_T(\F^{(1)},\,\k) \cong {\varprojlim}_{c \in \Cc(\Gamma)^{\rm op}} [H^{2i}_T(\cF_0(c),\,\k)] \,,\quad  H^{2i+1}_T(\F^{(1)},\,\k) \cong {\varprojlim}_{c \in \Cc(\Gamma)^{\rm op}}^1[H^{2i}_T(\cF_0(c),\,\k)]
\end{equation}
for all $i \ge 0 $. Now, for any small category $ \Cc$, the higher limits ${\varprojlim}_{\Cc}^p$ can be computed as cohomology of a standard cochain complex (see, e.g., \cite[Lemma 4.2]{JMO94}). In our case, this complex has only two nonzero terms, corresponding to the objects and morphisms
in $ \Cc(\Gamma)^{\rm op}$:
\begin{equation} \la{bk1} 
0 \to \bigoplus_{w \in W} \k[V]\ \oplus \bigoplus_{\alpha \in A \atop e(s_\alpha,w) \in E_\Gamma} \k[V]/\langle \alpha \rangle \ \xrightarrow{\quad \partial\quad } \bigoplus_{\alpha \in \A \atop (e(s_\alpha,w)\to w) \in \Cc(\Gamma)} \k[V]/\langle \alpha \rangle \to 0 \,.
\end{equation}
The differential $ \partial $ in \eqref{bk1} is given by
\begin{equation} \la{bk2} 
\big((p_w)_{w \in W},\ (q_{e(s_\alpha,w)})_{e(s_\alpha,w)\in E_\Gamma} \big) \ \mapsto \ (\bar{p}_w-q_w)_{(e(s_\alpha,w) \to w)\in \Cc(\Gamma)}\,,
\end{equation}
where $ \bar{p}_w \in \k[V]/\langle\alpha\rangle $ stands for the image of $ p_w \in \k[V]$ 
under the canonical projection. To see the isomorphisms of Corollary \ref{tgkm1sk} we simplify the complex \eqref{bk1} as follows. First,  we check that it is quasi-isomorphic to  
\begin{equation} \la{bk3} 
0 \to  \bigoplus_{w \in W} \k[V]\ \xrightarrow{\ \delta\ }\ \bigoplus_{\alpha \in \A} \ \bigoplus_{e(s_\alpha,w) \in E_\Gamma} \k[V]/\langle \alpha \rangle \to 0
\end{equation}
where $\, \delta\big((p_w)_{w \in W}) \,:=\, \pm (\bar{p}_w-\bar{p}_{s_\alpha w})_{e(s_\alpha,w) \in E_{\Gamma}}\,$, with the sign being positive if $w < s_\alpha w$ in the Bruhat order on $W$.
This immediately shows that
$$
\Ker(\partial)\,\cong\,\Ker(\delta)\,\cong \,\big\{ (p_w)\in \oplus_{w \in W} \, \k[V] \,:\, p_w-p_{s_\alpha w} \in \langle \alpha \rangle  \,\,\forall \,\alpha \in \A\,\,,\,\,\forall \ w \in W \big\}\,,
$$
which, by \eqref{bk0}, yields the first isomorphism of Corollary \ref{tgkm1sk}. Next, we observe that \eqref{bk3} is actually isomorphic to the two-term complex
\begin{equation} \la{bk4}
 0 \to \k[V \times W] \ \xrightarrow{\ \bar{\Delta}\ }\ \bigoplus_{\alpha \in \A} \,\k[V \times W]/\bQ_1(W_\alpha) \to 0
\end{equation}
This isomorphism identifies $\k[V \times W]$ with $\oplus_{w \in W} \k[V]$ as in \eqref{idenkwb}, and the component 
$$\bigoplus_{\alpha \in \A} \, \k[V \times W]/\bQ_1(W_\alpha) \to \bigoplus_{\alpha \in \A} \ \bigoplus_{e(s_\alpha,w) \in E_\Gamma} \k[V]/\langle \alpha \rangle $$
of this isomorphism sends the summand $ \k[V \times W]/\bQ_1(W_\alpha)$ to the summand $\oplus_{e(s_\alpha,w) \in E_\Gamma} \k[V]/\langle \alpha \rangle$ corresponding to the root $\alpha$ by the map
$$ [\sum p_w \otimes  w] \,\,{\rm mod}\,\, \bQ_1(W_\alpha) \mapsto \big(\pm (\bar{p}_w-\bar{p}_{s_\alpha w})_{e(s_\alpha,w)}\big)\,,$$
where the sign is positive if $w < s_\alpha w$ in the Bruhat order. Is is easy to see that the cokernel of $\bar{\Delta}$ is isomorphic to $\bigoplus_{\alpha \in \A }\,\k[V \times W]/R(\A)$, where
$$
R(\A) = \langle(f_{\alpha}-f)_{\alpha \in \A}\ : \ f\in \k[V \times W]\,,\ f_{\alpha} \in \bQ_1(W_{\alpha})\rangle\,.
$$
By \eqref{bk0}, this yields the second isomorphism of Corollary \ref{tgkm1sk} and completes the proof.
\end{proof}
Next, we observe that the map \eqref{eq:p-hocolim} factors through the equivalence \eqref{eq:we-1-skel}, inducing the  $G$-equivariant map 
\begin{equation}
\la{eq:G-act}
\bar{q}:\, G\times_N\bF^{(1)} \,\to\, \bF
\end{equation}
that extends the natural ($T$-equivariant) inclusion $\bF^{(1)}\into \bF$. Hence, to show that \eqref{eq:p-hocolim} is a mod-$p$ cohomology isomorphism it suffices to show that so is the map \eqref{eq:G-act}. To this end, we prove 
\begin{lemma}
\la{lem:act-fib}
    There is a homotopy fibration sequence 
    \begin{equation}\la{fibseq}
    \bigvee_{\alpha\in\R_+} G_{\alpha}/N_{\alpha} \,\to\, G\times_N\bF^{(1)} \,\xrightarrow{\bar{q}}\, \bF
    \end{equation}
    where 
    $G_{\alpha}\coloneqq C_G(T_{\alpha})\subseteq G$ is the centralizer of the singular torus
    $ T_{\alpha} \subseteq T$ in $G$ and $N_{\alpha} = G_{\alpha}\cap N $ is the normalizer of $T$ in $ G_{\alpha}$. 
\end{lemma}
\begin{proof}
For $ \alpha \in \R_+ $, let $ \bF^{(1)}_{\alpha} := \bF^{(1)}\setminus\bigcup_{\beta \not= \pm w(\alpha) } \O_{e(s_{\beta}, w)} $ denote the complement of all one-dimensional orbits $ \O_{e(s_{\beta}, w)}$ in $ \bF^{(1)} $ with $ s_\beta \not=  s_{w(\alpha)} $. Note that the fixed points $ \bF^{\bT} $ are contained in $ \bF^{(1)}_{\alpha} $ for all $ \alpha \in \R_+ $, 
and the space $\bF^{(1)} $ is obtained by gluing the subspaces $\, \bF^{(1)}_{\alpha}\,$ along their common intersection $ \bF^{\bT} $:
\begin{equation}
\la{decf1}
\bF^{(1)} \cong \bigvee_{\bF^{\bT}} \bF_{\alpha}^{(1)} \, :=\, \bigl(\coprod_{\alpha \in \R_+}\bF^{(1)}_{\alpha}\bigr)/\sim
\end{equation}
Next, observe that the action of $ N $ on $ \bF^{(1)}$  restricts to $ \bF^{(1)}_{\alpha}$ for each $\alpha \in \R_+ $, and this action restricts further to $\bF^{\bT} $, inducing a simply transitive action of $W$ on $\bF^{\bT}$. Hence,  \eqref{decf1} yields a $G$-equivariant decomposition
$$
G\times_{N} \bF^{(1)}\, \cong\, \bigvee_{G \times_N \bF^{\bT}} \left(G \times_N \bF_{\alpha}^{(1)}\right) \,\cong\,
\bigvee_{G/T} \left(G \times_{N_{\alpha}} \overline{\O}_{e(s_{\alpha}, 1)}\right)\,
$$
Since $ \overline{\O}_{e(s_{\alpha}, 1)} \cong G_{\alpha}/T $ for each $ \alpha \in \R_+$, we conclude that
\begin{equation}
    \la{decf2}
    G\times_N \bF^{(1)} \simeq \bigvee_{G/T} G \times_{G_{\alpha}}\left(G_{\alpha}\times_{N_{\alpha}} G_{\alpha}/T\right) \simeq \bigvee_{G/T} G \times_{G_{\alpha}}\left(G_{\alpha}/{N_{\alpha}} \times G_{\alpha}/T\right)\,,
\end{equation}
where the last isomorphism is induced by 
$$
G_{\alpha} \times_{N_{\alpha}} G_{\alpha}/T  \xrightarrow{\sim} G_{\alpha}/N_{\alpha} \times G_{\alpha}/T\ ,\quad
 [(g_1, g_2 T)]_{N_{\alpha}} \mapsto ([g_1 N_{\alpha}, g_1 g_2 T)\,.
$$
Now, under the decomposition \eqref{decf2}, the map $\bar{q}:  G\times_N \bF^{(1)} \to G/T $ contracts the factor $ G_{\alpha}/N_{\alpha} $ in each wedge summand
in \eqref{decf2} to a point, while identifying $ G \times_{G_{\alpha}} (G_{\alpha}/T) \cong G/T $ in the natural way. It follows that $ {\rm hofib}_{*}(\bar{q}) \simeq {\rm fib}_*(\bar{q}) \simeq \bigvee_{\alpha} G_{\alpha}/N_{\alpha}$, completing the proof of the lemma.
\end{proof}

\begin{proof}[Proof of \Cref{hDec}]
We have constructed above a `zig-zag' of morphisms of $W$-functors
\begin{equation}
\la{eq:F-zigzag}
    \begin{tikzcd}
        & G\times_T \wcF_0 \arrow[dl, "\tilde{p}"']\arrow[dr,"\tilde{q}"] & \\
        G\times_T \cF_0 && G/T
    \end{tikzcd}
\end{equation}
where $ G/T $ is regarded as a constant functor on $ \Cc(\Gamma)$ with trivial $W$-structure.
By Lemma~\ref{lem:W-fun-heq}, $\, \tilde{p}$ is a natural weak equivalence, while, by Lemma~\ref{LemmaL4}, $ \tilde{q}$ induces the map of spaces \eqref{eq:p-hocolim}, which is (weakly) equivalent to \eqref{eq:G-act}. Now, by Lemma~\ref{lem:act-fib}, the map \eqref{eq:G-act} fits in the 
fibration sequence \eqref{fibseq}. Associated to \eqref{fibseq} there is a Serre spectral sequence of the form
$$
E_2^{p,q}  \cong \bigoplus_{\alpha \in \R_+} H^p(\bF,\, H^{q}(G_{\alpha}/N_{\alpha},\,\k))\ \Rightarrow \ H^{p+q}(G \times_N \bF^{(1)},\,\k)\ ,
$$
where cohomology is taken with coefficients in an arbitrary field $ \k $.
Since $ G_{\alpha}/N_{\alpha} \cong (G_{\alpha}/T)/W_{\alpha} \cong \C\bP^1\!/\Z_2 \cong {\mathbb R}\bP^2 $ for all $ \alpha \in \R_+ $, the above spectral sequence collapses at $E_2$-term,  giving an isomorphism $ H^*(\bF,\,\k) \cong
H^{\ast}(G \times_N \bF^{(1)},\,\k) $, provided $ {\rm char}(\k)\not=2$. 

Thus, \eqref{eq:F-zigzag} induces a `zig-zag' of $G$-equivariant maps of spaces
$$
\hocolim_{\Cc(\Gamma)_{hW}}(G \times_T \cF_0) \, \xleftarrow{\sim}\, 
\hocolim_{\Cc(\Gamma)_{hW}}(G \times_T \wcF_0) \,\xrightarrow{q} \, G/T
$$
where the leftmost map is a (weak) homotopy equivalence, while the rightmost map is a mod-$p$ cohomology isomorphism for all primes $ p \not=2$. This completes the proof of Part $(2)$ of \Cref{hDec}. Part $(1)$ was established in Lemma~\ref{lem:W-fun-heq}.
\end{proof}

\begin{remark}
\la{pplus0}
In the next section, we will show that the map $ q $ defined by \eqref{eq:p-hocolim} is an example of a generalized plus construction (in the sense of \cite{BLO21}). Specifically, if we set $\, F_0(G,T) := \hocolim_{\Cc(\Gamma)_{hW}}(G \times_T \wcF_0)\,$, then $q: F_0(G,T) \to G/T $
exhibits $ G/T $ as a simultaneous $p$-plus construction $\, G/T \simeq F_0^{p+}(G,T)\,$ for  all primes $p\not=2$ (see Section~\ref{S5.3}).
\end{remark}

\begin{remark}
\la{2ndWact}
In addition to \eqref{Wact1} the moment category $ \Cc(\Gamma) $ carries another  $W$-action coming from the natural action of $W$ on $ \bF = G/T \,$ (see \eqref{WG/T}). Although  the action of $W$ on
$ \bF $ is not algebraic, it restricts to a continuous action on $ \bF^{(1)} $, preserving the $\bT$-orbits and the preorder relation \eqref{Lpreorder}. Under the isomorphism \eqref{quotientF}, it then induces a (strict) $W$-action on $ \Cc(\Gamma) $
given by
\begin{equation}
\la{Wact2} 
 w \mapsto wg^{-1}\ , \quad e(s_{\alpha}, w) \mapsto e(s_{\alpha},\, wg^{-1})\ , \quad \forall\,g \in W \,.
\end{equation}
The action \eqref{Wact2} obviously commutes with \eqref{Wact1}, making
$ \Cc(\Gamma) $ a $(W\times W)$-category (where we assume that the first factor in $W \times W$ acts as \eqref{Wact1} and the second as \eqref{Wact2}). In Section~\ref{S6.4}, we show that both  Lemma~\ref{lem:W-fun-heq}  and Theorem~\ref{hDec} can be refined to incorporate this $ (W \times W)$-action. In particular, we prove that the cohomology equivalence  \eqref{Gcohiso} is $W$-equivariant.

We should mention that various actions of $W$ on $T$-equivariant cohomology of flag manifolds $ \bF $ play an important role in Schubert calculus and have been studied extensively in the literature. In particular, the two $W$-actions  on the moment graph --- \eqref{Wact1} and \eqref{Wact2} --- first appeared in the work of A. Knutson \cite{Kn03}, who called them the left and right actions respectively, and were further studied by J. Tymosczko \cite{Ty08}, and more recently, by J. Carrell \cite{Ca20}.
\end{remark}

\section{Quasi-flag manifolds of compact Lie groups}
\la{S5}

In this section, we introduce our main objects: the quasi-flag manifolds $ F_m^{+}(G,T) $.
The construction of these spaces proceeds in two steps. First, motivated by Theorem~\ref{hDec},  we define the spaces $ F_m(G,T) $ as homotopy colimits of generalized $G$-orbit functors $ G \times_T \cF_m:\, \Cc^{(m)}(\Gamma) \to \Top^{G} $ defined 
over $W$-categories $ \Cc^{(m)}(\Gamma)$, depending  on a multiplicity  $m \in \M(W)$. For each $m \in \M(W)$, the category $ \Cc^{(m)}(\Gamma) = \Cc\left[\Gamma^{(m)}\right]$ is obtained by replacing the Bruhat moment graph $\Gamma $ with a higher dimensional CW-complex $\, \Gamma^{(m)}$, which we call the $m$-simplicial thickening of $ \Gamma $. 
The spaces $ F_m(G,T)$ thus defined
have correct mod-$p$ cohomology for all primes $p\not=2$; however, they are not simply connected. To remedy this problem, in our second step, we apply to $ F_m(G,T)$ a relative version of Quillen's plus construction introduced recently in \cite{BLO21}. Luckily, this construction can be performed in a universal functorial way (i.e., naturally for all primes $p \not=2 $ and for all $m \in \M(W) $ at once), and thus provides nice canonical models $ F_m^{+}(G,T) $ for quasi-flag manifolds.

\subsection{Simplicial thickening of moment graphs}
\la{S5.1}
We begin with a general construction of $\Gamma^{(m)}$ for an arbitrary graph $\Gamma$.
Let $ \Gamma = \Gamma(V,E)$ be a simple unoriented graph with vertex set $V_{\Gamma}$ and edge set $E_{\Gamma}$. Topologically, $\Gamma$ is a one-dimensional CW-complex obtained by `gluing together' the geometric $1$-simplices: $\,\Gamma = \Delta^1 \,\amalg_{\Delta^0} \Delta^1 \,\amalg_{\Delta^0} \ldots \,\amalg_{\Delta^0} \Delta^1\,$, one `$\Delta^1$' for each edge $ e \in E $. Given a function $ m: E_{\Gamma} \to \Z_{+}\,,\, e \mapsto m_e \,$ assigning multiplicities to the edges of $\Gamma$, we can construct a CW-complex $ \Gamma^{(m)} $ of dimension $ N = 1 + \max\{m_e\,:\, e \in E\}$) that `thickens' $ \Gamma$
in a natural way. To this end, we replace each $1$-dimensional simplex  $ \Delta^1 = S \Delta^{0} $ representing an edge $e$ in $\Gamma $ with the unreduced suspension $ S \Delta^{m_e} $ of the $ m_e $-dimensional simplex $\Delta^{m_e}$, where $ m_e $ is the multiplicity assigned to $e$.
As a result, we get the CW-complex
\begin{equation}
\la{Gm}
\Gamma^{(m)} \,:= \, S\Delta^{m_{e_1}}\,\amalg_{\Delta^0}\, S\Delta^{m_{e_2}}\,\amalg_{\Delta^0}\, \ldots \,\amalg_{\Delta^0}\, S\Delta^{m_{e_r}}   
\end{equation}
in which two  `fat' edges $ S \Delta^{m_{e_i}}$ and $ S\Delta^{m_{e_{i+1}}}\,$
are glued to each other at their suspension points according to the incidence relations of the original graph $\Gamma $ (see Figure~\ref{bicone-graph}). 

The CW-complex \eqref{Gm} described above has a simple categorical model. Recall that to each $m$-simplex $ \Delta^m $ there is a naturally associated poset  $ \sd(\Delta^m)$ --- called the  {\it barycentric subdivision} of $\Delta^m $ --- whose elements are all faces of $ \Delta^m $ ordered by inclusion. 
We identify this poset with the poset $\P_0([m])$ of all non-empty subsets of the ordinal $ [m] = \{0,1,\ldots,m\}$ by regarding the latter as the set of vertices of $ \Delta^m$. Consequently, we write the elements of $ \sd(\Delta^m)$ as subsets $ \sigma \subseteq [m] $. 

Now, given a graph $\Gamma $ with a multiplicity function $m: E_{\Gamma} \to \Z_+$, we consider the functor
\begin{equation}
\la{thickfun}
\sd^{(m)}(\Gamma):\ \Cc(\Gamma)\,\to\, \Cat\ , \quad e \mapsto 
{\sd}(\Delta^{m_e})^{\rm op} \ ,
\end{equation}
assigning to each edge object $ e \in \Cc(\Gamma) $ the opposite poset category of $ \sd(\Delta^m)$, while sending each vertex $ v \in \Cc(\Gamma) $ to the terminal (one-object) category  in $\Cat$. Then, we define the {\it $m$-simplicial thickening} of  $\Cc(\Gamma) $ to be the category
\begin{equation}
\la{thickcat}
\Cc^{(m)}(\Gamma)\, := \, \Cc(\Gamma) \smallint \sd^{(m)}(\Gamma)
\end{equation}
obtained by applying the Grothendieck construction to \eqref{thickfun} (see Definition~\ref{Grothcon}).
Note that $\Cc^{(m)}(\Gamma)$ comes with a canonical map $\, \pi^{(m)}:\, \Cc^{(m)}(\Gamma) \to \Cc(\Gamma)\,$, which is a cocartesian fibration 
(Grothendieck opfibration)  in 
${\sf Cat}$,
whose fibres over the edge objects $e \in \Cc(\Gamma) $ are the (opposite) subdivision posets $ {\rm sd}(\Delta^{m_e})^{\rm op} $. 

The category $\Cc^{(m)}(\Gamma)$ has the following explicit description.
The objects of $ \Cc^{(m)}(\Gamma)$ comprise the set of vertices of $\Gamma $ and the sets of ``fat'' edges $ \{e_{\sigma}\} $ labeled by the 
non-empty subsets $ \sigma \subseteq [m_e] $ for each edge $e \in E_{\Gamma} $:
$$
{\rm Ob}\,[\Cc^{(m)}(\Gamma)]\,=\, 
V_{\Gamma}\,\bigcup_{e \in E_{\Gamma}} \, \{e_{\sigma}\}_{\sigma \subseteq [m_e]} 
$$
The (non-identity) morphisms in $\Cc^{(m)}(\Gamma)$ are of two types: one arrow $\ e_{\sigma} \to v \,$ for each
$\,\sigma \subseteq [m_{e}]$, whenever the edge $ e $ is incident to the vertex $v$ in $\Gamma$, and one arrow $ e_{\sigma} \to e_{\tau} $ for each pair of subsets $ \sigma \supseteq \tau $ in $ [m_e] $ labeling the same edge $e$ in $\Gamma$. See Example~\ref{fatgr} below.

\vspace{1ex}

The next lemma shows that $ \Cc^{(m)}(\Gamma)$ is indeed a categorical model of  \eqref{Gm}.
\begin{lemma}
\la{sdlemma}
There is natural homotopy equivalence: $\, {B}[\Cc^{(m)}(\Gamma)]\, \simeq \,\Gamma^{(m)}\,$. 
 \end{lemma}
\begin{proof}
First, we recall that the classifying space $B\Cc$ of a small category $ \Cc $
is defined to be the geometric realization $ |N(\Cc)|$ of the simplicial nerve 
of $\Cc$. For $ \Cc = {\rm sd}(\Delta^m)^{\rm op} $, there is a natural cellular homeomorphism
$ B[{\rm sd}(\Delta^m)^{\rm op}] \,=\, B[{\rm sd}(\Delta^m)] \,\cong\, \Delta^m $ that takes a vertex $ \sigma = \{v_0, \ldots, v_k\}$ of the nerve of $ {\rm sd}(\Delta^m) $   
its ``barycenter'' $ \frac{1}{k+1}\,(v_0 + \ldots + v_k)\,$ (see \cite[Lemma III.4.1]{GJ99}). Then, by Thomason's Theorem (\cite[Theorem~1.2]{To79}), we conclude
\begin{eqnarray*}
B[\Cc^{(m)}(\Gamma)]\, = \, |\,N[\Cc^{(m)}(\Gamma)]\,| &\simeq &
|\,\hocolim^{\sset}_{\Cc(\Gamma)} (N \, \sd^{(m)}(\Gamma))\,| \\
&\simeq& \hocolim_{e \in \Cc(\Gamma)}[B({\rm sd}(\Delta^{m_e}))] \\
&\cong & \hocolim_{e \in \Cc(\Gamma)}[\Delta^{m_e}]\\
& \cong & \Gamma^{(m)}\,,
\end{eqnarray*}
where  $\sd^{(m)}(\Gamma) $ is the edge subdivision functor for $\Gamma $ defined by \eqref{thickfun}. 
\end{proof} 
Let us now specialize the above construction to the Bruhat moment graph $ \Gamma(\R_W) $. First, observe that $ \M(W)$ can be 
identified with the set of multiplicity functions on $\Gamma$, with
$ m \in \M(W) $ defining the function
$$
m:\, E_{\Gamma} \to \Z_+ \ ,\quad e(s_{\alpha}, w) \mapsto m_{\alpha}\,.
$$
Thus, for each $ m \in \M(W)$, we can define the CW-complex $\Gamma^{(m)}(\R_W) $ and the associated category $\Cc^{(m)}(\Gamma(\R_W))$, which we call the
$m$-thickened moment graph and moment category, respectively. 

For $ \Gamma = \Gamma(\R_W)$, the edge subdivision functor \eqref{thickfun}  has the canonical $W$-structure with respect to the $W$-action \eqref{Wact1}. Indeed, since $m$ is $W$-invariant, we can identify the posets $ \sd(\Delta^{m_{g(\alpha)}}) $ for all $ g \in W $ and define the natural transformations $ \vartheta_g:\, \sd^{(m)}(\Gamma) \to 
\sd^{(m)}(\Gamma) \circ \hat{g}\,$, specifying the $W$-structure on $\sd^{(m)}(\Gamma)$, to be the identity maps. The Grothendieck construction on a $W$-functor produces a $W$-category: hence, for each $m \in \M(W)$, the $m$-thickened moment category $\Cc^{(m)}(\Gamma)$ carries a natural $W$-action such that $ \pi^{(m)}: \Cc^{(m)}(\Gamma) \to \Cc(\Gamma)$ is a $W$-equivariant functor. Explicitly, this action is given by
\begin{equation}
\la{Wactm}
g\cdot w  = gw\ ,\quad  g\cdot e_{\sigma}(s_\alpha, w) = e_{\sigma}(gs_\alpha g^{-1}, gw)\,,\quad \forall\,\sigma \in \sd(\Delta^{m_{\alpha}})\,.
\end{equation}
The homotopy orbit category $ \Cc^{(m)}(\Gamma)_{hW} $ for this action fibers over
$ \Cc(\Gamma)_{hW}$: indeed, there is a canonical isomorphism of categories
\begin{equation}
\la{thickfunW}
 \Cc^{(m)}(\Gamma)_{hW} \,\cong\, 
 \Cc(\Gamma)_{hW} \smallint \sd^{(m)}(\Gamma)\,,
\end{equation}
defining a cocartesian fibration $  \pi^{(m)}_{hW}:\, \Cc^{(m)}(\Gamma)_{hW} \to \Cc(\Gamma)_{hW}$, whose classifying functor $\sd^{(m)}(\Gamma): \Cc(\Gamma)_{hW} \to \Cat $ is the extension of \eqref{thickfun} with respect to the $W$-structure described above.

\begin{example}
\la{fatgr}   
Let $ W = \langle s \rangle $ be a rank one Coxeter group generated by a reflection $s$. The Bruhat graph $\Gamma$ of $ W $ consists of a single edge $ e = e(s,1)$
connecting two vertices $1$ and $s$. The corresponding moment category is given by
the diagram $\, \Cc(\Gamma)\, = \,(\,1 \leftarrow e \rightarrow s\,) \,$. 
The simplicially thickened moment categories $ \Cc^{(m)}(\Gamma) $ for $m=1$ and $ m =2 $
are shown on Figure~\ref{CGamma}.
\end{example}
\begin{figure}[h!]
\centering
\begin{minipage}[c]{0.845\textwidth}\centering
\scalebox{0.82}{$\begin{tikzcd}[row sep=small, column sep=small, ampersand replacement=\&]
	1 \&\&\& 1 \&\&\&\& 1 \\
	\\
	\&\&\&\&\&\& {e_{\{0,1\}}} \&\& {e_{\{1\}}} \\
	{e_{\{0\}}} \&\& {e_{\{0\}}} \& {e_{\{0,1\}}} \& {e_{\{1\}}} \& {e_{\{0\}}} \&\& {e_{\{0,1,2\}}} \&\& {e_{\{1,2\}}} \\
	\&\&\&\&\&\& {e_{\{0,2\}}} \&\& {e_{\{2\}}} \\
	\\
	s \&\&\& s \&\&\&\& s
	\arrow[from=4-1, to=1-1]
	\arrow[from=4-1, to=7-1]
	\arrow[from=4-3, to=1-4]
	\arrow[from=4-3, to=7-4]
	\arrow[thmframe, from=4-4, to=1-4]
	\arrow[thmframe, from=4-4, to=4-3]
	\arrow[thmframe, from=4-4, to=4-5]
	\arrow[thmframe, from=4-4, to=7-4]
	\arrow[thmframe, from=4-5, to=1-4]
	\arrow[thmframe, from=4-5, to=7-4]
	\arrow[bend left=14, from=4-6, to=1-8]
	\arrow[bend right=14, from=4-6, to=7-8]
	\arrow[thmframe, from=3-7, to=1-8]
	\arrow[thmframe, from=3-7, to=3-9]
	\arrow[thmframe, from=3-7, to=4-6]
	\arrow[thmframe, from=3-7, to=7-8]
	\arrow[thmframe, from=3-9, to=1-8]
	\arrow[thmframe, from=3-9, to=7-8]
	\arrow[exframe!80!black, from=5-7, to=1-8]
	\arrow[exframe!80!black, from=5-7, to=4-6]
	\arrow[exframe!80!black, from=5-7, to=5-9]
	\arrow[exframe!80!black, from=5-7, to=7-8]
	\arrow[exframe!80!black, from=5-9, to=1-8]
	\arrow[exframe!80!black, from=5-9, to=7-8]
	\arrow[exframe!80!black, from=4-10, to=3-9]
	\arrow[exframe!80!black, from=4-10, to=5-9]
	\arrow[exframe!80!black, bend right=14, from=4-10, to=1-8]
	\arrow[exframe!80!black, bend left=14, from=4-10, to=7-8]
	\arrow[exframe!80!black, from=4-8, to=1-8]
	\arrow[exframe!80!black, from=4-8, to=3-7]
	\arrow[exframe!80!black, from=4-8, to=3-9]
	\arrow[exframe!80!black, from=4-8, to=4-6]
	\arrow[exframe!80!black, from=4-8, to=4-10]
	\arrow[exframe!80!black, from=4-8, to=5-7]
	\arrow[exframe!80!black, from=4-8, to=5-9]
	\arrow[exframe!80!black, from=4-8, to=7-8]
\end{tikzcd}$}
\end{minipage}\hfill
\begin{minipage}[c]{0.145\textwidth}\raggedright\footnotesize
\textcolor{black}{\rule{12pt}{1.5pt}}\;in $\Cc^{(0)}$\\[6pt]
\textcolor{thmframe}{\rule{12pt}{1.5pt}}\;new in $\Cc^{(1)}$\\[6pt]
\textcolor{exframe!80!black}{\rule{12pt}{1.5pt}}\;new in $\Cc^{(2)}$
\end{minipage}
\caption{The moment categories $\Cc^{(0)}(\Gamma)\,\subset\,\Cc^{(1)}(\Gamma)\,\subset\,\Cc^{(2)}(\Gamma)$ for $W=\Z/2\Z$.}
\label{CGamma}
\end{figure}
In Section~\ref{S4.1}, we defined the moment category $ \Cc(\Gamma) $ as the (opposite) face category of the graph $\Gamma $ viewed as a simplicial complex. In contrast, for $ m \ge 1 $, the `$m$-thickened' category $ \Cc^{(m)}(\Gamma) $ cannot be the face category of any simplicial complex (on the vertex set $V(\Gamma)$). Instead, in the next section, we will show that $ \Gamma^{(m)} $ carries the natural structure of a polyhedral complex of spherical type, and the category $ \Cc^{(m)}(\Gamma) $ can be then realized as the  scwol  associated to $ \Gamma^{(m)} $.

%

\subsection{$ \Gamma^{(m)}$ as a polyhedral complex}
\label{subsec:regCW}
Recall that a CW-complex $ K $ is called {\it regular} if for any $n$-dimensional cell $\sigma$, the characteristic map $ \Phi_{\sigma}: B^n \to K $ is a homeomorphism  onto its image. If $K$ is a regular CW-complex, the cells of $K$ form a poset $ \Cc(K)$ (called the cell category of $K$) with partial order $ \sigma \leq \sigma' $ defined by $\,\bar{\sigma} \supseteq \sigma'\,$, where $ \bar{\sigma}$ is the (CW) closure of the cell $ \sigma $ in $K$. 

We begin with the following simple observation.
\begin{lemma}\label{prop:regCW}
$\Gamma^{(m)}$ carries a canonical regular CW-structure whose cell category
is $\Cc^{(m)}(\Gamma)$. Consequently, there is a homeomorphism $B[\Cc^{(m)}(\Gamma)]\cong\Gamma^{(m)}$, refining the equivalence of Lemma~\ref{sdlemma}.
\end{lemma}

\begin{proof}
It suffices to specify a regular CW-structure on each `fat' edge 
$\, S\Delta^{m_e} $ of $ \Gamma^{(m)}$. We identify
$$
S\Delta^{m_e} = \left([-1,1] \times \Delta^{m_e}\right)/
\left(\{-1\} \times \Delta^{m_e}\right)/ \left(\{1\} \times \Delta^{m_e}\right)
$$
and put on it the cell structure comprising two $0$-cells $N_e,\, S_e$
(the suspension points of $S\Delta^{m_e}$), and the open cells $\,C_\sigma=\{(t,x)\in S\Delta^{m_e}:\ \operatorname{supp}(x)=\sigma,\,t\in(-1,1)\}\,$ of dimension $|\sigma|$, one for each non-empty subset $ \sigma\subseteq[m_e]$ (where $ {\rm supp}(x) := \{i \in [m_e]\,:\, x_i \not= 0\} \subseteq \Delta^{m_e} $).

Away from the suspension vertices, a point of $S\Delta^{m_e}$ has a well-defined simplex coordinate $x$ with support $\sigma(x)$, hence the cells $C_\sigma$ together with $\{N,S\}$ partition the bicone $S\Delta^{m_e}$. Each $C_\sigma\cong
\mathring\Delta_{\sigma}\times(-1,1)$ is an open $|\sigma|$-dimensional ball with embedded closure $ S\Delta_\sigma $, whose boundary is the cell subcomplex $\,\partial\overline{C_\sigma}=S(\partial\Delta_\sigma)\cup\{N_e,S_e\}
=\{N_e,S_e\}\cup\bigcup_{\tau\subsetneq\sigma}C_\tau$. Clearly, this 
CW-structure on $S\Delta^{m_e}$ is regular, and its cell poset is $\{N_e,S_e\}\sqcup\mathcal P_0([m_e])$. When glued along $\Gamma$
these cell posets yield exactly the poset category  $\Cc^{(m)}(\Gamma)$. 
The last claim follows from the fact that, for regular CW-complexes, the geometric realization of the cell poset is homeomorphic to the space.
\end{proof}

Next, we recall the notion of a $\kappa$-polyhedral complex, which is a natural generalization. Our main reference on polyhedral complexes is Chapter I.7 of the Bridson-Haefliger book \cite{BH99}.
\begin{defi}
\label{def:polyhedral}
Let $M_{\kappa}$ denote a complete, simply connected Riemannian manifold of constant sectional curvature $\kappa\in {\mathbb R}$.
An \emph{$\kappa$-polyhedral complex} $K$ is a finite-dimensional
CW-complex such that
\begin{itemize}
\item[(1)] each open cell of dimension $n$ is isometric to the interior of a
compact convex polyhedron in $M^n_\kappa$; 
\item[(2)] for each cell $\sigma$ of $K$, the restriction of the attaching
map to each open codimension-one face of $\sigma$ is an isometry onto an open
cell of $K$.
\end{itemize}
Such a complex is called \emph{Euclidean}, \emph{spherical}, or \emph{hyperbolic} if $\kappa=0$, $\kappa > 0$ or $\kappa < 0 $, respectively.
\end{defi}
Canonically associated to a polyhedral complex $K$ is its {\it scwol} $\Cc(K)$ --- a small category without loops\footnote{That is, a small category with only trivial endomorphisms.} --- whose objects are the cells of $K$  (see \cite[Chapter~III.C.1.4]{BH99}). If we take the (first) barycentric subdivision $\sd(K)$, the cells of $K$ can be identified with their barycenters, and the morphisms in $\Cc(K)$ are then defined to be the $1$-simplices in $\sd(K)$. By definition, each $1$-simplex corresponds to a pair of cells $ \tau \subset \sigma $, with (the barycenter of) the larger cell being the source of the morphism $ a: \sigma \to \tau $. 
Thus, the nonidentity morphisms in $\Cc(K)$ are the strict incidences $ \tau \subsetneq \sigma $ counted with multiplicity.
The compositions ($2$-composable maps) are determined by the $2$-simplices in $\sd(K)$, which correspond to triples of cells $ \tau \subset \sigma \subset \eta $ in $K$, etc. 

Note that a geodesic $m$-simplex (i.e., the convex hull of $(m+1)$  points  in $M_{\kappa}$ in general position) is naturally a polyhedral complex
(with ${m+1 \choose k+1}$ cells in dimension $k$); its scwol is the face poset of $\Delta^m$, with realization the barycentric subdivision. More generally, every simplicial complex can be realized as a polyhedral complex
(by gluing together geodesic simplices along isometric faces). On the other hand, a polyhedral complex need {\it not} be a simplicial complex, the simplest example being the `monogon' --- the circle
with one vertex (a single loop-cell) --- in which case the scwol is the parallel arrow category $ \{\,\bullet \rightrightarrows \bullet\,\}$.


Let $\Gamma=(V_\Gamma,E_\Gamma)$ be a simple unoriented graph with a multiplicity function $ m \colon E_\Gamma\to\Z_{\geq 0}$.
\begin{prop}\label{prop:noeuclid}
Assume that $m_e \ge 1$ for some edge $e \in E_{\Gamma}$.
Then there is no $\kappa$-polyhedral complex $K$ with $ \kappa\leq 0 $, whose scwol is isomorphic to $\Cc^{(m)}(\Gamma)$.
\end{prop}

\begin{proof}
We assign dimensions to the objects of $\Cc^{(m)}(\Gamma)$ by letting 
$ \dim(v) = 0 $ for every  $ v \in V_\Gamma $ and $ \dim(e_\sigma) =|\sigma| $ for a `fat' edge indexed by a subset $ \sigma \subseteq [m_e]$. 
%
%
For an edge $e \in E_{\Gamma}$ with $m_e \geq 1$, choose a subset $\sigma\subseteq[m_e]$ of cardinality $ d:=|\sigma|\ge 2 $ and consider the corresponding object $ e_{\sigma} $. Note that, since $ \dim(e_{\sigma}) = d \ge 2 $, there are exactly $d$  objects in $\Cc^{(m)}(\Gamma)$ of dimension $(d-1)$ under $e_\sigma\,$ (namely,  $ e_{\sigma\setminus\{i\}}$, one for each $i\in\sigma $), and hence, exactly $d$ `codimension-one morphisms'  with source $e_{\sigma} $.
Now, suppose that $ \Cc^{(m)}(\Gamma) $ is isomorphic to the scwol $ \Cc(K) $ of some $\kappa$-polyhedral complex $K$. Let $P_\sigma$
be the closed cell in $ \Cc(K) $ corresponding to $e_\sigma$.  Its dimension is $d$, and by definition, it is modeled on a compact, full-dimensional, convex polytope $P\subset M_\kappa^d$. Each facet of $P$ gives a codimension one incidence morphism with source
$P_\sigma$, and every such incidence arises from a facet. Hence\footnote{The facets are
counted with multiplicity: if distinct facets are identified with the same
cell of $K$, they still determine distinct arrows in the scwol $\Cc(K)$.}
\begin{equation}
\la{facecount}
 \#\{\text{facets of }P\}
 =\#\{\text{codimension-one arrows with source }e_\sigma\}\,
 =\,d\,.
\end{equation}
We show that this is impossible if $\kappa \le 0 $.

If $\kappa=0$, then $P$ is a compact, full-dimensional Euclidean
$d$-polytope; it admits an (irredundant) half-space presentation
$\, P=\bigcap_{i=1}^r  \{x\in {\mathbb R}^d:\langle u_i,x\rangle\leq c_i\} $, with one inequality for each facet.  Since $P$ is bounded, its outward
normal vectors $\{u_1,\ldots,u_r\}$ positively span ${\mathbb R}^d$.  A positive spanning
set in a $d$-dimensional vector space has at least $d+1$ elements.  Thus
$r\geq d+1$, contradicting \eqref{facecount}.
If $\kappa<0$, we rescale it to $\kappa=-1$ and use the classical Klein model $\,\mathbb H^d\cong B^d $, identifying the hyperbolic $d$-space with a $d$-dimensional ball in the Euclidean space $ {\mathbb R}^d $. The hyperbolic geodesics are then Euclidean line segments, the hyperbolic hyperplanes are intersections of Euclidean affine hyperplanes with $B^d$, and the hyperbolically convex sets are Euclidean convex. Because $P$ is compact in $\mathbb H^d$, its Klein image is contained in a smaller Euclidean ball.  It is therefore a compact Euclidean $d$-polytope, and its hyperbolic facets correspond exactly to its Euclidean facets.  It again has at least $d+1$ facets in 
contradiction to \eqref{facecount}.
\end{proof}
%
%
%
Thus, if $m\ge 1$, the category $\Cc^{(m)}(\Gamma)$ cannot be the scwol of any nonpositively curved polyhedral complex: in particular, it is not the face poset of any simplicial complex. This leaves the possibility of $ \kappa > 0 $, in which case the above `facet counting' argument does not work. In what follows, we will show that  $\Cc^{(m)}(\Gamma)$ can indeed be realized as the scwol of a spherical polyhedral complex $ \Gamma^{(m)}_{\rm sph}$, although defined in a slightly more general form
than in \cite{BH99}. 

Consider the $d$-dimensional sphere $S^d$ with two antipodal poles $N,\,S \in S^d$, and the equatorial hypersphere $S^{d-1} \subset S^d $.
Fix a spherical $(d{-}1)$-simplex $\tau$ contained in an open hemisphere of the equator $S^{d-1}$, and define 
\begin{equation*}
\Lambda_{\tau} :=\{N,S\}*\tau\;=\;\bigcup_{x\in\tau}
\bigl\{\text{meridian in $S^d$ passing through $x$ from $N$ to $S$}\}
\,.
\end{equation*}
If $ d=1 $, then $ \Lambda_\tau \subset S^1 $ is a semicircle with a single open $1$-cell and two vertices (facets) $N $ and $S$. In general, we have
\begin{lemma}\label{lem:lune}
Let $d\ge 2$. Then $\,\Lambda_\tau \subset S^d $ is a compact convex spherical polytope of dimension $d$. Its proper exposed\footnote{By an exposed face we mean the 
intersection of $\Lambda $ with its supporting great hypersphere in $S^d$.} faces are $\,\{N, \,S\}\,$ and $\Lambda_{\tau'} = \{N,S\} \ast \tau' $, where $ \tau' \subsetneq \tau $ ranges over the nonempty proper faces of $\tau $. There is a natural combinatorial face structure on $\Lambda_\tau$, in which $N$ and $S$ are considered as two distinct vertices, while the remaining proper faces are $\Lambda_{\tau'} = \{N,S\} \ast \tau' $ with $ \tau' \subsetneq \tau $. Equipped with this structure, $\Lambda_{\tau} $ has exactly two vertices and $d$ facets.
\end{lemma}
\begin{proof}We give an explicit model for $ \Lambda_{\tau}$. Let $ S^d $ be the unit sphere in $ {\mathbb R}^{d+1} $ with $\, N = e_0$, $\, S = - e_0 \,$, and  $\, S^{d-1} = S^d \,\cap\,e_0^{\perp} $. For $ \tau \subset S^{d-1}$, denote by $ C(\tau) := \{rx \in {\mathbb R}^{d+1}\,:\,r \ge 0,\,x\in \tau\} \subset e_0^{\perp} $ the Euclidean cone over $\tau$. Then
\begin{equation}
\la{lambdam}
\Lambda_{\tau} = S^d\,\cap\,\{{\mathbb R} e_0 + C(\tau)\} \,. 
\end{equation}
Indeed, every point of $\Lambda_{\tau} \setminus \{N,S\} $ can be written uniquely in the form
$$
y = t e_0 + \sqrt{1-t^2} \,x \ ,\quad 
-1 < t < 1 \ ,\quad  x \in \tau \ , 
$$
which gives (for a fixed $x$) a parametrization of the meridian 
passing through $x$ in $S^d$. A spherical $(d-1)$-simplex can be described as $\,\tau = \{x \in S^{d-1} : \langle a_i, x\rangle \ge 0\}\,$, where
$ \{a_1, a_2, \ldots, a_d\} \subset e_0^{\perp}$ are some linearly independent equatorial vectors. Then, by \eqref{lambdam},
$\,\Lambda_{\tau} = \{y \in S^{d} : \langle a_i, y\rangle \ge 0\,,\,\forall\,i=1,\ldots, d\}\,$, which shows that $ \Lambda_{\tau} $ is an intersection of closed hemispheres in $S^d$. Hence it is spherically convex and closed, and therefore compact (since $S^d$ is). As $ C(\tau)$ is a full-dimensional cone in the $d$-dimensional hyperplane $ e_0^{\perp}$, it also follows from \eqref{lambdam} that $ \dim(\Lambda_{\tau}) = d $ (see Figure~\ref{fig:lune}).
%
\begin{figure}[htbp]
\centering
\providecolor{thmframe}{RGB}{70,130,180}%
\providecolor{exframe}{RGB}{150,50,100}%
\def\Az{25}\def\El{22}                 
\begin{tikzpicture}[
  scale=2.2, line cap=round, line join=round,
  knot/.style={circle, fill=white, draw=black!65, semithick,
               inner sep=0pt, minimum size=4pt},
  declare function={
    ex(\t)    =  sin(\t-\Az);                              
    ey(\t)    = -sin(\El)*cos(\t-\Az);                     
    mx(\t,\u) =  cos(\u)*ex(\t);                           
    my(\t,\u) =  cos(\u)*ey(\t) + cos(\El)*sin(\u);        
    hz(\t)    = -atan(cos(\El)*cos(\t-\Az)/sin(\El));      
  }]
  \def\Ta{8}\def\Tb{64}\def\Tm{36}     
  \def\Tx{46}                          
  \def\Ed{1.8}                         
  \pgfmathsetmacro{\Ha}{hz(\Ta)}       
  \pgfmathsetmacro{\Hb}{hz(\Tb)}
  \pgfmathsetmacro{\Hv}{\Az-90}        

  \fill[white] (0,0) circle (1);

  \draw[black!45,dashed,thin]
       plot[domain=\Hv+180:\Hv+360,samples=120,variable=\t] ({ex(\t)},{ey(\t)});

  \fill[thmframe!16]
       plot[domain=90:-90,samples=100,variable=\u] ({mx(\Ta,\u)},{my(\Ta,\u)})
    -- plot[domain=-90:90,samples=100,variable=\u] ({mx(\Tb,\u)},{my(\Tb,\u)}) -- cycle;

  \fill[exframe!12] (0,0)
    -- plot[domain=\Ta:\Tb,samples=90,variable=\t] ({ex(\t)},{ey(\t)}) -- cycle;
  \foreach \t in {\Ta,\Tb}
    \draw[exframe!70,line width=.5pt] (0,0) -- ({\Ed*ex(\t)},{\Ed*ey(\t)});

  \draw[black!45,thin]
       plot[domain=\Hv:\Hv+180,samples=120,variable=\t] ({ex(\t)},{ey(\t)});

  \foreach \t/\h in {\Ta/\Ha,\Tb/\Hb}{
    \draw[thmframe!75!black,line width=1.1pt]
         plot[domain=\h:90,samples=100,variable=\u] ({mx(\t,\u)},{my(\t,\u)});
    \draw[thmframe!75!black,line width=1.1pt,dashed]
         plot[domain=-90:\h,samples=40,variable=\u] ({mx(\t,\u)},{my(\t,\u)});}

  \draw[exframe,line width=1.7pt]
       plot[domain=\Ta:\Tb,samples=120,variable=\t] ({ex(\t)},{ey(\t)});

  \draw[black,line width=.7pt] (0,0) circle (1);

  \foreach \a/\b/\c/\lab/\anc in
      {1/0/0/{x_0}/south, 0/1/0/{x_1}/{north east}, 0/0/1/{x_2}/west}{
    \pgfmathsetmacro{\px}{-sin(\Az)*\b + cos(\Az)*\c}
    \pgfmathsetmacro{\py}{cos(\El)*\a - sin(\El)*cos(\Az)*\b - sin(\El)*sin(\Az)*\c}
    \pgfmathsetmacro{\pn}{veclen(\px,\py)}
    \draw[black!45,dashed,line width=.45pt] (0,0) -- ({1.02*\px/\pn},{1.02*\py/\pn});
    \draw[->,black!55,line width=.5pt]
         ({1.05*\px/\pn},{1.05*\py/\pn}) -- ({1.30*\px/\pn},{1.30*\py/\pn});
    \node[black!55,font=\scriptsize,anchor=\anc]
         at ({1.34*\px/\pn},{1.34*\py/\pn}) {$\lab$};}

  \coordinate (N)  at (0,{cos(\El)});           \coordinate (S)  at (0,{-cos(\El)});
  \coordinate (V1) at ({ex(\Ta)},{ey(\Ta)});    \coordinate (V2) at ({ex(\Tb)},{ey(\Tb)});
  \coordinate (X)  at ({ex(\Tx)},{ey(\Tx)});
  \foreach \p in {N,S,V1,V2} \node[knot] at (\p) {};
  \fill[black!70] (0,0) circle (.7pt);   \fill[black!70] (X) circle (.7pt);

  \node[anchor=south west,font=\small]              at ([shift={(.035,.035)}]N)   {$N=e_0$};
  \node[anchor=north east,font=\small]              at ([shift={(-.035,-.035)}]S) {$S=-e_0$};
  \node[exframe,font=\scriptsize,anchor=south east] at ([shift={(-.025,.032)}]V1) {};
  \node[exframe,font=\scriptsize,anchor=south west] at ([shift={(.02,.045)}]V2)   {};
  \node[font=\scriptsize,anchor=north west]         at ([shift={(.02,-.025)}]X)   {$x$};
  \node[font=\scriptsize,anchor=east]               at (-.045,-.012)              {$0$};
  \node[exframe,font=\small,fill=white,inner sep=1.5pt]
       at ({.60*ex(\Tm)},{.60*ey(\Tm)}) {$C(\tau)$};

  \coordinate (TAU) at ({1.26*ex(30)},{1.26*ey(30)});
  \node[exframe,font=\small] at (TAU) {$\tau$};
  \draw[exframe,line width=.5pt,->]
       ([shift={(-.055,.02)}]TAU) -- ({1.03*ex(20)},{1.03*ey(20)});

  \node[thmframe!75!black,font=\small,anchor=west] at (1.02,.74) {$\Lambda_\tau=\{N,S\}*\tau$};
  \draw[black!55,line width=.45pt,->]
       (1.005,.71) .. controls (.589,.722) .. ({mx(\Tm,52)},{my(\Tm,52)});

  \coordinate (F1) at ({mx(\Ta,-40)},{my(\Ta,-40)});
  \coordinate (F2) at ({mx(\Tb,26)},{my(\Tb,26)});
  \node[thmframe!75!black,font=\scriptsize,anchor=east] at ([shift={(-.025,0)}]F1) {};
  \node[thmframe!75!black,font=\scriptsize,anchor=west] at ([shift={(.025,0)}]F2)  {};

  \node[font=\small,anchor=east]                at (-1.02,.62)  {$S^{d}$};
  \node[black!45,font=\scriptsize,anchor=east]  at (-1.04,-.02) {$S^{d-1}=S^{d}\cap e_0^{\perp}$};
\end{tikzpicture}
\caption{The lune polytope $\Lambda_\tau=\{N,S\}*\tau\subset S^{d}$ for $d=2$.
}
\label{fig:lune}
\end{figure}

Now, let $\, H_b =\{y \in S^d\,:\, \langle b,y\rangle = 0\,\}\,$ be a supporting great hypersphere of $\Lambda_{\tau} $ oriented so that 
$ \langle b,y\rangle \ge 0 $ for all $ y \in \Lambda_\tau$. Since
$ N = e_0 $ and $ S = - e_0 $ are both in $ \Lambda_\tau $, we have $ \langle b, e_0\rangle = 0 $. Thus every supporting great hypersphere of $\Lambda_{\tau} $ contains $ N $ and $S $. Furthermore, 
for any point $\,y = t e_0 + \sqrt{1-t^2} \,x $ in $ \Lambda_{\tau}\setminus\{N,S\}\,$ with $ x \in \tau $, we have
$ y \in H_b $ iff $ x \in H_b $. Hence 
$$
H_b\,\cap\,\Lambda_{\tau} \,=\,
\begin{cases} 
\{N,S\} \ast \tau_b  & \text{if } \tau_b \not= \varnothing \\ 
\{N,S\} & \text{if } \tau_b = \varnothing  
\end{cases}
$$
where $ \tau_b := \tau\,\cap\,H_b\,$. If $\tau_b$ is non-empty, it is an exposed face of the spherical simplex $\tau$.  This proves the desired classification of exposed faces; it also shows 
that a supporting great hypersphere  of $\Lambda_\tau$ cannot expose one pole without exposing the other, so neither singleton $\{N\}$ nor
$\{S\}$ is an exposed face.
For the combinatorial structure, we therefore choose a different
convention by defining a face of $\Lambda_{\tau}$ to be a  {\it connected component} of an exposed face. With this convention, the proper combinatorial faces of $\Lambda_{\tau}$ are $\{N\}$, $\{S\}$ and $\Lambda_{\tau'}$ for non-empty $ \tau'\subsetneq \tau $. In particular, when $ d\ge 2$, there are exactly two vertices and $d$ facets in $ \Lambda_\tau$.
\end{proof}
We call the spherical polytope $ \Lambda_{\tau} $ of  Lemma~\ref{lem:lune} the {\it lune polytope} of $\tau $. We then define a \emph{spherical lune complex} to be a cell complex obtained by gluing spherical lune polytopes along their pole vertices, where cells of diameter $\pi$ are allowed and the two antipodal poles of each lune are considered to be distinct combinatorial vertices. We remark that, strictly speaking, 
the spherical lune complexes are {\it not} spherical polyhedral complexes in the sense of Bridson-Haefliger (Definition~\ref{def:polyhedral}) as they violate some of the geometric conventions adopted in \cite{BH99}\footnote{Specifically, in \cite{BH99}, the faces of all polyhedra are required to be exposed faces, and the cells of a complex of curvature $ \kappa > 0 $ must lie in an open ball of radius $ \pi/(2\sqrt{\kappa})$. These conditions are needed for the proof of \cite[Theorem~I.7.19]{BH99} on geodesic completeness of polyhedral complexes, which is not used in this paper.}. However, the metric space structure and the notion of a scwol extend naturally to spherical lune complexes.

Let $\Gamma$ be a graph with multiplicity function $m = \{m_e\}_{e \in E_{\Gamma}}$. We realize each `fat' edge of $ \Gamma^{(m)} $ as a lune polytope $\Lambda_{e} := \{N_e,S_e\}*\tau_e \subset S^{m_e+1}$ by 
choosing a spherical $m_e$-simplex $\tau_e \subset S^{m_e} $ contained in an open hemisphere, index the vertices of $\tau_e $ by $[m_e]$, and for each non-empty $ \sigma \subseteq [m_e] $, put $\,\Lambda_{e,\sigma}:=\{N_e,S_e\}*\tau_{e,\sigma} \subset \Lambda_e \,$, where $\tau_{e,\sigma} $ be the face of $ \tau_e$ corresponding to $ \sigma $. We then equip $\Lambda_e$ with the combinatorial lune structure of Lemma~\ref{lem:lune}, the vertices being $N_e,S_e$ and the remaining closed cells being $\Lambda_{e,\sigma}$. 
Gluing the lune polytopes $\Lambda_{e} $ at their poles according to the incidence relations of $\Gamma$, identifying $N_e,S_e$ with the endpoints of $e$, we get a spherical lune complex 
$\Gamma_{\mathrm{sph}}^{(m)}$.
\begin{theorem}
\label{thm:lune}
There is a cellular homeomorphism $\,\Gamma_{\mathrm{sph}}^{(m)}\, \cong \,\Gamma^{(m)}$, where $\Gamma^{(m)}$ is equipped with the regular CW-structure of Lemma~\ref{prop:regCW}. Consequently,
the scwol $\Cc[\Gamma_{\mathrm{sph}}^{(m)}]$ of the lune complex $\Gamma_{\mathrm{sph}}^{(m)}$ is naturally isomorphic to the category $\, \Cc^{(m)}(\Gamma)\,$. 
If $ \Gamma = \Gamma(\cR_W) $ is the Bruhat graph and $m$ is $W$-invariant, the $W$-action \eqref{Wactm}  is realized by cellular isometries of $\Gamma^{(m)}_{\mathrm{sph}}\,$: in particular, each $s_\alpha$ acts on a lune corresponding to a `fat' edge $e_{\sigma}(s_{\alpha}, w)$
by reflection in its equatorial great hypersphere, interchanging its two poles.
\end{theorem}
%
%
\begin{proof}
By Lemma~\ref{lem:lune}, each $\Lambda_e$ is compact, spherically convex polytope of dimension $m_e+1$. For $\varnothing\neq\sigma\subseteq[m_e]$, we consider the cells $\Lambda_{e,\sigma}$ and their interiors
$\,
 \Lambda_{e,\sigma}^{\circ}
 :=\Lambda_{e,\sigma}\!\setminus\!
 \bigl(
   \{N_e,S_e\}
   \cup
   \bigcup_{\varnothing\neq\tau\subsetneq\sigma}
          \Lambda_{e,\tau}
 \bigr)\,$.
Using the meridian coordinates $(t,x)$ introduced in the proof of Lemma~\ref{lem:lune}, we can identify 
$$
 \Lambda_{e,\sigma}^{\circ}
 \cong(-1,1)\times\mathring\Delta^{|\sigma|-1}
 \cong {\mathbb R}^{|\sigma|}\,.
$$
This shows that each $\Lambda_{e,\sigma}^{\circ}$ is an open cell of dimension $|\sigma|$, with closure
$\Lambda_{e,\sigma}$, and the boundary 
\begin{equation}\label{eq:lune-boundary}
 \partial\Lambda_{e,\sigma}
 =\{N_e,S_e\}
  \cup
  \bigcup_{\varnothing\neq\tau\subsetneq\sigma}
       \Lambda_{e,\tau}^{\circ}.
\end{equation}
Hence, the objects of $\Cc[\Gamma_{\mathrm{sph}}^{(m)}]$ that lie strictly under
$\Lambda_{e,\sigma}^{\circ}$ are the two pole vertices and all cells
$\Lambda_{e,\tau}^{\circ}$ with $\varnothing\neq\tau\subsetneq\sigma$.
If $|\sigma|\geq 2$, the codimension-one incidences under
$\Lambda_{e,\sigma}^{\circ}$ are precisely the morphisms
$\, \Lambda_{e,\sigma}^{\circ} \rightarrow \Lambda_{e,\sigma\setminus\{i\}}^{\circ} \,$ with $ i \in \sigma $, while in the case $|\sigma|=1$, the cell is one-dimensional and its codimension-one incidences are the pole vertices themselves.

Now, if we glue the lunes at their poles according to the incidence relations of $\Gamma$, different lunes will meet at pole vertices only.  
The resulting complex therefore has one vertex for every $v\in V_\Gamma$ and one open cell $\Lambda_{e,\sigma}^{\circ}$ for every $ e\in E_\Gamma $ and $\,\varnothing\neq\sigma\subseteq[m_e]\,$. By \eqref{eq:lune-boundary}, the strict incidences in $ \Lambda_{e} $ are $\, \Lambda_{e,\tau}^{\circ} \subsetneq \Lambda_{e,\sigma}^{\circ}\,$, one for each pair of non-empty subsets $\tau\subsetneq \sigma \subset [m_e] $. Hence, for each $e \in E_{\Gamma}$, there is a cellular homeomorphism $ \Lambda_{e} \xrightarrow{\sim} S\Delta^{m_e} $, identifying 
$  \Lambda_{e,\sigma}^{\circ} \cong C_\sigma $, where $C_\sigma $ are the open cells of $ S\Delta^{m_e} \subseteq \Gamma^{(m)}$  defined in (the proof of) Lemma~\ref{prop:regCW}. This extends to a homeomorphism of CW-complexes $\,\Gamma_{\rm sph}^{(m)} \xrightarrow{\sim} \Gamma^{(m)}\,$, which, in turn, induces an isomorphism of categories 
\begin{equation}
\la{isocats}
 \Cc\bigl[\Gamma_{\mathrm{sph}}^{(m)}\bigr]\,
\xrightarrow{\sim}\,\Cc^{(m)}(\Gamma)\,.
\end{equation}
The non-identity morphisms in 
$ \Cc[\Gamma_{\mathrm{sph}}^{(m)}] $ are exactly
$\,\Lambda_{e,\sigma}^{\circ}\rightarrow v\,$, when $v$ is an endpoint of $e$, and $\, \Lambda_{e,\sigma}^{\circ} \rightarrow\Lambda_{e,\tau}^{\circ}\,$ when $\varnothing\neq\tau\subsetneq\sigma $; thus, \eqref{isocats} is given explicitly by $\, v\mapsto v\,$ and $\, \Lambda_{e,\sigma}^{\circ} \mapsto e_\sigma\,$.

Suppose now that $\Gamma=\Gamma(\R_W)$ and $m$ is $W$-invariant.  Choose the
simplices $\tau_e$ so that lunes associated to edges in the same $W$-orbit
are isometric, with the face labels $\sigma$ preserved.  For $g\in W$, map
the lune associated to $e(s_\alpha,w)$ isometrically onto the lune
associated to
\[
 g\cdot e(s_\alpha,w)
 =e(gs_\alpha g^{-1},gw),
\]
sending the cell indexed by $\sigma$ to the cell indexed by the same
subset.  On poles this sends the endpoints $w,s_\alpha w$ to
$gw,gs_\alpha w$, so these isometries respect the gluing and define a
cellular $W$-action.

For $g=s_\alpha$, the underlying unoriented edge of type $\alpha$ is fixed
and its endpoints are interchanged.  In terms of meridian coordinates on 
the lune, the corresponding map is given by
\[
 t e_0+\sqrt{1-t^2}\,x
 \,\longmapsto
\, -t e_0+\sqrt{1-t^2}\,x.
\]
It is reflection in the equatorial great hypersphere $e_0^\perp$; it fixes
$\tau_e$ pointwise and interchanges $N_e$ and $S_e$.  This proves the last
assertion.
\end{proof}

\begin{remark}
\la{Rsimpcom}
The `$m$-thickening' construction of Section~\ref{S5.1} extends naturally to polyhedral complexes. Given such a complex $K$, we define a {\it multiplicity functor} on $K$ to be a contravariant functor $\, m: \Cc(K) \to 
(\Z_+,\, \le) \,$, assigning to each cell $ \sigma \in K $ a nonnegative integer $ m_{\sigma} $ in such a way that $ m_\tau \le m_{\sigma} $ whenever $\tau \subsetneq \sigma $, with convention that $ m_{v} = 0 $ for all vertices ($0$-cells) in $K$. (In addition, if a group $W$ acts on $K$ by cellular isometries, then $\Cc(K)$ carries a strict $W$-action by scwol automorphisms (see \cite[Definition III.${\mathcal C}$.1.11]{BH99}),
we also require $m$ to be $W$-invariant.) 
Next, we consider the subdivision functor $\, \sd^{\rm op}: \Z^{\rm op}_+ \to \Cat $, $\,m \mapsto \sd(\Delta^m)^{\rm op} $ that assigns to $\, n \le m \,$ the min-truncation map $\, \sd(\Delta^m) \to  \sd(\Delta^n) \,$ induced by the natural projection $ \Delta^m \onto \Delta^n $ collapsing all extra vertices of $ \Delta^m $ to the last vertex of $ \Delta^n $. Combining this with the multiplicity functor, we get the (covariant) functor $\,\sd^{(m)}(K): \ \Cc(K) \xrightarrow{m} \Z^{\rm op}_+ \xrightarrow{\sd^{\rm op}} \Cat\,$, and we define $\Cc^{(m)}(K)$ to be the Grothendieck construction on $ \sd^{(m)}(K)$:
$$
\Cc^{(m)}(K) := \Cc(K) \smallint \sd^{(m)}(K)
$$
When $ K = \Gamma $ is a graph (viewed a $1$-dimensional simplicial complex), the above definition reduces precisely to  \eqref{thickcat}. 
\end{remark}

\subsection{Generalized orbit functors}
\la{S5.2}
The next step of our construction is to define $W$-functors on $\Cc^{(m)}(\Gamma)\,$, generalizing the orbit functors \eqref{F0} and \eqref{GF0}. 
The $T$-obit functor \eqref{F0} extends to $\Cc^{(m)}(\Gamma)\,$ by
\begin{equation}
\la{Fm}
\cF_m\,:\ \Cc^{(m)}(\Gamma)\,\to\, \Top^T\ ,\quad w \mapsto \{w\}\ ,\quad 
e_{\sigma}(s_\alpha, w) \,\mapsto\,\O^{\times \sigma}_{e(s_\alpha,w)}\ ,
\end{equation}
where $\,\O^{\times \sigma}_{e(s_\alpha,w)} := \prod_{\sigma}\O_{e(s_\alpha,w)} $ denotes the product 
of copies of the $\bT$-orbit $\O_{e(s_\alpha,w)}$ indexed by elements of the subset $ \sigma \subseteq [m_{\alpha}]$, with  $T \subset \bT $ acting diagonally. The arrows $ e_{\sigma}(s_{\alpha}, w) \to e_{\tau}(s_{\alpha}, w')$ are mapped by $ \cF_m $ to the canonical projections 
$\,\O^{\times \sigma}_{e(s_\alpha,w)}\to \O^{\times \tau}_{e(s_\alpha,w)}\,$ 
indexed by $ \tau \subseteq \sigma $. We note that the above construction of
generalized orbit functors \eqref{Fm} is reminiscent of polyhedral products of spaces over simplicial complexes in toric topology (see, e.g., \cite{BP15}, \cite{BBC20}), one difference being that the categories $ \Cc^{(m)}(\Gamma)$ are not simplicial complexes in general (see Remark~\ref{Rsimpcom}).

Just as in the case $m=0$, the functor \eqref{Fm} itself is  not a $W$-functor, but 
its extension to $G$-spaces is:
\begin{lemma} \la{WstrT}
For each multiplicity $ m \in \M(W) $, the functor 
\begin{equation}
\la{FFGm}
G \times_T \cF_m\,:\ \Cc^{(m)}(\Gamma)\,\to\, \Top^G\ ,\quad
e_{\sigma}(s_\alpha, w) \,\mapsto\, G \times_T \O^{\times \sigma}_{e(s_\alpha,w)} \,,
\end{equation}
admits a natural $W$-structure  with respect to the action \eqref{Wactm} of $W$ on $\Cc^{(m)}(\Gamma)$.
\end{lemma}
\begin{proof}
The proof is similar to that of Lemma~\ref{lem:W-fun-heq}.
The $W$-structure on  \eqref{FFGm} is defined by the family of natural transformations $\, \{\,\vartheta^{(m)}_w:\, G \times_T \cF_m \to (G \times_T \cF_m) \circ \hat{w}\,\}_{w \in W} \,$, with components 
\begin{equation}
\la{eq:tildeFm-ver1}
\vartheta^{(m)}_{w,v}:\ G \times_T \{v\} \to G \times_T \{wv\}\,,\quad
[(g,v)]_T \mapsto [(g \dot{w}^{-1}, wv)]_T\,,
\end{equation}
where $ \dot{w} \in N $ denotes a representative of $ w \in W $ and $ v \in V_{\Gamma} = W $, 
\begin{equation}
\la{eq:tildeFm-edge1}
\vartheta^{(m)}_{w,e}:\ G \times_T \O^{\times \sigma}_e \to G \times_T \O^{\times \sigma}_{we}\,,\quad
[(g, (x_i)_{i \in \sigma})]_T \mapsto [(g \dot{w}^{-1}, \, (\dot{w} \cdot x_i)_{i \in \sigma})]_T\,,
\end{equation}
where $ e \in E_{\Gamma} $, $ (x_i)_{i \in \sigma} \in \O^{\times \sigma}_e $, and 
$we$ denotes the action of $ w\in W $ on $e \in E_{\Gamma}$ as in \eqref{Wact1}. The `dot' product $ \dot{w} \cdot x_i $ on $\bT$-orbits  is obtained by restricting to $N$ the natural (left) $G$-action on $ \bF = G/T $ (see \cite{Kn03}). 
\end{proof}
Using the $W$-structure of Lemma~\ref{WstrT}, we extend $G \times_T \cF_m$ to a 
functor $\,\Cc^{(m)}(\Gamma)_{hW}\,\to\, \Top^G\,$ and define
\begin{equation}
\la{FmG}
F_m(G,T)\,:=\, \hocolim_{\Cc^{(m)}(\Gamma)_{hW}}\,(G \times_T \cF_{m})\,,
\end{equation}
where the homotopy colimit is taken in the category of $G$-spaces.
Note, for $ m = 0 $, by Lemma~\ref{LemmaL4}, we have
$$
F_0(G,T) = \hocolim_{\Cc(\Gamma)_{hW}}(G\times_T \cF_0) \,\simeq\,
\hocolim_{\Cc(\Gamma)_{hW}}(G\times_T \wcF_0) \,\simeq\, G\times_N \bF^{(1)}\,,
$$
where $ \bF^{(1)}$ is the $T$-equivariant $1$-skeleton of $ \bF = G/T$, see \eqref{bF1}.

The spaces \eqref{FmG} come with natural $G$-equivariant maps
\begin{equation}
\la{phimm}
 \varphi_{m, m'}:\, F_m(G,T)\,\to \, F_{m'}(G,T) 
\end{equation}
defined for each pair of multiplicities $ m \leq m' \,$ in $\M(W)$.
The maps \eqref{phimm} are compatible in the sense that $ \varphi_{m,m''} =  
\varphi_{m',m''} \circ \varphi_{m,m'} $
for all $\, m \leq m' \leq m'' \,$ in $ \M(W)$. To see this it suffices to 
observe that the $W$-functors \eqref{FFGm} defining the spaces $ F_m(G,T)$
are related to each other by 
\begin{equation}
\la{GTmm}
G \times_T \cF_{m} = (G \times_T \cF_{m'}) \circ \Phi_{m,m'}\ ,\quad\forall\,m \leq m'\ ,
\end{equation}
where
\begin{equation}
\la{Phmm}
\Phi_{m,m'}:\ \Cc^{(m)}(\Gamma) \into  \Cc^{(m')}(\Gamma) \ ,\quad e_{\sigma}(s_{\alpha}, w) \mapsto e_{\sigma}(s_{\alpha}, w)\ ,
\end{equation}
are the full embeddings of categories defined by the natural inclusions of posets $ \sd(\Delta^{m_{\alpha}}) \subseteq \sd(\Delta^{m_{\alpha}'})$ for  $ m_\alpha \leq m_{\alpha}' $. In view of \eqref{GTmm}, we have canonical maps
of homotopy colimits
$$
\Phi^*_{m,m'}:\ \hocolim_{\Cc^{(m)}(\Gamma)_{hW}}\,(G \times_T \cF_{m}) \,\to\,\hocolim_{\Cc^{(m')}(\Gamma)_{hW}}\,(G \times_T \cF_{m'})
$$
which define the maps \eqref{phimm}. By definition, the functors $ \Phi_{m,m'}$ satisfy the cocycle conditions: $ \Phi_{m,m''} = \Phi_{m',m''} \circ \Phi_{m,m'}  $, inducing the corresponding conditions for \eqref{phimm}. Thus, the spaces \eqref{FmG} assemble to a $\M(W)$-diagram of $G$-spaces which we denote by
\begin{equation}
\la{towFm}
F_*(G,T):\ \M(W) \to \Top^G\ ,\quad m \mapsto F_m(G,T)\ ,\quad
(m \leq m') \,\mapsto \,\varphi_{m,m'}\,.
\end{equation}

To describe the properties of the diagram \eqref{towFm} we construct  --- 
in addition to \eqref{FmG} --- two other homotopy decompositions for the spaces $F_m(G,T)$: one is over the homotopy orbit category $ \Cc(\Gamma)_{hW} $, or rather its skeletal subcategory $\overline{\Cc(\Gamma)}_{hW}$ (see Lemma~\ref{CGdec}) and the other over the basic reflection category $ \Sc(W) $ of $W$ (see Proposition~\ref{natdec}). 
Each of these decompositions has its own advantages and will be used repeatedly in our calculations. 

We begin by giving a more explicit description of the homotopy orbit category $ \Cc(\Gamma)_{hW} $. 
\begin{lemma}
\la{skCG}
For any finite Coxeter group $W$, the category $  \Cc(\Gamma)_{hW} $ associated to the Bruhat graph $\Gamma = \Gamma(\R_W) $
is equivalent to its full subcategory $ \overline{\Cc(\Gamma)}_{hW} $ 
spanned by the objects $ 1 \in W $ and $ e(s_{\alpha}) := 
e(s_{\alpha}, 1) $ for all $\,\alpha \in \R_+$. The non-identity
morphisms in $ \overline{\Cc(\Gamma)}_{hW} $ are depicted by arrows on Figure~\ref{skC}; the generating relations are $\, s_{\alpha}^2 = \id \,$ for all $ \alpha \in \R_+ $.
\end{lemma}
\begin{figure}[h!]
\[\begin{tikzcd}[scale cd= 0.9]
	&& {e(s_{\beta})} \\
	{e(s_{\alpha})} &&&& {e(s_{\gamma})} \\
	&& 1 \\
	{e(s_{\nu})} &&&& {e(s_{\delta})} \\
	&& {e(s_{\lambda})}
	\arrow["{s_{\beta}}"{description}, from=1-3, to=1-3, loop, in=55, out=125, distance=10mm]
	\arrow["{f_\beta}"{description}, curve={height=6pt}, from=1-3, to=3-3]
	\arrow["{f_\beta s_\beta}"{description}, shift left, curve={height=-6pt}, from=1-3, to=3-3]
	\arrow["{s_{\alpha}}"{description}, from=2-1, to=2-1, loop, in=100, out=170, distance=10mm]
	\arrow["{f_\alpha}"{description}, shift left, curve={height=6pt}, from=2-1, to=3-3]
	\arrow["{f_\alpha s_\alpha}"{description}, shift left, curve={height=-6pt}, from=2-1, to=3-3]
	\arrow["{s_\gamma}"{description}, from=2-5, to=2-5, loop, in=10, out=80, distance=10mm]
	\arrow["{f_\gamma}"{description}, curve={height=6pt}, from=2-5, to=3-3]
	\arrow["{f_\gamma s_\gamma}"{description}, curve={height=-6pt}, from=2-5, to=3-3]
	\arrow["{f_\nu}"{description}, curve={height=6pt}, from=4-1, to=3-3]
	\arrow["{f_\nu s_\nu}"{description}, curve={height=-6pt}, from=4-1, to=3-3]
	\arrow["{s_\nu}"{description}, from=4-1, to=4-1, loop, in=260, out=190, distance=10mm]
	\arrow["{f_\delta}"{description}, shift left, curve={height=6pt}, from=4-5, to=3-3]
	\arrow["{f_\delta s_\delta}"{description}, shift left, curve={height=-6pt}, from=4-5, to=3-3]
	\arrow["{s_{\delta}}"{description}, from=4-5, to=4-5, loop, in=350, out=280, distance=10mm]
	\arrow["{f_\lambda}"{description}, curve={height=6pt}, from=5-3, to=3-3]
	\arrow["{f_\lambda s_\lambda}"{description}, shift left, curve={height=-6pt}, from=5-3, to=3-3]
	\arrow["{s_{\lambda}}"{description}, from=5-3, to=5-3, loop, in=305, out=235, distance=10mm]
\end{tikzcd}\]
\caption{The skeleton of  $ \Cc(\Gamma)_{hW} $}
\label{skC}
\end{figure}
\begin{proof}
Recall the category $\Cc(\Gamma)_{hW}$ is defined by \eqref{ChW}; the non-identity morphisms in this category are listed explicitly in \eqref{morChW}. It is immediate from the first equation in \eqref{morChW} that $\, w \cong 1 \,$ in $ \Cc(\Gamma)_{hW}$ for all objects $w$, while from the second equation it follows that
$$
e(s_\alpha, w) \,\cong\, e(s_{\alpha'}, w')\ \Leftrightarrow\ w^{-1} s_{\alpha} w = (w')^{-1} s_{\alpha'} w' \ \Leftrightarrow\ s_{w^{-1}(\alpha)} = s_{(w')^{-1}(\alpha')}
$$
This shows that $\, e(s_{\alpha}, w) \cong e(s_{w^{-1}(\alpha)})\,$ for all objects
$e(s_{\alpha}, w) \in \Cc(\Gamma)_{hW}$ and $\, e(s_{\alpha}) \not\cong  e(s_{\alpha'})\,$ unless $ \alpha = \alpha' $ in $ \R_+$; moreover, we have
$ \End_{\Cc(\Gamma)_{hW}}[e(s_{\alpha})] = W_{\alpha} $ for all $ \alpha \in \R_+$
Finally, the last equations on the list \eqref{morChW} show that the Hom-sets
$ \Hom_{\Cc(\Gamma)_{hW}}(e(s_{\alpha}),\,1) $ can be identified with the right cosets $ 1 \cdot W_{\alpha} $ of $W_\alpha$, whose elements we denote by $\{f_{\alpha}, \, f_{\alpha} s_{\alpha}\}$.  Thus, the full subcategory $ \overline{\Cc(\Gamma)}_{hW} $ with objects 
$ \{1\} \cup  \{e(s_{\alpha})\,:\,\alpha \in \R_+\} $ depicted on 
Figure~\ref{skC} is skeletal for $\Cc(\Gamma)_{hW}$.
\end{proof}
\begin{remark}
\la{comp_groups}
The skeletal category $ \overline{\Cc(\Gamma)}_{hW} $ has a  natural geometric interpretation. Let us regard the moment graph $ \Gamma  $ as a polyhedral complex and identify the moment category $\Cc(\Gamma)$ with its scwol. An easy check shows that the action of $W$ on $\Gamma$ induces an action on the {\it opposite} scwol $ \Cc(\Gamma)^{\rm op}$ by scwol automorphisms which are admissible in the sense of \cite[Def. III.C.1.11]{BH99}. The (na\"ive) categorical quotient $\Cc(\Gamma)^{\rm op}\!/W $ is then a scwol that comes equipped with a canonical functor ${\mathscr G}: \,\Cc(\Gamma)^{\rm op}\!/W\, \to\, {\rm Gr} \,$, $\, [e(s_{\alpha}, w)] \mapsto W_{\alpha} \,$,
encoding the isotropy groups of the quotient orbifold $ [\Gamma/W] $.
This functor defines a {\it complex of groups} in the sense of A. Haefliger \cite{Ha92}, and it turns out that its classifying category $ B{\mathscr G} $ is dual (contravariantly isomorphic) to $ \overline{\Cc(\Gamma)}_{hW} $. With Theorem~\ref{thm:lune}, this duality extends to arbitrary multiplicities:
$$
\overline{\Cc^{(m)}(\Gamma)}_{hW}\, \cong 
\, (B{\mathcal G}^{(m)})^{\rm op}
$$
where $ {\mathcal G}^{(m)} $ is a complex of groups associated 
to the quotient spherical orbifold $\, [\Gamma_{\rm sph}^{(m)}/W]\,$.
\end{remark}

Next, for each $m \in \M(W)$, we define the functor 
\begin{equation}
\la{Fm*}
\cF_{m}^{\,*}:\ \Cc(\Gamma) \to \Top^T\ ,\quad w \mapsto \{w\}\ ,\quad 
e(s_{\alpha}, w) \mapsto  \O_{e(s_{\alpha}, w)}^{\ast (m_{\alpha} + 1)}\,,
\end{equation}
and consider its extension to $G$-spaces parallel to \eqref{GF0}:
\begin{equation}
\la{FmG*}
G \times_T \cF_{m}^{\,*}:\ \Cc(\Gamma)_{hW} \to \Top^G\ ,\quad
e(s_{\alpha}, w) \mapsto G \times_{T} \O_{e(s_{\alpha}, w)}^{\ast (m_{\alpha} + 1)}\,,
\end{equation}
where $ \O_{e(s_{\alpha}, w)}^{\ast (m_{\alpha} + 1)} $ denotes the iterated ($(m_{\alpha}+1)$ times) topological join of the orbit spaces $\O_{e(s_{\alpha},w)}$
equipped with diagonal $T$-action (see, for example, \cite[Appendix A]{BR1}).
The same construction as in Lemma~\ref{lem:W-fun-heq} shows that
\eqref{FmG*} is indeed a well-defined $W$-functor with $W$-structure induced by the diagonal action of $ N_{\alpha}$ on  $ \O_{e(s_{\alpha}, w)}^{\ast (m_{\alpha} + 1)} $ extending the $T$-action. 
The restriction of \eqref{FmG*} to the
subcategory $ \overline{\Cc(\Gamma)}_{hW}$ is represented
by the diagram \eqref{dGTFm}, where the action of the reflections 
$ s_{\alpha} \in W_{\alpha} $, $\, s_{\beta} \in W_{\beta} $, $\ldots$
are induced by the diagonal action of $ N_{\alpha}, \,N_{\beta}, \ldots $ 
on  $ \O_{e(s_{\alpha})}^{\ast (m_{\alpha} + 1)}$,
$ \,\O_{e(s_{\beta})}^{\ast (m_{\beta} + 1)}, \ldots $ while the maps $ f_{\alpha},\,f_{\beta}, \ldots $ are defined by contracting 
these joins to the basepoint $[1]_T $ in $ G/T$.
\begin{equation}
\la{dGTFm}
\begin{tikzcd}[scale cd= .9, column sep=small]
	&& {G \times_T {\O}^{\ast(m_{\beta}+1)}_{e(s_{\beta})}} \\
	{G \times_T {\O}^{\ast(m_{\alpha}+1)}_{e(s_{\alpha})}} &&&& {G \times_T {\O}^{\ast(m_{\gamma}+1)}_{e(s_{\gamma})}} \\
	&& {G/T} \\
	{G \times_T {\O}^{\ast(m_{\nu}+1)}_{e(s_{\nu})}} &&&& {G \times_T {\O}^{\ast(m_{\delta}+1)}_{e(s_{\delta})}} \\
	&& {G \times_T {\O}^{\ast(m_{\lambda}+1)}_{e(s_{\lambda})}}
    \arrow["{s_{\beta}}", from=1-3, to=1-3, loop, in=55, out=125, distance=6mm]
	\arrow["{f_\beta}"{description}, curve={height=6pt}, from=1-3, to=3-3]
	\arrow["{f_\beta s_\beta}"{description}, shift left, curve={height=-6pt}, from=1-3, to=3-3]
    \arrow["{s_{\alpha}}", from=2-1, to=2-1, loop, in=100, out=170, distance=6mm]
	\arrow["{f_\alpha}"{description}, shift left, curve={height=6pt}, from=2-1, to=3-3]
	\arrow["{f_\alpha s_\alpha}"{description}, shift left=2, curve={height=-6pt}, from=2-1, to=3-3]
    \arrow["{s_{\gamma}}", from=2-5, to=2-5, loop, in=10, out=80, distance=6mm]
	\arrow["{f_\gamma}"{description}, curve={height=6pt}, from=2-5, to=3-3]
	\arrow["{f_\gamma s_\gamma}"{description}, shift left, curve={height=-6pt}, from=2-5, to=3-3]
	\arrow["{f_\nu}"{description}, shift right, curve={height=6pt}, from=4-1, to=3-3]
	\arrow["{f_\nu s_\nu}"{description}, curve={height=-6pt}, from=4-1, to=3-3]
    \arrow["{s_{\nu}}", from=4-1, to=4-1, loop, in=250, out=200, distance=6mm]
	\arrow["{f_\delta}"{description}, shift left=2, curve={height=6pt}, from=4-5, to=3-3]
	\arrow["{f_\delta s_\delta}"{description}, shift left=2, curve={height=-6pt}, from=4-5, to=3-3]
    \arrow["{s_{\delta}}", from=4-5, to=4-5, loop, in=340, out=290, distance=6mm]
	\arrow["{f_\lambda}"{description}, curve={height=6pt}, from=5-3, to=3-3]
	\arrow["{f_\lambda s_\lambda}"{description}, shift left, curve={height=-6pt}, from=5-3, to=3-3]
    \arrow["{s_{\lambda}}", from=5-3, to=5-3, loop, in=305, out=235, distance=6mm]
\end{tikzcd}
\end{equation}
\begin{lemma}
\la{CGdec}  
For all $m \in \M(W)$, there are natural $($weak$)$ homotopy equivalences of $G$-spaces
\begin{equation}
\la{FGmdec}
F_m(G,T) \,\simeq\, \hocolim_{\Cc(\Gamma)_{hW}}(G\times_T \cF_m^{\,*}) 
\,\simeq\,  \hocolim_{\,\overline{\Cc(\Gamma)}_{hW}}(G\times_T \cF_m^{\,*})
\end{equation}
Under \eqref{FGmdec}, the maps $\varphi_{m,m'}: F_m \to F_{m'} $ in  \eqref{towFm}
are induced by the morphisms of functors $ G\times_T \cF_m^{\,*} \to G\times_T \cF_{m'}^{\,*}$ defined  by the natural inclusions
of iterated joins $ \,\O^{\ast(m_{\alpha} +1)}_{e(s_{\alpha})} \into 
\O^{\ast(m'_{\alpha} +1)}_{e(s_{\alpha})}$ for $\, m_{\alpha} \leq m_{\alpha}'$.
\end{lemma}
\begin{proof}
The second  equivalence in \eqref{FGmdec} is immediate from Lemma~\ref{skCG}, while the first follows from a classical theorem of R. Thomason on homotopy colimits over cofibered categories (see, e.g., \cite[Section 26]{ChS02}). To apply this theorem 
we recall that the categories $ \Cc^{(m)}(\Gamma)_{hW} $ are cofibered over 
$ \Cc(\Gamma)_{hW}$, with fibres being the subdivision categories ${\sd}(\Delta^{m_\alpha})^{\rm op}$ (see \eqref{thickfunW}). Thus, by \cite[Proposition 26.8]{ChS02}, we can decompose
\begin{eqnarray*}
  F_m(G,T) &=& \hocolim_{\Cc^{(m)}(\Gamma)_{hW}}\,(G \times_T \cF_{m}) \\
  &\simeq& \hocolim_{e_{\sigma}(s_{\alpha}, w) \in \Cc(\Gamma)_{hW}}\,\hocolim_{\sigma \in {\sd}(\Delta^{m_\alpha})^{\rm op}}\,[\,G \times_T \cF_{m}(e_{\sigma}(s_{\alpha}, w))]\\
  &\simeq& \hocolim_{e_{\sigma}(s_{\alpha}, w) \in \Cc(\Gamma)_{hW}}\,\hocolim_{\sigma \in {\sd}(\Delta^{m_\alpha})^{\rm op}}\,
  [\,G \times_T \O^{\times \sigma}_{e(s_\alpha,w)}]\\
  &\simeq& \hocolim_{e_{\sigma}(s_{\alpha}, w) \in \Cc(\Gamma)_{hW}}\,
  [\,G \times_T \hocolim_{\sigma \in {\sd}(\Delta^{m_\alpha})^{\rm op}}\,
  (\, \O^{\times \sigma}_{e(s_\alpha,w)})]\\
  &\simeq& \hocolim_{e_{\sigma}(s_{\alpha}, w) \in \Cc(\Gamma)_{hW}}\,
  [\,G \times_T \O_{e(s_{\alpha}, w)}^{\,\ast (m_{\alpha} + 1)}]\\
  & = & \hocolim_{\Cc(\Gamma)_{hW}}(G\times_T \cF_m^{\,*})
\end{eqnarray*}
where on the last step we use the well-known decomposition formula 
for iterated joins of spaces $\, X_0 \ast X_1\ast\ldots \ast X_m 
\simeq \hocolim_{\sigma \in \sd(\Delta^m)^{\rm op}} 
\prod_{i \in \sigma} X_i\,$ (see, e.g., \cite[Example 5.8.8]{MV15}).
\end{proof}
As a first topological application of Lemma~\ref{skCG} and 
Lemma~\ref{CGdec}, we prove
\begin{prop}
\la{fgrProp}
Assume that the Lie group $G$ is simply connected. Then, for all $ m \in \M(W)$,
the spaces $F_m(G,T)$ are connected with fundamental group
\begin{equation*}\la{freeprd}
\pi_1[F_m(G,T)]\,\cong \Ast_{\alpha \in \R_+}\!\!\! W_{\alpha}
\end{equation*}
Moreover, for all $ m \leq m' $, the maps $\,\varphi_{m,m'}: F_m \to F_{m'} $ in \eqref{towFm} induce isomorphisms on $\pi_1$.
\end{prop}
\begin{proof}
We will use a general version of the Seifert-Van Kampen Theorem for fundamental groupoids of diagrams of spaces (as proven, e.g., in \cite{DF04}). 
We apply this theorem to the diagram \eqref{dGTFm}: specifically, we consider the canonical map of homotopy colimits
\begin{equation}
\la{fundgrm}
\hocolim_{\,\overline{\Cc(\Gamma)}_{hW}}(G\times_T \cF_m^{\,*}) \,\to\,
\hocolim_{\,\overline{\Cc(\Gamma)}_{hW}}({\rm pt}) 
\end{equation}
induced by the morphism of $\overline{\Cc(\Gamma)}_{hW}$-diagrams collapsing  \eqref{dGTFm} to a point. This map fits naturally in the commutative diagram of groupoids
\begin{equation}
\la{diagrpd}
\begin{diagram}[small]
{\Pi}_1 [\,\hocolim_{\,\overline{\Cc(\Gamma)}_{hW}}(G\times_T \cF_m^{\,*})\,] & \rTo_{} & {\Pi}_1[\hocolim_{\,\overline{\Cc(\Gamma)}_{hW}}({\rm pt})] \\
\dTo^{\simeq}\ & & \dTo_{\simeq}\\
\hocolim_{\,\overline{\Cc(\Gamma)}_{hW}}\! \Pi_1[(G\times_T \cF_m^{\,*})\,]  & \rTo & 
\hocolim_{\,\overline{\Cc(\Gamma)}_{hW}}\! \Pi_1({\rm pt})
\end{diagram}
\end{equation}
where the homotopy colimits in the bottom are taken in the category of groupoids.
Now, by \cite[Theorem 1.1]{DF04}, the two vertical arrows in \eqref{diagrpd} are equivalences of groupoids. We claim that the bottom arrow is also an equivalence.
To prove this it suffices to show that all spaces appearing in the diagram \eqref{dGTFm} are $1$-connected. It is well-known that $G/T$ is $1$-connected for any compact connected Lie group $G$ (see \cite[Theorem V.5.16]{MT78}).
Then, it follows from the fibration sequence
$$
\O_{e(s_{\alpha})}^{\,\ast (m_{\alpha} + 1)} \to G \times_{T} \O_{e(s_{\alpha})}^{\ast (m_{\alpha} + 1)} \to G/T\,,
$$
that the spaces $ G \times_{T} \O_{e(s_{\alpha})}^{\,\ast (m_{\alpha} + 1)} $ are $1$-connected for all $ m_\alpha \ge 1 $. On the other hand, for $m_\alpha = 0$, the space $\, G \times_{T} \O_{e(s_{\alpha})} \simeq G/T_{\alpha}$ is $1$-connected, provided $G$ is $1$-connected\footnote{If $G$ is not simply connected, the homogeneous spaces $ G/T_{\alpha}$ may fail to be simply connected for some $\alpha$ (for a simple example, take $ G = {\rm SO}(3) \times {\rm SU}(2)$).} (which follows from the Borel fibration sequence $ G/T_{\alpha} \to BT_{\alpha} \to BG $). Thus, if $G$ is simply connected, the map on the bottom of the diagram \eqref{diagrpd} is an equivalence of groupoids, and hence so is the map on the top. Now, by Lemma~\ref{CGdec}, the map \eqref{fundgrm} is (weakly) homotopy equivalent to the map
\begin{equation*}
 F_m(G,T) \to B[\overline{\Cc(\Gamma)}_{hW}]    
\end{equation*}
where $B[\overline{\Cc(\Gamma)}_{hW}]$ is the classifying space of the category $\overline{\Cc(\Gamma)}_{hW}$.  It follows that 
\begin{equation}
\la{Pi1}
\Pi_1[F_m(G,T)] \,\simeq\,  \Pi_1(B[\overline{\Cc(\Gamma)}_{hW}])\, \cong\, \Pi_1[\overline{\Cc(\Gamma)}_{hW}] 
\end{equation}
for all $m \in \M(W)$. 
From the presentation of $\overline{\Cc(\Gamma)}_{hW}$ given in Lemma~\ref{skCG} it is obvious that the category $ \overline{\Cc(\Gamma)}_{hW} $ is connected, while its $ \pi_1 $ can be easily calculated. Indeed, by \cite[Lemma~IV.3.4]{Wei13}, 
for every small connected category $ \Cc $, the fundamental group has the 
presentation
\begin{equation}
\la{pi1C}
\pi_1(\Cc) \,\cong\, \F\langle \, [f] : f \in {\rm Mor}(\Cc)\rangle/R\ ,
\end{equation}
where the relations $R$ are determined by the choice of a maximal tree $T$ in 
${\rm Mor}(\Cc)$: namely, $[t] = 1 $ for every $t \in T $, $ [\id_c] = 1 $ for every
$ c \in {\rm Ob}(\Cc) $, and $[f \circ g] = [f] \cdot [g] $ for every pair $(f,g) $ of composable morphisms in $\Cc$. In our case, for $ \Cc = \overline{\Cc(\Gamma)}_{hW} $, there is an obvious choice for $T$: namely,  
$ T= \{f_{\alpha}\}_{\alpha \in \R_+} =\{f_{\alpha}, f_{\beta}, f_{\delta}, \ldots\} $ (see Figure~\ref{skC}). With this choice, it follows from \eqref{pi1C}
immediately that
$$
\pi_{1}[\overline{\Cc(\Gamma)}_{hW}]\,\cong\, \langle\,[s_{\alpha}]\ :\ s_{\alpha}^2 =1 \,,\ \forall\,\alpha \in \R_+\,\rangle\,\cong \Ast_{\alpha \in \R_+}\!\!\! W_{\alpha}\,.
$$
To prove the last statement it suffices to note that, by naturality, the maps \eqref{FmG} fit in the  commutative diagram
%
\[\begin{tikzcd}[scale cd=0.85]
	{F_m(G,T)} && {F_{m'}(G,T)} \\
	& { B[\overline{\Cc(\Gamma)}_{hW}]}
	\arrow["{\varphi_{m,m'}}", from=1-1, to=1-3]
	\arrow[from=1-1, to=2-2]
	\arrow[from=1-3, to=2-2]
\end{tikzcd}
\]
where, as shown above, the diagonal maps are $\pi_1$-isomorphisms.
\end{proof}
\subsection{$p\,$-plus construction}\la{S5.3}
Proposition~\ref{fgrProp} shows that the spaces $F_m(G,T)$, while having correct cohomology, are not simply connected, and thus cannot serve as integral models for quasi-flag manifolds\footnote{Recall that, up to rational equivalence, the quasi-flag manifolds are to be represented by coaffine stacks which {\it are} simply connected (see Section~\ref{Sullivan}).}. In homotopy theory, there is a well-known functorial construction --- called the {\it Quillen plus construction} \cite{Q70}---  that, by attaching $2$- and $3$-cells to a path-connected space $X$, eliminates its fundamental group $\pi_1(X)$  without changing cohomology. Unfortunately, this construction applies to $X$ (if and) only if $\pi_1(X)$ is a perfect group (see, e.g., \cite[Ch. IV, \S1]{Wei13}): i.e. $\pi_1(X) = [\pi_1(X),\,\pi_1(X)]\,$, which obviously does not hold for 
$ \pi_1[F_m(G,T)] \cong (\Z/2\Z) \ast \ldots \ast (\Z/2\Z) $. Thus,  we cannot
perform the Quillen plus construction on our spaces $ F_m(G,T)$.

There is, however, a natural generalization of the plus construction that applies in our situation.
Given a commutative ring $R$, we define an {\it $R$-plus construction}\footnote{We refer the reader to the interesting recent papers \cite{BLO21}, where the $R$-plus construction is defined and its basic properties are established, and \cite{CSS23}, where the questions of uniqueness and functoriality are studied.}
for a connected space $X$ to be a $1$-connected space $ X^{R+}$  together with a map $ q: X \to X^{R+} $ that induces  an isomorphism on homology with coefficients in any $R$-module $M$. 
To state precise conditions for the existence of this construction we recall that a (discrete) group $ \pi $ is {\it $ R$-perfect} if $\, H_1(\pi, R) = 0 \,$, and {\it strongly $R$-perfect} if it is $R$-perfect and $ {\rm Tor}^1_{\Z}(R, \,H_1(\pi, \Z)) = 0 \,$. When $ R = \F_p $ (for $p$ a prime) or $ R = \Q $ (for $p=0$), we call the $R$-perfect groups $p$-perfect and refer to the corresponding $R$-plus constructions as {\it $p\,$-plus constructions}; in this case, we will also abuse the notation writing $X^{p+}$ --- or simply $X^+$ --- instead of $ X^{R+}$. 
The following result is due to Broto, Levi and Oliver \cite{BLO21}.
\begin{theorem}[\cite{BLO21}, Proposition A.5]
\la{plusthm}
Let $R$ be a commutative ring with $1$, and let $X$ be a connected CW-complex. An $R$-plus construction for $X$ exists if and only if
either $ {\rm char}(R) \not = 0 $ and $ \pi_1(X) $ is $R$-perfect, 
or $ {\rm char}(R) = 0 $ and $ \pi_1(X) $ is strongly $R$-perfect.
\end{theorem}
To apply this theorem in our situation, we observe that, by Proposition~\ref{fgrProp}, for any  ring $R$,
$$
H_1(\pi_1(F_m), R)\,\cong\, R \otimes_{\Z} \pi_1(F_m)_{\rm ab} 
\cong R \otimes_{\Z} \bigl(\oplus_{\alpha \in \R_{+}}\! W_{\alpha}\bigr) \cong (R/2R)^{|\R_+|}\,.
$$
Hence the groups $ \pi_1(F_m)$ are $R$-perfect whenever $ \frac{1}{2} \in R \,$; in particular, they are $p$-perfect
for all primes $ p>2$ and strongly $p$-perfect for $p=0$. By Theorem~\ref{plusthm}, the spaces $ F_m $ thus admit $p$-plus constructions for all values of $ p $, except $\,p=2$. In what follows, we give a functorial model for these plus constructions: namely, 
we construct a diagram of spaces $ F_*^{+}(G,T): \M(W) \to \Top^G $ together with a morphism of diagrams
$q_*: F_\ast(G,T) \to  F_*^{+}(G,T) $ such that, for all $ m \in \M(W) $,
the maps $ q_m:\, F_m(G,T) \to  F_m^{+}(G,T)\,$ provide (simultaneous) $p$-plus constructions for $p=0$ and all primes $p > 2$ (see Corollary~\ref{corpplus}).

As a first step, we assemble the diagrams \eqref{FFGm} into a single diagram of spaces. To this end we consider the functor
$$
\Phi(\Gamma):\,\M(W) \to \Cat\ ,\quad 
m \mapsto \Cc^{(m)}(\Gamma)\ ,\quad (m \leq m') \mapsto \Phi_{m,m'}\,,
$$
where $\Phi_{m,m'}$ are the full embeddings of categories $\Cc^{(m)}(\Gamma) \into \Cc^{(m')}(\Gamma)$ defined by \eqref{Phmm}. Using the Grothendieck construction
(Definition~\ref{Grothcon}), we then define
\begin{equation}\la{CMG}
\Cc^{(\M)}(\Gamma) := \M(W) \smallint\Phi(\Gamma)    
\end{equation}
The objects of this category are all pairs $ (m, w)$ and $(m, e_{\sigma}(s_{\alpha}, w)) $, where $ m = (m_{\alpha})\in \M(W) $, $\,w \in W\,$ and
$\, e_{\sigma}(s_{\alpha}, w) \in \Cc^{(m)}(\Gamma)\,$ with $\sigma \subseteq [m_{\alpha}]$. The morphisms $ (m,c) \to (m', c') $ are morphisms $\, \Phi_{m,m'}(c) \to c' \,$ in $ \Cc^{(m')}(\Gamma)$ defined whenever $ m \leq m' $. The category \eqref{CMG} is fibered over the poset $ \M(W)$: i.e., there is a natural map (forgetful functor)
\begin{equation}\la{piCM}
\Cc^{(\M)}(\Gamma) \to \M(W)   
\end{equation}
whose fibers over $ m  $ are canonically isomorphic to $ \Cc^{(m)}(\Gamma)\,$.
The functor \eqref{piCM} is a $W$-functor with respect to the trivial action of $W$ on $ \M(W)$; hence $\, \Cc^{(\M)}(\Gamma)\,$ is a $W$-category, with $W$-action preserving the fibers of \eqref{piCM}. We write $\Cc^{(\M)}(\Gamma)_{hW} $ for the corresponding homotopy orbit category. 
Then \eqref{piCM} factors through the natural inclusion $ \Cc^{(\M)}(\Gamma) \into \Cc^{(\M)}(\Gamma)_{hW} $ inducing a map
\begin{equation}
\la{piCMh}
\pi: \Cc^{(\M)}(\Gamma)_{hW} \to \M(W)    
\end{equation}
This last map is again a cocartesian fibration in $\Cat$
with fibers  $\,\pi^{-1}(m) \cong \Cc^{(m)}(\Gamma)_{hW}\,$ for all $ m \in \M(W)$. Thus we have a canonical isomorphism of categories
\begin{equation}
  \Cc^{(\M)}(\Gamma)_{hW} \,\cong \, \M(W)\smallint \Phi(\Gamma)_{hW} 
\end{equation}
where $\Phi(\Gamma)_{hW}:\, \M(W) \to \Cat $ is the classifying (fiber) functor for \eqref{piCMh} given by $\,m \mapsto  \Cc^{(m)}(\Gamma)_{hW}\,$.
We will refer to $ \Cc^{(\M)}(\Gamma)$ and $ \Cc^{(\M)}(\Gamma)_{hW}$ as the {\it universal $\M$-thickened categories} associated to $\Gamma $.

Next, we observe that the orbit functors \eqref{Fm} extend naturally to a functor $ \cF $ on $  \Cc^{(\M)}(\Gamma)$, which, in turn, extends to
the $W$-functor
\begin{equation}
\la{FFGM}
G \times_T \cF\,:\ \Cc^{(\M)}(\Gamma)_{hW}\,\to\, \Top^G\ ,\quad
(m, \, e_{\sigma}(s_\alpha, w)) \,\mapsto\, G \times_T \O^{\times \sigma}_{e(s_\alpha,w)} \  \, (\sigma \subseteq [m_\alpha])\, .
\end{equation}
The functor \eqref{FFGM} assembles the diagrams of spaces \eqref{FFGm} together in the sense that
\begin{equation}
\la{restm}
G \times_T \cF_m = i_m^* (G \times_T \cF) \quad \mbox{for all} \ \,  
m\in\M(W)\,,
\end{equation}
where $ i_m^* $ denotes the restriction  to a fiber  of \eqref{piCMh} via the fiber inclusion functor
\begin{equation}
\la{fibin}
i_{m}:\, \Cc^{(m)}(\Gamma)_{hW} \into  \Cc^{(\M)}(\Gamma)_{hW}\ ,\quad c \mapsto (m,c)\,.
\end{equation}
By formal properties of cocartesian fibrations, our basic diagram \eqref{towFm} can be then recovered from  \eqref{FFGM} by the  formula
\begin{equation}
\la{Kanext}
F_*(G,T) \,\simeq \, \pi_{!}\,(G \times_T \cF)\,,
\end{equation}
where $\pi_!: \Top^{\Cc^{(\M)}(\Gamma)_{hW}} \to \Top^{\M(W)} $ is 
the (homotopy) left Kan extension along the map \eqref{piCMh}.
Indeed, for all $m \in \M(W) $, we have 
%
$$
\pi_{!}\,(G \times_T \cF)(m) \,\simeq\, \hocolim_{\pi^{-1}(m)}\,[\,i_m^*(G\times_T\cF)\,]\,=\,\hocolim_{\Cc^{(m)}(\Gamma)_{hW}}(G\times_T\cF_m) 
= F_m(G,T)\,.
$$
As a consequence, we get a natural (weak) homotopy equivalence of $G$-spaces
$$
\hocolim_{\Cc^{(\M)}(\Gamma)_{hW}}(G\times_T \cF)\,\simeq\, \hocolim_{m \in \M(W)} [\,F_m(G,T)\,]
$$

Now, for the next step, observe that, by Proposition~\ref{fgrProp}, the fundamental groups of all spaces $ F_m = F_m(G,T) $ are  naturally isomorphic to the fundamental group of the category $ \Cc(\Gamma)_{hW} $ (see \eqref{Pi1}).  To make these spaces simply connected we modify the definition\footnote{Recall that the spaces $F_m $ are defined as homotopy 
colimits over the categories $ \Cc^{(m)}(\Gamma)_{hW} $  (see \eqref{FmG}), and thus these categories provide the `gluing data' for $F_m$.} of the $m$-thickened moment categories
 $ \Cc^{(m)}(\Gamma)_{hW} $  by homotopically contracting $ \Cc(\Gamma)_{hW} $ in
 $ \Cc^{(m)}(\Gamma)_{hW} $ for each $ m $. Just as in the case of spaces, 
 such a contraction amounts to taking the homotopy cofiber of the inclusion
 map $ \Cc(\Gamma)_{hW} \into \Cc^{(m)}(\Gamma)_{hW} $, which, in turn,
 amounts to attaching a cone over  $ \Cc(\Gamma)_{hW} $ inside $ \Cc^{(m)}(\Gamma)_{hW} $ for each $ m $. 
 The above reformulation allows us to perform this procedure as a universal construction (for all $m$ at once):
 \begin{equation}
 \la{CG-}
 \Cc^{(\M)}(\Gamma)^{+}_{hW} \,:=\,{\rm hocofib}\,[\,\Cc(\Gamma)_{hW} \xrightarrow{i}  \Cc^{(\M)}(\Gamma)_{hW}\,]\,,
 \end{equation}
where $ i = i_0 $ is the $0$-fiber inclusion functor defined by \eqref{fibin}. The above homotopy cofiber is taken in $\Cat$ and thus represented by the Grothendieck construction
 \begin{equation}
 \la{CGG-}
 \Cc^{(\M)}(\Gamma)^{+}_{hW} \,=\,\I \smallint \{\ast \leftarrow \Cc(\Gamma)_{hW} \xrightarrow{i}  \Cc^{(\M)}(\Gamma)_{hW}\} \,,
 \end{equation}
where $\, \I := \{1 \leftarrow 0 \rightarrow 2 \}\,$ is the category parametrizing pushouts. The objects in $\Cc^{(\M)}(\Gamma)^{+}_{hW}$ is the disjoint union of the objects of
$\Cc(\Gamma)_{hW}$, $\,\Cc^{(\M)}(\Gamma)_{hW} $
and the (unique) object of $\ast$. Instead of listing the morphisms, we describe the data defining
a functor on such a category. 
\begin{lemma}
\la{funIG}
Let $\,\{\Cc_1 \xleftarrow{F_1} \Cc_0 \xrightarrow{F_2} \Cc_2\} \,$ be a pushout diagram in $\Cat$. A functor 
$$
G:\
\I\smallint\, \{\Cc_1 \xleftarrow{F_1} \Cc_0 \xrightarrow{F_2} \Cc_2\}\,\to\, \Dc
$$
is determined by  three functors $\,G_i: \Cc_i \to \Dc \,$ for $\,i=0,1,2\,$ given together with two natural transformations $\,\varphi_1: G_0 \to G_1 \circ F_1 \,$ and $\,\varphi_2: G_0 \to G_2 \circ F_2\,$. 
We  represent $G$ diagrammatically as
%
\begin{equation}
\la{fungrc}
\begin{tikzcd}[scale cd= 0.9]
	{{\mathscr C}_1} && {{\mathscr C}_0} && {{\mathscr C}_2} \\
	\\
	&& {{\mathscr D}}
	\arrow["{G_1}"', curve={height=6pt}, from=1-1, to=3-3]
	\arrow["{F_1}"', from=1-3, to=1-1]
	\arrow["{F_0}", from=1-3, to=1-5]
	\arrow["{G_0}"{description}, from=1-3, to=3-3]
	\arrow["{G_2}", curve={height=-6pt}, from=1-5, to=3-3]
	\arrow["{\varphi_1}"', shift right=5, shorten <=22pt, shorten >=22pt, Rightarrow, from=3-3, to=1-1]
	\arrow["{\varphi_2}", shift left=5, shorten <=22pt, shorten >=22pt, Rightarrow, from=3-3, to=1-5]
\end{tikzcd}
\end{equation}
%
\end{lemma}
\begin{proof}
Immediate from the universal mapping property of Grothendieck construction (see, e.g., \cite[Prop. 1.3.1]{To79}). 
\end{proof}

Now, recall that, in \Cref{S4.3}, 
we defined the `geometric' functor $G \times_T \tilde{\cF}_0: \, \Cc(\Gamma)_{hW} \to \Top^G $:
$$
 w \,\mapsto\,  G\times_T \biggl(\,\bigvee_{e(s_{\alpha}, w) \in E_{\Gamma}(w)}\!\! \!\overline{\O}_{e(s_{\alpha}, w)}\!\setminus\!\{s_{\alpha} w\}\biggl)\ ,
\quad e(s_{\alpha}, w) \mapsto G \times_T \O_{e(s_{\alpha}, w)}\,,
$$
together with two natural transformations \eqref{GTF} and \eqref{GTFq}:
\begin{equation*}
G/T \,\xleftarrow{\tilde{q}}\,G \times_T \tilde{\cF}_0 \,\xrightarrow{\tilde{p}} \,G \times_T \cF_0\,,
\end{equation*}
which we used in the proof of Theorem~\ref{hDec}. By \eqref{restm}, the functor $ G \times_T \cF_0 $ factors as $ (G \times_T \cF) \circ i $, where $ G \times_T \cF $ is the universal $W$-functor \eqref{FFGM}. Hence, by Lemma~\ref{funIG}, we can assemble these data into the diagram 
\begin{equation}
\la{diagF-}
%
\begin{tikzcd}
	\ast && {{\mathscr C}(\Gamma)_{hW}} && {{\mathscr C}^{(\M)}(\Gamma)_{hW}} \\
	\\
	&& {\Top^G}
	\arrow["{G/T}"', curve={height=6pt}, from=1-1, to=3-3]
	\arrow[from=1-3, to=1-1]
	\arrow["i", from=1-3, to=1-5]
	\arrow["{G\times_T\widetilde{\mathcal F}_0}"{description}, from=1-3, to=3-3]
	\arrow["{G \times_T\mathcal{F}}", curve={height=-6pt}, from=1-5, to=3-3]
	\arrow["{\tilde{q}}"', shift right=5, shorten <=24pt, shorten >=24pt, Rightarrow, from=3-3, to=1-1]
	\arrow["{\tilde{p}}", shift left=5, shorten <=27pt, shorten >=27pt, Rightarrow, from=3-3, to=1-5]
\end{tikzcd}
\end{equation}
which defines a $W$-functor $\,G\times_T \cF^{+}:\, \Cc^{(\M)}(\Gamma)_{hW}^{+} \to \Top^G$.

On the other hand, the category \eqref{CGG-} comes with a canonical map $\,\Cc^{(\M)}(\Gamma)_{hW} \to  \Cc^{(\M)}(\Gamma)_{hW}^{+} \,$ that exhibits it as a cofiber, and the fibration \eqref{piCMh} factors through this map inducing 
\begin{equation*}
\pi^{+}: \Cc^{(\M)}(\Gamma)^{+}_{hW} \to \M(W)  
\end{equation*}
Diagrammatically, this last functor can be represented by 
\begin{equation}
\la{piCMh-}
   \begin{tikzcd}[scale cd=0.9]
	{\ast} && { {\Cc(\Gamma)_{hW}} } && {{\Cc^{(\M)}(\Gamma)_{hW}}} \\
	\\
	&& {\M(W)}
	\arrow["{0}"', curve={height=6pt}, from=1-1, to=3-3]
	\arrow["{}"', from=1-3, to=1-1]
	\arrow["{i}", from=1-3, to=1-5]
	\arrow["{0}"', from=1-3, to=3-3]
	\arrow["{\pi}", curve={height=-6pt}, from=1-5, to=3-3]
	\arrow["{\id}"', shift right=5, shorten <=22pt, shorten >=22pt, Rightarrow, from=3-3, to=1-1]
	\arrow["{\id}", shift left=5, shorten <=22pt, shorten >=22pt, Rightarrow, from=3-3, to=1-5]
\end{tikzcd} 
\end{equation}
Thus, we have constructed a universal category $\Cc^{(\M)}(\Gamma)^{+}_{hW}$  together with two natural functors
$G \times_T \cF^{+}$ and $ \pi^{+}$ replacing the functors
$G \times_T \cF $ and $\pi$ on $ \Cc^{(\M)}(\Gamma)_{hW}$.
Recall that, in terms of these last two functors, the spaces $ F_m(G,T)$ --- or rather the diagram \eqref{towFm} ---  can be concisely defined by \eqref{Kanext}. This motivates us to make the following definition:
\begin{equation}
\la{towFm-}
F^{+}_{*}(G,T) :=   \pi^{+}_{!}(G \times_T \cF^{+}) \,,
\end{equation}
where $\pi^{+}_{!}$ denotes the homotopy left Kan extension along \eqref{piCMh-}. Formula \eqref{towFm-} defines a functor $ F^{+}_{*}(G,T): \M(W) \to \Top^G $ that assigns to each multiplicity $ m \in \M(W) $ a $G$-space $ F_m^{+}(G,T) $. Our next goal is to describe these spaces explicitly.

To this end, for each $ m \in \M(W)$, we introduce the category
 \begin{equation}
 \la{CGGm-}
 \Cc^{(m)}(\Gamma)^{+}_{hW} \,:=\,\I \smallint \{\ast \leftarrow \Cc(\Gamma)_{hW} \xrightarrow{\Phi_{0,m}}  \Cc^{(m)}(\Gamma)_{hW}\} \,,
 \end{equation}
and define the functor
$$ 
G\times_T \cF_m^{+}: \, \Cc^{(m)}(\Gamma)^{+}_{hW} \to \Top^G 
$$ 
by
\begin{equation}
\la{diagFm-}
%
\begin{tikzcd}[scale cd= 0.9]
	\ast && {{\mathscr C}(\Gamma)_{hW}} && {{\mathscr C}^{(m)}(\Gamma)_{hW}} \\
	\\
	&& {\Top^G}
	\arrow["{G/T}"', curve={height=6pt}, from=1-1, to=3-3]
	\arrow[from=1-3, to=1-1]
	\arrow["\Phi_{0,m}", from=1-3, to=1-5]
	\arrow["{G\times_T\widetilde{\mathcal F}_0}"{description}, from=1-3, to=3-3]
	\arrow["{G \times_T\mathcal{F}_m}", curve={height=-6pt}, from=1-5, to=3-3]
	\arrow["{\tilde{q}}"', shift right=5, shorten <=24pt, shorten >=24pt, Rightarrow, from=3-3, to=1-1]
	\arrow["{\tilde{p}}", shift left=5, shorten <=27pt, shorten >=27pt, Rightarrow, from=3-3, to=1-5]
\end{tikzcd}
\end{equation}
where the map $ \Phi_{0,m}$ is given by \eqref{Phmm}. Note that  the diagram \eqref{diagFm-} does represent a functor on $\Cc^{(m)}(\Gamma)^{+}_{hW}$ in view of \eqref{GTmm}.
The next proposition provides a natural decomposition of the 
spaces $ F_m^{+}(G,T) $ similar to \eqref{FmG}.
\begin{prop}
\la{Prop5}
For all $m \in \M(W) $, there are natural $G$-homotopy equivalences of $G$-spaces    
\begin{equation}
\la{F-m}
F^{+}_m(G,T) \simeq \hocolim_{ \Cc^{(m)}(\Gamma)^{+}_{hW}}(G \times_T{\cF_m^{+}}) 
\end{equation}
\end{prop}
\begin{proof}
Before proving \eqref{F-m}, we remark that the left Kan extension \eqref{towFm-} cannot be computed directly by
restricting the functor $ G \times_T \cF^{+} $ to fibers of  \eqref{piCMh-} and taking homotopy colimits of such restrictions.  
This is because the map \eqref{piCMh-} is not a cocartesian fibration in $\Cat$, i.e. unlike $ \Cc(\Gamma)_{hW}$, the category $\Cc^{(\M)}(\Gamma)^{+}_{hW}$ is not fibered over $ \M(W)$. To remedy this problem
we construct a fully faithful embedding of  $\Cc^{(\M)}(\Gamma)^{+}_{hW}$ in a larger category that {\it is} fibered over $\M(W)$:
\begin{equation}
\la{embGr}
\begin{tikzcd}[scale cd= 0.9]
	{{\mathscr C}^{(\M)}(\Gamma)^{+}_{hW}\!} && {\widetilde{{\mathscr C}}^{(\M)}(\Gamma)^{+}_{hW}} \\
	& {\M(W)}
	\arrow["f", hook, from=1-1, to=1-3]
	\arrow["{\pi^{+}}"', from=1-1, to=2-2]
	\arrow["{\tilde{\pi}}", two heads, from=1-3, to=2-2]
\end{tikzcd}
\end{equation}
The category $ \widetilde{\Cc}^{(\M)}(\Gamma)^{+}_{hW} $ in \eqref{embGr} is defined by $\,\M(W) \smallint \Phi^{+}(\Gamma)\,$, where the classifying functor is
\begin{equation*}
\Phi^{+}(\Gamma):\, \M(W) \to \Cat\ , \quad m \mapsto \Cc^{(m)}(\Gamma)^{+}_{hW} \ ,\quad (m \leq m') \mapsto \Phi^{+}_{m,m'}\ ,
\end{equation*}
with  $\Cc^{(m)}(\Gamma)^{+}_{hW}$ given by \eqref{CGGm-}
and the $\Phi^{+}_{m,m'}$  induced by the following maps of $\I$-diagrams:
%
\[\begin{tikzcd}[scale cd= 0.9]
	\ast && {{\mathscr C}(\Gamma)_{hW}} && {{\mathscr C}^{(m)}(\Gamma)_{hW}} \\
	\ast && {{\mathscr C}(\Gamma)_{hW}} && {{\mathscr C}^{(m')}(\Gamma)_{hW}}
	\arrow[shorten <=1pt, shorten >=1pt, equals, from=1-1, to=2-1]
	\arrow[from=1-3, to=1-1]
	\arrow["{\Phi_{0,m}}", from=1-3, to=1-5]
	\arrow[equals, from=1-3, to=2-3]
	\arrow["{\Phi_{m,m'}}", from=1-5, to=2-5]
	\arrow[from=2-3, to=2-1]
	\arrow["{\Phi_{0,m'}}"', from=2-3, to=2-5]
\end{tikzcd}\]
Now, to define the functor $ f $ in \eqref{embGr} we observe that, by commutativity of Grothendieck constructions, there is a natural equivalence of categories
\begin{equation}
\la{eqnceC}
\widetilde{\Cc}^{(\M)}(\Gamma)^{+}_{hW}\ \simeq\ 
\I \,\smallint\, \{\M(W) \xleftarrow{\,\rm pr} \M(W) \times \Cc(\Gamma)_{hW} \xrightarrow{ i_* \Phi_{0,\ast}}  \Cc^{(\M)}(\Gamma)_{hW}\} \,,
\end{equation}
where the map $i_* \Phi_{0,\ast}$ is given by 
$\,(m,c) \mapsto (m,\, i_m \Phi_{0,m}(c))\,$, see \eqref{Phmm} and \eqref{fibin}. Identifying  $\widetilde{\Cc}^{(\M)}(\Gamma)^{+}_{hW}$ as in \eqref{eqnceC}, we can define the functor
$f$ by the following natural map of $\I$-diagrams
%
\[\begin{tikzcd}[scale cd= 0.9]
	\ast && {{\mathscr C}(\Gamma)_{hW}} && {{\mathscr C}^{(\M)}(\Gamma)_{hW}} \\
	{\M(W)} && {\M(W) \times {\mathscr C}(\Gamma)_{hW}} && {{\mathscr C}^{(\M)}(\Gamma)_{hW}}
	\arrow["0"', hook, from=1-1, to=2-1]
	\arrow[from=1-3, to=1-1]
	\arrow["{i=i_0}", from=1-3, to=1-5]
	\arrow["{0 \times \id}"', hook, from=1-3, to=2-3]
	\arrow[equals, from=1-5, to=2-5]
	\arrow[two heads, from=2-3, to=2-1]
	\arrow["{i_*\Phi_{0,\ast}}"', from=2-3, to=2-5]
\end{tikzcd}\]
To check the commutativity of \eqref{embGr} we observe that
with identification \eqref{eqnceC}, the functor 
$ \tilde{\pi}: \ \widetilde{\Cc}^{(\M)}(\Gamma)^{+}_{hW} \to \M(W) $ is represented by
%
\[\begin{tikzcd}[scale cd= 0.9]
	{\M(W)} && {\M(W)\times {\mathscr C}(\Gamma)_{hW}} && {{\mathscr C}^{(\M)}(\Gamma)_{hW}} \\
	\\
	&& {\M(W)}
	\arrow["\id"', curve={height=6pt}, from=1-1, to=3-3]
	\arrow["{{\rm pr}}"', two heads, from=1-3, to=1-1]
	\arrow["{i_* \Phi_{0,*}}", from=1-3, to=1-5]
	\arrow["{{\rm pr}}"{description}, two heads, from=1-3, to=3-3]
	\arrow["\pi", curve={height=-6pt}, from=1-5, to=3-3]
	\arrow["\id"', shift right=3, shorten <=28pt, shorten >=28pt, Rightarrow, from=3-3, to=1-1]
	\arrow["\id", shift left=3, shorten <=30pt, shorten >=30pt, Rightarrow, from=3-3, to=1-5]
\end{tikzcd}\]
Hence, the composition $ \tilde{\pi} \circ f $ is given by the
diagram \eqref{piCMh-} which represents the functor $\pi^{+}$.

Now, by Lemma~\ref{funIG}, the functor $f: \Cc^{(\M)}(\Gamma)^{+}_{hW} \into \widetilde{\Cc}^{(\M)}(\Gamma)^{+}_{hW} $ is fully faithful\footnote{In fact, it is easy to show that  $f$ admits a 
left inverse, so that $\Cc^{(\M)}(\Gamma)^{+}_{hW}$ is a deformation retract of $ \widetilde{\Cc}^{(\M)}(\Gamma)^{+}_{hW} $.}, and the left Kan extension of \eqref{diagF-} along $f$ 
is easy to compute: specifically, with identification \eqref{eqnceC},
\begin{equation}
\la{fshr}
f_!\,(G\times_T \cF^{+}):\  \widetilde{\Cc}^{(\M)}(\Gamma)^{+}_{hW}\,\to\, \Top^G
\end{equation}
is represented by the diagram
\begin{equation}
\la{diafshr}
\begin{tikzcd}[scale cd= 0.9]
	{\M(W)} && {\M(W) \times {\mathscr C}(\Gamma)_{hW}} && {{\mathscr C}^{(\M)}(\Gamma)_{hW}} \\
	\\
	&& {\Top^G}
	\arrow["{{\rm const}(G/T)}"', curve={height=6pt}, from=1-1, to=3-3]
	\arrow["{{\rm pr}_1}"', two heads, from=1-3, to=1-1]
	\arrow["{i_* \Phi_{0,*}}", from=1-3, to=1-5]
	\arrow["{(G\times_T\widetilde{\mathcal F}_0)\,{\rm pr}_2}"{description}, from=1-3, to=3-3]
	\arrow["{G \times_T\mathcal{F}}", curve={height=-6pt}, from=1-5, to=3-3]
	\arrow["{\tilde{q}}"', shift right=4, shorten <=29pt, shorten >=29pt, Rightarrow, from=3-3, to=1-1]
	\arrow["{\tilde{p}}", shift left=4, shorten <=31pt, shorten >=31pt, Rightarrow, from=3-3, to=1-5]
\end{tikzcd}
\end{equation}
Since $\, \pi^{+} \cong \tilde{\pi} \circ f\,$ by \eqref{embGr} and the map $\tilde{\pi} $ is a cocartesian
fibration,  we conclude 
\begin{eqnarray*}
F^{+}_m(G,T) &:=& \pi^{+}_{!}(G \times_T \cF^{+})(m) \cong
(\tilde{\pi} \circ f)_{!} \,(G \times_T \cF^{+})(m) 
\cong \tilde{\pi}_!\,[\,f_{!} \,(G \times_T \cF^{+})](m)\\
&\simeq&
\hocolim_{\,\tilde{\pi}^{-1}(m)}\,[\,i_m^* f_{!}\,(G \times_T \cF^{+})\,]
\cong
\hocolim_{\Cc^{(m)}(\Gamma)^{+}_{hW}}(G\times_T\cF^{+}_m) 
\end{eqnarray*}
for all $ m \in \M(W)$. This completes the proof of the proposition.
\end{proof}
As a consequence of Proposition~\ref{Prop5}, we get
\begin{cor}
\la{corpplus}
$(a)$ For all $m \in \M(W)$, there are homotopy cocartesian squares in $\Top^G:$
\begin{equation}\la{diagrm-}
\begin{diagram}[small]
F_0(G,T) & \rTo^{\ \varphi_{0,m}\ } & F_m(G,T)\\
\dTo^{q}\ &  & \dTo_{q_m}\\
G/T   & \rTo & {\NWpbk\!\!}F^{+}_m(G,T)
\end{diagram}  
\end{equation}
where the maps $ q $ and $ \varphi_{0,m}  $ are defined by 
\eqref{eq:p-hocolim} and \eqref{phimm}. Thus,
$$
F^{+}_m(G,T) \,\simeq\,\hocolim\,\{\,G/T \xleftarrow{q}  F_0(G,T) \xrightarrow{\varphi_{0,m}}F_m(G,T)\,\}\,.
$$
\noindent
$(b)$ Assume that the group $G$ is simply connected. Then, for all $ m \in \M(W) $,
\begin{enumerate}
\item[$(1)$]
the  spaces $ F^{+}_m(G,T) $
are simply connected, and
\item[$(2)$]  
the maps $ q_m $ in \eqref{diagrm-} are mod $p$ cohomology isomorphisms for all $p\not=2$.
\end{enumerate}
Thus, the  spaces $ F^{+}_m(G,T) $ together with the maps $\, q_m: F_m(G,T) \to F^{+}_m(G,T) $ provide a $($simultaneous$)$ $p$-plus construction  for $F_m(G,T)$ for all $\, p \not = 2\,$.
\end{cor}
\begin{proof}
Part $(a)$ follows from Proposition~\ref{Prop5} and Thomason's Theorem (\cite[Prop. 26.8]{ChS02}): since $\Cc^{(m)}(\Gamma)^{+}_{hW}$ is fibered over $\I$ (see \eqref{CGGm-}), we have
\begin{eqnarray*}
&&\!\!\!\!\!\!\!\!\!\! F^{+}_m(G,T)\ \simeq \ \hocolim_{\Cc^{(m)}(\Gamma)^{+}_{hW}}(G\times_T\cF^{+}_m) \\
&&\simeq
\ \hocolim\,\{\hocolim_{\ast}(G/T) \xleftarrow{q} \hocolim_{\Cc(\Gamma)_{hW}}(G \times_T \tilde{\cF}_0) \xrightarrow{\Phi^*_{0,m}} \hocolim_{\Cc^{(m)}(\Gamma)_{hW}}(G \times_T {\cF}_m) \}\\
&& \simeq \ \hocolim\,\{\,G/T \xleftarrow{q}  F_0(G,T) \xrightarrow{\varphi_{0,m}}F_m(G,T)\,\}\,.
\end{eqnarray*}
Part $(b)$ is a consequence of Part $(a)$ and Proposition~\ref{fgrProp}.
Indeed, since all spaces $ F_m $ are connected and the diagram
\eqref{diagrm-} is  cocartesian, we can use the Seifert-van Kampen Theorem to represent the fundamental groups of the spaces $F_m^{+} = F_m^{+}(G,T) $ as amalgamated products
\begin{equation}\la{amalgpr}
\pi_1(F_m^{+}) \,\cong\, \pi_1(G/T) \Ast_{\pi_1(F_0)} \pi_1(F_m)
\end{equation}
But, by Proposition~\ref{fgrProp}, the maps $ \varphi_{0,m} $ induce
isomorphisms on $ \pi_1$'s, hence \eqref{amalgpr} reduces to
$\pi_1(F_m^{+}) \cong\ \pi_1(G/T) = \{1\}$ for all $ m \in \M(W)$.
This proves $(1)$. Now, by (the proof of) Theorem~\ref{hDec}, the map
$q$ in \eqref{diagrm-} is a mod-$p$ cohomology isomorphism. From the Mayer-Vietoris exact sequence associated to \eqref{diagrm-} it then follows that so are the maps $q_m$. This proves $(2)$.
\end{proof}

We are now in position to state our main definition.
\begin{defi}
\la{maindef}
For $m\in \M(W) $, we call the $G$-space $ F_m^{+}(G,T) $ defined by \eqref{F-m}:
\begin{equation}
\la{F+m}
F^{+}_m(G,T) \,:=\, \hocolim_{ \Cc^{(m)}(\Gamma)^{+}_{hW}}(G \times_T{\cF_m^{+}})\, 
\end{equation}
the $m$-{\it quasi-flag manifolds} associated to $ (G,T)$. 
\end{defi}
Abusing terminology, we will also refer to the $G$-space $ F_m(G,T) $ defined 
by \eqref{FmG}:
\begin{equation*}
\la{FmG2}
F_m(G,T) \,:=\, \hocolim_{ \Cc^{(m)}(\Gamma)_{hW}}(G \times_T{\cF_m})
\end{equation*}
as a (non-simply connected) {\it model} of the $m$-quasi-flag manifold. 
In view of Corollary~\ref{corpplus}, the spaces $F_m(G,T)$ and $F^+_m(G,T)$ are 
cohomologically equivalent and can be used interchangeably in many situations
(see Section~\ref{S6} below). 

In addition to the above two basic spaces, 
we introduce the $T$-space
\begin{equation}
\la{tilFmT}
\tilde{F}_m(T)\,:=\,
\hocolim_{ \Cc^{(m)}(\Gamma)}(\cF_m)\,\simeq\,\hocolim_{\Cc(\Gamma)}(\cF_{m}^{\,*})\,.
\end{equation}
%
We observe that the corresponding $G$-space $\,\tilde{F}_m(G,T) := 
G \times_T  \tilde{F}_m(T)\,$
%
%
carries a natural {\it free} $W$-action, such that
\begin{equation}
\la{quotW}
\tilde{F}_m(G, T)/W \,\simeq\, \tilde{F}_m(G, T)_{hW} \,\simeq\, F_m(G,T)\,,
\end{equation}
where the last equivalence follows from Lemma~\ref{lem:Thom}.
Thus, the canonical restriction map induced by the functor $ \Cc^{(m)}(\Gamma) \into \Cc^{(m)}(\Gamma)_{hW}$  defines a $W$-bundle of $G$-spaces
\begin{equation}
\la{Wbund}
\tilde{F}_m(G, T) \to F_m(G,T)
\end{equation}
In the next two sections, we will study the cohomological and $K$-theoretic properties of the spaces $F_m(G,T)$, $\,F^+_m(G,T) $, $ \tilde{F}_m(G,T)$ and the associated $W$-bundle \eqref{Wbund}. Here we conclude by constructing some explicit models for $F_{m}(G,T)$
and $ F^{+}_m(G,T)$.
\subsection{Geometric models}
\la{S5.4}
Note that, by \eqref{quotW}, there is a $G$-equivariant homotopy equivalence
$$
F_m(G,T)\, \simeq\, \tilde{F}_m(G, T)/W \,\simeq\, G \times_N \tilde{F}_m(T)
$$
Thus, to construct a geometric model for the $G$-space $F_m(G,T)$, it suffices to construct a geometric model for the $T$-space $\tilde{F}_m(T)$ defined by \eqref{tilFmT}, which, in turn, requires constructing an appropriate cofibrant resolution for the $\Cc(\Gamma)$-diagram $ \cF_m^{\,*}$ (see \eqref{Fm*}). To this end we will use the following general lemma that describes cofibrant diagrams over the face categories of  finite graphs.
\begin{lemma}
\la{Lcofres}    
Let $ \Gamma = (V_{\Gamma}, E_{\Gamma}) $ be a finite simple graph, and let
 $ \cF: \Cc(\Gamma) \to \M $ be a functor defined on the $($opposite$)$ face category of $\Gamma $ with values in a model category $\M$. Assume that the natural maps determined by the incidence relations of $\Gamma$:
 \begin{equation}
 \la{cofcond}
 i_{v}: \ \coprod_{e \in E_{\Gamma}(v)}\cF(e)\,\into\, \cF(v)
 \end{equation}
 are cofibrations in $\M$ for all  $ v \in V_{\Gamma} $. Then $ \cF $ is  cofibrant in the projective model structure on the category of $\Cc(\Gamma)$-diagrams  $ \M^{\Cc(\Gamma)} $ in $\M$. Consequently, 
 \begin{equation}
 \la{colimfor}
 \hocolim_{\Cc(\Gamma)}(\cF)\,\simeq\,\colim_{\Cc(\Gamma)}(\cF) \,\cong\, 
 \frac{\left(\amalg_{e \in E_{\Gamma}} \cF(e)\right) \,\amalg\, \left(\amalg_{v \in V_{\Gamma}} \cF(v)\right)}{\langle (x \in \cF(e)) \sim (i_v(x) \in \cF(v))\,,\ \forall\,e \in E_{\Gamma}(v) \rangle}\ ,
 \end{equation}
 where $ E_{\Gamma}(v) = \{e \in E_{\Gamma}\,:\, v \leq e\,\} $ is the subset of edges of $\Gamma$ adjacent to a given vertex $v \in V_{\Gamma}$.
\end{lemma}
\begin{proof}
For any finite simple graph $\Gamma $, the category $ \Cc(\Gamma)$ is `very small' (see Definition~\ref{vscat});  by  \cite[Section 10.13]{DS95}, a map of diagrams $f: {\mathcal X} \to {\mathcal Y} $ is then a cofibration in $ \M^{\Cc(\Gamma)} $ iff for all 
$c \in \Cc(\Gamma)$, the maps $ \partial_{c}(f) \into {\mathcal Y}(c) $ are cofibrations in $ \M$. Here $\partial_{c}(f) := {\mathcal X}(c) \amalg_{\partial_c({\mathcal X})}\partial_c({\mathcal Y})$, where $ \partial_{c}({\mathcal X})$ denotes the colimit of the composite functor
$\,
 \Cc(\Gamma)^{\circ}_{/c} \into \Cc(\Gamma)_{/c} \to \Cc(\Gamma) \xrightarrow{{\mathcal X}} \M 
\,$
defined on the full subcategory $ \Cc(\Gamma)^{\circ}_{/c} $ of the slice category $ \Cc(\Gamma)_{/c} $ obtained by discarding the 
identity object $ \id_c $ from $ \Cc(\Gamma)_{/c} $. It follows that $ \cF \in  \M^{\Cc(\Gamma)}$ is cofibrant iff for all objects $ c \in \Cc(\Gamma) $,
\begin{equation}
\la{cofibr}
\partial_c(\cF) \into \cF(c) \ \, \mbox{is a cofibration in}\ \M
\end{equation}
Note that, for all edge objects $ e \in E_{\Gamma} $, we have $ \partial_e(\cF) = \varnothing $; hence, \ref{cofibr} holds automatically. On the other hand, for $\, v \in V_{\Gamma}\,$, we have $ \partial_v(\cF) \cong \amalg_{e \in E_{\Gamma}(v)} \cF(e) $, so that \eqref{cofibr} is equivalent to \eqref{cofcond}. Thus $\cF $ is cofibrant in 
$ \M^{\Cc(\Gamma)}$ iff \eqref{cofcond} holds for all $ v \in V_{\Gamma}$. 
\end{proof}
We apply Lemma~\ref{Lcofres} to the Bruhat moment graph $ \Gamma = \Gamma(\R_W)$.
Given $ m \in \M(W)$, we define the functor 
\begin{equation}
\la{wFm*}
\wcF_{m}^{\,*}:\ \Cc(\Gamma) \to \Top^T\ ,\quad w \mapsto \bigvee_{e \in E_{\Gamma}(w)} C_{\{w\}}[\O_{e(s_{\alpha}, w)}^{\ast (m_{\alpha} + 1)}]\,,\quad 
e(s_{\alpha}, w) \mapsto  \O_{e(s_{\alpha}, w)}^{\ast (m_{\alpha} + 1)}\,,
\end{equation}
where $ C_{\{w\}}(X) = \{w\} \ast X $ denotes the cone over a space $X$ with vertex at 
$\{w\} $. The condition~\ref{cofcond} of Lemma~\ref{Lcofres} obviously holds, hence 
\eqref{wFm*} represents a cofibrant diagram in $ \Top^{\Cc(\Gamma)}$. On the other hand,
there is a natural projection $ \wcF_m^{\,*} \to \cF_m^{\,*} $ contracting each cone in
$\wcF_m^{\,*}(w)$ to its vertex $\{w\}$. This is obviously a homotopy equivalence in $ \Top^{\Cc(\Gamma)}$ that we can use as a cofibrant replacement for  \eqref{Fm*}. Then, by formula \eqref{colimfor}, we get
\begin{equation}
\la{bFwm}
\hocolim_{\Cc(\Gamma)}(\cF_m^{\,*})\,\simeq\,\colim_{\Cc(\Gamma)}(\wcF_m^{\,*})\,\cong\,
\coprod_{E_{\Gamma}} S\bigl(\O_{e(s_{\alpha}, w)}^{\ast (m_{\alpha} + 1)}\bigr)\big/\!\sim\ ,
\end{equation}
where $S(X) := C_{\{w\}}(X) \,\cup_X\,C_{\{s_{\alpha}w\}}(X) = \{W_{\alpha}w\} \ast X $
is the unreduced suspension of a space $X$ with suspension points $ W_{\alpha}w = \{w,\,s_{\alpha} w\}$ (see Fig~\ref{SO}). The spaces  $ S\bigl(\O_{e(s_{\alpha}, w)}^{\ast (m_{\alpha} + 1)}\bigr) $  are glued to each other 
at their suspension points according to the incidence relations of the graph $\Gamma $. We denote the resulting topological space in \eqref{bFwm} by $ \widetilde{\bF}^{(m)} $.

\begin{figure}[h!]
\centering
\begin{tikzpicture}[
  knot/.style={circle, fill=white, draw=black!65, semithick,
               inner sep=0pt, minimum size=3.2pt},
  scale=1
]
  \coordinate (T) at (0, 2);
  \coordinate (B) at (0,-2);
 
  \shade[left color=thmframe!40, right color=thmframe!14]
    (T) -- (-1.2,0) arc[start angle=180, end angle=360,
                        x radius=1.2cm, y radius=0.3cm] -- cycle;
  \shade[left color=thmframe!40, right color=thmframe!14]
    (B) -- (-1.2,0) arc[start angle=180, end angle=360,
                        x radius=1.2cm, y radius=0.3cm] -- cycle;
  \shade[left color=thmframe!26, right color=thmframe!14]
    (0,0) ellipse (1.2cm and 0.3cm);
  \draw[thmframe!50!black, densely dotted, thin]
    (-1.2,0) arc[start angle=180, end angle=0,
                 x radius=1.2cm, y radius=0.3cm];
  \draw[thmframe!65!black, semithick]
    (-1.2,0) arc[start angle=180, end angle=360,
                 x radius=1.2cm, y radius=0.3cm];
  \draw[thmframe!65!black, semithick] (T) -- (-1.2,0);
  \draw[thmframe!65!black, semithick] (T) -- ( 1.2,0);
  \draw[thmframe!65!black, semithick] (B) -- (-1.2,0);
  \draw[thmframe!65!black, semithick] (B) -- ( 1.2,0);
 
  \node at (0,0) {${\O}_{e(s_{\alpha}, w)}^{\ast (m_{\alpha}+1)}$};
  \node[knot] at (T) {};   \node[knot] at (B) {};
  \node[anchor=south] at (0.3, 2.05) {$s_{\alpha} w$};
  \node[anchor=north] at (0.3,-2.05) {$w$};
 
  \node[font=\small, anchor=west] at (1.45,-1.45)
    {\textcolor{thmframe}{\rule{12pt}{1.5pt}}\;$\alpha$};
\end{tikzpicture}
\caption{The space $S\bigl(\O_{e(s_{\alpha}, w)}^{\ast (m_{\alpha}+1)}\bigr)$.}
\label{SO}
\end{figure}
%
%
Observe that $ \widetilde{\bF}^{(m)}$ carries a natural $N$-action making  
the equivalence \eqref{bFwm} $N$-equivariant. Specifically, on each suspension component, an element $ n \in N $ acts as the map
\begin{equation}
\la{Nact}
 \{W_{\alpha}w\} \ast  \O_{e(s_{\alpha}, w)}^{\ast (m_{\alpha} + 1)}\ \to\ 
 \{g W_{\alpha}w\} \ast  \O_{e(gs_{\alpha}g^{-1},\,gw)}^{\ast (m_{\alpha} + 1)}\ ,\quad
 x \mapsto n \cdot x\,,
\end{equation}
where $ g = nT \in W $ and the `dot' denotes the diagonal action on iterated joins.
Thus, we have a $G$-equivariant equivalence
\begin{equation}
\la{cmodel}
F_m(G,T) \,\simeq \,G \times_{N} \widetilde{\bF}^{(m)}
\end{equation}
that provides an explicit model for the space $ F_m(G,T) $. While this model is geometric, it is not compact as the spaces $ S\bigl(\O_{e(s_{\alpha}, w)}^{\ast (m_{\alpha} + 1)}\bigr) $ are built from the complex orbits $ \O_{e(s_{\alpha}, w)} $ (which are homeomorphic to $\c^* $). To remedy this problem we will replace these complex orbits with their real compact counterparts by constructing deformation retracts in an $N$-equivariant way. To this end, for each $ \alpha \in \R_+ $, we consider the orbit  $ \O_{e(s_{\alpha}, 1)} \subset \bF $ and choose a basepoint $ x_\alpha \in \O_{e(s_{\alpha}, 1)}$ in such a way that $ \dot{s}_{\alpha} \cdot x_\alpha = x_\alpha $ for some $ \dot{s}_{\alpha} \in N_{\alpha}$ representing $ s_{\alpha} \in W_{\alpha}$. Clearly, such $x_{\alpha}$'s exist for all $ \alpha  \in \R_+  $, and they define the $N_{\alpha}$-equivariant maps
\begin{equation}
\la{NOTa}
i_{\alpha}:\ T/T_{\alpha} \into \O_{e(s_{\alpha},1)}\,,\quad t T_{\alpha} \mapsto t \cdot  x_\alpha\,,
\end{equation}
which are homotopy equivalences. Next, we consider the union of $N$-orbits of these points $\, N \cdot \{\,x_{\beta}\,\}_{\beta \in \R_+}\, = \,\bigcup_{\beta \in \R_+} N \cdot x_{\beta} $ in $ \bF $ and take its intersection with the complex $\bT$-orbits corresponding to the edges of the moment 
graph $\Gamma $. It easy to check that, for each $ e(s_{\alpha}, w) \in E_{\Gamma} $, the subspace
\begin{equation}
\la{realorb}
\cO_{e(s_{\alpha}, w)}\, := \, 
N \cdot \{\,x_{\beta}\,\}_{\beta \in \R_+}\,\cap\,\O_{e(s_{\alpha}, w)} 
\end{equation}
is a connected retract of $\O_{e(s_{\alpha}, w)}$ which is ($T$-equivariantly) homeomorphic to $T/T_{\alpha} $.
In particular, we have $\,\cO_{e(s_{\alpha}, 1)} = \im(i_{\alpha})\,$ for all $ \alpha \in \R_+ $, where $ i_{\alpha}$ is defined by \eqref{NOTa}. 

Now, replacing in \eqref{wFm*}  the complex orbits $\O_{e(s_{\alpha}, w)} $ with  their compact retracts \eqref{realorb}, we define the functor
\begin{equation}
\la{RwFm*}
\wcF_{m, {\mathbb R}}^{\,*}:\ \Cc(\Gamma) \to \Top^T\ ,\quad w \mapsto \bigvee_{e \in E_{\Gamma}(w)} C_{\{w\}}[\cO_{e(s_{\alpha}, w)}^{\ast (m_{\alpha} + 1)}]\,,\quad 
e(s_{\alpha}, w) \mapsto  \cO_{e(s_{\alpha}, w)}^{\ast (m_{\alpha} + 1)}\,,
\end{equation}
and the space
\begin{equation}
\la{realFm}
\widetilde{\bF}_{\mathbb R}^{(m)} \,:=\, \coprod_{E_{\Gamma}}S \big(\cO^{\ast (m_{\alpha} + 1)}_{e(s_{\alpha}, w)}\bigr)\,\big/\!\sim
\end{equation}
which,  by Lemma~\ref{Lcofres}, represents the homotopy colimit of $ \wcF_{m, {\mathbb R}}^{\,*} $. Since $\,T/T_{\alpha} \cong S^1 \,$, we have
$$
S \big(\cO^{\ast (m_{\alpha} + 1)}_{e(s_{\alpha}, w)}\bigr) \,\cong\,
S\bigl[(T/T_{\alpha})^{\ast (m_{\alpha} + 1)}\bigr]\,\cong\,S^{2m_{\alpha}+2}
$$
Thus, topologically, the space \eqref{realFm} is the union of even-dimensional spheres, the dimensions depending on $m$, attached to each other at their poles `the same way' as the edges of the moment graph $\Gamma $ (see Fig.~\ref{modelFm}).

\usetikzlibrary{calc,shadings}
\begin{figure}[htbp]
\centering
\begin{tikzpicture}[
  knot/.style={circle, fill=white, draw=black!65, semithick,
               inner sep=0pt, minimum size=4pt},
  lead/.style={thin, draw=black!45},
  elabel/.style={font=\tiny, fill=white, inner sep=1.5pt}
]
  \def\coneR{7.5pt}\def\waistR{2.6pt}
  \def\bicone#1#2#3{%
    \coordinate (bcA) at (#2);
    \coordinate (bcB) at (#3);
    \coordinate (bcM) at ($(bcA)!.5!(bcB)$);
    \pgfmathanglebetweenpoints{\pgfpointanchor{bcA}{center}}%
                              {\pgfpointanchor{bcB}{center}}%
    \edef\bcang{\pgfmathresult}%
    \pgfpointdiff{\pgfpointanchor{bcA}{center}}%
                 {\pgfpointanchor{bcB}{center}}%
    \pgfgetlastxy{\bcx}{\bcy}%
    \pgfmathsetlengthmacro{\bcL}{veclen(\bcx,\bcy)/2}%
    \begin{scope}[shift={(bcM)}, rotate=\bcang]
      \shade[left color=#1!40, right color=#1!14]
         (-\bcL,0pt) -- (0pt,\coneR) -- (0pt,-\coneR) -- cycle;
      \shade[left color=#1!14, right color=#1!40]
         ( \bcL,0pt) -- (0pt,\coneR) -- (0pt,-\coneR) -- cycle;
      \shade[left color=#1!26, right color=#1!14]
         (0pt,0pt) ellipse [x radius=\waistR, y radius=\coneR];
      \draw[#1!50!black, densely dotted, thin]
         (0pt,\coneR) arc[start angle=90, end angle=270,
                          x radius=\waistR, y radius=\coneR];
      \draw[#1!65!black, semithick]
         (0pt,-\coneR) arc[start angle=-90, end angle=90,
                           x radius=\waistR, y radius=\coneR];
      \draw[#1!65!black, semithick]
         (-\bcL,0pt) -- (0pt,\coneR) -- (\bcL,0pt) -- (0pt,-\coneR) -- cycle;
    \end{scope}}
 
  \def\ball#1#2#3{\shade[ball color=#1!42] (#2) circle (#3);
    \draw[draw=#1!60!black, thin, opacity=.55] (#2) circle (#3);}
  \pgfmathsetmacro{\Rb}{0.84}
  \pgfmathsetmacro{\rb}{0.56}
  \pgfmathsetmacro{\rr}{\Rb+\rb}
  \pgfmathsetmacro{\fB}{\Rb/\rr}
  \pgfmathsetmacro{\fS}{\rb/\rr}
 
\begin{scope}[shift={(-5.35,0)}, scale=1.06]
  \foreach \k in {0,...,5}{\coordinate (P\k) at ({60*\k}:1.62);}
  \bicone{thmframe}{P0}{P1}   \bicone{defframe}{P1}{P2}
  \bicone{thmframe}{P2}{P3}   \bicone{defframe}{P3}{P4}
  \bicone{thmframe}{P4}{P5}   \bicone{defframe}{P5}{P0}
  \foreach \k in {0,...,5}{\node[knot] at (P\k) {};}
  \draw[lead] (P0) -- ++(0.70,0.44)
     node[anchor=south west, inner sep=1pt, font=\scriptsize, text=black!60]
     {$w$};
\end{scope}
 
  \node[font=\Large, text=black!40] at (-2.05,0) {$\ \rightsquigarrow\ $};
 
\begin{scope}[shift={(1.8,0)}, scale=1.06]
  \foreach \k in {0,...,5}{\coordinate (K\k) at ({60*\k+30}:\rr);}
  \foreach \k in {0,2,4}{\ball{thmframe}{K\k}{\Rb}}
  \foreach \k in {1,3,5}{\ball{defframe}{K\k}{\rb}}
  \foreach \i/\j in {0/1, 2/3, 4/5}{\coordinate (Q\i) at ($(K\i)!\fB!(K\j)$);}
  \foreach \i/\j in {1/2, 3/4, 5/0}{\coordinate (Q\i) at ($(K\i)!\fS!(K\j)$);}
  \foreach \k in {0,...,5}{\node[knot] at (Q\k) {};}
  \draw[lead] ($(K0)+(30:\Rb)$) -- ++(0.76,0.34)
     node[anchor=west, font=\scriptsize, text=black!60] {$S^{2m_\alpha+2}$};
  \draw[lead] ($(K1)+(90:\rb)$) -- ++(0.10,0.62)
     node[anchor=south, font=\scriptsize, text=black!60] {$S^{2m_\beta+2}$};
  \draw[lead] (Q5) -- ++(1.34,0.30)
     node[anchor=south west, inner sep=1pt, font=\scriptsize, text=black!60]
     {$w$};
\end{scope}
 
  \node[font=\small, anchor=west] at (5.05,-1.25)
    {\textcolor{thmframe}{\rule{12pt}{1.5pt}}\;$\alpha$};
  \node[font=\small, anchor=west] at (5.05,-1.85)
    {\textcolor{defframe}{\rule{12pt}{1.5pt}}\;$\beta$};
\end{tikzpicture}
\caption{The geometric model of $\widetilde{F}_m(T)$.}
\label{modelFm}
\end{figure}

The inclusions $ \cO_{e(s_{\alpha}, w)} 
\subset \O_{e(s_{\alpha}, w)} $ induce a morphism of functors $ \wcF_{m, {\mathbb R}}^{\,*} \into  \wcF_{m}^{\,*}$ and hence a natural $T$-equivariant map
$\,\widetilde{\bF}_{\mathbb R}^{(m)} \into \widetilde{\bF}^{(m)} $, which is a homotopy equivalence. Identifying $\,\widetilde{\bF}_{\mathbb R}^{(m)} $ with the image of this last map, we can easily check that the subspace $ \widetilde{\bF}_{\mathbb R}^{(m)} \subset \widetilde{\bF}^{(m)} $
is stable under the $N$-action defined by \eqref{Nact}. Thus, if we replace the space $\widetilde{\bF}^{(m)}$ in \eqref{cmodel} with its compact subspace $ \,\widetilde{\bF}_{\mathbb R}^{(m)} $ we obtain the desired compact model for the space $F_m(G,T)$:
\begin{equation}
\la{rmodel}
F_m(G,T) \,\simeq \,G \times_{N} \widetilde{\bF}_{\mathbb R}^{(m)}
\end{equation}
Note that, for all $ m \leq m' $ in $ \M(W) $,
the natural inclusions of joins 
$ \O^{\ast(m_{\alpha} +1)}_{e(s_{\alpha}, w)} \into \O^{\ast(m'_{\alpha} +1)}_{e(s_{\alpha}, w)} $ restrict to $ \cO^{\ast(m_{\alpha} +1)}_{e(s_{\alpha}, w)} \into \cO^{\ast(m'_{\alpha} +1)}_{e(s_{\alpha}, w)} $, inducing the $G$-equivariant cofibrations
\begin{equation}
\la{rmodelmaps}
G \times_{N} \widetilde{\bF}_{\mathbb R}^{(m)}\,\into \,G \times_{N} \widetilde{\bF}_{\mathbb R}^{(m')}
\end{equation}
that model the maps $\varphi_{m,m'}: F_m(G,T) \to F_{m'}(G,T)\,$, see \eqref{phimm}.

Next, recall that in Section~\ref{S4},  
we introduced a functor $ \wcF_0: \Cc(\Gamma) \to \Top^T $ (see \eqref{FFP0})
and constructed two natural maps $p: \,\wcF_0 \to  \cF_0 $
and $\tilde{q}:\,G \times_T \wcF_0 \to G/T $, the first being a natural weak equivalence 
(see \eqref{GTF} and \eqref{GTFq}). By Lemma~\ref{Lcofres}, $\wcF_0$ is actually cofibrant  in $ \Top^{\Cc(\Gamma)}$, and the map $ p $ is a cofibrant resolution of  $\cF_0 $. Now, observe that, for all $ \alpha \in \R_+ $, there are natural maps 
$$
 C_{\{w\}}[\cO_{e(s_{\alpha}, w)}]\,\to\, \overline{\O}_{e(s_{\alpha}, w)}\!\setminus\!\{s_{\alpha} w\}
$$
extending the inclusions $ \cO_{e(s_{\alpha}, w)} \subset
\O_{e(s_{\alpha}, w)} \subset \overline{\O}_{e(s_{\alpha}, w)}\!\setminus\!\{s_{\alpha} w\} $.
These maps induce an equivalence of diagrams $ \wcF_{0, {\mathbb R}}^{\,*} \,\xrightarrow{\sim} \,\wcF_0 $, and hence, by Lemma~\ref{LemmaL4}, an $N$-equivariant map 
\begin{equation*}
\widetilde{\bF}_{\mathbb R}^{(0)} \,\xrightarrow{\sim}\,\bF^{(1)}
\end{equation*}
which is a homotopy equivalence of spaces. Whence we get the map of $G$-spaces
\begin{equation}
\la{rmap}
G \times_{N} \widetilde{\bF}_{\mathbb R}^{(0)} \,\xrightarrow{\sim}\,
G \times_{N} \bF^{(1)} \xrightarrow{\bar{q}} G/T
\end{equation}
where $\bar{q} $ is the canonical $G$-map extending the inclusion $ \bF^{(1)} \subset G/T $, see \eqref{eq:G-act}. Now, combining \eqref{rmodel}, \eqref{rmodelmaps},
\eqref{rmap} with the result of Corollary~\ref{corpplus}, we obtain the $G$-equivariant equivalence
\begin{equation}
\la{rmodel+}
F^+_m(G,T) \,\simeq\, 
G/T \bigcup_{G \times_{N} \widetilde{\bF}_{\mathbb R}^{(0)}} (G \times_{N} \widetilde{\bF}_{\mathbb R}^{(m)})
\end{equation}
that provides us with a compact model for the space $ F_m^+(G,T)$. 
We will use this model for computing equivariant $K$-theory of quasi-flag manifolds
in Section~\ref{S7}.

\section{Cohomology of quasi-flag manifolds}\la{S6}
In this section, we prove our main theorem relating quasi-flag manifolds to coaffine stacks of quasi-invariants constructed in Section~\ref{S3.2}. Recall that the quasi-flag manifolds are defined by the (relative) plus construction applied to the spaces $ F_m(G,T)$, which are given by two (canonically equivalent) homotopy decompositions \eqref{FmG} and \eqref{FGmdec} over the categories $ \Cc^{(m)}(\Gamma)_{hW}$ and $ \Cc(\Gamma)_{hW} $, respectively. On the other hand,
the coaffine stacks of quasi-invariants are defined as colimits over the poset $ \Sc(W) $ of elementary reflection subgroups of $W$.  To relate these objects we thus need to compare the homotopy colimits of diagrams over the categories  $ \Cc(\Gamma)_{hW} $ and $ \Sc(W) $. The corresponding comparison theorem (Theorem~\ref{3hdec}) is non-trivial: it requires  technical tools and results from abstract homotopy theory that we prove in Appendix~\ref{AC}.

\subsection{Homotopy decomposition over $\Sc(W)$}
\la{S6.2}
Given $ X: \Cc(\Gamma)_{hW} \to \Top $, we define a functor $ X^{\natural}: \Sc(W) \to \Top $ by
\begin{equation}
\la{Fnat}
X^{\natural}(W_0) = X(1)\ ,\quad X^{\natural}(W_{\alpha}) = X(1) \,\amalg_{X[e(s_{\alpha})]} X[e(s_{\alpha})]_{hW_{\alpha}}
\end{equation}
with $ W_0 \leq W_{\alpha} $ mapping to the canonical inclusions 
$ X(1) \into X(1) \,\amalg_{X(e(s_{\alpha}))} X[e(s_{\alpha})]_{hW_{\alpha}} $.
We recall that $ e(s_\alpha) $ denotes the object $ e(s_{\alpha},1) \in \Cc(\Gamma)_{hW} $, which comes with the `edge' map $ f_{\alpha}: e(s_\alpha) \to 1 $ (see Fig.~\ref{skC}). This map is used to define the attaching map for the first summand in \eqref{Fnat}. The second summand in \eqref{Fnat} is attached via the canonical fiber inclusion $ X[e(s_\alpha)] \into X[e(s_{\alpha})]_{hW_{\alpha}} $ into the space of Borel homotopy orbits
$\,X[e(s_{\alpha})]_{hW_{\alpha}} := EW_{\alpha} \times_{W_{\alpha}} X[e(s_\alpha)]\,$. 
Note that this inclusion is always a cofibration in $\Top $, hence the space $X^{\natural}(W_{\alpha})$ represents the homotopy colimit
\begin{equation}
\la{hocnat}
X^{\natural}(W_{\alpha})\, \simeq \, \hocolim\,\{X(1) \leftarrow X[e(s_{\alpha})] \to X[e(s_{\alpha})]_{hW_{\alpha}}\}
\end{equation}
for all $ \alpha \in \R_+$.
\begin{theorem}
\la{3hdec}
For any $ X: \Cc(\Gamma)_{hW} \to \Top $, there is a natural weak homotopy equivalence
\begin{equation}
\la{CGS}
\hocolim_{\Cc(\Gamma)_{hW}}(X) \, \simeq \, \hocolim_{\Sc(W)}(X^{\natural})\ ,
\end{equation}
where the functor $ X^{\natural}: \Sc(W) \to \Top $ is defined by \eqref{Fnat}.
\end{theorem}
\begin{proof}
First, we note that the category $\Cc(\Gamma)_{hW}$ itself is not an EIA-category (it contains isomorphisms which are not automorphisms); however, its skeletal subcategory 
$ \overline{\Cc(\Gamma)}_{hW} $ described in Lemma~\ref{skCG} is EIA. 
Thus, to apply the result of Corollary~\ref{EIA} we first restrict the functor $X$ to $ \overline{\Cc(\Gamma)}_{hW} $. Now, the universal poset $\Pc[\overline{\Cc(\Gamma)}_{hW}]$ under $ \overline{\Cc(\Gamma)}_{hW} $ is simply the set of its objects, and (the opposite of) its barycentric subdivision is shown on Fig.~\ref{sdP}. 
\begin{figure}[h!]
\[\begin{tikzcd}
	&& {(e(s_{\alpha}))} \\
	{(e(s_{\omega}))} &&&& {(e(s_{\beta}))} \\
	&& {(1)} \\
	{(e(s_{\nu}))} &&&& {(e(s_{\gamma}))} \\
	&& {(e(s_{\lambda}))}
	\arrow["{(e(s_{\alpha}), 1)}"{description}, tail reversed, from=1-3, to=3-3]
	\arrow["{(e(s_{\omega}), 1)}"{description}, tail reversed, from=2-1, to=3-3]
	\arrow["{(e(s_{\beta}), 1)}"{description}, tail reversed, from=2-5, to=3-3]
	\arrow["{(e(s_{\nu}), 1)}"{description}, tail reversed, from=4-1, to=3-3]
	\arrow["{(e(s_{\gamma}), 1)}"{description}, tail reversed, from=4-5, to=3-3]
	\arrow["{(e(s_{\lambda}), 1)}"{description}, tail reversed, from=5-3, to=3-3]
\end{tikzcd}\]
    \caption{The poset $\, \sd[\Pc(\overline{\Cc(\Gamma)}_{hW})]^{\rm op}$}
    \label{sdP}
\end{figure}

\noindent
The corresponding functor \eqref{Xnat} is given by
\begin{eqnarray}
\la{Xnatob}
X^{\natural}(1) &=& X(1)\ ,\nonumber\\
X^{\natural}(e(s_\alpha)) &=& X[e(s_{\alpha})]_{hW_{\alpha}}\ ,\quad\\
X^{\natural}(e(s_\alpha),1) &=& 
[W_{\alpha} \times X(e(s_{\alpha}))]_{hW_{\alpha}} 
\cong X(e(s_{\alpha}))\ ,\nonumber
\end{eqnarray}
with morphisms $\,X^{\natural}(1) \leftarrow X^{\natural}(e(s_\alpha),1) \to X^{\natural}(e(s_\alpha)) \,$ induced by the natural maps
\begin{equation}\la{Xnatmor}
X(1)\, \xleftarrow{X(f_{\alpha})}\, X(e(s_{\alpha})) \,\into \, X[e(s_{\alpha})]_{hW_{\alpha}}
\end{equation}

To complete the proof it remains to relate the poset categories $\sd[\Pc(\overline{\Cc(\Gamma)}_{hW})]^{\rm op} $ and $\Sc(W) $. To this end, 
we notice that the former category can be expressed naturally as a colimit  over the latter. Namely,  $ \sd[\Pc(\overline{\Cc(\Gamma)}_{hW})]^{\rm op} $
is obtained by amalgamating the categories $\, \I_{\alpha} = \{\,(1) \leftarrow (e(s_{\alpha}), 1) \to (e(s_{\alpha}))\,\} \,$   along 
their common object $(1) $ --- a diagram in $\Cat$ indexed by $ \Sc(W)$.
Thus, there is an equivalence of categories
$$ 
\sd[\Pc(\overline{\Cc(\Gamma)}_{hW})]^{\rm op} \,\simeq\,\colim_{\Sc(W)}(\Phi)  \,,
$$
where $\Phi$ is defined by 
\begin{equation}
\la{glufun}
\Phi:\, \Sc(W) \to \Cat \ ,\quad  W_{0} \mapsto (1) \ ,\quad
 W_{\alpha} \mapsto \I_{\alpha} \,,
\end{equation}
This suggests that replacing $\colim_{\Sc(W)}(\Phi)$ with $ \hocolim_{\Sc(W)}(\Phi) $
and representing the latter by the Grothendieck construction, we get a `zig-zag' of categories 
\begin{equation}
\la{cat-zigzag}
    \begin{tikzcd}
        & \Sc(W)\smallint \Phi \arrow[dl, "p"']\arrow[dr,"\phi"] & \\
        \Sc(W) && \sd[\Pc(\overline{\Cc(\Gamma)}_{hW})]^{\rm op}
    \end{tikzcd}
\end{equation}
where the functor $\phi $ represents the canonical map: 
$\,\hocolim_{\Sc(W)}(\Phi) \to \colim_{\Sc(W)}(\Phi)\,$.
The category $\,\Sc(W)\smallint \Phi\,$ is obtained from $\sd[\Pc(\overline{\Cc(\Gamma)}_{hW})]^{\rm op}$ by `blowing up' the object $(1)$
to a full subcategory $ \Sc(W) \times \{(1)\}$: thus,
the objects of $\,\Sc(W)\smallint \Phi\,$ comprise the union of sets
$ \{(W_0, 1),\,(W_\alpha,1)\} $ and $ \{(W_\alpha, e(s_{\alpha})),\, (W_\alpha, e(s_{\alpha}), 1)\}$  for all $\alpha \in \R_+ $, while the (non-identity) morphisms are 
$$
(W_{\alpha}, e(s_{\alpha})) \leftarrow (W_\alpha, e(s_{\alpha}), 1) \rightarrow 
(W_\alpha,1) \leftarrow (W_0, 1)\,.
$$
The functor $ \phi $ is then given by the projection
$$
(W_0,1) \mapsto (1)\ ,\quad (W_\alpha,1) \mapsto (1)\ ,\quad (W_\alpha, e(s_{\alpha}), 1) \mapsto (e(s_{\alpha}), 1)\ ,\quad (W_\alpha, e(s_{\alpha})) \mapsto (e(s_{\alpha}))
$$
contracting $\,\Sc(W) \times \{(1)\} \,$ in $\, \Sc(W)\smallint \Phi \,$ to $\,(1)$. A simple verification shows that the comma categories $ \phi_{c/} $ under all objects of $\,\sd[\Pc(\overline{\Cc(\Gamma)}_{hW})]^{\rm op}$ are contractible. Hence the functor $ \phi $ is homotopy cofinal. Combining this fact with the result of Corollary~\ref{EIA}, we conclude
\begin{eqnarray*}
  \hocolim_{\Cc(\Gamma)_{hW}} (X) & \cong & \hocolim_{\overline{\Cc(\Gamma)}_{hW}} (X) \\
  &\simeq& \hocolim_{\sd[\Pc(\overline{\Cc(\Gamma)}_{hW})]^{\rm op}}(X^{\natural})\\
  &\simeq& \hocolim_{\Sc(W)}(p_! \,\phi^* X^{\natural})\\
  &\simeq& \hocolim_{W_\alpha \in \Sc(W)}\,[\,\hocolim_{\Phi(W_\alpha)} (\phi^* X^\natural)\,]
\end{eqnarray*}
It remains to note that the diagram $ \phi^*X^{\natural} $ is represented by \eqref{Xnatmor}, and hence (see \eqref{hocnat}),
$$
\hocolim_{\Phi(W_\alpha)} (\phi^* X^\natural) \,\simeq \, 
\colim_{\Phi(W_\alpha)}(\phi^* X^\natural) \,=\, X^{\natural}(W_{\alpha}) \,.
$$
This completes the proof of the theorem. 
\end{proof}
We now apply Theorem~\ref{3hdec} to our orbit functors \eqref{FmG*}. Recall our basic notation:
$\,T_{\alpha} := \Ker(\alpha) \subseteq T\,$ denotes the singular torus corresponding to the root $\,\alpha \in \R_+ \subset \Lambda(T)$; $\, G_{\alpha} := C_{G}(T_\alpha) $ is the centralizer of $T_\alpha $ in $G$ (which is a compact connected Lie group containing $T$ as a maximal torus of corank 1), and $\, N_\alpha  =  N_G(T) \cap G_{\alpha}\,$ is the  normalizer of $ T $ in $ G_{\alpha}$.

Given  $ m \in \M(W)$, we define a functor $ (G \times_T \cF_m)^\natural: \Sc(W) \to \Top^G $ by 
\begin{equation}
\la{GNW}
W_0 \,\mapsto\, G/T\ ,\quad W_\alpha \,\mapsto\, 
G \times_{N_{\alpha}}\bigl[(G_{\alpha}/T) \ast (T/T_{\alpha})^{\ast\, m_{\alpha}}\bigr]\,,
\end{equation}
with $W_0 \leq W_{\alpha}$ corresponding to the canonical map 
$$ 
G/T = G \times_{N_{\alpha}} (N_{\alpha}/T)\,\to\, G \times_{N_{\alpha}}\bigl[(G_{\alpha}/T) \ast (T/T_{\alpha})^{\ast\, m_{\alpha}}\bigr] 
$$ 
that extends the inclusion $\,N_{\alpha}/T  \into G_{\alpha}/T $ into the first join factor.
The action of $N_{\alpha}$ on the iterated join
$ (G_{\alpha}/T) \ast (T/T_{\alpha})^{\ast\, m_{\alpha}}$ is diagonal: by left translation
on the first factor and by conjugation on the other factors. We also consider the composite $ (G \times_T \cF_m)_{hG}^\natural $ of the above functor with the homotopy quotient $ (\,\mbox{--}\,)_{hG} = EG \times_G (\,\mbox{--}\,): \Top^G \to \Top $, which is explicitly given by
\begin{equation}
\la{GNWG}
W_0 \,\mapsto\, BT\ ,\quad 
W_\alpha \,\mapsto\,  \bigl[(G_{\alpha}/T) \ast (T/T_{\alpha})^{\ast\, m_{\alpha}}\bigr]_{hN_{\alpha}}\,.
\end{equation}
\begin{prop}
\la{natdec}    
For all $m \in \M(W) $, there are natural homotopy decompositions
\begin{eqnarray}
F_m(G,T) & \simeq & \hocolim_{\Sc(W)}\,(G \times_T \cF_m)^\natural\, , \la{Fmnat}\\
F_m(G,T)_{hG} & \simeq & \hocolim_{\Sc(W)}\,(G \times_T \cF_m)_{hG}^\natural \,,\la{Xmnat}
\end{eqnarray}
where the functors $ (G \times_T \cF_m)^\natural $  and $ 
(G \times_T \cF_m)_{hG}^\natural $ are defined by \eqref{GNW} and \eqref{GNWG}, respectively.
\end{prop}
\begin{proof}
By Lemma~\ref{CGdec}, $\,F_m(G,T) \simeq \hocolim_{\Cc(\Gamma)_{hW}}(G \times_T \cF_m^{\,*})\,$, where the functor $\,(G \times_T \cF_m)^*\,$ is defined by \eqref{FmG*} (see also \eqref{dGTFm}). Applying  Theorem~\ref{3hdec} to this last functor, we compute
\begin{eqnarray*}
(G \times_T \cF_m)^{\natural}(W_0) & = & (G \times_T \cF_m)^{\ast}(1) \,= \, G/T\ ,\\*[2ex]
(G \times_T \cF_m)^{\natural}(W_\alpha) &\simeq & 
\hocolim\,\bigl[\,G/T \leftarrow G \times_T \O_{e(s_{\alpha})}^{\ast (m_\alpha + 1)} \to G \times_{N_{\alpha}} \O_{e(s_{\alpha})}^{\ast (m_\alpha + 1)}\bigr]\\
&\simeq & 
G \times_{N_{\alpha}} 
\hocolim\,\bigl[\,N_{\alpha}/T \leftarrow  N_{\alpha} \times_T \O_{e(s_{\alpha})}^{\ast (m_\alpha + 1)} \to \O_{e(s_{\alpha})}^{\ast (m_\alpha + 1)}\bigr]\\
& \cong & G \times_{N_{\alpha}} 
\hocolim\,\bigl[\,N_{\alpha}/T \leftarrow  (N_{\alpha}/T) \times \O_{e(s_{\alpha})}^{\ast (m_\alpha + 1)} \to \O_{e(s_{\alpha})}^{\ast (m_\alpha + 1)}\bigr]\\
&\simeq& 
G \times_{N_{\alpha}}\bigl[(N_{\alpha}/T) \ast \O_{e(s_{\alpha})}^{\ast (m_\alpha + 1)} \bigr]\\
&=&
G \times_{N_{\alpha}}\bigl[W_{\alpha} \ast \O_{e(s_{\alpha})}^{\ast (m_\alpha + 1)}\bigr]\ ,
\end{eqnarray*}
where on the first step we use formula \eqref{hocnat}, and on the third --- the fact that each iterated join
$ \O_{e(s_{\alpha})}^{\ast (m_\alpha + 1)} $ carries a natural (diagonal) action
of $N_{\alpha}$. 
Now, observe that,  for each $\alpha \in \R_+$, there is an $N_{\alpha}$-equivariant homotopy equivalence
\begin{equation} 
\la{hoeqO}
W_{\alpha} \ast \O_{e(s_{\alpha})}\,\xrightarrow{\sim} \, G_{\alpha}/T\,
\end{equation}
that restricts to the natural inclusion $ \O_{e(s_{\alpha})} \into  G_{\alpha}/T $ on the second factor of the join, while mapping the first factor $ W_{\alpha} \subset W_{\alpha} \ast \O_{e(s_{\alpha})} $ to the $T$-fixed points $(G_{\alpha}/T)^T $ in $G_{\alpha}/T$. To construct \eqref{hoeqO} we  identify 
\begin{equation} 
\la{hoeq1}
G_{\alpha}/T \,\cong \,\bP_{\alpha}/\bB\,\cong\,\c\bP^1\,,
\end{equation}
where $\, \bP_{\alpha} = \bB W_{\alpha} \bB \, $ is the minimal 
parabolic subgroup of $ \bG $ corresponding to the root $\alpha$.
The first map in \eqref{hoeq1} is a homeomorphism induced by the canonical inclusion
$ G_{\alpha} \subseteq G \into \bG $, while the second is an isomorphism of complex varieties induced by a root homomorphism $ SL(2, \c) \to \bG $ corresponding to $\alpha$. Under \eqref{hoeq1}, the points $[1]_T$ and $ [s_{\alpha}]_T $ in $G_\alpha/T $ correspond to  $[1:0]$ and 
$[0:1]$, while the orbit $ \O_{e(s_{\alpha})} \cong \{[\lambda:\lambda^{-1}] \in \c\bP^1\ :\ \lambda \in \c^*\} $. The map \eqref{hoeqO} can then be written explicitly (in terms of homogeneous coordinates in $ \c\bP^1$) by
$$
(1-t) e + t [\lambda:\lambda^{-1}] \mapsto [\lambda: t\lambda^{-1}]\ ,
\quad (1-t)s_{\alpha} + t [\lambda:\lambda^{-1}] \mapsto [t\lambda: \lambda^{-1}]\ ,
$$
where $ t \in \Delta^1 $ is the join parameter. It is easy to check that this map is
equivariant with respect to the action of $N_{\alpha}$. 
Now, combining \eqref{hoeqO} with the maps \eqref{NOTa} constructed in Section~\ref{S5.4}, we identify
$$
G \times_{N_{\alpha}}\bigl[W_{\alpha} \ast \O_{e(s_{\alpha})}^{\ast (m_\alpha + 1)}\bigr] \,\simeq \,
G \times_{N_{\alpha}}\bigl[(G_{\alpha}/T) 
\ast \O_{e(s_{\alpha})}^{\ast m_\alpha }\bigr]
\,\simeq \,G \times_{N_{\alpha}}\bigl[(G_{\alpha}/T) \ast (T/T_{\alpha})^{\ast\, m_{\alpha}}\bigr]
$$
that brings the functor $ (G \times_T \cF_m)^{\natural} $ to the required form \eqref{GNW}.
Thus, we get the first decomposition \eqref{Fmnat}; the second, \eqref{Xmnat}, is a formal consequence of the first since homotopy quotients commute with homotopy colimits.
\end{proof}
\begin{remark}
Note that, for each $ \alpha \in \R_+$, there is a natural fibration sequence
$$
G_{\alpha}/N_{\alpha} \to G \times_{N_{\alpha}} (G_{\alpha}/T) \xrightarrow{p_{\alpha}} G/T
$$
where the map $p_{\alpha} $ takes the $N_{\alpha}$-orbit of 
$ (g, [g_{\alpha}]_T) \in G \times (G_\alpha/T)$ to $ [gg_{\alpha}]_T \in G/T$. 
These sequences assemble together into a fibration sequence of 
$\Sc(W)$-diagrams of spaces over the constant $\Sc(W)$-diagram $G/T$. 
Hence, by Puppe's Theorem \cite{Pu74}, taking homotopy colimits, we get the homotopy fibration sequence over the same base:
\begin{equation}
\la{fibseq1}
    \bigvee_{\alpha\in\R_+} G_{\alpha}/N_{\alpha} \,\to\, 
    \hocolim_{\Sc(W)}(G \times_T \cF_0)^{\natural}\, \to \,G/T
    \end{equation}
By Proposition~\ref{natdec}, the total space of \eqref{fibseq1} is equivalent 
to $ F_0(G,T)$, and the sequence \eqref{fibseq1} itself is equivalent 
to \eqref{fibseq}. Thus, for $m=0$, Proposition~\ref{natdec} implies 
the result of Lemma~\ref{lem:act-fib}.
\end{remark}

We now use Proposition~\ref{natdec} to construct an $\Sc(W)$-decomposition
for the quasi-flag manifolds $ F_m^+(G,T) $. We will do this under the additional 
assumption that $m$ is even: namely, for $\, m =  2k \,$, we construct a model for a $G$-equivariant map $ F_m(G,T) \to F^+_m(G,T) $ that represents
a $p$-plus construction on $F_m(G,T) $ for all primes $p \not= 2$. 

We begin by defining a functor $ (G \times_T \cF_{2k})^{\natural,+}: 
\Sc(W) \to \Top^G $:
\begin{equation}
\la{GGW}
W_0 \,\mapsto\, G/T\ ,\quad W_\alpha \,\mapsto\, 
G \times_{G_{\alpha}}\bigl[(G_{\alpha}/T) \ast (G_{\alpha}/T_{\alpha})^{\ast\, k_{\alpha}}\bigr]\,,
\end{equation}
with $W_0 \leq W_{\alpha}$ corresponding to the canonical map 
$$ 
G/T = G \times_{G_{\alpha}} (G_{\alpha}/T)\,\to\, G \times_{G_{\alpha}}\bigl[(G_{\alpha}/T) \ast (G_{\alpha}/T_{\alpha})^{\ast\, k_{\alpha}}\bigr] 
$$ 
extending the natural inclusion  $ G_{\alpha}/T \into (G_{\alpha}/T) \ast (G_{\alpha}/T_{\alpha})^{\ast\, k_{\alpha}}$.
To describe the properties of this functor --- in particular, to relate it to
the functor \eqref{GNW} --- we recall a few basic topological facts about compact
Lie groups. By definition, for each $ \alpha \in \R_+ $, the singular torus $ T_{\alpha} $
is a central (and hence, normal) subgroup of $ G_{\alpha} $. Hence $G_{\alpha}/T_{\alpha} $
is a compact connected Lie group of rank 1 with maximal torus $T/T_{\alpha} $
(see \cite[Lemma~V.3.25]{MT78}). By classification of such Lie groups
(see {\it loc. cit.}, Corollary~V.3.27 and Remark~V.3.28), there are two possibilities: $G_{\alpha}/T_{\alpha} \,\cong\,SU(2) $ or $G_{\alpha}/T_{\alpha} \,\cong\,SO(3) $, depending on whether $ G_{\alpha}/T_{\alpha} $ is simply connected or 
not. In either case, there is an $ N_{\alpha}$-equivariant map
\begin{equation}
\la{TTG}
(T/T_{\alpha})\ast (T/T_{\alpha})\,\to\, G_{\alpha}/T_{\alpha}
\end{equation}
which is a homeomorphism when $ \pi_1(G_{\alpha}/T_{\alpha}) = 0 $, and
a universal covering  when $ \pi_1(G_{\alpha}/T_{\alpha}) = \Z/2 $
(see \cite[(3.36)]{BR1} for an explicit formula for this map). Iterating \eqref{TTG}, we can define a morphism of functors $\,(G \times_T \cF_{2k})^{\natural} \to (G \times_T \cF_{2k})^{\natural,+}\,$ by
\begin{equation}
\la{NGmaps}
G \times_{N_{\alpha}}\bigl[(G_{\alpha}/T) \ast (T/T_{\alpha})^{\ast\, 2k_{\alpha}}\bigr] 
\ \to\ G \times_{G_{\alpha}}\bigl[(G_{\alpha}/T) \ast (G_{\alpha}/T_{\alpha})^{\ast\, k_{\alpha}}\bigr] 
\end{equation}
This induces a map of homotopy colimits:
\begin{equation}
\la{NGcol}
\hocolim_{\Sc(W)} \,(G \times_T \cF_{2k})^\natural\,\to
\,\hocolim_{\Sc(W)} \,(G \times_T \cF_{2k})^{\natural,+}\,
\end{equation}
By Proposition~\ref{natdec}, the domain of the map \eqref{NGcol} represents the space $ F_{2k}(G,T)$. We claim that its target is a $p$-plus construction $F_{2k}^{+}(G,T)$. Since the homotopy fiber of \eqref{NGmaps} is $ G_{\alpha}/N_{\alpha} \cong {\mathbb R}{\bP}^2 $ if
$ \pi_1(G_{\alpha}/T_{\alpha}) = 0 $, or in any case $2$-local, the map \eqref{NGcol} is a mod-$p$ cohomology isomorphism for all $ p \not=2$. To show that it represents the $p$-construction we need to check that the target of \eqref{NGcol} is simply connected. For this, by Seifert-Van Kampen Theorem, it suffices to show that
\begin{equation*}
 \pi_{1}\bigl(G \times_{G_{\alpha}}\bigl[(G_{\alpha}/T) \ast (G_{\alpha}/T_{\alpha})^{\ast\, k_{\alpha}}\bigr] \bigr) \ = \ 0  
\end{equation*}
for all $ \alpha \in \R_+ $ and all $ k_{\alpha} \ge 0 $. This follows from the long homotopy exact sequence associated to the natural fibration sequence
$$
(G_{\alpha}/T) \ast (G_{\alpha}/T_{\alpha})^{\ast\, k_{\alpha}}\, \to\,
G \times_{G_{\alpha}}\bigl[(G_{\alpha}/T) \ast (G_{\alpha}/T_{\alpha})^{\ast\, k_{\alpha}}\bigr]\,\to\, G/G_{\alpha}
$$
where both the basespace and the homotopy fibre are simply connected. Indeed,
the fact that $ \pi_1(G/G_{\alpha}) = 0 $ follows from the fibration sequence
$ G_{\alpha}/T \to G/T \to G/G_{\alpha} $, where $\,G_{\alpha}/T \,$ and $G/T$ are simply connected (since $T$ is a maximal torus in both $G$ and $G_{\alpha}$). The fact
$(G_{\alpha}/T) \ast (G_{\alpha}/T_{\alpha})^{\ast\, k_{\alpha}}$ is simply connected
for all $k_{\alpha}$ is a consequence of Milnor's Lemma (see, e.g., \cite[Lemma~A.2]{BR1}).
Thus, we conclude that, for $ m = 2k $,
\begin{equation}
\la{F2k+}
F^+_{2k}(G,T)\,\simeq\,
\,\hocolim_{\Sc(W)} \,(G \times_T \cF_{2k})^{\natural,+}\,
\end{equation}
where the functor $ (G \times_T \cF_{2k})^{\natural,+} $ is given by \eqref{GGW}.
Since the Borel quotient construction commutes with homotopy colimits, the equivalence \eqref{F2k+} induces the following decomposition for the space of homotopy $G$-orbits in $F^+_{2k}(G,T)$:
\begin{equation}
\la{C2k+}
F^+_{2k}(G,T)_{hG}\,\simeq\,
\,\hocolim_{\Sc(W)} \,(G \times_T \cF_{2k})_{hG}^{\natural,+}\,,
\end{equation}
where $ (G \times_T \cF_{2k})^{\natural,+}_{hG}: \, \Sc(W) \to \Top\, $ is defined by
\begin{equation}
\la{XGGW}
W_0 \mapsto BT\ ,
\quad
W_{\alpha} \mapsto [(G_{\alpha}/T) \ast (G_{\alpha}/T_{\alpha})^{\ast k_{\alpha}}]_{h G_{\alpha}}
\end{equation}
We remark that the functors \eqref{GGW} and \eqref{XGGW} can be obtained by a relative version of the classical Ganea construction studied in our forthcoming paper.

\subsection{Algebraic models of quasi-flag manifolds}\la{S6.3}
Recall that if $X$ is a topological space, the cochain complex
$C_{\k}^*(X)$ of $X$ with coefficients in a field $\k$ has the natural structure 
of an $\bE_{\infty}$-algebra over $\k$ (see Section~\ref{A1}). Furthermore, if $\k = \Q $ and $X$ is a simply connected space of finite type over $\Q$, then $ C_{\Q}^*(X)$ defines a coaffine stack $\, \cSpec\,[C_{\Q}^*(X)] \in \cAff_{\Q} \,$ 
that determines the rational homotopy type of $X$. As we explained in Appendix~\ref{AB},
this is a natural generalization of the classical theorem of D. Sullivan \cite{Su77}  (see Theorem~\ref{SuTh}) stated and proved in the context of derived algebraic geometry by B. To\"en \cite{To06} and J. Lurie \cite{DAGXIII}. We apply this result to our basic 
$G$-spaces $ F_m^+ = F_m^+(G,T) $ defined by \eqref{maindef} and their homotopy $G$- and 
$T$-quotients:
\begin{eqnarray}
X_m^+(G,T) &:=& F_m^{+}(G,T)_{hG} \,=\, EG \times_G F^+_m(G,T) \,,   \la{XmG}\\*[1ex]
\bX_m^+(G,T) &:=& F_m^{+}(G,T)_{hT} \,=\, ET \times_T F^+_m(G,T) \,.   \la{XmT}
\end{eqnarray}
We note that, by Corollary~\ref{corpplus}, all these spaces are simply connected, provided so is the Lie group $G$. The next theorem is the main result of the present paper.

\vspace*{1ex}

\begin{theorem}
\la{MainT}
Let $G$ be a compact simply connected Lie group with maximal torus $T \subseteq G $
and Weyl group $ W = W_G(T) $. Then, for any $ m \in \M(W) $, there are natural equivalences in $\cAff_{\Q} $:
\begin{eqnarray}
\cSpec \,[\,C_{\Q}^*(F^+_m)\,] & \simeq & F^c_{k}(W)\,, \la{coaffF}\\
\cSpec \,[\,C_{\Q}^*(X^+_m)\,] & \simeq & V^c_{k}(W)\,,\la{coaffX}\\
\cSpec \,[\,C_{\Q}^*(\bX^+_m)\,] & \simeq & \bU^c_{\!k}(W)\,, \la{coaffbX}
\end{eqnarray}
where $ \,k = [\frac{m+1}{2}]\,$ and the coaffine stacks $V^c_{k}(W)$, $\,F^c_{k}(W)$ and 
$ \bU^c_{\!k}(W) $ are defined in Section~\ref{S3.2} $($see \eqref{defvmc}, \eqref{Fck} and \eqref{bvck} respectively$)$.
\end{theorem}      
\begin{proof}[Proof of \Cref{MainT}]
We first establish the equivalence \eqref{coaffX} for the spaces $ X_m^+ = X^+_m(G,T)$. As explained in  Appendix~\ref{AB}, we can replace the classical cochain functor $C_{\Q}^*(\,\mbox{--}\,) $ with the Bousfield-Gughenheim functor $A_{\rm PL} $ (see \eqref{CAPL}). This last functor has good functorial properties: in particular, 
being a left Quillen contravariant functor, it maps arbitrary (small) homotopy colimits to homotopy limits: i.e.
\begin{equation}
\la{APLcolim}    
A_{\rm PL}[{\rm hocolim}_{i \in \I} \,X(i)] \,\simeq\,{\rm holim}_{i\in \I^{\rm op}}\,A_{\rm PL}(X(i))
\end{equation}
for any small $\I$-diagram of spaces. Moreover, it also preserves homotopy pullbacks:
\begin{equation}
\la{APLpull}  
A_{\rm PL}(X \times^{h}_Z Y) \,\simeq\, A_{\rm PL}(X) \,\otimes^{\bL}_{A_{\rm PL}(Z)} \,A_{\rm PL}(Y)\,,
\end{equation}
provided $Z$ is simply connected. 

Now, by definition, the spaces $X_m^+(G,T)$ fit in the homotopy fibration sequence
\begin{equation}\la{fibseq+}
F_m^+(G,T)\,\to\, X_m^+(G,T)\,\xrightarrow{p^+_m}\, BG   
\end{equation}
which shows (in combination with part $(b)$ of Corollary~\ref{corpplus}) that 
all $X_m^+(G,T)$ are simply connected. Further, by part $(a)$ of Corollary~\ref{corpplus}, the natural map $ q_m: F_m(G,T) \to F^+_m(G,T) $ is
a rational cohomology isomorphism. Defining $\, X_m(G,T) := F_m(G,T)_{hG} \,$,
by functoriality of the Borel construction, we get the induced map
$ \tilde{q}_{m}:\, X_m(G,T)\,\to \, X^+_m(G,T)$, which is also a rational cohomology isomorphism. Thus, for all $m \in \M(W)$, $ \tilde{q}_m $ induces an
equivalence
\begin{equation}
\la{APLX}
A_{\rm PL}[X^+_m(G,T)]\,\simeq\, A_{\rm PL}[X_m(G,T)]  
\end{equation}
Next, since the Borel quotients commute with homotopy colimits, 
Proposition~\ref{natdec} yields
\begin{equation}
\la{XmGT}
X_m(G,T)\,\simeq\,\hocolim_{\Sc(W)}\,(EG \times_T \cF_m)^\natural
\end{equation}
where $ (EG \times_T \cF_m)^\natural: \Sc(W) \to \Top $ is given by
\begin{eqnarray*}
\la{XGNW}
(EG \times_T \cF_m)^\natural(W_0) &=& EG \times_G (G/T) \simeq BT \ ,\\
(EG\times_T \cF_m)^\natural(W_\alpha) &=& 
EG \times_{N_{\alpha}}\bigl[(G_{\alpha}/T) \ast (T/T_{\alpha})^{\ast\, m_{\alpha}}\bigr]\,\simeq \,\bigl[(G_{\alpha}/T) \ast (T/T_{\alpha})^{\ast\, m_{\alpha}}\bigr]_{hN_{\alpha}}
\end{eqnarray*}
Since $ N_{\alpha} \cong T \rtimes W_{\alpha}\,$, there is an equivalence of homotopy quotients
$$
\bigl[(G_{\alpha}/T) \ast (T/T_{\alpha})^{\ast\, m_{\alpha}}\bigr]_{hN_{\alpha}} \,\simeq\, \bigl[\bigl((G_{\alpha}/T) \ast (T/T_{\alpha})^{\ast\, m_{\alpha}}\bigr)_{hT}\,\bigr]_{hW_{\alpha}}
$$
which induces
$$
A_{\rm PL}\bigl[\bigl((G_{\alpha}/T) \ast(T/T_{\alpha})^{\ast\, m_{\alpha}}\bigr)_{hN_{\alpha}}\bigr]\ \simeq\ 
A_{\rm PL}\bigl[\bigl((G_{\alpha}/T) \ast(T/T_{\alpha})^{\ast\, m_{\alpha}}\bigr)_{hT}\bigr]^{W_{\alpha}}
$$
for all $ \alpha \in \R_+ $. Then, thanks to \eqref{APLcolim}, $ \, A_{\rm PL}[X_m(G,T)]$ can be represented as the homotopy limit of the following diagram of dg algebras in
$ \ccdga_\k^{\ge 0} $ indexed by $\Sc(W)^{\rm op}$:

\begin{equation}
\la{APLdia}
\begin{tikzcd}[scale cd=0.8]
	{A_{\rm PL}\bigl[\bigl((G_{\alpha_1}/T) \ast (T/T_{\alpha_1})^{\ast\, m_{\alpha_1}}\bigr)_{hT}\bigr]^{W_{\alpha_1}}} & \ldots\ldots & {A_{\rm PL}\bigl[\bigl((G_{\alpha_r}/T) \ast (T/T_{\alpha_r})^{\ast\, m_{\alpha_r}}\bigr)_{hT}\bigr]^{W_{\alpha_r}}} \\
	& {A_{\rm PL}(BT)}
	\arrow[from=1-1, to=2-2]
	\arrow[from=1-3, to=2-2]
\end{tikzcd}
\end{equation}
Clearly, $\,A_{\rm PL}(BT) \cong \Q[V] \,$, while, by Lemma~\ref{APLjoin} below, 
$$
A_{\rm PL}\bigl[\bigl((G_{\alpha}/T) \ast (T/T_{\alpha})^{\ast\, m_{\alpha}}\bigr)_{hT}\bigr]^{W_{\alpha}}\ \simeq\ 
\bigl(\,\Q[V]\,\times_{\Q[V]/(\alpha)^{m_{\alpha}+1}}\,\Q[V]\bigr)^{W_\alpha} \,\cong\, Q_{[\frac{m_{\alpha}+1}{2}]}(W_{\alpha})
$$
Thus, we see that the diagram \eqref{APLdia} is equivalent in 
$ \ccdga_\k^{\ge 0} $ to \eqref{Fgpull} for $ k = [(m+1)/2] $. Hence,
by \eqref{hocvmc}, the homotopy limit of this diagram represents the coaffine
stack $ V^c_{[\frac{m+1}{2}]}(W)$. In combination with \eqref{APLX} and \eqref{CAPL}, this proves the equivalence \eqref{coaffX}.

Now, the equivalences \eqref{coaffF} and \eqref{coaffbX} follow formally from \eqref{coaffX}. Indeed, since 
$ A_{\rm PL}(BG) \simeq A_{\rm PL}(BT)^W $,  there is an equivalence in $\cAff_\Q$:
\begin{equation}
\la{coaffBG}
\cSpec\,[C_{\Q}^*(BG)]\,\simeq \, \cSpec [A_{\rm PL}(BG)]  \,\simeq\, V^c\!/\!/W\,.
\end{equation}
Then, since $ BG $ is simply connected, by \eqref{APLpull}, it follows  from \eqref{fibseq+} that
$$
A_{\rm PL}[F_m^+(G,T)] \,\simeq\,A_{\rm PL}(X_m^+(G,T) \times^h_{BG} \{\ast\}) \,
\simeq\, A_{\rm PL}[X_m^+(G,T)] \otimes_{A_{\rm PL}(BG)}^{\bL} \Q
$$
Hence, by formula \eqref{pullb}, the coaffine stack $ \cSpec\,A_{\rm PL}[F_m^+(G,T)] $ represents the (homotopy) fibre of the natural map 
$$
\cSpec\,A_{\rm PL}[X_m^+(G,T)]\,\to\, \cSpec [A_{\rm PL}(BG)]
$$
By \eqref{coaffX} and \eqref{coaffBG}, this last map is equivalent in $ \cAff_{\Q} $ to the projection $\, V^c_{[\frac{m+1}{2}]}(W) \to V^c\!/\!/W$, and hence, by \eqref{Fck}, its fibre is equivalent to $F^c_{[\frac{m+1}{2}]}(W)$. This proves \eqref{coaffF}, and the argument 
for \eqref{coaffbX} is completely similar.
\end{proof}
\begin{lemma}
\la{APLjoin}
For all $ m \ge 0 $, there is a  $W_{\alpha}$-equivariant  quasi-isomorphism of dg algebras
\begin{equation*}
A_{\rm PL}\bigl[\bigl((G_{\alpha}/T) \ast (T/T_{\alpha})^{\ast\, m}\bigr)_{hT}\bigr]\ \simeq\ \Q[V]\,\times_{\Q[V]/(\alpha)^{m+1}}\,\Q[V] \,.
\end{equation*}
Consequently, 
$$
A_{\rm PL}\bigl[\bigl((G_{\alpha}/T) \ast (T/T_{\alpha})^{\ast\, m}\bigr)_{hT}\bigr]^{W_{\alpha}}\ \simeq\  Q_{[\frac{m+1}{2}]}(W_{\alpha})
$$
\end{lemma}
\begin{proof}
First, we note that 
$$
\bigl[(G_{\alpha}/T) \ast (T/T_{\alpha})^{\ast\, m}\bigr]_{hT}\,\simeq\,
(G_{\alpha}/T)_{hT} \ast_{BT} \bigl[(T/T_{\alpha})^{\ast\, m}\bigr]_{hT}\ ,
$$
where $ X \ast_Z Y $ stands for the relative join of two maps $ X \to Z $ and $ Y \to Z $
in $\Top$ defined by 
$$ 
X \ast_Z Y := \hocolim(X \leftarrow X \times^h_Z Y \to Y ).
$$
Now, we can identify
\begin{equation}
\la{GaT}
A_{\rm PL}\left[(G_{\alpha}/T)_{hT}\right]\,\simeq\, H^*_T(G_{\alpha}/T, \Q)
\,\simeq\,  \Q[V]\,\times_{\Q[V]/(\alpha)}\,\Q[V]\,\cong\, \Q[V] \oplus (\alpha)\ ,
\end{equation}
where
the last isomorphism is an isomorphism of graded $\Q[V]$-modules given by $(p,q) \mapsto (\frac{1}{2}(p+q),\,\frac{1}{2}(p-q))$. Furthermore, by induction on $m$, it is
easy to check that
\begin{equation}
\la{TaT}
A_{\rm PL}\left[\bigl((T/T_{\alpha})^{\ast\, m}\bigr)_{hT}\right] \,\simeq\,
H_T^*((T/T_{\alpha})^{\ast\, m},\,\Q) \,\simeq\, \Q[V]/(\alpha)^{m}.
\end{equation}
Then, using \eqref{APLpull}, we conclude from \eqref{GaT} and \eqref{TaT} that
\begin{eqnarray*}
A_{\rm PL}\left[(G_{\alpha}/T)_{hT} \times^h_{BT} \bigl((T/T_{\alpha})^{\ast\, m}\bigr)_{hT}\right] &\simeq & A_{\rm PL}\left[(G_{\alpha}/T)_{hT}\right] \otimes^{\bL}_{A_{\rm PL}(BT)} A_{\rm PL}\left[\bigl((T/T_{\alpha})^{\ast\, m}\bigr)_{hT}\right]\\*[2ex]
&\simeq& \left(\Q[V]\,\times_{\Q[V]/(\alpha)}\,\Q[V]\right) \otimes_{\Q[V]} \Q[V]/(\alpha)^{m}\\*[2ex]
&\cong& \Q[V]/(\alpha)^{m}\,\oplus\,(\alpha)/(\alpha)^{m+1}
\end{eqnarray*}
and hence, by \eqref{APLcolim}, the dg algebra
$A_{\rm PL}\bigl[(G_{\alpha}/T)_{hT} \ast_{BT} \bigl((T/T_{\alpha})^{\ast\, m}\bigr)_{hT} \bigr]$ is equivalent to the pullback of the following diagram in $\ccdga^{\ge 0}_{\Q}\,$:
$$
\Q[V] \oplus (\alpha)\,\onto\,\Q[V]/(\alpha)^{m}\,\oplus\,(\alpha)/(\alpha)^{m+1} \,\hookleftarrow\, \Q[V]/(\alpha)^{m}\,,
$$
where the first map is the natural projection and the second is the inclusion into the first factor. Thus, we conclude
$$
A_{\rm PL}\bigl[(G_{\alpha}/T)_{hT} \ast_{BT} \bigl((T/T_{\alpha})^{\ast\, m}\bigr)_{hT} \bigr]\,\simeq\, \Q[V] \oplus (\alpha)^{m+1}\,\cong\, 
 \Q[V]\,\times_{\Q[V]/(\alpha)^{m+1}}\,\Q[V] \,, 
$$
proving the first claim of the lemma. Now, note that the $W_{\alpha}$-action on the above fiber product is given by
$\,(p,q) \mapsto (s_{\alpha}(q),\,s_{\alpha}(p))\,$. Hence the composite 
map
$$
\left(\Q[V]\,\times_{\Q[V]/(\alpha)^{m+1}}\,\Q[V]\right)^{W_\alpha}\,\into\,
 \Q[V]\,\times_{\Q[V]/(\alpha)^{m+1}}\,\Q[V]\,\xrightarrow{{\rm pr}_1}\,\Q[V]
$$
is injective, and its image consists of the polynomials $ p \in \Q[V] $ satisfying
$ s_{\alpha}(p) \equiv p\ \mbox{mod}\ (\alpha)^{m+1}$. The  second claim
of the lemma follows.
\end{proof}

We now state a few interesting consequences of Theorem~\ref{MainT}. Recall our basic algebras of
quasi-invariants, quasi-coinvariants and quasi-covariants introduced in Section~\ref{S2}:
\begin{eqnarray*}
Q_k(W) &:= & \{p \in \k[V]\ :\ s_{\alpha}(p) \equiv p\ \, \mod\,\langle \alpha \rangle^{2 k_{\alpha}}\ ,\ \forall\,\alpha \in \R_+ \}\ ,  \\*[1ex]
\H_k(W) & := & Q_{k}(W)/\langle Q_{k}(W)^W_{+}\rangle\ =\ Q_{k}(W)/\langle \k[V]^W_{+}\rangle\ ,\\
\bQ_m(W) & := & \{(p_w)_{w\in W} \in \prod_{w \in W}\k[V]\ :\ p_{s_{\alpha} w} \equiv p_w\ \, \mod\,\langle \alpha \rangle^{m_{\alpha}}\ ,\ \forall\,\alpha \in \R_+\,,\ \forall\,w \in W \}\ .
\end{eqnarray*}

\begin{cor} 
\la{CorHev}
Let $\k$ be a field of characteristic $0$. Then, for any $ m \in \M(W) $, there are natural isomorphisms of graded algebras
    \begin{eqnarray}
    \la{evencohom}
H^{\rm ev}(F^+_m(G,T),\,\k) \,\cong\, H^{\rm ev}(F_m(G,T),\,\k) & \cong & \H_{k}(W)\ ,\nonumber\\*[1ex]
H_G^{\rm ev}(F^+_m(G,T),\,\k) \,\cong\, H_G^{\rm ev}(F_m(G,T),\,\k) & \cong & Q_{k}(W)\ ,\\*[1ex]
    H_T^{\rm ev}(F^+_m(G,T),\,\k) \,\cong\, H_T^{\rm ev}(F_m(G,T),\,\k) & \cong & \bQ_{2k+1}(W)\nonumber\ ,
\end{eqnarray}
where $\,k = [\frac{m+1}{2}]\,$.
\end{cor}
\begin{proof}
This follows from Theorem~\ref{MainT} combined with results of Section~\ref{S3.2}: namely,  Theorem~\ref{Tcc}, Theorem~\ref{Hcomp} and Corollary~\ref{CHcomp}.
\end{proof}
\begin{remark}
By Corollary~\ref{CorHev}, the spaces $X_{2k}(G,T) = F_{2k}(G,T)_{hG}$ and $ X^+_{2k}(G,T) = F_{2k}^+(G,T)_{hG}$ provide topological realizations for 
the algebras $Q_k(W)$ of classical $W$-quasi-invariants. This generalizes the result of \cite{BR1}  for $ G = SU(2) $: in fact, the decomposition \eqref{C2k+} for 
$  X^+_{2k}(G,T) $ shows that these spaces are equivalent to those constructed in \cite{BR1}
in the rank one case.
\end{remark}

Note that the last isomorphism in \eqref{evencohom} can be viewed as a {\it GKM-type presentation} for $T$-equivariant cohmomology of quasi-flag manifolds that generalizes the well-known GKM presentation for $ H^*_{T}(G/T,\,\k) $ due to Kumar-Kostant \cite{KK87} and Goresky-Kottwitz-MacPherson \cite{GKM98}. We have also a natural generalization of the Borel presentation. As in the classical case, this arises from the canonical map between two Borel fibration sequences:
\begin{equation} \la{Borseqs}
\begin{diagram}[small]
G/T & \rTo & F^{(+)}_m(G,T)_{hT} & \rTo & F^{(+)}_m(G,T)_{hG}\\
\| & & \dTo & & \dTo\\
G/T & \rTo & BT & \rTo & BG
\end{diagram}
\end{equation}
where $ F^{(+)}_m(G,T) $ is either $ F_m(G,T) $ or $ F^{+}_m(G,T) $.
In fact, \eqref{Borseqs} induces the map of spaces $\,F^{(+)}_m(G,T)_{hT} \,\to\, BT \times_{BG} F^{(+)}_m(G,T)_{hG} \,$, which, in turn, induces a map on even-dimensional cohomology:
$$
H_T^{*}(BT,\,\k) \,\otimes_{ H_T^{*}(BG,\,\k)} \, H_G^{\rm ev}(F^{(+)}_m(G,T),\,\k) \ \to \  H_T^{\rm ev}(F^{(+)}_m(G,T),\,\k)
$$
Now, by Corollary~\ref{CorHev} and Proposition~\ref{bQodd} (see \eqref{oddisoalg}), this last 
map must be an isomorphism of algebras which yields a {\it Borel-type presentation}
for $T$-equivariant cohomology of $ F_m(G,T)$:
\begin{equation}
  H_T^{\rm ev}(F^{(+)}_m(G,T),\,\k)\,\cong\, \k[V] \otimes_{\k[V]^W} Q_{[\frac{m+1}{2}]}(W) 
\end{equation}
As a consequence of Theorem~\ref{ThFree} and Theorem~\ref{Qfatodd}, we have
\begin{cor}\la{CorHev2}
For all $ m \in \M(W) $, $\, H_T^{\rm ev}(F^{(+)}_m(G,T),\,\k) \,$ is a free $ H^*(BT)$-module of rank $ |W|$. This module structure extends to a left module structure over the nil-Hecke algebra $ \NH(W) $.
\end{cor}
\begin{remark}
\la{leftDem}
When $m=0$, the nil-Hecke action of Corollary~\ref{CorHev2}    
specializes to the so-called {\it left} nil-Hecke action on $\, H_T^{\rm ev}(F^{+}_0(G,T),\,\k) = H^*_T(G/T,\,\k) \,$ originally appeared in \cite{KK87} and \cite{Br97} (see \cite[Section 16.5]{AF23}).
\end{remark}
Finally, recall the $T$-space  $\tilde{F}_m(T)$ defined by \eqref{tilFmT}.
The corresponding $G$-space  $ \tilde{F}_m(G,T) := G \times_T \tilde{F}_m(T) $ carries a natural free $W$-action defining a $W$-bundle \eqref{Wbund} over $ F_m(G,T)$:
\begin{equation*}\la{Wbundle}
W \to \tilde{F}_m(G,T) \to F_m(G,T)
\end{equation*}
\begin{prop}
\la{tilProp} 
For all $ m \in \M(W) $, there are natural isomorphisms of algebras
\begin{equation}
\la{bQtop}
H_G^{\rm ev}(\tilde{F}_m(G,T),\,\k) \,\cong\, H_T^{\rm ev}(\tilde{F}_m(T),\,\k)\,\cong\, \bQ_{m+1}(W)\,,
\end{equation}
which fit in the commutative diagram
\begin{equation}
    \begin{diagram}[small]
H_G^{\rm ev}(F_m(G,T),\,\k) & \rTo^{\cong} & H_G^{\rm ev}(\tilde{F}_m(G,T),\,\k)^W\\
\dTo^{\eqref{evencohom}} & & \dTo_{\eqref{bQtop}} \\
Q_{[\frac{m+1}{2}]} & \rTo^{\eqref{isoqq}} & \bQ_{m+1}(W)^W 
\end{diagram}
\end{equation}
where the top isomorphism is induced by the quotient map of the $W$-bundle \eqref{Wbund} and the one on the bottom is defined in $($the proof of$)$ Lemma~\ref{bVV}.
\end{prop}
\begin{proof}
We claim that   
\begin{equation} \la{apltfm} A_{\rm PL}[\tilde{F}_m(T)_{hT}]\,\simeq\,\bQ^{\ast}_{m+1}(W)\ .\end{equation}
In order to prove \eqref{apltfm}, we replace $\tilde{F}_m(T)$ by its compact model $\widetilde{\bF}_{\mathbb{R}}^{(m)}$ (see \eqref{realFm}). Observe that 
\begin{equation} \la{tfmwedgeh} \widetilde{\bF}_{\mathbb{R}}^{(m)}\,=\, \bigvee_{\substack{{\bF}(W_0)\\\alpha \in \A}} \bF^{(m_\alpha)}(W_{\alpha})\,,\end{equation}
where
\begin{equation} \la{defbf}
\bF(W_0) \,=\,\coprod_{w \in W} \{w\}\,, \qquad
\bF^{(m_\alpha)}(W_\alpha) \,=\,\coprod_{e(s_\alpha,w) \in E_{\Gamma}} S(\cO_{e(s_\alpha,w)}^{\ast (m_\alpha+1)})\,, \end{equation}
and for each $\alpha \in \A$, the map $f_\alpha\,:\,\bF(W_0) \into \bF^{(m_\alpha)}(W_\alpha)$ maps the elements of each pair $\{w, s_\alpha w\} \subseteq \bF(W_0)$ to the poles of $S(\cO_{e(s_\alpha,w)}^{\ast (m_\alpha+1)})$. Since  each $f_\alpha$ is injective (i.e., a cofibration), \eqref{tfmwedgeh} represents the {\it homotopy colimit} of the functor $\bF^{(m)}\,:\,\Sc(W) \to \Top^T$ sending $W_0$ to $\bF(W_0)$ and $W_\alpha$ to $\bF^{(m_\alpha)}(W_\alpha)$. Hence, 
\begin{align} \la{apltilfm} A_{\rm PL}[\tilde{F}_m(T)_{hT}] & \simeq\, A_{\rm PL}[(\widetilde{\bF}_{\mathbb{R}}^{(m)})_{hT}] \,\simeq\, A_{\rm PL}[(\hocolim_{\Sc(W)} \bF^{(m)})_{hT}]\\ 
& \simeq\,A_{\rm PL}[\hocolim_{\Sc(W)} (\bF^{(m)}_{hT})]\,\simeq\, \holim_{\alpha \in \Sc(W)^{\rm op}} A_{\rm PL}[\bF^{(m_\alpha)}(W_\alpha)_{hT}]\,,\nonumber\end{align}
where the last (weak-)equivalence is because $A_{\rm PL}$ is a left Quillen contravariant functor. Next, by \eqref{defbf}, and by (the proof of) Lemma \ref{APLjoin},
\begin{equation} A_{\rm PL}[(\bF^{(m_\alpha)}(W_\alpha))_{hT}] \,\simeq\, \prod_{\{w,s_\alpha w\} \in W_\alpha \backslash W} \Q[V] \times_{\Q[V]/(\alpha^{(m_\alpha+1)})} \Q[V] \,\cong\, \bQ_{m_\alpha+1}(W_\alpha)\,,\ \end{equation}
where the last isomorphism is the identification \eqref{idenkwb}. Clearly, $A_{\rm PL}(\bF(W_0)) \simeq \Q[V] \otimes \Q[W]$. In other words, $A_{\rm PL}[\tilde{F}_m(T)_{hT}]$ is the homotopy limit of the contravariant functor $\Sc(W) \to \ccdga^{\geq 0}$ sending $W_0$ to $\Q[V] \otimes \Q[W]$ and $W_\alpha$ to $\bQ_{m_\alpha+1}(W_\alpha)$. By \eqref{QcW}, this is precisely $\bQ^{\ast}_{m+1}(W)$ (the fiber product therein is an explicit model of the above homotopy limit). 
This proves \eqref{apltfm}. The first assertion of Proposition \ref{tilProp} then follows immediately from Theorem \ref{Hcomp}. To verify the second assertion, note that the map $\coprod_{w \in W} \{w\} \into \bF^{(m)}_{\mathbb R}$ is $N(T)$-equivariant. Hence the $W$-action on $H^{\rm ev}_T(\bF^{(m)}_{\mathbb R}, \k) \subseteq  \k[V] \otimes \k[W]$ is the restriction of the $W$-action on $  \k[V] \otimes \k[W]$ induced by the $N(T)$-action on $\coprod_{w \in W} \{w\}$. This is indeed the left action \eqref{act1}. The second assertion of Proposition \ref{tilProp} then follows from Lemma \ref{bVV}.
\end{proof}
As a consequence of Proposition~\ref{tilProp} and Theorem~\ref{Qfat}, we get
\begin{cor}
\la{cohCh}
For $ m \in \M(W) $ odd, the $T$-equivariant cohomology $H_T^{\rm ev}(\tilde{F}_m(T),\,\k)$ carries
a natural action of the rational Cherednik algebra $ \bH_{\frac{m+1}{2}}(W) $ inducing 
the action of the spherical algebra $ \e \bH_{\frac{m+1}{2}}(W) \e $ on $ H^{\rm ev}_G(F_m(G,T),\,\k)$.
\end{cor}
Corollary~\ref{cohCh} provides a topological realization for the (differential) action of the rational Cherednik algebra $ \bH_k(W) $ on quasi-covariants constructed in \cite{BC11}, which, in turn, descends to a natural action of the spherical subalgebra $ \e\bH_k(W)\e $ on the $G$-equivariant cohomology $ H^{{\rm ev}}_G(F_m(G,T), \k) \cong Q_k(W) $.
\subsection{$W$-action on quasi-flag manifolds}
\la{S6.4}
In this section, we construct an action of $W$ on the spaces $ F_m(G,T) $ and $ F^+_m(G,T)$ that is natural in $m$ and commutes with the given $G$-action: in other words, we enrich the $\M(W)$-diagrams of $G$-spaces \eqref{FmG} and \eqref{F-m} to those of $(G \times W)$-spaces:
$$
F^{(+)}_\ast(G,T):\ \M(W) \to \Top^{G \times W}\,,\ m \mapsto F^{(+)}_m(G,T)\,.
$$
We will focus on the spaces $F_m(G,T)$ but our construction works {\it mutatis mutandis} for $F^+_m(G,T)$ as well. Recall that we defined the spaces $ F_m $ by gluing the products of orbits $ G \times_T \cO^{\times \sigma}_{e(s_{\alpha}, w)} $ via the $ \Cc^{(m)}(\Gamma)$-diagram \eqref{FFGm}. In Lemma~\ref{WstrT}, we put a $W$-structure on these diagrams related to the $W$-action \eqref{Wactm} on the category $\Cc^{(m)}(\Gamma)$. To avoid confusion we point out that
this $W$-structure does {\it not} induce any $W$-action on $F_m\,$ as it gets `factored out' in the process of taking the homotopy colimit \eqref{FmG}. However, apart from \eqref{Wactm}, the category $\Cc^{(m)}(\Gamma)$ carries another $W$-action induced by \eqref{Wact2}. This second
$W$-action commutes with the first, making $\Cc^{(m)}(\Gamma)$ a $(W\times W)$-category: 
explicitly, $W \times W$  acts on $\Cc^{(m)}(\Gamma)$  by
\begin{equation}
\la{WWact}
(w_1, w_2)\cdot v  = w_1 v w_2^{-1}\ ,\quad  
(w_1, w_2)\cdot e_{\sigma}(s_\alpha, v) = e_{\sigma}(w_1 s_\alpha w_1^{-1}, w_1 v w_2^{-1})\,,
\end{equation}
where the first and the second actions correspond to the inclusions
$ W \times 1 \into W \times W $ and $1 \times W \into W \times W $.
The next lemma refines the result of Lemma~\ref{WstrT}.
\begin{lemma} \la{WWstrT}
The functor  \eqref{FFGm} admits a $(W\times W)$-structure  with respect to the action \eqref{WWact} on $\Cc^{(m)}(\Gamma)$ that extends the $W$-structure of Lemma~\ref{WstrT}. Consequently, the spaces $ F_m(G,T)$ carry a natural $W$-action induced by \eqref{Wact2}, the corresponding homotopy quotient being represented by
%
\begin{equation}
\la{FmWW}
F_m(G,T)_{hW} \,\simeq\, \hocolim_{\Cc^{(m)}(\Gamma)_{h(W\times W)}}\,(G \times_T \cF_{m})
\,\simeq\, \hocolim_{\Cc(\Gamma)_{h(W\times W)}}\,(G \times_T \cF^*_{m})\,,
\end{equation}
where the functor $ \cF_m^* $ is defined by \eqref{FmG*}
\end{lemma}
\begin{proof} Recall, the first $W$-structure on  \eqref{FFGm} given in Lemma~\ref{WstrT} is defined by the family of natural transformations $\, \{\,{'\vartheta}^{(m)}_w:\, G \times_T \cF_m \to (G \times_T \cF_m) \circ (w,1)\,\}_{w \in W} \,$, with components ({\it cf} \eqref{eq:tildeFm-ver1},\eqref{eq:tildeFm-edge1})
\begin{equation*}
\la{eq:tildeF0-ver1}
'\vartheta^{(m)}_{w,v}:\ G \times_T \{v\} \to G \times_T \{wv\}\,,\quad
[(g,v)]_T \mapsto [(g \dot{w}^{-1}, wv)]_T\,,
\end{equation*}
\begin{equation*}
\la{eq:tildeF0-edge1}
{'\vartheta}^{(m)}_{w,e}:\ G \times_T \O^{\times \sigma}_e \to G \times_T \O^{\times \sigma}_{(w,1)e}\,,\quad
[(g, (x_i)_{i \in \sigma})]_T \mapsto [(g \dot{w}^{-1}, \, (\dot{w} \cdot x_i)_{i \in \sigma})]_T\,,
\end{equation*}
where $ e \in E_{\Gamma} $, $ (x_i)_{i \in \sigma} \in \O^{\times \sigma}_e $, and 
$(w,1)e$ denotes the action of $ w\in W $ on $e \in E_{\Gamma}$ as in \eqref{Wact1}. The `dot' product $ \dot{w} \cdot x_i $ on $\bT$-orbits  is obtained by restricting the natural (left) $G$-action on $ \bF $. 

Next, recall that the flag manifold $\bF$ carries the natural (free)
$W$-action \eqref{WG/T}, which we denote by $\,  w \ast (gT) := g \dot{w}^{-1} T\,$, 
where $ w = \dot{w}T \in W $. As mentioned in Remark~\ref{2ndWact}, this action restricts to a continuous action on $ \bF^{(1)} $, mapping the $\bT$-orbits to $\bT$-orbits and inducing the $W$-action \eqref{Wact2} on $ \Cc(\Gamma) $. We define the second $W$-structure 
$\,\{\,{''\vartheta}^{(m)}_w:\, G \times_T \cF_m \to (G \times_T \cF_m) \circ (1,w)\,\}_{w \in W}\,$
on \eqref{FFGm} by 
\begin{equation*}
\la{eq:tildeF0-ver2}
''\vartheta^{(m)}_{w,v}:\ G \times_T \{v\} \to G \times_T \{vw^{-1}\}\,,\quad
[(g,v)]_T \mapsto [(g,\, vw^{-1})]_T\,,
\end{equation*}
\begin{equation*}
\la{eq:tildeF0-edge2}
''\vartheta^{(m)}_{w,e}:\ G \times_T \O^{\times \sigma}_e \to G \times_T \O^{\times \sigma}_{(1,w)e}\,,\quad
[(g, (x_i)_{i \in \sigma})]_T \mapsto [(g, \, (w \ast x_i)_{i \in \sigma})]_T\,,
\end{equation*}
where $(1,w)e$ denotes the action of $ w\in W $ on $e \in E_{\Gamma}$ as in \eqref{Wact2}.
This second $W$-structure commutes with the first, and we can compose the two to give
 \eqref{FFGm} the $(W\times W)$-structure  
\begin{equation}
\la{eq:tildeF0-edge12}
\vartheta^{(m)}_{(w_1,w_2)}\, := \, {'\vartheta}^{(m)}_{w_1} \circ {''\vartheta}^{(m)}_{w_2}:\ 
G \times_T \cF_m \to (G \times_T \cF_m) \circ (w_1,w_2)
\end{equation}
with respect to the $(W\times W)$-action \eqref{WWact}.

The equivalence \eqref{FmWW} follows by formal properties of homotopy colimits:
\begin{eqnarray*}
F_m(G,T)_{hW} & = & F_m(G,T)_{h(1\times W)} \, =\, [\,\hocolim_{\Cc^{(m)}(\Gamma)_{h(W \times 1)}} (G \times_T \cF_m)]_{h(1\times W)} \\*[1ex]
  & \simeq & [\,(\hocolim_{\Cc^{(m)}(\Gamma)} (G \times_T \cF_m))_{h(W \times 1)}]_{h(1\times W)} 
\qquad (\mbox{by Lemma~\ref{lem:Thom}})  
  \\*[1ex]
  & \simeq & [\,\hocolim_{\Cc^{(m)}(\Gamma)} (G \times_T \cF_m)]_{h(W \times W)}\\*[1ex]
   & \simeq & \,\hocolim_{\Cc^{(m)}(\Gamma)_{h(W \times W)}} (G \times_T \cF_m)
   \qquad (\mbox{by Lemma~\ref{lem:Thom}})  \\*[1ex]
    & \simeq & \hocolim_{\Cc(\Gamma)_{h(W \times W)}} (G \times_T \cF^*_m) 
    \qquad (\mbox{by Lemma~\ref{CGdec}})\,. 
\end{eqnarray*}
This finishes the proof of the lemma.
\end{proof}
Next, we describe the `double' homotopy quotient category $ \Cc(\Gamma)_{h(W\times W)} $ for the action \eqref{WWact} by exhibiting explicitly its skeleton. The following lemma should be compared to Lemma~\ref{skCG}.
\begin{lemma}
\la{skGT2W}
The category $\Cc(\Gamma)_{h(W\times W)}$ is equivalent to its skeletal subcategory $\overline{\Cc(\Gamma)}_{h(W\times W)}$ spanned by
the objects $ \{1\} \,\cup\,\{e(\bar{\alpha}) = e([s_{\alpha}])\}_{\bar{\alpha} \in \A/W} $, 
where $ [s_{\alpha}] $ denotes the conjugacy class of the reflection $s_{\alpha}$ in $W$. The morphisms in $\overline{\Cc(\Gamma)}_{h(W\times W)}$ are represented by the diagram
\begin{equation}   
\la{skWW}
\begin{tikzcd}[scale cd=0.8]
	&& {e(\bar{\beta})} \\
	{e(\bar{\alpha})} &&&& {e(\bar{\gamma})} \\
	&& 1
	\arrow["{W_{\beta} \times C_{\beta}}", from=1-3, to=1-3, loop, in=55, out=125, distance=10mm]
	\arrow["{S_{\bar{\beta}}}"{description}, shift left, Rightarrow, from=1-3, to=3-3]
	\arrow["{W_{\alpha}} \times C_{\alpha}", from=2-1, to=2-1, loop, in=100, out=170, distance=10mm]
	\arrow["{S_{\bar{\alpha}}}"{description}, shift left, Rightarrow, from=2-1, to=3-3]
	\arrow["{W_{\gamma} \times C_{\gamma}}", from=2-5, to=2-5, loop, in=10, out=80, distance=10mm]
	\arrow["{S_{\bar{\gamma}}}"{description}, Rightarrow, from=2-5, to=3-3]
	\arrow["W"', from=3-3, to=3-3, loop, in=310, out=230, distance=15mm]
\end{tikzcd}
\end{equation}
where $C_{\alpha} := C_W(W_{\alpha}) $ denotes the centralizer of $W_{\alpha} $ in $W$, and
$\, S_{\bar{\alpha}} := \Hom(e(\bar{\alpha}), 1) \cong W \times W_{\alpha} $.
The compositions in \eqref{skWW} are given by the action maps
\begin{equation}
\la{rightcomp}
S_{\bar{\alpha}} \times (W_{\alpha} \times C_{\alpha}) 
\to  S_{\bar{\alpha}}\ ,\quad (w, s_{\alpha}) \circ (g_{\alpha}, c_{\alpha}) = 
(wc_{\alpha}, s_{\alpha}g_{\alpha})\,,
\end{equation}
\begin{equation}
\la{leftcomp}
W \times S_{\bar{\alpha}} \to  S_{\bar{\alpha}}\ ,\quad g \circ (w, s_{\alpha}) = (gw, s_{\alpha})\,,
\end{equation}
where $w,g \in W$, $\,s_{\alpha}, g_{\alpha} \in W_{\alpha} \,$, $\, c_{\alpha} \in C_{\alpha} $.
\end{lemma}
\begin{proof} 
Recall, by definition, the objects in $\Cc(\Gamma)_{h(W\times W)}$ are the same as in $ \Cc(\Gamma)$,
and since $ \Cc(\Gamma)$ is a poset category, the morphisms can be represented by elements of $ W \times W $. Specifically, the non-empty Hom-sets  in $\Cc(\Gamma)_{h(W\times W)}$ are
\begin{eqnarray}
\Hom(w_1, w_2) &=& \{\,(g, \,w_2^{-1} g w_1) \in W \times W\ : \ g \in W\,\} \,, \la{Homww}\\*[1ex]
\Hom(e(s_{\alpha}, w_1),\, w_2) &=& \{\,(gW_{\alpha}, \,w_2^{-1} g w_1) \in W \times W\ : 
\ g \in W\,\}\,,\la{Homew}\\*[1ex]
\Hom(e(s_{\alpha}, w_1),\, e(s_{\beta}, w_2)) &=& 
\{\,(gW_{\alpha}, \,w_2^{-1} g w_1) \in W \times W\ : 
\ g \in N_W(W_{\alpha}, W_{\beta})\,\}\,\la{Homee}\ ,
\end{eqnarray}
where $\,gW_{\alpha} = \{g,\,gs_{\alpha}\} \,$ and $\, N_W(W_{\alpha}, W_{\beta}) :=\{g \in W\ :\ gW_{\alpha}g^{-1} \subseteq W_{\beta} \} \,$. It follows from \eqref{Homww} that every vertex 
object $ w \in W $ is isomorphic to the object $1 \in W $ in $\Cc(\Gamma))_{h(W \times W)} $.
Hence, to get a skeleton we can discard all objects except $1 \in W $ and $ e(s_{\alpha},1) $
for $ \alpha \in \A $. Then,  from \eqref{Homee} we see that $ e(s_{\alpha}, 1) \cong e(s_{\beta}, 1) $ in $\Cc(\Gamma))_{h(W \times W)} $ iff $ s_{\alpha} $ and $ s_{\beta} $ are conjugate  reflections  in $W$. We therefore choose a representative in each conjugacy class $ [s_{\alpha}] $ and denote the corresponding object by $ e(\bar{\alpha}) $. Thus, the object set of $\overline{\Cc(\Gamma)}_{h(W\times W)}$ is the union $ \{1\} \,\cup\, \{ e(\bar{\alpha})\}_{\bar{\alpha} \in \A/W}$. The morphism sets shown in the diagram \eqref{skWW} are then easily computed from  \eqref{Homww}-\eqref{Homee}. These morphism sets are embedded into $ W \times W $ as
\begin{eqnarray}
\nonumber \End(1) &\cong&  W \into W \times W\ , \quad g \mapsto (g,\,g)\,, \\*[1ex]
\la{wwaemb} \End[e(\bar{\alpha})] &\cong& W_{\alpha} \times C_\alpha \into W \times W \ ,\quad 
(g_{\alpha}, c_{\alpha}) \mapsto (g_{\alpha} c_{\alpha}, \, c_{\alpha})\,,\\*[1ex]
\nonumber S_{\bar{\alpha}} := \Hom(e(\bar{\alpha}), 1) &\cong& W \times W_{\alpha} \into W \times W\ ,\quad (w,s_{\alpha}) \mapsto (w s_{\alpha}, w)\,,
\end{eqnarray}
and the maps \eqref{rightcomp} and \eqref{leftcomp} under these embeddings correspond to multiplication in $W \times W$.
\end{proof}
Now, let $ \overline{\Sc}(W) := \Sc(W)/W $ denote the quotient poset of $\Sc(W) $ modulo the  conjugation-action of $W$. Viewed as a category, $\overline{\Sc}(W)$
has the congugacy classes $\,  W_{0} = [{W}_{0}] $ and $
W_{\bar{\alpha}} := [{W}_{\alpha}] $ as its objects, and the arrows $\, W_{0} \to W_{\bar{\alpha}} $, one each $ \bar{\alpha} \in \A/W $, as its (non-identity) morphisms.  
Given a functor $ X: \Cc(\Gamma)_{h(W\times W)} \to \Top $, we define $\,X^{\nn}: \ \overline{\Sc}(W) \to \Top\,$ by
\begin{eqnarray}
\la{Fnn}
X^{\nn}(W_0) &=& X(1)_{hW}\,, \la{Fnn1}\\
X^{\nn}(W_{\bar{\alpha}}) &=& 
\hocolim\,\{X(1)_{hW} \xleftarrow{(f_{\alpha})_*} X[e(s_{\alpha})]_{hC_{\alpha}} \xrightarrow{(i_{\alpha})_*}  X[e(s_{\alpha})]_{h(W_{\alpha} \times C_{\alpha})}\}\,,
\la{Fnn2}
\end{eqnarray}
where the map  $ (i_{\alpha})_* $ is induced on homotopy quotients by the natural inclusion
$i_{\alpha}: C_{\alpha} \into W_{\alpha} \times C_{\alpha} $, and the map $ (f_{\alpha})_* $ is defined as the composition 
$$
EW \times_{C_{\alpha}} X[e(\alpha)] \xrightarrow{\id \times X(f_{\alpha})} EW \times_{C_{\alpha}} X[e(1)] \onto EW \times_{W} X[e(1)]\,
$$ 
The next result should be compared to Theorem~\ref{3hdec}.
\begin{theorem}
\la{WWhdec}
For any $(W\times W)$-diagram  $ X$, there is a natural weak homotopy equivalence
\begin{equation}
\la{WWCGS}
\hocolim_{\Cc(\Gamma)_{h(W\times W)}}(X) \, \simeq \, \hocolim_{\,\overline{\Sc}(W)}(X^{\nn})\,.
\end{equation}
where the $\overline{\Sc}(W)$-diagram $\,X^{\nn}\,$ is defined by \eqref{Fnn1} and \eqref{Fnn2}.
\end{theorem}
\begin{proof} The proof is parallel to that of Theorem~\ref{3hdec}. First, we restrict the given
functor  $ X: \Cc(\Gamma)_{h(W\times W)} \to \Top $ to the skeletal subcategory of $ \Cc(\Gamma)_{h(W\times W)} $ described in Lemma~\ref{skGT2W}. This does not change the 
homotopy colimit, and since $ \overline{\Cc(\Gamma)}_{h(W\times W)} $ is an EIA category, we can apply Corollary~\ref{EIA}. Arguing as in the proof of Theorem~\ref{3hdec}, we then get the functor $\,X^{\nn}: \ \overline{\Sc}(W) \to \Top\,$ in the form: 
$\, X^{\nn}(W_0) = X(1)_{hW} \,$ and 
\begin{equation}
\la{XnnS} 
X^{\nn}(W_{\bar{\alpha}})\,=\,\hocolim\,\{ X(1)_{hW} \xleftarrow{\rm ev} [S_{\bar{\alpha}} \times X(e(\bar{\alpha})]_{h(W \times W_{\alpha} \times C_{\alpha} )} \xrightarrow{\rm pr} X[e(\bar{\alpha})]_{h(W_{\alpha} \times C_{\alpha})}\}\,. 
\end{equation}
We need only to show that \eqref{XnnS} is equivalent to \eqref{Fnn2}. First,
note that the action of $W \times W_{\alpha} \times C_{\alpha} $ on $ S_{\bar{\alpha}} \times X[e(\bar{\alpha})] $ defining the homotopy quotient space in the middle of \eqref{XnnS} is given explicitly by ({\it cf.}  formula \eqref{Autact}):
%
$$
(g, g_{\alpha}, c_{\alpha}) \cdot ((w, w_{\alpha}),\,x) \,=\, (g \circ (w, w_{\alpha}) \circ (g_{\alpha}, c_{\alpha})^{-1}, \,(g_{\alpha}, c_{\alpha}) \cdot x)
\,=\, ((gwc_{\alpha}^{-1}, w_{\alpha} g_{\alpha}^{-1}), \,(g_{\alpha},c_{\alpha}) \cdot x)\,,
$$
%
We see that the subgroup $\, W \subset W \times W_{\alpha} \times C_{\alpha}  \,$ acts on $ S_{\bar{\alpha}} \times X[e(\bar{\alpha})]$
freely on the first factor and trivially on the second factor. Hence,
$$
[S_{\bar{\alpha}} \times X(e(\bar{\alpha})]_{h(W \times W_{\alpha} \times C_{\alpha} )} \,
\simeq\, [(S_{\bar{\alpha}}/W) \times X(e(\bar{\alpha})]_{h(W_{\alpha} \times C_{\alpha} )} \,
\cong\, [(W_{\alpha} \times X(e(\bar{\alpha})]_{h(W_{\alpha} \times C_{\alpha} )}\,,
$$
where the induced action of $ W_{\alpha} \times C_{\alpha}  $ on 
$\,W_{\alpha} \times X(e(\bar{\alpha}))\,$ in the last homotopy quotient
is given by $\,(g_\alpha, c_{\alpha}) \cdot (w_{\alpha}, x) = (w_{\alpha} g_{\alpha}^{-1}, (g_{\alpha},c_{\alpha}) \cdot x)\,$. Again, $\, W_{\alpha} \subset W_{\alpha} \times C_{\alpha}  $ acts freely on the first factor, hence 
$\,[(W_{\alpha} \times X(e(\bar{\alpha}))]_{h(W_{\alpha} \times C_{\alpha} )} \simeq X[e(\bar{\alpha})]_{hC_{\alpha}} $, which brings \eqref{XnnS} to the required form \eqref{Fnn2}.
\end{proof}
We now apply Theorem~\ref{WWhdec} to our generalized orbit functors defining quasi-flag manifolds: specifically, we will use the functor $G \times_T \cF^{\ast}_m:\Cc(\Gamma) \to \Top^G$ defined by \eqref{FmG*}. By Lemma \ref{WWstrT}, this functor admits a $(W \times W)$-structure:
$$ 
\vartheta^{(m)}_{(w_1,w_2)}\,:\,G \times_T \cF_m^\ast \to G \times_T \cF_m^\ast \circ (w_1,w_2)
$$
induced by \eqref{eq:tildeF0-edge12}. The following theorem should be viewed as a generalization to quasi-flag manifolds of a well-known result of A. Borel ({\it cf.} \cite[Theorem 20.3]{Bo67}):
\begin{theorem} 
\la{WFmTh} 
For all $m \in \M(W)$, the $G$-spaces $ F_m(G,T)$ and $ F^+_m(G,T) $ carry a natural $W$-action induced by \eqref{WG/T}. The corresponding homotopy quotients $F_m(G,T)_{hW}$ and $F_m^+(G,T)_{hW}$ are $p$-acyclic for all primes $\, p \not|\ |W|\,$.
Consequently, there are natural mod $p$ cohomology isomorphisms
\begin{equation}
\la{XmWs6}
X_m(G,T)_{hW} \,\xrightarrow{\sim} \,X_m^+(G,T)_{hW} \,\xrightarrow{\sim} \, BG\,,
\end{equation}
where $\, X_m(G,T) := EG \times_G F_m(G,T)\,$ and $\, X^+_m(G,T) := EG \times_G F^+_m(G,T) \,$.
\end{theorem}

The proof of Theorem \ref{WFmTh} is based on the following  lemma. 

\begin{lemma} \la{wwactcoh}
For any field $ \k $ of characteristic $ p \not= 2$, there are natural isomorphisms
\begin{equation}
\la{orbitiso}
H_G^*(G \times_T \O_{e(s_\alpha,w)}^{\ast (m_\alpha+1)},\,\k)\,\cong\,
H_T^*(\O_{e(s_\alpha,w)}^{\ast (m_\alpha+1)},\,\k)\,\cong\,\k[V]/\langle \alpha^{m_\alpha+1}\rangle\,,
\end{equation}
such that the $(W\times W)$-structure maps
$$
\vartheta^{(m)}_{(w_1,w_2)}(e(s_{\alpha},w)): \ 
G \times_T \O_{e(s_\alpha,w)}^{\ast (m_\alpha+1)}\,\to\,
G \times_T \O_{e(s_{w_1(\alpha)}, \,w_1 w w_2^{-1})}^{\ast (m_\alpha+1)}
$$
induce on the $G$-equivariant cohomology the following maps
\begin{equation*}
\k[V]/\langle \alpha^{m_\alpha+1}\rangle \,\to\,
\k[V]/\langle \alpha^{m_\alpha+1}\rangle\,,\quad
(w_1\cdot f)\, {\rm mod}\, \langle (w_1\cdot\alpha)^{m_\alpha+1}\rangle \,\mapsto \, 
f\, {\rm mod}\, \langle \alpha^{m_\alpha+1} \rangle \,.
\end{equation*}
\end{lemma}
\begin{proof} 
The first isomorphism in \eqref{orbitiso} is formal; the second can be constructed as follows. For each $\alpha \in \R_+$, we pick a basepoint  $x_{\alpha} \in \O_{e(s_\alpha,1)}$ as in \eqref{NOTa}, and for each $ e(s_\alpha,w)$, we denote by $\cO_{e(s_\alpha,w)}$ a connected ($T$-equivariant) retract of $\O_{e(s_\alpha,w)}$ defined as in \eqref{realorb}. Fixing a representative  $\dot{w} \in N$ for the element $ w \in W $, we pick the basepoint $\dot{w} \cdot x_{w^{-1}\alpha} \in \cO_{e(s_\alpha,w)}$ to obtain the $T$-equivariant isomorphism 
$$i_{\alpha,w}\,:\,T/T_\alpha \stackrel{\sim}{\to} \cO_{e(s_\alpha,w)}\,,\qquad tT_\alpha \mapsto t\dot{w}x_{w^{-1}\alpha}\ .$$
Since our choice of a basepoint in $\cO_{e(s_\alpha,w)}$ is unique up to left action by an element of $T$, the pullback map $\,i_{\alpha,w}^{\ast}\,$ determines a unique isomorphism on $T$-equivariant cohomology:
\begin{equation}\la{orisom1}
H^\ast_T(\cO_{e(s_\alpha,w)},\,\k) \,\cong \, H^\ast_T(T/T_\alpha,\,\k)\,\cong\,
H^\ast(BT_\alpha,\,\k)\,\cong\, \k[V]/\langle \alpha \rangle\,,
\end{equation}
where the last isomorphism is induced by the canonical map $ BT_{\alpha} \to BT $.
Composing \eqref{orisom1} with the isomorphism induced by the natural inclusion $\cO_{e(s_\alpha,w)} \hookrightarrow \O_{e(s_\alpha,w)}$ yields  
\begin{equation} \la{orisom0} H^\ast_T(\O_{e(s_\alpha,w)},\k) \,\cong\, \k[V]/\langle \alpha \rangle\,, \end{equation}
which is our isomorphism \eqref{orbitiso} for $m_\alpha=0$. Then, representing the iterated joins as
$$
\O^{\ast (m_\alpha+1)}_{e(s_\alpha,w)} = 
\hocolim\{\O_{e(s_\alpha,w)} \leftarrow \O_{e(s_\alpha,w)} \times \O^{\ast\, m_\alpha}_{e(s_\alpha,w)} \to \O^{\ast \,m_\alpha}_{e(s_\alpha,w)}\}
$$
it is easy to construct (by induction on $m_{\alpha}$) the required 
isomorphisms  \eqref{orbitiso} for all $ m_\alpha \ge 0 $:
\begin{equation} \la{orisomal}  H^\ast_T(\O^{\ast (m_\alpha+1)}_{e(s_\alpha,w)},\k) \,\cong\, \k[V]/\langle \alpha^{m_\alpha+1} \rangle\ \end{equation}
This inductive construction shows that the embedding $\cO_{e(s_\alpha,w)} \hookrightarrow 
\cO^{\ast (m_\alpha+1)}_{e(s_\alpha,w)}$ induces the canonical projection $\k[V]/\langle \alpha^{m_\alpha+1}\rangle \onto \k[V]/\langle \alpha \rangle$ on $T$-equivariant cohomologies. 

We then note that there is a commutative diagram of $G$-spaces
\begin{equation*}
\begin{diagram}[small]
G \times_T \O_{e(s_\alpha,w)} &\rInto &   G \times_T \F^{(1)} \\
  \dTo^{\vartheta_{(w_1,w_2)}} & & \dTo_{(w_1,w_2)}\\
  G \times_T \O_{e(s_{w_1\alpha},w_1ww_2^{-1})} & \rInto & G \times_T \F^{(1)}\\
  \end{diagram}
\end{equation*}
where $\F^{(1)}$ is the $1$-skeleton of $\bF $ given in \eqref{bF1}. The horizontal arrows in the above diagram are induced by the inclusions $\O_{e} \hookrightarrow \bF^{(1)}$, where $e=e(s_\alpha,w) $ or $ e(w_1s_\alpha w_1^{-1}, w_1ww_2^{-1})$. 
On (even-dimesnional) $G$-equivariant cohomology, this yields the commutative diagram 
\begin{equation} \la{of1eq}
\begin{diagram}[small]
H^{\rm ev}_G(G \times_T \F^{(1)}, \k)  & \rTo^{(w_1^{-1},w_2^{-1}) \cdot (\mbox{--})} & H^{\rm ev}_G(G \times_T \F^{(1)},\k) \\
\dTo & & \dTo\\
H^{\rm ev}_G(G \times_T \O_{e(s_{w_1\alpha},w_1ww_2^{-1})},\k) & \rTo^{\vartheta^{\ast}_{(w_1,w_2)}} & H^{\rm ev}_G(G \times_T \O_{e(s_\alpha,w)},\k)\\
\end{diagram}\ .
\end{equation}
By Corollary \ref{tgkm1sk}, we have isomorphisms
\begin{eqnarray}  \nonumber H^{\rm ev}_G(G \times_T \F^{(1)},\k) & \cong & H^{\rm ev}_T(\F^{(1)},\k)\\ 
\la{gkm1sk} & \cong &  \{\sum_{w \in W} p_w \otimes w \in \k[V] \otimes \k[W]\,:\,p_{s_\alpha w} \equiv p_w\ {\rm mod}\ \langle \alpha\rangle \,, \forall \,\alpha\in \R_+\,,\,\forall \,w\in W\} \ .\end{eqnarray}
Since the inclusion $G \times_T (G/T)^T = G \times_T \coprod_{w \in W} \{w\} \into G \times_T \F^{(1)}$ is $(W \times W)$-equivariant, the $(W \times W)$-action on $ H^{\rm ev}_G(G \times_T \F^{(1)},\k)$ is the restriction of the $(W \times W)$-action on $H^\ast_G(G \times_T (G/T)^T,\k) \cong \k[V] \otimes \k[W]$. This action is easily seen to be given by the formula\footnote{
When $\k$ has characteristic $0$, this $(W \times W)$-action restricts to the so-called
(`dot' $ \times $ `star')-action on the GKM presentation of $ H^\ast_G(G \times_T G/T,\,\k) = H^\ast_T(G/T,\,\k)$ described in \cite[Sect. 3.1]{Ty08} (see also \cite[Theorem 7]{Ca20}).}
\begin{equation} \la{wwactkvw} (w_1,w_2) \cdot \sum_{w \in W} p_w \otimes w \,=\, \sum_{w \in W} w_1 \cdot p_w \otimes w_1ww_2^{-1}\ . \end{equation}
Under the identifications \eqref{gkm1sk} and \eqref{orisom0}, the map on (even) $G$-equivariant cohomology induced by $G \times_T \O_{e(s_\alpha,w)} \into G \times_T \F^{(1)}$ is the map
$$ \sum_{w' \in W} p_{w'} \otimes w' \mapsto p_w\, {\rm mod}\ \langle\alpha \rangle \ .$$
Therefore,
\begin{equation}  \la{emul1oe}
\vartheta^{\ast}_{(w_1,w_2)}((w_1\cdot f)\ {\rm mod}\ (w_1\cdot\alpha)) \,=\, f\ {\rm mod}\ \alpha\ .\end{equation}
Finally, note that the diagram below commutes:
\begin{equation*}
\begin{diagram}
 G \times_T \O_{e(s_{\alpha},w)} & \rTo &  G \times_T \O^{\ast (m_\alpha+1)}_{e(s_{\alpha},w)}\\
 \dTo^{\vartheta_{(w_1,w_2)}} & & \dTo_{\vartheta^{(m)}_{(w_1,w_2)}}\\
 G \times_T \O_{e(s_{w_1\alpha},w_1ww_2^{-1})} & \rTo &  G \times_T \O^{\ast (m_\alpha+1)}_{e(s_{w_1\alpha},w_1ww_2^{-1})}\\
\end{diagram}
\end{equation*}
Since the map on $G$-equivariant cohomology induced by $G \times_T \O_{e(s_\alpha,w)} \to G \times_T \O^{\ast (m_\alpha+1)}_{e(s_\alpha,w)}$ is the canonical projection $\k[V]/ \langle\alpha^{m_\alpha+1}\rangle \onto \k[V]/ \langle \alpha \rangle $, the last claim of Lemma \ref{wwactcoh} follows from \eqref{emul1oe}.
\end{proof}
\begin{proof}[Proof of Theorem \ref{WFmTh}]
By Corollary \ref{corpplus}, the ($W$-equivariant) maps $F_m(G,T) \to F_m^+(G,T)$ are mod $p$ cohomology isomorphisms for all $p \neq 2$. It therefore suffices to prove Theorem \ref{WFmTh} for the spaces $F_m(G,T)$. For these spaces, the first assertion of Theorem \ref{WFmTh} follows immediately from Lemma \ref{WWstrT}. Indeed, as clarified by the proof thereof, the natural $W$-action on $F_m(G,T)$ induced by \eqref{Wact2} coincides with that induced by \eqref{WG/T}. Next, we claim that for $X=G \times_T \cF_m:\Cc(\Gamma)_{h(W\times W)} \to \Top$, the map (see \eqref{Fnn2})
\begin{equation} \la{ial} (i_{\alpha})_\ast \,:\, X[e(s_\alpha)]_{hC_{\alpha}} \to X[e(s_\alpha)]_{h(W_\alpha \times C_{\alpha})} \end{equation}
is a mod $p$ cohomology isomorphism for any prime $p \not|\ |W|$. Our claim implies that for any prime $p \not|\ |W|$, there are natural mod $p$ cohomology isomorphisms
$$X^{\nn}(W_{\bar{\alpha}}) \stackrel{\sim}{\to} X(1)_{hW}$$
for all $\alpha \in \A/W$. By Theorem \ref{WWhdec}, it follows that for $p \not|\ |W|$,
$$ F_m(G,T)_{hW} \simeq \hocolim_{\Cc(\Gamma)_{h(W\times W)}} X \simeq_p \hocolim_{\bar{\Sc}(W)}(\underline{X(1)}_{hW}) \simeq X(1)_{hW} = (G/T)_{hW} \simeq G/N\ .$$
Here, $\simeq_p$ stands for mod $p$ cohomology isomorphism and $\underline{\mbox{--}}$ stands for a constant functor. Since $G/N \simeq_p {\rm pt}$ for $p \not|\ |W|$, Theorem \ref{WFmTh} follows for the $F_m(G,T), m \in \M(W)$. It remains to prove our claim. Since $BG$ is simply connected, and $(i_\alpha)_\ast$ (see \eqref{ial}) is a morphism of $G$-spaces, the (natural in $Y$) Eilenberg-Moore spectral sequence 
$$E_2^{\ast,\ast}={\rm Tor}^{\ast,\ast}_{H^\ast(BG)}(H^*_G(Y), \k) \implies H^\ast(Y)$$
corresponding to the Borel fibration sequence $Y \to Y_{hG} \to BG$ implies that it is enough to show that $[(i_\alpha)_\ast]_{hG}$ is a mod $p$ cohomology isomorphism for any prime $p \not|\ |W|$. Since the $G$-action on $X[e(s_\alpha)]$ commutes with the $W_\alpha \times C_\alpha$-action, for $H=C_{\alpha}, W_\alpha \times C_\alpha$ we have 
$$ (X[e(s_\alpha)]_{hH})_{hG}\,\cong \, (X[e(s_\alpha)]_{hG})_{hH}\ . $$
Since $|H| \,|\,|W|^2$, the Serre spectral sequence associated with the fibration sequence
$X[e(s_\alpha)]_{hG} \to (X[e(s_\alpha)]_{hG})_{hH} \to BH$ collapses for $\k = \F_p\,,\, p\ \not|\ |W|$ to yield an isomorphism
$$ H^\ast(X[e(s_\alpha)]_{hG}, \k)^{H} \,\cong\, H^\ast((X[e(s_\alpha)]_{hG})_{hH}, \k)\ .$$
Hence, on cohomologies for $\k = \F_p\,,\, p\ \not|\ |W|$, $[(i_\alpha)_\ast]_{hG}$ induces the inclusion 
\begin{equation} \la{cohinc}  H^\ast_G(X[e(s_\alpha)], \k)^{W_\alpha \times C_{\alpha}} \cong \big(\k[V]/\langle\alpha^{m_\alpha+1}\rangle\big)^{W_\alpha \times C_{\alpha}} \hookrightarrow \big(\k[V]/\langle\alpha^{m_\alpha+1}\rangle\big)^{C_{\alpha}} \,\cong\,  H^\ast_G(X[e(s_\alpha)], \k)^{C_{\alpha}}\ .\end{equation}
By Lemma \ref{wwactcoh}, for $(w_1,w_2) \in \End[e(\bar{\alpha})] \subset W \times W$, $(w_1,w_2)$ acts on $H^\ast_G(X[e(s_\alpha)], \k) \cong \k[V]/\langle\alpha^{m_\alpha+1}\rangle$ by
$$(w_1,w_2) \cdot (f\ {\rm mod}\ \alpha^{m_\alpha+1}) \,=\,(w_1\cdot f) \ {\rm mod} \ \alpha^{m_\alpha+1}\ .$$
It follows from \eqref{wwaemb} that 
$$\big(\k[V]/\langle\alpha^{m_\alpha+1}\rangle\big)^{W_\alpha \times C_{\alpha}} \,=\, \big(\k[V]/\langle \alpha^{m_\alpha+1}\rangle\big)^{W_\alpha C_{\alpha}}\,=\, \big(\k[V]/\langle \alpha^{m_\alpha+1}\rangle\big)^{C_{\alpha}}\,,$$
where the last equality is because $W_\alpha \subseteq C_\alpha$. Hence, the map \eqref{cohinc} is an isomorphism for $\k=\F_p\,,\,p \ \not|\ |W|$. This completes the proof of our claim about $(i_\alpha)_\ast$, in turn completing the verification of Theorem \ref{WFmTh}.
\end{proof}
As a consequence of Theorem \ref{WFmTh}, we obtain
\begin{cor} \la{Coroddcoh}
Let $ \k = \Q$ or $ \k=\F_p $ with $ p\ \not|\ |W|$. Then, for all $m \in \M(W)$, 
$$
H^\ast_G(F_m^+(G,T), \k)^W \,\cong\,H^\ast_G(F_m(G,T), \k)^W \,\cong\,H^\ast(BG,\k) \,\cong\, \k[V]^W\,.
$$
In particular,  we have
$\,H^{\rm odd}_G(F_m^+(G,T),\k)^W \cong H^{\rm odd}_G(F_m(G,T),\k)^W = 0\,$.
\end{cor}
\begin{proof}
Associated to the Borel fibration  $ X_m(G,T) \to X_m(G,T)_{hW} \to BW $ there is the Serre spectral sequence 
$$
E_2^{pq} = H^{p}(W,H^q(X_m(G,T),\k)) \,\Rightarrow\, H^\ast(X_m(G,T)_{hW},\k)\,.
$$
If $ \k = \Q$ or $ \k=\F_p $ with $ p \not| \ |W| $, this spectral sequence
collapses, giving an isomorphism $\,H^\ast(X_m(G,T),\k)^W \cong H^\ast(X_m(G,T)_{hW},\k)\,$. The claim now follows from
 \eqref{XmWs6}.
\end{proof}

\section{Equivariant \texorpdfstring{$K$}{K}-theory}
\la{S7}
In this section, we first define exponential analogs of rings of quasi-invariants and rings of quasi-covariants and prove results analogous to those of Sections \ref{S2.3} and \ref{S2.4} for these rings (in Section \ref{S7.1}). Next, in Section \ref{S7.3}, we compute the $G$-equivariant and $T$-equivariant $K$-theory of quasi-flag manifolds in terms of the rings of quasi-invariants and quasi-covariants defined in Section \ref{S7.1}. These results provide a topological motivation for the results in Section \ref{S7.1}.

\subsection{Exponential quasi-invariants}
\la{S7.1}
Multiplicative analogues of the rings of quasi-invariants $ Q_k(W) $  and their $q$-generalizations were introduced and studied in the work of O. Chalykh \cite{Ch00, Ch02}. 
In \cite[Section 5.2]{BR1}, motivated by $K$-theoretic computations (see Section~\ref{S7.3} below), we proposed a variant of these rings that we called the ring $ \cQ_k(W)$ of {\it exponential quasi-invariants} of a Weyl group $W$. In this section, we review our definition of exponential quasi-invariants and introduce the related notion of exponential quasi-covariants of $W$.

Let $G$ be a compact connected Lie group with maximal torus $T$ and associated Weyl group $W$. Let $ \hat{T} := \Hom(T, U(1)) $ denote the character lattice and $R(T)$ the representation ring of $T$. It is well known that $ R(T) \cong \Z[\hat{T}] $ via the canonical map induced by taking characters of representations, and $ R(T)^W \cong R(G) $ via the restriction map $ i^*: R(G) \to R(T) $ induced by the inclusion $ i: T \into G $ (see, e.g., \cite[Chap. IX, Sect. 3]{Bou82}). Using the first isomorphism, we identify $ R(T) = \Z[\hat{T}] $ and write $ e^{\lambda} $ for the elements of $ R(T) $ corresponding to characters $ \lambda \in \hat{T} $. Next, we let $ \R \subseteq \hat{T} $ denote the root system of $W$ determined by $(G,T)$ and choose a subset $ \R_+ \subset \R $ of positive roots in $ \R $. If $ s_\alpha \in W $ is the reflection in $W$ corresponding to $ \alpha \in \R_+ $, then the difference $ e^{\lambda} - e^{s_{\alpha}(\lambda)} $ in $ R(T)$ is uniquely divisible by $ 1 - e^{\alpha} $ for any $\lambda \in \hat{T}$. Following \cite{D74},  we define a linear endomorphism $\, \Lambda_{\alpha}: R(T) \to R(T)\,$ for each $\alpha \in \R_+ $, such that
\begin{equation}
\la{extdem}
(1 - s_\alpha) f = \Lambda_{\alpha}(f) \cdot (1 - e^{\alpha})\ .
\end{equation}
The operator $ \Lambda_\alpha $ is an exponential analogue of the divided difference operator. Note that the conditions \eqref{In9} defining the usual quasi-invariant polynomials can be written in terms of the divided difference operators as $ \nabla_\alpha(p) \equiv 0 \ \mbox{mod}\, (\alpha)^{2 k_\alpha} $.
This motivates the following definition of quasi-invariants in the exponential case. Throughout this section, we use the notation
$$
\bar\alpha := \begin{cases} \frac{\alpha}{2} & \text{if}\ \, \frac{1}{2}\alpha \in \hat{T}\\*[1ex] \alpha & \text{otherwise} \end{cases}
$$
\begin{defi}
\la{Expqi}
An element $ f \in R(T)$ is called an {\it exponential quasi-invariant of $W$ of multiplicity} $ k \in \M(W) $ if
\begin{equation}
\la{EQinv}
\Lambda_{\alpha}(f) \,\equiv\, 0 \ \,\mbox{mod}\ (1 - e^{\bar{\alpha}})^{2 k_\alpha} \ ,\quad 
\forall\,  \alpha \in \R_+\ .
\end{equation}
We write  $\cQ_k(W)$ (resp., $\cQ_{k}(W_\alpha)$) for the set of all  $ f \in R(T) $ satisfying \eqref{Expqi} for a fixed multiplicity $k$ (resp., satisfying the congruence \eqref{EQinv} for a {\it fixed} $\alpha \in \R_+$, and a fixed nonnegative integer $k$). 
\end{defi}
It is easy to see that $\cQ_k(W)$ and $\cQ_k(W_{\alpha})$ are closed under addition and multiplication in $R(T)$, i.e. $\cQ_k(W)$ is a commutative subring of $ R(T)$. 
(That $\cQ_k(W)$ and $\cQ_k(W_{\alpha})$ are closed under multiplication follows from the twisted derivation property of Demazure operators 
$\, 
\Lambda_{\alpha}(f_1 f_2) = \Lambda_\alpha(f_1)\cdot f_2 + s_{\alpha}(f_1) \cdot \Lambda_{\alpha}(f_2) 
\,$ 
that holds for all $  \alpha \in \R $, see \cite[Section 5.5]{D74}.) 
\begin{remark}
\la{R202}
As mentioned above, the exponential analogues of the rings of polynomial quasi-invariants were defined by O. Chalykh in his work on bispectral problem for quantum integrable systems of 
Ruijsenaars type \cite{Ch00} and its application to Macdonald conjectures \cite{Ch02}.
Specifically, if we ignore the difference between $\alpha$ and $\bar{\alpha}$, our formula \eqref{EQinv} appears as the limit $ q \to 1 $ of Equation (2.15) defining $q$-deformations of 
exponential quasi-invariants in \cite{Ch02}. In a different form, this definition also appears 
in \cite{Ch00} (see Equation (45) and Theorem 6.4).
\end{remark}
\begin{remark}
\la{R22}
In \cite[Section 5.2]{BR1}, $f \in R(T)$ was called an exponential quasi-invariant of multiplicity $k$ if 
\begin{equation}
\la{EQinv1}
\Lambda_{\alpha}(f) \,\equiv\, 0 \ \,\mbox{mod}\ (1 - e^{\frac{\alpha}{2}})^{2 k_\alpha} \ ,\quad 
\forall\,  \alpha \in \R_+\,,
\end{equation}
In general, it may happen that $ \frac{\alpha}{2} \not\in \hat{T} $ for some $ \alpha \in \R_+ $, so that 
$\, e^{\frac{\alpha}{2}} \not\in R(T) $. In this case, we view
\eqref{EQinv1} as a congruence in the extended ring $ R(T)[e^{\frac{\alpha}{2}}] \subset \Z[\frac{1}{2} \hat{T}] $ that naturally contains $ R(T) $. Note that  there is an involution $\sigma$ of $R(T)[e^{\frac{\alpha}{2}}]$ fixing $R(T)$ and mapping $e^{\frac{\alpha}{2}}$ to $-e^{\frac{\alpha}{2}}$. If $g \in R(T)$, then $g \equiv 0\ \, \mbox{mod}\ (1 - e^{\frac{\alpha}{2}})^{2 k_\alpha}$ implies that $g=\sigma(g) \equiv 0\ \,\mbox{mod}\ (1 + e^{\frac{\alpha}{2}})^{2 k_\alpha}$ which, in turn, implies that $g \equiv 0\ \,\mbox{mod}\ (1 - e^{\alpha})^{2 k_\alpha}$. This shows that Definition \ref{Expqi} above is equivalent to the one proposed in \cite{BR1}.  
\end{remark}

Next, we define two exponential analogues of the ring $\bQ_m(W)$ of quasi-covariants of $W$.
\begin{defi} 
\la{expqmring}
For $m \in \M(W)$, we let
\begin{equation} 
\la{bcqm2} \bcQ_m(W):=\bigg\{ \sum_{w \in W} p_w \otimes w \in R(T) \otimes \Z[W]\,\,:\,\, p_{s_{\alpha}w} \equiv p_w\ \mod\, \langle(1 - e^{\bar{\alpha}})^{m_\alpha}\rangle \, ,\
\forall\,  \alpha \in \R_+ \bigg\}
\end{equation}
and
\begin{equation} 
\la{bcqm1}
\bcQ'_m(W)\,:=\,\bigg\{\sum_{w \in W} p_w \otimes w \in R(T) \otimes \Z[W]\,\,:\,\, p_{s_{\alpha}w} \equiv p_w \ \mod\, \langle(1-e^{\alpha})(1 - e^{\bar{\alpha}})^{m_\alpha-1}\rangle  \, ,\
\forall\,  \alpha \in \R_+ \bigg\}
\end{equation}
\end{defi}
We note that both $ \bcQ_m(W) $ and $ \bcQ_m'(W) $ are commutative rings, 
$ \bcQ_m'(W) $ being a subring of $ \bcQ_m(W) $. In
fact, it is easy to check that
\begin{equation*} 
\bcQ'_m(W) = \bcQ_m(W) \cap \bcQ_1'(W) = \bigg\{ \sum_{w \in W} p_w \otimes w \in \bcQ_m(W)\ :\ p_{s_\alpha w} \equiv p_w \ \mod\,\langle(1-e^{{\alpha}}) \rangle \, ,\ \forall\,  \alpha \in \R_+ \bigg\}
\end{equation*}
for all $ m \in \M(W)$.
Moreover,  $ \bcQ_m'(W) $ contains $ R(T) $ as a subring via the canonical ring homomorphism
\begin{equation} 
\la{rtmodstr} 
R(T) \hookrightarrow \bcQ'_{m}(W)\,,\qquad f \mapsto \sum\, f \otimes w \,,
\end{equation}
Thus, for all $ m\in\M(W) $, we have the inclusions
\begin{equation} 
\la{rtmodstr1}
R(T) \,\subseteq \, \bcQ'_{m}(W) \,\subseteq\, \bcQ_m(W)\,\subseteq\,R(T) \otimes \Z[W]
\,.
\end{equation}

Observe that the ring $R(T) \otimes \Z[W]$, viewed as the ring of $R(T)$-valued functions on $W$, carries the following three $W$-actions:
 \begin{eqnarray}
   g\cdot (f \otimes w)& = & g(f) \otimes gw\,, \la{expact1}\\ 
 g\cdot (f \otimes w) &=& f \otimes gw \,, \la{expact2}\\
g\cdot (f \otimes w) &=& f \otimes wg^{-1} \,, \la{expact3}
 \end{eqnarray}
  which are analogous to the corresponding actions of $W$ on $\k[V] \otimes \k[W]$ (see Section \ref{S2.3}) and which we denote by $g \otimes g, 1 \otimes g$ (acting on the left), and $1 \otimes g^{-1}$ (acting on the right), respectively. It is not difficult to check that, for all $m \in \M(W)$, the rings $\bcQ_m(Q)$ and $\bcQ'_m(W)$ are stable under the diagonal action \eqref{expact1} and the right action \eqref{expact3}. In particular, the inclusions \eqref{rtmodstr1} are equivariant under the diagonal action of $W$.
  In this context, we have the following analogue of Lemma \ref{bVV}:
\begin{lemma}
\la{QmulHalf}
 The natural inclusions \eqref{rtmodstr1} restrict to  
 $$
 R(T)^W \,\subseteq \, \bcQ'_{m}(W)^W \,=\, \bcQ_m(W)^W\,\subseteq\,R(T)\,,
 $$
 where the $W$-invariants are taken with respect to the diagonal action \eqref{expact1}.
 Moreover, there are $W$-equivariant isomorphisms of commutative rings
 \begin{equation} \la{bcqwinv} 
 \bcQ'_m(W)^{W} \,=\, \bcQ_m(W)^W\,\cong\, \cQ_{[\frac{m}{2}]}(W)\,,
 \end{equation}
 where the $W$-structures are induced by the right action \eqref{expact3}. 
\end{lemma}
\begin{proof}
 Note that there is an algebra homomorphism 
 $$R(T) \to R(T) \otimes \Z[W]\,,\qquad f \mapsto \sum_{w \in W} (w\cdot f) \otimes w\,,$$
 whose image is the subring of $W$-invariants in $R(T) \otimes \Z[W]$ (with respect to the diagonal action). Using \eqref{bcqm1} and \eqref{bcqm2}, the invariant subrings $\bcQ_m(W)^W\,=\,\bcQ_m(W) \cap R(T)^W$ and $\bcQ'_m(W)^W\,=\,\bcQ'_m(W) \cap R(T)^W$ can then be described as
 \begin{eqnarray*}
 \bcQ_m(W)^W &\cong & \left\{f \in R(T)\,:\, w(f)-s_{\alpha}w(f)\,\in\, \langle(1 - e^{\bar{\alpha}})^{m_\alpha}\rangle \,,\,
\forall w \in W\,,\,\forall  \alpha \in \A\right\}\\
 \bcQ'_m(W)^W &\cong & \left\{f \in R(T)\,:\, w(f)-s_{\alpha}w(f)\,\in\, \langle(1-e^{\alpha})(1 - e^{\bar{\alpha}})^{m_\alpha-1}\rangle \, ,
\,\forall w \in W\,,\,\forall  \alpha \in \A\right\}\ .
  \end{eqnarray*}
Since $m$ is $W$-invariant, we can simplify the above conditions on $f$ to obtain
 \begin{eqnarray*}
 \bcQ_m(W)^W &\cong & \left\{f \in R(T)\,:\, f-s_{\alpha}(f)\,\in\, \langle(1 - e^{\bar{\alpha}})^{m_\alpha}\rangle \,,\,\forall  \alpha \in \A\right\}\\
 \bcQ'_m(W)^W &\cong & \left\{f \in R(T)\,:\, f-s_{\alpha}(f)\,\in\,  \langle(1-e^{\alpha})(1 - e^{\bar{\alpha}})^{m_\alpha-1}\rangle\,,\,\forall  \alpha \in \A\right\}\ .
  \end{eqnarray*}
In order to prove the desired lemma, it therefore suffices to show that for $f \in R(T)$, and for any $\alpha \in \A$, and for any non-negatve integer $m$, the following congruence conditions are equivalent:
\begin{eqnarray}
f-s_\alpha(f) & \in &\langle(1-e^{\alpha})(1 - e^{\bar{\alpha}})^{m-1}\rangle \la{congc1} \\
 f-s_\alpha(f) &\in & \langle(1 - e^{\bar{\alpha}})^{m}\rangle \la{congc2}\\
\Lambda_{\alpha}(f)  &\in & \langle(1 - e^{\bar{\alpha}})^{2[\frac{m}{2}]}\rangle\la{congc3}\ .\end{eqnarray}
The \eqref{congc1} and \eqref{congc2} are identical if $\frac{\alpha}{2} \not\in \hat{T}$. If $\frac{\alpha}{2} \in \hat{T}$, note that \eqref{congc1} clearly implies \eqref{congc2}. Conversely, if $(1-e^{\frac{\alpha}{2}})^m\,|\,f-s_\alpha(f)$, then since $(1-e^{\alpha})\,|\,f-s_\alpha(f)$, $(1+e^{\frac{\alpha}{2}})\,|\,f-s_\alpha(f)$. Hence, 
$(1-e^{\alpha})(1 - e^{\bar{\alpha}})^{m-1}=(1+e^{\frac{\alpha}{2}})(1-e^{\frac{\alpha}{2}})^m\,|\,f-s_\alpha(f)$, showing that \eqref{congc2} implies \eqref{congc1}. If $m$ is odd, \eqref{congc3} is clearly equivalent to \eqref{congc1}. It remains to verify the same for $m$ even. 

For this, note that the exact sequence of torii $1 \to T_{\alpha} \to T \to T/T_{\alpha} \to 1$ splits (albeit, non-canonically), yielding an isomorphism of representation rings 
\begin{equation} \la{rtsp} R(T) \,\cong \,R(T_{\alpha}) \otimes R(T/T_{\alpha}) \,\cong\, R(T_{\alpha}) \otimes \Z[e^{\pm \bar{\alpha}}] \,,\end{equation}
with $s_{\alpha}$ acting trivially on the first factor. While the isomorphism \eqref{rtsp} is not canonical, the inclusion of rings $\Z[e^{\pm \bar{\alpha}}] \hookrightarrow R(T),\,f \mapsto 1 \otimes f$ is indeed canonical. Hence, to complete the proof of the desired lemma, it suffices to show that \eqref{congc3} is equivalent to \eqref{congc2} for $f \in \Z[e^{\pm \bar{\alpha}}]$ for $m$ odd. This follows immediately from the fact (see, e.g., \cite[Example 5.5]{BR1}) that for $f \in \Z[e^{\pm \bar{\alpha}}]$,
$$\Lambda_{\alpha}(f) \,\equiv\, 0\, \mbox{mod}\ (1 - e^{\bar{\alpha}})^{2k} \,\iff\, f \in \Z \oplus \Z \cdot (e^{\frac{\bar{\alpha}}{2}}-e^{-\frac{\bar{\alpha}}{2}})^2 \oplus \ldots \oplus \Z \cdot (e^{\frac{\bar{\alpha}}{2}}-e^{-\frac{\bar{\alpha}}{2}})^{2k-2} \oplus \Z \cdot (e^{\frac{\bar{\alpha}}{2}}-e^{-\frac{\bar{\alpha}}{2}})^{2k}\Z[e^{\pm \bar{\alpha}}]\, ,$$
completing the proof of the desired lemma.
\end{proof}
Next, we have the following analogue of Lemma \ref{fQGKM}, whose proof is a simple modification (which we leave to the reader) of the proof thereof. 
For all $m \in \M(W)$, we have 
\begin{align}
 \bcQ_m(W) &= \big\{p \in R(T) \otimes \Z[W]\,\,|\,\, (1 \otimes (1-s_\alpha))p\,\in \langle(1 - e^{\bar{\alpha}})^{m_\alpha}\rangle \otimes \Z[W] \, ,
\forall\,  \alpha \in \A \big\} \la{fQexpGKM}\\
 \bcQ'_m(W) &=\big\{p \in R(T) \otimes \Z[W]\,\,|\,\, (1 \otimes (1-s_\alpha))p \in \langle(1-e^{\alpha})(1 - e^{\bar{\alpha}})^{m_\alpha-1}\rangle \otimes \Z[W] \, ,
\forall\,  \alpha \in \A \big\}\ . \la{fQexpGKM2}
\end{align}
Consider the homomorphism of commutative rings
\begin{equation} \la{phimul} \varphi: R(T) \otimes_{R(T)^W} R(T) \to R(T) \otimes \Z[W]\,,\qquad f \otimes g \mapsto \sum_{w \in W} fw(g) \otimes w\ .\end{equation}
We then have the following exponential analog of Proposition \ref{bQodd}. Its proof is, however, a little more subtle than that of Proposition \ref{bQodd}.
\begin{prop} \la{bcQodd}
For all $k \in \M(W)$, the map \eqref{phimul} restricts to  a ring isomorphism 
\begin{equation} \la{phimulk} \varphi_k\,:\,R(T) \otimes_{R(T)^W} \cQ_k(W) \stackrel{\sim}{\to} \bcQ'_{2k+1}(W)\ . \end{equation}
\end{prop}
\begin{proof}
 The map \eqref{phimul} being injective with image $\bcQ'_1(W)$ is a well known result on the $T$-equivariant $K$-theory $K^\ast_T(G/T)$ of the flag manifold (see \cite{HHH05,KK87}). Hence, for $k=0$, the map \eqref{phimulk} is an isomorphism, whose inverse we denote by 
 $$ \psi:= \varphi_0^{-1}\,:\,\bcQ'_1(W) \to R(T) \otimes_{R(T)^W} R(T)\ .$$
 As in the proof of Proposition \ref{bQodd} freeness of $R(T)$ as an $R(T)^W \cong R(G)$-module (see \cite{P72,St75}) implies that $\varphi$ restricts to injective maps $\varphi_k:R(T) \otimes_{R(T)^W} \cQ_k(W) \to \bcQ'_1(W)$. The verification that ${\rm Im}(\varphi_k) \subset \bcQ'_{2k+1}(W)$, which we leave to the reader, is similar to the corresponding verification in the proof of Proposition \ref{bQodd} (the only extra ingredient being the equivalence of conditions \eqref{congc1} and \eqref{congc3}). 
 
 The proof of the opposite inclusion, which we present below, is similar to the corresponding part of the proof of Proposition \ref{bQodd}. First, we note that for all $s \in W$, $f,g \in R(T)$
 $$ \varphi[(1 \otimes s)\cdot (f \otimes g)] = \varphi(f \otimes g)\cdot (1 \otimes s^{-1})\ .$$
 Hence, 
 \begin{equation} \la{2.54mul} (1 \otimes (1-s_{\alpha})) \cdot \psi(\bcQ'_{2k+1}(W)) \,=\, \psi[\bcQ'_{2k+1}(W) \cdot (1 \otimes (1-s_{\alpha}))]\,,\qquad \forall\,\alpha \in \A\ . \end{equation}
 Next, we note that for all $p \in \bcQ'_{2k+1}(W)$,
 \begin{equation} \la{2.55mul} p \cdot (1 \otimes (1-s_{\alpha}))\,=\, \varphi(1 \otimes \lambda_{k_\alpha}(e^{\bar{\alpha}})) \cdot \sum_{w \in W} (q_w \otimes w)\,,\end{equation}
 where $\lambda_k(e^{\bar{\alpha}}):=(1-e^{\alpha})(1-e^{\bar{\alpha}})^{2k}$ and the `dot' product on the right denotes the commutative multiplication in $R(T) \otimes \Z[W]$. We then claim that
 $$\sum_{w \in W} q_w \otimes w  \in \bcQ'_1(W)\ .$$
 Indeed, for any root $\beta$, we have 
 \begin{equation} \la{2.57mul} \lambda_{k_\alpha}(e^{\overline{w\alpha}}) \lambda_{k_\alpha}(e^{\overline{s_{\beta}w\alpha}})(q_w -q_{s_\beta w})\,=\,  \lambda_{k_\alpha}(e^{\overline{s_{\beta}w\alpha}})(p_w -p_{ws_{\alpha}})-  \lambda_{k_\alpha}(e^{\overline{w\alpha}})(p_{s_{\beta}w}-p_{s_{\beta}ws_{\alpha}})\ .\end{equation}
 As in the proof of Proposition \ref{bQodd}, it is easy to check that the right hand side of \eqref{2.57mul} is divisible by $1-e^{\beta}$ for all $\beta \in \A$. If $\beta \neq \pm w\alpha$, then $1-e^{\beta}$ has no common irreducible factor with the polynomials $\lambda_{k_\alpha}(e^{\overline{w\alpha}})$ (whose irreducible factors are $1-e^{\overline{w\alpha}}$ and, if $\frac{\alpha}{2} \in \hat{T}$, $1+e^{\overline{w\alpha}}$) and  $\lambda_{k_\alpha}(e^{\overline{s_{\beta}w\alpha}})$. Hence, if $\beta \neq \pm w\alpha$, $(1-e^{\beta})|(q_w-q_{s_{\beta}w})$. On the other hand, if $\beta=\pm w\alpha$, the right hand side of \eqref{2.57mul} does {\it not} vanish, unlike in the proof of Proposition \ref{bQodd}. In this case, we note that $s_{\beta}w\alpha=-w\alpha$ and $s_{\beta}w=ws_{\alpha}$. By \eqref{2.57mul}, we have 
 \begin{eqnarray*}
 q_w-q_{s_\beta w} &=& \frac{(\lambda_{k_\beta}(e^{\mp \bar{\beta}})+\lambda_{k_\beta}(e^{\pm \bar{\beta}}))(p_w-p_{s_\beta w})}{\lambda_{k_\beta}(e^{\mp \bar{\beta}})\lambda_{k_\beta}(e^{\pm \bar{\beta}})}\\*[2ex]
 &=& \left\{\begin{array}{lr}
      \frac{(1-e^{\pm (k_\beta+1)\beta})\lambda_{k_\beta}(p_w-p_{s_\beta w})}{\lambda_{k_\beta}(e^{\pm \bar{\beta}})} & \text{if } \frac{\beta}{2} \in \hat{T}\\*[2ex]
      \frac{(1-e^{\pm (2k_\beta+1)\beta})\lambda_{k_\beta}(p_w-p_{s_\beta w})}{\lambda_{k_\beta}(e^{\pm \bar{\beta}})} & \text{if } \frac{\beta}{2} \not\in \hat{T}
    \end{array}\right.
 \end{eqnarray*}
 Since $\sum_w p_w \otimes w \in \bcQ'_{2k+1}(W)$, $\lambda_{k_\beta}(e^{\pm \bar{\beta}}) | (p_w-p_{s_\beta}w)$. Since $(1-e^{\beta})|(1-e^{l\beta})$ for any integer $l$, $(1-e^{\beta})|(q_w-q_{s_\beta w})$. Hence, $\sum_{w \in W} q_w \otimes w \in \bcQ'_1(W)$. As in the proof of Proposition \ref{bQodd}, it then follows from \eqref{2.54mul} and \eqref{2.55mul} that 
 $$(1 \otimes (1-s_\alpha)) \cdot \psi(\bcQ'_{2k+1}(W)) \subseteq R(T) \otimes_{R(T)^W} \langle \lambda_{k_\alpha}(e^{\bar{\alpha}}) \rangle\,,\,\,\forall\,\alpha \in \A\ .$$
 From this, and freeness of $R(T)$ as a $R(T)^W=R(G)$-module, it follows that $\bcQ'_{2k+1}(W) \subseteq {\rm Im}(\varphi_k)$, completing the proof of the desired proposition.
\end{proof}

\subsubsection{Actions of the Hecke algebra $\D(W)$ and the trigonometric DAHA}
\la{S7.1.1}
The Hecke algebra $\D(W)$ associated with $R(G)$ (see \cite[Sec. 2]{HLS10}) is the subalgebra of $\End_{R(G)}(R(T))$ generated by the elements of $R(T)$ and the Demazure operators $\delta_{\alpha}:=\frac{1-e^{-\alpha}s_{\alpha}}{1-e^{-\alpha}}, \alpha \in \R_+$. Note that the Demazure operators are related to the divided difference operators $\Lambda_\alpha$ in \eqref{extdem} by the formula $\delta_\alpha=1-\Lambda_\alpha$. Let 
$$ R(T)_{\rm reg}:= R(T)[(1-e^{\alpha})^{-1}\,:\,\alpha \in \R_+]\ .$$
Since $R(T)$ is an integral domain, $R(T) \hookrightarrow R(T)_{\rm reg}$. To state an analog of Theorem \ref{Qfatodd}, we consider the action of $\D(W)$ on $R(T)_{\rm reg} \otimes \Z[W]$ obtained by restricting the natural left multiplication action of $R(T)_{\rm reg} \rtimes W$. Explicitly, the Demazure operators act on $R(T)_{\rm reg} \otimes \Z[W]$ by 
\begin{equation} \la{muldemact} \Lambda_\alpha[\sum_{w \in W} f_w \otimes w]\,=\, \sum_{w \in W} \frac{(f_w-s_{\alpha}f_{s_\alpha w})}{(1-e^{\alpha})} \otimes w \ .\end{equation}
We claim
\begin{theorem} \la{cQfatodd}
The abelian subgroup $\bcQ'_{2k+1}(W)$ defined by \eqref{bcqm1} for $m=2k+1$ is stable under the action \eqref{muldemact} and hence carries a left module structure over the Hecke algebra $\D(W)$ extending the natural $R(T)$-module structure on $\bcQ'_{2k+1}(W)$ defined via \eqref{rtmodstr}.
\end{theorem}
\begin{proof}
 Let $u_w$ denote the coefficient of $w \in W$ in $u \in R(T) \otimes \Z[W]$. As in the proof of Theorem \ref{Qfatodd}, it is easy to check that for any $\alpha \in \R_+$, and $u = \sum_{w \in W}f_w \otimes w \in \bcQ'_{2k+1}(W)$, $\Lambda_{\alpha}(u) \in R(T) \otimes \Z[W]$ and that 
 \begin{equation} \la{checkdiffodd} (\Lambda_\alpha(u))_w- (\Lambda_\alpha(u))_{s_{\beta}w} \,=\, 
 \frac{(f_w-f_{s_\beta w})-s_\alpha(f_{s_\alpha w}-f_{s_\alpha s_\beta w})}{1-e^{\alpha}}\ .\end{equation}
 For $\beta \neq \pm \alpha$, the corresponding steps in the proof of Theorem \ref{Qfatodd} go through with trivial modifications to show that $ (\Lambda_\alpha(u))_w- (\Lambda_\alpha(u))_{s_{\beta}w} \in \langle \lambda_{k_\beta}(e^{\bar{\beta}})\rangle$, where $\lambda_k(e^{\bar{\beta}}):=(1-e^{\beta})(1-e^{\bar{\beta}})^{2k}$. For $\beta= \pm \alpha$, \eqref{checkdiffodd} equals $\frac{(1+s_\alpha)(f_w-f_{s_\alpha w})}{1-e^{\alpha}}$. Since $f_w-f_{s_\alpha w}=\lambda_{k_\alpha}(e^{\bar{\alpha}})g_w$ for some $g_w \in R(T)$, we have 
 $$(1+s_\alpha)(f_w -f_{s_\alpha w}) = \left\{\begin{array}{lr}
     \frac{\lambda_{k_\alpha}(e^{\bar{\alpha}})e^{-(k+1)\bar{\alpha}}\left(e^{(k+1)\bar{\alpha}}g_w-s_\alpha(e^{(k+1)\bar{\alpha}}g_w)\right)}{1-e^{\alpha}} & \text{if }\ \frac{\alpha}{2} \in \hat{T}\\*[2ex]
     \frac{\lambda_{k_\alpha}(e^{\bar{\alpha}})e^{-(k+1)\bar{\alpha}}\left(e^{(k+1)\bar{\alpha}}g_w-e^{\alpha}s_\alpha(e^{(k+1)\bar{\alpha}}g_w)\right)}{1-e^{\alpha}} &\ \text{otherwise.}
 \end{array}\right.$$
Since $1-e^{\alpha}|h-s_\alpha(h)$, the above formulas imply that $\lambda_{k_\alpha}(e^{\bar{\alpha}})|(\Lambda_\alpha(u))_w- (\Lambda_\alpha(u))_{s_{\alpha}w}$ as well. This proves the first assertion of Theorem \ref{cQfatodd}. The second follows from \eqref{fQexpGKM2} (specialized to the case $m=2k+1$).
\end{proof}

Theorem \ref{Qfat} has an exponential analog as well. For $k \in \M(W)$, define the trigonometric double affine Hecke algebra (a.k.a trigonometric Cherednik algebra)  $\bcH_k(W)$ as the $\k$-subalgebra of $\End_{\k}(R(T)_{\rm reg} \otimes \k[W])$ generated by the elements of $R(T)$ (viewed as multiplication operators acting by multiplication on the first factor), elements of $W$ (acting diagonally as in \eqref{expact1}) and the {\it trigonometric Dunkl operators}
\begin{equation} \la{TrigDunk} T^{\rm trig}_{\xi}:= \partial_{\xi}+\frac{1}{2} \sum_{\alpha \in \R_+} k_\alpha \bar{\alpha}(\xi) \left(\frac{1+e^{\bar{\alpha}}}{1-e^{\bar{\alpha}}}\right)(1-s_\alpha)\,,\,\,\xi \in \mathfrak{h}\, ,\end{equation}
where for $\xi \in \h$, $\partial_{\xi}$ stands for the first order differential operator on $R(T)_{\rm reg} \otimes \k$ sending $e^{\nu}$ to $\nu(\xi)e^{\nu}$ for all $\nu \in \hat{T} \subset \h^{\ast}$. Explicitly, the trigonometric Dunkl operators act on $R(T)_{\rm reg} \otimes \k[W] \cong (R(T)_{\rm reg} \otimes \k) \otimes_{\k} \k[W]$ as follows:
\begin{equation} \la{TrigDact} T^{\rm trig}_{\xi}\,\bigl[\sum_{w \in W} f_w \otimes w\,\bigr] \,=\, \sum_{w \in W} \big( \partial_\xi f_w +\frac{1}{2} \sum_{\alpha \in \R_+} k_\alpha \bar{\alpha}(\xi) \left(\frac{1+e^{\bar{\alpha}}}{1-e^{\bar{\alpha}}}\right)(f_w-s_\alpha f_{s_\alpha w})\big) \otimes w\ .\end{equation}

\begin{theorem} \la{cQfat}
The subspace $\bcQ_{2k}(W) \otimes \k$ defined by \eqref{bcqm2} is stable under the action \eqref{TrigDact} and hence carries a left module structure over the trigonometric Cherednik algebra  $\bcH_k(W)$ extending the natural $R(T)$-module structure on $\bcQ_{2k}(W) \otimes \k$ defined via \eqref{rtmodstr}.
\end{theorem}
\begin{proof} For $u \in R(T)_{\rm reg} \otimes \k[W]$, let $u_w \in R(T)_{\rm reg} \otimes \k$ denote the coefficient of $w \in W$ in $u$. For any $u=\sum_{w \in W} f_w \otimes w \in \bcQ_{2k}(W)$, \eqref{TrigDact} implies that $T^{\rm trig}_\xi(u) \in R(T) \otimes \k[W]$. Indeed, in this case, whenever $k_\alpha >0$, $1-e^{\bar{\alpha}}$ divides $f_w-f_{s_\alpha w}$ as well as $f_{s_\alpha w}-s_\alpha f_{s_\alpha w}$, implying that the coefficient of each $w \in W$ on the right hand side of \eqref{TrigDact} lies in $R(T) \otimes \k$. Further, for any $\alpha \in \R_+$ and $u=\sum_{w \in W} f_w \otimes w \in \bcQ_{2k}(W)$, \eqref{TrigDact} implies that for any $\xi \in \h$,
\begin{align} (T^{\rm trig}_\xi(u))_w -(T^{\rm trig}_\xi(u))_{s_\alpha w} &\, =\, \partial_\xi(f_w-f_{s_\alpha w})+\frac{k_\alpha}{2}\left(\frac{1+e^{\bar{\alpha}}}{1-e^{\bar{\alpha}}}\right)(1+s_\alpha)(f_w-f_{s_\alpha w})  \la{TDactDiff} \\
&+ \sum_{\substack{\beta \in \R_+\\ \beta \neq \alpha}} \frac{k_\beta}{2}\left(\frac{1+e^{\bar{\beta}}}{1-e^{\bar{\beta}}}\right)\big((f_w-f_{s_\alpha w})-s_\beta(f_{s_\beta w}-f_{s_{s_\beta \alpha}s_\beta w})\big) \ . \nonumber \end{align}
Note that $(1-e^{\bar{\alpha}})^{2k_\alpha}$ divides (each bracketed summand of) the expression $(f_w-f_{s_\alpha w})-s_\beta(f_{s_\beta w}-f_{s_{s_\beta \alpha}s_\beta w})$ for any $\beta \in \R_+, \beta \neq \alpha$. Since for such $\beta$, $1-e^{\bar{\beta}}$ does not divide $(1-e^{\bar{\alpha}})^{2k_\alpha}$, the second summand of \eqref{TDactDiff} is divisible by $(1-e^{\bar{\alpha}})^{2k_\alpha}$. 

It remains to study the first summand of \eqref{TDactDiff}. Putting $f_w-f_{s_\alpha w}=(1-e^{\bar{\alpha}})^{2k_\alpha}h_w$, we have 
\begin{equation} \la{rank1c1} \partial_\xi(f_w-f_{s_\alpha w})\,\equiv\, -2k_\alpha\bar{\alpha}(\xi)(1-e^{\bar{\alpha}})^{2k_\alpha-1}e^{\bar{\alpha}}h_w\,\equiv\,-2k_\alpha\bar{\alpha}(\xi)(1-e^{\bar{\alpha}})^{2k_\alpha-1}h_w\,,  \end{equation}
where the congruences are modulo $\langle (1-e^{\bar{\alpha}})^{2k_\alpha}\rangle$. On the other hand, 
\begin{eqnarray}
\frac{k_\alpha}{2}\left(\frac{1+e^{\bar{\alpha}}}{1-e^{\bar{\alpha}}}\right)(1+s_\alpha)(f_w-f_{s_\alpha w}) &=& \frac{k_\alpha}{2}\left(\frac{1+e^{\bar{\alpha}}}{1-e^{\bar{\alpha}}}\right)\big((1-e^{\bar{\alpha}})^{2k_\alpha}h_w+(1-e^{-\bar{\alpha}})^{2k_\alpha}s_\alpha(h_w)\big) \nonumber \\
&=& -k_\alpha \bar{\alpha}(\xi)\big((1-e^{\bar{\alpha}})^{2k_\alpha}h_w+(1-e^{-\bar{\alpha}})^{2k_\alpha}s_\alpha(h_w)\big) \nonumber\\
&+& \frac{k_\alpha \bar{\alpha}(\xi)}{(1-e^{\bar{\alpha}})}\big((1-e^{\bar{\alpha}})^{2k_\alpha}h_w+e^{-2k_\alpha \bar{\alpha}}(1-e^{\bar{\alpha}})^{2k_\alpha}s_\alpha(h_w)\big) \nonumber\\
&\equiv & k_\alpha \bar{\alpha}(\xi)(1-e^{\bar{\alpha}})^{2k_\alpha -1}\big(h_w+e^{-2k_\alpha \bar{\alpha}}s_\alpha(h_w)\big) \nonumber\\
&=& k_\alpha \bar{\alpha}(\xi)(1-e^{\bar{\alpha}})^{2k_\alpha -1}e^{-k_\alpha \bar{\alpha}}\big(e^{k_\alpha \bar{\alpha}}h_w+s_\alpha(e^{k_\alpha \bar{\alpha}}h_w)\big) \nonumber \\
& \equiv & 2k_\alpha \bar{\alpha}(\xi)(1-e^{\bar{\alpha}})^{2k_\alpha -1}h_w\,,\la{rank1c2}
\end{eqnarray}
where all congruences are modulo $\langle (1-e^{\bar{\alpha}})^{2k_\alpha}\rangle$ and the last congruence is because for all $g \in R(T)$, $g+s_{\alpha}g \equiv 2g$ modulo $1-e^{\bar{\alpha}}$. By \eqref{rank1c1} and \eqref{rank1c2}, the first summand of \eqref{TDactDiff} is also divisible by $(1-e^{\bar{\alpha}})^{2k_\alpha}$. This verifies the first assertion in the desired theorem. The second follows from \eqref{fQexpGKM} (specialized to the case $m=2k$).
\end{proof}
We close this section with the following conjecture motivated by Theorem \ref{ThFree}.
\begin{conjecture} 
\la{freeQmul}
For all $k \in \M(W)$, after $($possibly$)$ inverting $ |W|$, the following is true:
\begin{enumerate}
\item[$(i)$] $\bcQ_{2k}(W)$ is a free module of rank $|W|$ over $R(T)$,
\item[$(ii)$] $\cQ_{k}(W)$ is a free module of rank $|W|$ over $R(T)^W$,
\item[$(iii)$] $\bcQ'_{2k+1}(W)$ is a free module of rank $|W|$ over $R(T)$,
\end{enumerate}
where the $R(T)$-module structure on $\bcQ_{2k}(W) $ and $ \bcQ_{2k+1}'(W)$ is defined via \eqref{rtmodstr}. 
\end{conjecture} 
Note that we have implications: $\,(i)\,\Rightarrow\,(ii)\,\Rightarrow\,(iii)\,$, where the first one follows from Lemma~\ref{QmulHalf} and the second from Proposition~\ref{bcQodd}. Thus, to prove 
Conjecture~\ref{freeQmul} it suffices to prove Claim $(i)$. Recall that the analogous claim in the additive case --- the fact that $\bQ_{2k}(W)$ is a free module over $\k[V]$ --- is proved in \cite{{BC11}} by studying the module structure on $ \bQ_{2k}(W) $ over the rational Cherednik algebra $\bH_k(W)$ (see \cite[Proposition~8.1]{BC11}). We expect that a similar approach 
can be developed to prove Claim $(i)$ (at least, over a field $\k$ of characteristic zero). 

\subsection{Equivariant $K$-theory of quasi-flag manifolds} 
\la{S7.3}
The main result of this section is the following description of $K^\ast_G(F_m(G,T))$. Let $\Sc_k(\A)$ denote the submodule of $R(T)^{\A}$ of tuples of the form 
    $$
    \bigl\{ (f-f_\alpha)_{\alpha\in \R_+} \, \vert \, f\in R(T)\,,\, f_\alpha \in \cQ_{k_\alpha}(W_\alpha)\bigr\}\ .
    $$
\begin{theorem}
\la{KGFm}
    For any $m \in \M(W)$, let $k:=[\frac{m+1}{2}]$. There are isomorphisms (over $\Z[\frac{1}{2}]$)
    \begin{align} 
        K^0_G(F^{+}_m(G,T)) & \cong\, K^0_G(F_m(G,T))\,\cong\,  \cQ_{k}(W)\, \nonumber\\ 
     \la{kfmgt}   K^1_G(F^{+}_m(G,T)) & \cong\, K^1_G(F_m(G,T)) \,\cong\, \bigl(R(T)^{\A}/\Sc_k(\A)\bigr)\,,
    \end{align}
    where the first isomorphism is an isomorphism of rings and the second isomorphism is an isomorphism of $\cQ_k(W)$-modules. 
    If, moreover, $\pi_1(G)$ is torsion free, there are isomorphism of rings (over $\Z[\frac{1}{2}]$)
    $$
    K^0_T(F_m^{+}(G,T))\,\cong\, K^0_T(F_m(G,T))\, \cong \,\bcQ'_{2k+1}(W)\ .
    $$
\end{theorem}
We start with some computations that we need for our proof.
\begin{lemma}
\la{KTTam}
    $K^\ast_T((T/T_\alpha)^{\ast m})\cong \frac{R(T)}{(1-e^{\bar{\alpha}})^m}$.
\end{lemma}
\begin{proof}
    Note that for any $\alpha$, the exact sequence $1\to T_\alpha \to T \to T/T_{\alpha} \to 1$ of tori splits, and a choice of such splitting yields an isomorphism 
    $$
    K^\ast_T(\text{pt})\cong R(T) \cong R(T_\alpha) \otimes_\Z R(T/T_\alpha) \cong R(T_\alpha) \otimes_\Z \Z[e^{\pm \bar{\alpha}}].
    $$
    Note that the $T$-action on $(T/T_\alpha)^{\ast m}\cong S^{2m-1}$ factors through the surjection $T\twoheadrightarrow T/T_\alpha$, and the diagonal $T/T_\alpha$-action on $(T/T_\alpha)^{\ast m}\cong S^{2m-1}$ is free, hence 
    $$
    K^\ast_T((T/T_\alpha)^{\ast m}/(T/T_\alpha)) \cong K^\ast(S^{2m-1}/S^1) \cong K^\ast(\C \mathbb{P}^{m-1}) \cong \frac{\Z[\gamma]}{(\gamma-1)^m} \cong \frac{\Z[\gamma^{\pm1}]}{((\gamma-1)^m,}
    $$
    where $\gamma$ corresponds to the Hopf bundle $S^{2m-1}\to \C\mathbb{P}^{m-1}$, which in turn, corresponds to the generating representation $e^{\bar{\alpha}}$ of $R(T/T_\alpha)$. Therefore we have
    $$
    K^\ast_T((T/T_\alpha)^{\ast m}) \cong K^\ast_{T_\alpha}(\text{pt})\otimes_\Z K^\ast_{T/T_\alpha}\left((T/T_\alpha)^{\ast m}\right) \cong R(T_\alpha) \otimes_\Z \frac{\Z[e^{\pm \bar{\alpha}}]}{(1-e^{\bar{\alpha}})^m} \cong \frac{R(T)}{(1-e^{\bar{\alpha}})^m}\ .
    $$
    %
    %
\end{proof}
The following result is an analog of \Cref{APLjoin} for equivariant $K$-theory.
\begin{lemma}
\la{KTGam}
   For all $m\geq 0$, there is a $W_\alpha$-equivariant isomorphism of $\Z/2\Z$-graded rings 
    \begin{eqnarray*}
    K^\ast_T((G_\alpha/T)\ast(T/T_\alpha)^{\ast m}) \cong R(T) \times_{\frac{R(T)}{(\theta_m(e^{\bar{\alpha}}))}} R(T)\,,
    \end{eqnarray*}
    where $\theta_m(e^{\bar{\alpha}})=(1-e^{\alpha})(1-e^{\bar{\alpha}})^m$. Consequently, there is an isomorphism of $\Z/2\Z$-graded rings over $\Z[\frac{1}{2}]$
    \begin{eqnarray*}
    K^\ast_G(G\times_{N_\alpha}\left(G_\alpha/T\ast(T/T_\alpha)^{\ast m}\right)) \cong \cQ_{k}(W_\alpha) \,,
    \end{eqnarray*}
    where $k=[\frac{m+1}{2}]$.
\end{lemma}
\begin{proof}
    By \cite[Theorem 3.13]{KK87} (see also \cite[Remark 3.4]{HHH05}), 
    $$
    K^\ast_T(G_\alpha/T) \cong R(T)\times_{R(T)/(1-e^{\alpha})} R(T) \cong R(T) \oplus (1-e^{\alpha})\,,
    $$
    where the second isomorphism above is given by $(p,q)\mapsto (p,q-p)$.  
    Since $K^\ast_T(G_\alpha/T)$ is a free module over $R(T)$, combining with \Cref{KTTam}, the Hodgkin spectral sequence (see e.g., \cite[Theorem 2.3]{BrZ00} or \cite[Theorem 5.1]{BR1}) collapses to give
    \begin{eqnarray*}
         K^\ast_T((G_\alpha/T)\times (T/T_\alpha)^{\ast m}) & \cong & K^\ast_T(G_\alpha/T)\otimes_{R(T)} K^\ast_T((T/T_\alpha)^{\ast m}) \\
         & \cong & \left(R(T) \times_{\frac{R(T)}{(1-e^{\alpha})}} R(T)\right) \otimes_{R(T)} \frac{R(T)}{(1-e^{\bar{\alpha}})^m} \\
         & \cong & \frac{R(T)}{(1-e^{\bar{\alpha}})^m} \oplus \frac{(1-e^{\alpha})}{(\theta_m(\alpha))}\ .
    \end{eqnarray*}
    By \cite[Lemma 5.2]{BR1}, there is a six-term exact sequence 
    \begin{equation*}
    \begin{diagram}[small]
        K^0_T((G_\alpha/T)\ast(T/T_\alpha)^{\ast m}) & \rTo & K^0_T(G_\alpha/T) \oplus K^0_T((T/T_\alpha)^{\ast m}) & \rTo & K^0_T(G_\alpha/T\times (T/T_\alpha)^{\ast m})\\
        \uTo^{\partial} & &  & & \dTo_{\partial}\\
        K^1_T(G_\alpha/T\times (T/T_\alpha)^{\ast m}) & \lTo^{} & K^1_T(G_\alpha/T) \oplus K^1_T((T/T_\alpha)^{\ast m}) & \lTo &  K^1_T((G_\alpha/T)\ast(T/T_\alpha)^{\ast m})
    \end{diagram}
    \end{equation*}
    where $K^1_T(G_\alpha/T) = K^1_T((T/T_\alpha)^{\ast m}) = K^1_T(G_\alpha/T\times (T/T_\alpha)^{\ast m})=0$. 
    Since $K^0_T(G_\alpha/T) \twoheadrightarrow K^0_T(G_\alpha/T\times (T/T_\alpha)^{\ast m})$ is surjective, $K^1_T((G_\alpha/T)\ast(T/T_\alpha)^{\ast m})=0$ and $K^0_T((G_\alpha/T)\ast(T/T_\alpha)^{\ast m})$ is given by the pullback
    $$
    \lim \, \left\{R(T)\oplus (1-e^{\alpha}) \,\onto \,\frac{R(T)}{(1-e^{\bar{\alpha}})^m} \oplus \frac{(1-e^{\alpha})}{(\theta_m(e^{\bar{\alpha}}))}  \,\hookleftarrow \,\frac{R(T)}{(1-e^{\bar{\alpha}})^m} \,\right\}\, ,
    $$
    where the second map is the inclusion into the first factor. 
    Hence 
    $$
    K^0_T((G_\alpha/T)\ast(T/T_\alpha)^{\ast m}) \cong R(T) \oplus (\theta_m(e^{\bar{\alpha}})) \cong R(T) \times_{\frac{R(T)}{(\theta_m(e^{\bar{\alpha}}))}} R(T)\ .
    $$
    Note that the $W_\alpha$-action on the fiber product is $(p,q)\mapsto (s_\alpha(q),s_\alpha(p))$.
    Hence the composite map 
    $$
    \big(R(T) \times_{\frac{R(T)}{(\theta_m(e^{\bar{\alpha}}))}} R(T)\big)^{W_\alpha} \,\hookrightarrow \,R(T) \times_{\frac{R(T)}{(\theta_m(e^{\bar{\alpha}}))}} R(T) \,\xrightarrow{\text{pr}_1}\, R(T)
    $$
    is injective, with its image being 
    $$
    \cQ_{k}(W_\alpha) = \{p\in R(T) \,\vert\,\, p\,\equiv\, s_\alpha(p) \,\, \mod \,\,(\theta_m(e^{\bar{\alpha}}))\}\,,
    $$
    where $k=[\frac{m+1}{2}]$ by (the proof of) \Cref{QmulHalf} (more precisely, the equivalence of the congruence conditions \eqref{congc1} and \eqref{congc3}. 
    Next, observe that
    $$
    G\times_{N_\alpha} \left((G_\alpha/T)\ast(T/T_\alpha)^{\ast m}\right) \cong G\times_T \left((G_\alpha/T)\ast(T/T_\alpha)^{\ast m}\right)\,/\,W_\alpha\ .
    $$
    As the $W_\alpha$-action on $G\times_T \left((G_\alpha/T)\ast(T/T_\alpha)^{\ast m}\right)$ is free and commutes with the $G$-action, there are isomorphisms of $\Z/2\Z$-graded rings over $\Z[1/|W_\alpha|]$
    \begin{eqnarray*}
        K^\ast_G\left(G\times_{N_\alpha} \left((G_\alpha/T)\ast(T/T_\alpha)^{\ast m}\right)\right) & \cong & \left(K^\ast_G\left(G\times_T \left((G_\alpha/T)\ast(T/T_\alpha)^{\ast m}\right)\right)\right)^{W_\alpha} \\
        & \cong &  \left(K^\ast_T \left((G_\alpha/T)\ast(T/T_\alpha)^{\ast m}\right)\right)^{W_\alpha}\\
        & \cong & \cQ_{k}(W_\alpha) \ .
    \end{eqnarray*}
  This proves the desired lemma. 
\end{proof}

\begin{lemma}
\la{KGwedge}
    Let $G$ be a compact connected Lie group. Given compact $G$-spaces $X_0,X_1,\cdots,X_r$ and closed embeddings $f_i:X_0\to X_i$ of compact $G$-spaces for $1\leq i\leq r$, we denote ${X\coloneqq \left(\bigvee_{X_0} X_i\right)}$. 
    If $K^1_G(X_i)=0$ are injective for $0\leq i\leq r$, and $f^\ast_i:K^0_G(X_i)\to K^0_G(X_0)$ are injective for $1\leq i \leq r$, then 
    \begin{eqnarray*}
        K^0_G\left(X\right) \cong \Ker \left(\overline{\Delta}\right), \quad K^1_G\left(X\right) \cong \mathrm{Coker} \left(\overline{\Delta}\right).
    \end{eqnarray*}
    where
    \begin{equation}
    \la{KGmap}
        \overline{\Delta}: K^0_G(X_0) \to \bigoplus_{i=1}^r \frac{K^0_G(X_0)}{K^0_G(X_i)}
    \end{equation}
    is induced from the diagonal map and $K^0_G(X_i)$ are viewed as submodules in $K^0_G(X_0)$ via $f^\ast_i$.
\end{lemma}
\begin{proof}
    By \cite{Segal68} Proposition 2.6, Definition 2.7 and Proposition 3.5, there are six-term exact sequences for $1\leq i\leq r$:
    \[
    \begin{diagram}[small]
    K^0_G(X_i,X_0) & \rTo & K^0_G(X_i) & \rTo & K^0_G(X_0)\\
    \uTo^{\partial} & &  & & \dTo_{\partial}\\
    K^1_G(X_0) & \lTo & K^1_G(X_i) & \lTo & K^1_G(X_i,X_0) 
    \end{diagram}
    \]
    Hence $K^0_G(X_i,X_0)=0$ and 
    $$
    K^1_G(X_i,X_0)\cong \frac{K^0_G(X_0)}{K^0_G(X_i)}.
    $$ 
    By \cite{Segal68} Proposition 2.9, there are natural isomorphisms $K^\ast(X_i\backslash X_0) \cong K^\ast_G(X_i,X_0)$, therefore
    \[
    K^\ast_G(X\backslash X_0) \cong K^\ast_G\left(\bigsqcup_{i=1}^n (X_i\backslash X_0)\right) \cong \bigoplus_{i=1}^n K^\ast_G(X_i,X_0).
    \]
    Therefore $K^0_G(X,X_0)=K^0_G(X\backslash X_0)=0$ and 
    \[
    K^1_G(X,X_0)=K^1_G(X\backslash X_0) = \bigoplus_{i=1}^r K^1_G(Xi,X_0) = \bigoplus_{i=1}^r \frac{K^0_G(X_0)}{K^0_G(X_i)}.
    \]
    We have another six-term exact sequence
    \[
    \begin{diagram}[small]
    K^0_G(X,X_0) & \rTo & K^0_G(X) & \rTo & K^0_G(X_0)\\
    \uTo^{\partial} & &  & & \dTo_{\partial}\\
    K^1_G(X_0) & \lTo & K^1_G(X) & \lTo & K^1_G(X,X_0) 
    \end{diagram}
    \]
    where $\partial:K^0_G(X_0) \to K^1_G(X,X_0)$ can be identified with \eqref{KGmap}. This proves the desired lemma.
\end{proof}
We are now ready to prove the main result in this section.
\begin{proof}[Proof of \Cref{KGFm}]
    By \Cref{natdec}, 
    $F_m(G,T) \cong \hocolim_{\Sc(W)}\left( G\times_T \cF_m\right)^\natural$ 
    where $\left( G\times_T \cF_m\right)^\natural:\Sc(W)\to \Top^G$ is given by \eqref{GNWG}. Since for any $\alpha\in \cA$, we have a natural inclusion $N_\alpha/T\hookrightarrow G_\alpha/T\ast (T/T_\alpha)^{\ast m_\alpha}$, we have a natural inclusion 
    \begin{equation}
    \la{GTinc}
        G/T =G\times_{N_\alpha}N_\alpha/T \hookrightarrow G\times_{N_\alpha}\left(G_\alpha/T\ast(T/T_\alpha)^{\ast m_\alpha}\right)\ .
    \end{equation}
    It follows that 
    \[
    F_m(G,T) \cong \bigvee_{\substack{G/T \\\alpha\in \A}} G\times_{N_\alpha}\left(G_\alpha/T\ast (T/T_\alpha)^{\ast m_\alpha}\right)\ .
    \]
    In addition, the spaces $G/T$ and $G\times_{N_\alpha}\left(G_\alpha/T\ast (T/T_\alpha)^{\ast m_\alpha}\right)$ are compact $G$-spaces, and therefore the inclusions \Cref{GTinc} are proper. By \Cref{KTGam}, since $K^{\ast}_G(G/T)=R(T)$, the induced maps (over $\Z[\frac{1}{2}]$)
    \[
    K^\ast_G\left(G\times_{N_\alpha}(G_\alpha/T\ast (T/T_\alpha)^{\ast m_\alpha})\right) \,\cong\, \cQ_{k_\alpha}(W_\alpha) \hookrightarrow R(T) = K^\ast_G(G/T)
    \]
    are injective, where $k_\alpha=[\frac{m_\alpha+1}{2}]$. Hence, by Lemma \ref{KGwedge}, there are isomorphisms over $\Z[\frac{1}{2}]$
    \begin{equation*}
        K^0_G(F_m(G,T) \,\cong\,\Ker \left(R(T) \xrightarrow{\overline{\Delta}} \bigoplus_{\alpha\in\R_+}\frac{R(T)}{\cQ_{k_\alpha}(W_\alpha)}\right)
        \,\cong\, \bigl(\bigcap_{\alpha\in\R_+} \cQ_{k_\alpha}(W_\alpha)\bigr) 
        \,=\, \cQ_{k}(W)\,,
    \end{equation*}
    and 
    \begin{equation*}
        K^1_G(F_m(G,T) \cong \mathrm{Coker} \left(R(T) \xrightarrow{\overline{\Delta}} \bigoplus_{\alpha\in\R_+}\frac{R(T)}{\cQ_{k_\alpha}(W_\alpha)} \right)\cong \left(R(T)^{\A}/\Sc_k(A)\right)\,,
    \end{equation*}
    where the second isomorphism is a direct computation similar to the computation of \eqref{HQodd}: indeed, the two term complex 
    $$\bigl(\prod_{\alpha \in \A} \cQ_{k_\alpha}(W_\alpha) \bigr) \oplus R(T) \to R(T)^{\A}\,,\qquad ((f_\alpha)_{\alpha \in \A},f) \mapsto (f_{\alpha}-f)_{\alpha \in \A} $$
    has an acyclic subcomplex $\prod_{\alpha \in \A} \cQ_{k_\alpha}(W_\alpha) \stackrel{\id}{\to} \prod_{\alpha \in \A} \cQ_{k_\alpha}(W_\alpha) \subset R(T)^{\A}$, and the quotient complex is isomorphic to $R(T) \xrightarrow{\overline{\Delta}} \oplus_{\alpha \in\A}\frac{R(T)}{\cQ_{k_\alpha}(W_\alpha)}$. This proves the first two assertions of Theorem \ref{KGFm} regarding $K^{\ast}_G(F_m(G,T))$. In particular, we have isomorphisms over $\Z[\frac{1}{2}]$
    \begin{equation} \la{kf0gt} K^0_G(F_0(G,T))\,\cong\,\cQ_0(W) = R(T)\,,\qquad K^1_G(F_0(G,T))\,=\,0 \ .\end{equation}
    In particular, if $q:F_0(G,T) \to G/T$ is the map defined by \eqref{eq:p-hocolim} then $q^{\ast}\,:\,R(T)\,\cong\,K^{\ast}_G(G/T) \to K^{\ast}_G(F_0(G,T))$ is an isomorphism of ($\Z/2\Z$-graded algebras). By the homotopy cocartesian square \eqref{diagrm-} in $\Top^G$ from Corollary \ref{corpplus}, whose compact $G$-equivariant model is given by \eqref{rmodel+}, and by \cite[Lemma 5.2]{BR1}, there is a six term exact sequence
    $$
    \begin{diagram}[small]
    K^0_G(F_m^{+}(G,T)) & \rTo & K^0_G(F_m(G,T)) \oplus K^0_G(G/T) & \rTo & K^0_G(F_0(G,T))\\
    \uTo^{\partial} & &  & & \dTo_{\partial}\\
    K^1_G(F_0(G,T)) & \lTo & K^1_G(F_m(G,T)) \oplus K^1_G(G/T) & \lTo & K^1_G(F_m^+(G,T)) 
    \end{diagram}
    $$
    It then follows from \eqref{kf0gt} that $K^\ast(F_m^+(G,T))\,\cong\,K^{\ast}(F_m(G,T))$, completing the proof of the first two assertions of Theorem \ref{KGFm}.

    If $\pi_1(G)$ is torsion free, then by \cite[Theorem 4.4]{HLS10}, there is an isomorphism of $\Z/2\Z$-graded $R(T)$-algebras
    $$ R(T) \otimes_{R(G)} K^{\ast}_G(F^+_m(G,T))\,\cong\,K^{\ast}_T(F^+_m(G,T))\,, \qquad R(T) \otimes_{R(G)} K^{\ast}_G(F_m(G,T))\,\cong\,K^{\ast}_T(F_m(G,T))\ . $$
    The final assertion of Theorem \ref{KGFm} then follows immediately from the first assertion of Theorem \ref{KGFm} and Proposition \ref{bcQodd}.
\end{proof}

As a consequence of Theorem \ref{KGFm} and Theorem \ref{cQfatodd} we have 
\begin{cor} \la{mulHeckeact}
 If $\pi_1(G)$ is torsion free, then the $R(T)$-module structure on $K^0_T(F_m(G,T))$ extends to a left module structure over the Hecke algebra $\D(W)$ (over $\Z[\frac{1}{2}]$).  
\end{cor}
\begin{remark} When $m=0$, the $\D(W)$-action on $K^0_T(F_0(G,T))$ specializes to the localization (inverting the prime $2$) of the left Hecke action on $K^0_T(G/T)$ that appears in \cite{KK87}.  
\end{remark}
Finally, recall the $T$-space $\tilde{F}_m(T)$ defined by \eqref{tilFmT}. The corresponding $G$-space $\tilde{F}_m(G,T):=G \times_T \tilde{F}_m(T)$ carries a natural free $W$-action defining a $W$-bundle \eqref{Wbund} over $F_m(G,T)$: $\tilde{F}_m(G,T) \to F_m(G,T$. The following proposition is a multiplicative analog of Proposition \ref{tilProp}.
\begin{prop} \la{kgfmtilde}
For all $m \in \M(W)$, there are natural isomorphisms of rings
\begin{equation} \la{bcQptop} K^0_G(\tilde{F}_m(G,T)) \,\cong\, K^0_T(\tilde{F}_m(T))\,\cong\, \bcQ_{m+1}(W) \end{equation}
which fit in the commutative diagram (over $\Z[1/|W|]$)
\begin{equation} \la{}
\begin{diagram}[small]
 K^0_G(F_m(G,T)) &\rTo^{\cong} & K^0_G(\tilde{F}_m(G,T))^{W}\\
 \dTo^{\eqref{kfmgt}} & & \dTo_{\eqref{bcQptop}}\\
 \cQ_{[\frac{m+1}{2}]}(W) &\rTo^{\eqref{bcqwinv}} & \bcQ_{m+1}(W)^{W}\\
 \end{diagram}
\end{equation}
where the top isomorphism is induced by the $W$-bundle \eqref{Wbund} and the bottom is defined in (the proof of) Lemma \ref{QmulHalf}.
\end{prop}
\begin{proof} Clearly, 
$$K^\ast_G(\tilde{F}_m(G,T))\,=\, K^{\ast}_G(G \times_{T} \tilde{F}_m(T))\,\cong\, K^{\ast}_T(\tilde{F}_m(T))\ .$$
In order to compute the $T$-equivariant $K$-theory $K^{\ast}_T(\tilde{F}_m
(T))$, we consider the compact model $\widetilde{\bF}_{\mathbb{R}}^{(m)}$ (see \eqref{realFm}) of $\tilde{F}_m$. Recall that 
\begin{equation} \la{tfmwedge} \widetilde{\bF}_{\mathbb{R}}^{(m)}\,=\, \bigvee_{\substack{{\bF}(W_0)\\\alpha \in \A}} \bF^{(m_\alpha)}(W_{\alpha})\,,\end{equation}
where
$$\bF(W_0)= \coprod_{w \in W} \{w\}\,,\qquad \bF^{(m_\alpha)}(W_\alpha)=\coprod_{e(s_\alpha,w) \in E_{\Gamma}} S(\cO_{e(s_\alpha,w)}^{\ast (m_\alpha+1)})\,, $$
and for each $\alpha \in \A$, the map $f_\alpha\,:\,\bF(W_0) \into \bF^{(m_\alpha)}(W_\alpha)$ maps the elements of each pair $\{w, s_\alpha w\} \subseteq \bF(W_0)$ to the poles of $S(\cO_{e(s_\alpha,w)}^{\ast (m_\alpha+1)})$. Note that there are $T$-equivariant homotopy equivalences
$$S(\cO_{e(s_\alpha,w)}^{\ast (m_\alpha+1)}) \simeq C_{w}(\cO_{e(s_\alpha,w)}^{\ast (m_\alpha+1)}) \cup C_{s_\alpha w}(\cO_{e(s_\alpha,w)}^{\ast (m_\alpha+1)})\,,\qquad  C_{w}(\cO_{e(s_\alpha,w)}^{\ast (m_\alpha+1)}) \cap C_{s_\alpha w}(\cO_{e(s_\alpha,w)}^{\ast (m_\alpha+1)}) \simeq \cO_{e(s_\alpha,w)}^{\ast (m_\alpha+1)}\,,$$
where $C_{w}$ stands for the cone with vertex $w$. Since the cones are $T$-equivariantly contractible, $K^\ast_T(C_{w}(\cO_{e(s_\alpha,w)}^{\ast (m_\alpha+1)})) \cong K^\ast_T(C_{s_\alpha w}(\cO_{e(s_\alpha,w)}^{\ast (m_\alpha+1)})) \cong R(T)$. By Lemma \ref{KTTam}, $K^{\ast}_T(\cO_{e(s_\alpha,w)}^{\ast (m_\alpha+1)}) \cong R(T)/((1-e^{\bar{\alpha}})^{(m_\alpha+1)})$. It follows from the six term exact sequence (see \cite[Lemma 5.2]{BR1})
 $$
    \begin{diagram}[small]
    K^0_T(S(\cO_{e(s_\alpha,w)}^{\ast (m_\alpha+1)})) & \rTo & K^0_T(C_{w}(\cO_{e(s_\alpha,w)}^{\ast (m_\alpha+1)})) \oplus K^0_T(C_{s_\alpha w}(\cO_{e(s_\alpha,w)}^{\ast (m_\alpha+1)})) & \rTo & K^0_T(\cO_{e(s_\alpha,w)}^{\ast (m_\alpha+1)})\\
    \uTo^{\partial} & &  & & \dTo_{\partial}\\
    K^1_T(\cO_{e(s_\alpha,w)}^{\ast (m_\alpha+1)})) & \lTo &  K^1_T(C_{w}(\cO_{e(s_\alpha,w)}^{\ast (m_\alpha+1)})) \oplus K^1_T(C_{s_\alpha w}(\cO_{e(s_\alpha,w)}^{\ast (m_\alpha+1)})) & \lTo & K^1_T(S(\cO_{e(s_\alpha,w)}^{\ast (m_\alpha+1)})) 
    \end{diagram}
    $$

that $K^{\ast}_T(S(\cO_{e(s_\alpha,w)}^{\ast (m_\alpha+1)}))\,\cong\, R(T) \times_{R(T)/((1-e^{\bar{\alpha}})^{(m_\alpha+1)})} R(T)$. Hence,
\begin{equation} \la{ktfmwa} K^{\ast}_T(\bF^{(m_\alpha)}(W_\alpha))\,\cong\, \bcQ_{m_\alpha+1}(W_{\alpha})\,,\end{equation}
where $\bcQ_{m}(W_{\alpha}):= \{\sum_{w \in W} f_w \otimes w \in R(T) \otimes \Z[W]\,|\, f_w-f_{s_\alpha w} \in \langle (1-e^{\bar{\alpha}})^{m}\rangle \,\forall \,w \in W\}$ and the maps $f_\alpha^{\ast}:R(T)/((1-e^{\bar{\alpha}})^{(m_\alpha+1)}) \to K^{\ast}_T(\bF(W_0)) \cong R(T) \otimes \Z[W]$ are identified with the obvious inclusions 
$\bcQ_{m_\alpha+1}(W_\alpha) \into R(T) \otimes \Z[W]$. It follows from Lemma \ref{KGwedge} that 
$$K^0_T(\tilde{F}_m(T)) \,\cong\, K^0_T(\bF^{(m)}_{\mathbb R})\,\cong\, \bigcap_{\alpha \in \A} \bcQ_{m_\alpha+1}(W_{\alpha}) \,=\,\bcQ_{m+1}(W)\ . $$
This proves the first assertion of Proposition \ref{kgfmtilde}. The verify the second assertion, note that the map $\coprod_{w \in W} \{w\} \into \bF^{(m)}_{\mathbb R}$ is $N(T)$-equivariant. Hence the $W$-action on the $W$-action on $K^0_T(\bF^{(m)}_{\mathbb R}) \subseteq  R(T) \otimes \Z[W]$ is the restriction of the $W$-action on $ R(T) \otimes \Z[W]$ induced by the $N(T)$-action on $\coprod_{w \in W} \{w\}$. This is indeed the left action \eqref{expact1}. The second assertion of Proposition \ref{kgfmtilde} then follows from Lemma \ref{QmulHalf}.
\end{proof}

As a consequence of Proposition \ref{kgfmtilde} and Theorem \ref{cQfat}, we get

\begin{cor} \la{trigDAHAact}
 For $m \in \M(W)$ odd, the $T$-equivariant $K$-theory $K^0_T(\tilde{F}_m(T)) \otimes \k$ carries a natural action of the trigonometric Cherednik algebra $\bcH_{\frac{m+1}{2}}(W)$ inducing an action of the spherical algebra $\boldsymbol{e}\bcH_{\frac{m+1}{2}}(W)\boldsymbol{e}$ on $K^0_G(F_m(G,T)) \otimes \k$.
\end{cor}

Corollary \ref{trigDAHAact} provides a topological realization for the (differential) action of the trigonometric Cherednik algebra $\bcH_k(W)$ on $\bcQ_{2k}(W)$ constructed in Theorem \ref{cQfat}, which, in turn, is a multiplicative analog of the action of the rational Cherednik algebras $\bH_k(W)$ on quasi-covariants constructed in \cite{BC11}.

\appendix
\section{\texorpdfstring{$W$}{W}-model categories and \texorpdfstring{$W$}{W}-functors}
\la{AA}
For a (discrete) group $W$, the notion of a $W$-category --- a small category equipped with a strict or weak action of $W$ --- is standard: it arises in many areas of mathematics and has been studied extensively in the literature (see, e.g., \cite{De97}). Less standard
(and apparently less known) is the related notion of a $W$-functor (or $W$-diagram). This  notion was introduced in topology in \cite{JS01} (and independently in \cite{VF99}), and in recent years, has been developed systematically as part of a general model-categorical framework for equivariant homotopy theory (see \cite{DM16, Dotto16, Dotto17, Moi20}). In this Appendix, we briefly review definitions and basic constructions related to $W$-functors, following \cite{DM16}, and prove a few technical results needed for the present paper.  

\vspace*{1ex}

Throughout this section, $W$ stands for a fixed finite group.

\subsection{\texorpdfstring{$W$}{W}-categories and \texorpdfstring{$W$}{W}-functors}
\la{AA1}
If $ \Cc $ a (small) category, a strict action of $W$ on $ \Cc $ is 
defined to be a functor $ W \to \Cat $ sending $ \ast $ to $ \Cc$, where $W$ is viewed as a category with the single object `$\ast$'.  Abusing the notation, we continue to denote this functor by $\Cc$.
Explicitly, giving such an action  amounts to giving a family of endofunctors $ \hat{g}: \Cc \to \Cc $, one for each element
$g \in W $, satisfying the natural conditions: $ \hat{e} = {\rm Id}_{\Cc}$ and $ \widehat{gh} = \hat{g} \circ \hat{h} $ for all $g,h \in W$. We  refer to a category equipped with a strict $W$-action as a {\it $W$-category}.  

Next, we introduce our basic notion of a $W$-functor.
\begin{defi}[\cite{JS01}]
\la{Wfun}
Let $\Cc$ be a $W$-category, and let $F: \Cc \to \Dc $ be a functor from $\Cc$ to an
arbitrary category $\Dc$. A {\it $W$-action} on $F$ (with respect to the given action of $W$ on $\Cc$)
is defined by a family of natural transformations $ \vartheta_g: F \to F \circ \hat{g} $, one for each
element $g \in W $ with $ \vartheta_e := {\rm Id}_F $, such that the following diagrams  commute for all elements $ g,h \in W \,$:
\begin{equation}
\la{diagact}
\begin{diagram}[small] 
                F & \rTo^{\vartheta_h} & F \circ \hat{h} \\
                  & \rdTo_{\vartheta_{gh}} & \dTo_{(\vartheta_g)|_{h}}  \\            
                  &                    &      F \circ \hat{g} \circ \hat{h}  \\
               \end{diagram}  
\end{equation}
where $ (\vartheta_g)|_h $ denotes the restriction of  $ \vartheta_g: F \to F \circ \hat{g}$ along the endofunctor $\hat{h}: \Cc \to \Cc$. We will refer to a functor equipped with a $W$-action as a {\it $W$-functor} (or a {\it $W$-diagram} in $\Dc$ of shape $\Cc$).
\end{defi}
\begin{remark}
If all objects of a category $ \Dc$ are equipped with $W$-action and all morphisms are $W$-equivariant (e.g., $ \Dc = \Top^W$, the category of topological $W$-spaces), then any functor $ F: \Cc \to \Dc $ is a $W$-functor (with respect to the trivial action of $W$ on $\Cc$).
In this way, the notion of a $W$-functor generalizes the usual diagrams of $W$-spaces and $W$-equivariant maps. 
\end{remark}

There is a useful alternative characterization of $W$-functors in terms of the classical Grothendieck construction. Recall
\begin{defi}
\la{Grothcon}
Given a functor $\,F: \I \to \Cat\,$, the (covariant) {\it Grothendieck construction on $F$} is the category $\,\I\smallint F\,$ defined by the following data: 
\begin{itemize}
\item The objects of $\,\I\smallint F\,$ are the pairs $\,(i,\,x)\,$, where $ i \in \rm Ob(\I) $ 
and $x \in {\rm Ob}[F(i)] $.
\item  For a pair of objects $(i,\,x)$ and $(i',\,x')$, a morphism
$(i,\,x) \to (i',\,x') $ in $\,\I\smallint F\,$ is represented by a pair $(f,\, \varphi)$, where
$f: i\to i' $ is a morphism in $ \I $ and $\varphi: F(f)(x) \to x' $ is a morphism in the 
category $ F(i')$.
\item Composition in $\,\I\smallint F\,$ is defined by $\,(g, \psi) \circ (f,\varphi)=(gf, \psi \circ F(g)(\varphi))\,$.
\end{itemize}
\end{defi}
The category $\,\I\smallint F\,$ comes equipped with the canonical functor $\, p: \, I\smallint F \to \I \,$ defined by $\, (i,x) \mapsto i $ and $(f,\varphi) \mapsto f $,
that represents a cocartesian fibration in $\Cat$. The functor $\,F: \I \to \Cat\,$ 
plays the role of (and is often referred to as) a  classifying functor for this fibration.

The meaning of the Grothendieck construction from a homotopy-theoretic point of view is clarified
by Thomason's Homotopy Colimit Theorem \cite{To79}. This theorem asserts that there is a natural weak equivalence of simplicial sets:
\begin{equation}
\la{ThomForm}
\hocolim_{\I}(NF)\,\simeq\, N(\I\smallint F) \,,  
\end{equation}
where $N: \Cat \to \sset $ is the classical nerve functor. Thus, $ \I\int F $ is a canonical model for the homotopy colimit $\hocolim_{\I}(F) $ of the functor $F: \I \to \Cat $ in the category $\Cat$ of small categories. 

In particular, we can apply the above observation to a $W$-category $\Cc$, viewed as a functor $ W \to \Cat\,$. By definition, the homotopy colimit $\, \Cc_{hW} := \hocolim_{W}(\Cc) \,$ in $\Cat $
is the homotopy quotient of $\Cc$ with respect to the given $W$-action, and the corresponding Grothendieck construction $\,W \smallint \Cc\,$ is thus a canonical model for this quotient. We will identify $\,\Cc_{hW} = W \smallint \Cc \,$ so that the category $ \Cc_{hW} $ has the following explicit description:
\begin{itemize}
\item $\,{\rm Ob}(\Cc_{hW})={\rm Ob}(\Cc)\,$.
\item A morphism $\,c \to c'\,$ in $ \Cc_{hW} $ is a pair $\,(g,\ \varphi: \,g(c) \to c')\,$,
where $ g \in W $ and $ \varphi \in {\rm Mor}(\Cc)$.
\item Composition is given by 
$$
(h,\ \psi:\, h(c') \to c'') \circ 
(g,\ \varphi: \, g(c) \to c') = (hg,\ \psi \circ h(\varphi):\, hg(c) \to c'')\,.
$$
\end{itemize}
\noindent
If $ \Cc $ is a $W$-category, there is a canonical inclusion functor $\,\iota:\Cc \into \Cc_{hW}$ that acts as the identity on objects and sends $\varphi \in {\rm Mor}(\Cc) $ to 
$(e, \,\varphi) \in {\rm Mor}(\Cc_{hW})$. Now, any $W$-functor 
$\,F:\,\Cc \to \Dc\,$ extends naturally through $ \iota \,$: indeed, the corresponding
functor $\tilde{F}:\,\Cc_{hW} \to \Dc$  on $\Cc_{hW}$ is defined on objects
by $\,\tilde{F}(c) := F(c)\,$, while on morphisms $\tilde{F}(g, \varphi) := F(\varphi) \circ \vartheta_g(c)$, where $\vartheta $ is a $W$-structure on $F$. It turns out that this extension property characterizes $W$-functors: namely, we have (see \cite[Lemma 1.9]{DM16})
\begin{lemma}
\la{LA4}
For any $W$-category $\Cc$, the assignment $F \mapsto \tilde{F}$ defines an equivalence from the category $\Dc^{\Cc}_W$ of $W$-diagrams of shape $\Cc$ in $\Dc$ to the category $\Dc^{\Cc_{hW}}$ of diagrams of shape $\Cc_{hW}$ in $\Dc$, the inverse being the restriction functor $\iota^\ast$. 
\end{lemma}

The key property of a $W$-functor taking values in a homotopical category (e.g., the category of
topological spaces) is that its homotopy limit and homotopy colimit --- whichever exists --- carries a natural $W$-action. We discuss this property in full generality in the next section, after developing a necessary model-categorical background, but at the moment we state its elementary consequence (see {\cite[Corollary 1.5]{DM16}}).

\begin{lemma}\la{Wob}
 Let $ \Cc $ be a $W$-category, and $ F: \Cc \to \Dc $ a $W$-functor with $W$-action $\vartheta$. If $ \colim_{\Cc}(F) $ exists in $\Dc $, then $X := \colim_{\Cc}(F) $ is a $W$-object in $\Dc$ with 
$W$-action induced  by $\vartheta$. Similarly, if $\lim_{\Cc}(F)$ exists in $\Dc$, then $Y:=\lim_{\Cc}(F)$ is a $W$-object in $\Dc$ with 
$W$-action induced  by $\vartheta$.
\end{lemma}
\begin{proof}
By definition, the colimit of $ F $ is represented by a collection of morphisms 
$\{\varphi_{c}: F(c) \to X \}$ in $\Dc$ indexed by the objects $ c \in \Cc $ that form a universal (initial) cone on $F$. Up to re-indexing, the same collection of morphisms represents the colimit of each functor $ F \circ \hat{g} $. Hence, we can identify $ \colim_{\Cc}(F \circ \hat{g}) \cong \colim_{\Cc}(F) = X $ for all $ g \in W $. With this identification, the natural
transformations $ \vartheta_g: F \to F \circ \hat{g} $ defining the $W$-action on $F$ induce  
automorphisms $ \tilde{g}: X \to  X $ defining a $W$-action on $X$. The proof that $Y:=\lim_{\Cc}(F)$ carries a $W$-action induced by $\vartheta$ is similar and left to the reader.
\end{proof}
With a little more effort, we can prove the following
\begin{cor} \la{Wcolim}
Let $\Cc$ be a $W$-category, and $\Dc$ a category with arbitrary $($small$)$ colimits and limits. Then the adjunction $\colim_{\Cc}:\Dc^{\Cc} \rightleftarrows \Dc\,:\,\Delta$ enriches to the adjunction
$$\colim_{\Dc,W}\,:\,\Dc^{\Cc}_W \rightleftarrows \Dc^W\,:\,\Delta\,,$$
where $\Delta$ is the constant diagram functor. Dually, the adjunction 
$\Delta\,:\,\Dc \rightleftarrows \Dc^{\Cc}\,:\,\lim_{\Cc}$ enriches to 
$$ \Delta\,:\,\Dc^{W} \rightleftarrows \Dc^{\Cc}_W\,:\,{\rm lim}_{\Cc,W}\ .$$
The functors $ \lim_{\Cc,W} $ and $ \colim_{\Cc,W} $ are called
{\rm $W$-limit} and {\rm $W$-colimit}, respectively.
\end{cor}

\subsection{$W$-model categories and homotopy colimits} 
\la{AA2}
We assume that the reader is familiar with basics of Quillen's theory of model categories. A comprehensive reference is \cite{Hir03}, but we also recommend \cite{DS95} for a gentle introduction
to the subject and \cite[Appendix~A.2]{HTT} for an excellent technical overview. 

Given a category $\Mc$ and a subgroup $ H \leq  W$, we write $\Mc^H $ for the category of $H$-objects in $\Mc$, which is simply the functor category $ \Mc^H := \Fun(H, \Mc)\,$.
The next definition is due to \cite{DM16}.
\begin{defi} 
\la{WMCat}
A {\it $W$-model category} is a cofibrantly generated simplicial model category $\Mc$  together with the data of a cofibrantly generated model structure on $\Mc^{H}$ for every subgroup $H$ of $W$, satisfying:
\begin{enumerate}
\item[$(1)$] The model structure on $\Mc^H$, together with the 
$\sset^{H}$-enrichment (cotensor and tensor)  induced from $\Mc$ is a cofibrantly generated $\sset^H$-enriched model structure on $\Mc^H$.

\item[$(2)$] For every pair of subgroups $H,H'$ of $W$ and any finite set $K$ equipped with commuting free left $H'$-action and free right $H$-action, the natural adjunction
$$ K \otimes_H (\mbox{--})\,:\ \Mc^H \rightleftarrows \Mc^{H'}\ :\,\Hom_{H'}(K,\,\mbox{--}\,)  $$
is a Quillen pair. 
\end{enumerate}
\end{defi}
\begin{remark} \la{res}
Let $ H' \leqslant H $ be subgroups of $W$. Take $ K = H \,$ to be the underlying set of $H$. If $K=H$ is given the $H'$-action by left multiplication and the $H$-action by right multiplication, then $H \otimes_H \mbox{--}$ coincides with the restriction functor $\Mc^H \to \Mc^{H'}$. One the other hand, if $K=H$ is given with left $H$-action and right $H$-action by multiplication, then the functor $\Hom_{H}(H, \mbox{--}):\Mc^H \to \Mc^{H'}$ also coincides with the restriction $\Mc^H \to \Mc^{H'}$. Hence, the restriction $\Mc^H \to \Mc^{H'}$ is left, as well as right, Quillen. Taking $H'$ to be the trivial subgroup, we see that the forgetful functor $\Mc^H \to \Mc$ is left --- as well as right --- Quillen. 
\end{remark}

\begin{example} \la{naive} For any cofibrantly generated $\sset$-enriched model category $\Mc$, the collection of projective model structures on the $\Mc^H$ for $H \leqslant W$ defines a $W$-model structure on $\Mc$ (see \cite[Example 2.3]{DM16}). We refer to this $W$-model structure as the {\it coarse} $W$-model structure on $\Mc$. In this paper, we consider the coarse $W$-model structure on $\Mc$ in the following cases:
\begin{itemize}
\item $\Comm_\k$ equipped with the trivial model structure (Section \ref{S2}).
\item The category $\ccdga_{\k}^{\leq 0}$ of non-positively graded cochain CDGAs with standard model structure induced from non-negatively graded chain complexes (see Section \ref{S3.1} and Appendix \ref{DGModel}).
\item The category $\ccdga_{\k}^{\geq 0}$ of non-negatively graded cochain CDGAs with Bousfield-Gugenheim model structure (see Section \ref{S3.1} and Appendix \ref{DGModel}.)
\end{itemize}
\end{example}

\begin{example} \la{fine}
Let $\Mc$ be a cofibrantly generated $\sset$-enriched model category. For any pair of subgroups $L \leqslant H$ of $W$, consider the fixed point functor $(\mbox{--})^L:\Mc^W \to \Mc$. If these functors are cellular in the sense of (\cite{GuM20},\cite[Prop. 2.6]{Ste16}) then $\Mc^H$ admits a (fine) model structure in which a weak equivalence (resp., fibration) is precisely a map that is sent by $(\mbox{--})^L$ to a weak equivalence (resp. fibration) for every $L \leqslant H$. The fine model structures on $\Mc^H, H\leqslant W$ assemble to yield a $W$-model structure on $\Mc$, which we shall refer to as the {\it fine} $W$-model structure (see \cite{GuM20}, \cite[Prop. 2.6]{Ste16} and \cite[Example 2.4]{DM16}). 
\end{example}

Let $\Cc$ be a $W$-category, and let $F:\Cc \to \Mc$ be a $W$-functor. Then, for each object $c \in \Cc$, $F(c)$ inherits an action of the stabilizer $W_c$ of $c$ in $W$. This defines an evaluation functor ${\rm ev}_c: \Mc^{\Cc}_W \to \Mc^{W_c}$. 

\begin{theorem}[{\cite[Theorem 2.6]{DM16}}] \la{Wmodel}
Let $\Mc$ be a $W$-model category. There is a cofibrantly generated model structure on the category $\Mc^{\Cc}_W$ of $W$-diagrams of shape $\Cc$ in $\Mc$, with
\begin{enumerate}
\item[$(1)$] weak equivalences the maps $f:X \to Y$ of $W$-diagrams whose evaluations ${\rm ev}_c(f)$ are weak equivalences in $\Mc^{W_c}$ for every $c \in \Cc$.

\item[$(2)$] fibrations the maps $f:X \to Y$ of $W$-diagrams whose evaluations ${\rm ev}_c(f)$ are fibrations in $\Mc^{W_c}$ for every $c \in \Cc$\,,
\end{enumerate}
\end{theorem}

By Remark \ref{res} above, the forgetful functor $\Mc^{\Cc}_W \to \Mc^{\Cc}$ preserves weak equivalences (resp., fibrations). Hence,
\begin{cor} \la{identification} 
The equivalence of categories $($hence, adjunction$)$
$$ 
\iota^\ast\,:\,\Mc^{\Cc_{hW}} \rightleftarrows \Mc^{\Cc}_W \,:\,\bar{(\,\mbox{--}\,)} 
$$
is a Quillen pair, where $\Mc^{\Cc_{hW}}$ is equipped with the projective model structure.
\end{cor}

\begin{remark} \la{naive1} If the categories $\Mc^H$ are equipped with the projective model structure (see Remark \ref{naive}), then the weak equivalences (resp., fibrations) in $\Mc^{\Cc}_W$ are indeed the morphisms that are taken by the forgetful functor $\Mc^{\Cc}_W \to \Mc^{\Cc}$ to objectwise weak equivalences (resp., fibrations). 
\end{remark}

Let $K\,:\,\Cc \to \sset$ and $L:\Cc^{\rm op} \to \sset$ be $W$-diagrams. Let $X\,:\,\Cc \to \Mc$ be a $W$-diagram. As explained in \cite[Sec. 1.2]{DM16}, the object ${\rm Map}_{\Cc}(K,X) \in \Mc$, defined as the equalizer 
\begin{equation*}
    \begin{tikzcd}[column sep=1cm] 
 {\rm Map}_{{\Cc}}(K,X) \rightarrowtail  \prod_{c \in \Cc}{\rm Map}_{\Mc}(K_c,X_c) \arrow[r,"s",yshift=2pt] \arrow[r,"t"',yshift=-2pt] & \prod_{\alpha:c \to d}{\rm Map}_{\Mc}(K_c,X_d)
    \end{tikzcd}
\end{equation*} 
carries a (explicit) $W$-action. We denote the resulting object in $\Mc^W$ by ${\rm Map}^a_{{\Cc}}(K,X)$. Similarly, the object 
$L \otimes_{\Cc} X$, defined as the coequalizer
\begin{equation*}
    \begin{tikzcd}[column sep=1cm] 
 \coprod_{\alpha:c \to d} L(d) \otimes X(c) \arrow[r,"s",yshift=2pt] \arrow[r,"t"',yshift=-2pt] & \coprod_{c \in \Cc} L(c) \otimes X(c) \twoheadrightarrow L \otimes_{\Cc} X
    \end{tikzcd}
\end{equation*} 
carries an explicit $W$-action. We denote the resulting object in $\Mc^W$ by $L \otimes^a_{\Cc} X$. Note that the $\sset$-valued functors 
$$N(\Cc/\mbox{--})\,:\,\Cc \to \sset\,,\qquad N(\mbox{--}/\Cc):\Cc^{\rm op} \to \sset$$
are $W$-functors for any $W$-category $\Cc$. 

\begin{defi} \la{Whocolim}
$(i)$ The {\it $W$-homotopy colimit}: $\,\hocolim_{\Cc,W}: \Mc^{\Cc}_W \to \Mc^W$ is the 
functor
$$X \mapsto  N(\mbox{--}/\Cc) \otimes_{\Cc}^a (Q \circ X)\,,$$
where $Q:\Mc \to \Mc$ is a(ny) cofibrant replacement functor.\\

$(ii)$ The {\it $W$-homotopy limit}: $\,\holim_{\Cc,W}: \Mc^{\Cc}_W \to \Mc^W$ is the 
functor
$$X \mapsto  {\rm Map}_{\Cc}^a\big(N(\Cc/\,\mbox{--}), (R\circ X)\big)\,,$$
where $R:\Mc \to \Mc$ is a(ny) fibrant replacement functor.
\end{defi}

By \cite[Corollary 2.22]{DM16}, the functors $\,\hocolim_{\Cc,W}: \Mc^{\Cc}_W \to \Mc^W$ and  $\,\holim_{\Cc,W}: \Mc^{\Cc}_W \to \Mc^W$ preserve weak equivalences. They therefore, descend to functors $\hocolim_{\Cc,W}:\Ho(\Mc^{\Cc}_W) \to \Ho(\Mc^W)$ and $\holim_{\Cc,W} :\Ho(\Mc^{\Cc}_W) \to \Ho(\Mc^W)$ on the corresponding homotopy categories. 
\begin{prop} \la{Whob} 
The $W$-homotopy colimit and limit functors restrict to the ordinary homotopy colimit and limit respectively: i.e., there are
a commutative diagrams
$$
\begin{diagram}[small]
\Ho(\Mc^{\Cc}_W) &\rTo^{\hocolim_{\Cc,W}}& \Ho(\Mc^W) & & & \Ho(\Mc^{\Cc}_W) &\rTo^{\holim_{\Cc,W}}& \Ho(\Mc^W)\\
 \dTo & & \dTo & & &\dTo & & \dTo\\
 \Ho(\Mc^{\Cc}) & \rTo^{\hocolim_{\Cc}} & \Ho(\Mc) & & & \Ho(\Mc^{\Cc}) & \rTo^{\holim_{\Cc}} & \Ho(\Mc)\\
\end{diagram} 
$$
where the vertical arrows are the forgetful functors.
Consequently, for any $W$-functor $F:\Cc \to \Mc$, $\hocolim_{\Cc}(F) $ and $\holim_{\Cc}(F)$ come equipped with  $($homotopy$)$ $W$-actions.
\end{prop}
\begin{proof}
This follows immediately from Definition \ref{Whocolim} and the fact that the Bousfield-Kan constructions represent the ordinary homotopy colimit and limit respectively (see \cite[Corollary 5.13]{Rie14}). 
\end{proof}
\begin{lemma} 
\la{lem:Thom}
 Let  $\Mc$ be a $W$-model category, and let $F:\Cc \to \Mc$ be a $W$-functor. Then, 
 $$ {\rm hocolim}_{\Cc_{hW}}(\tilde{F}) \,\cong\, ({\rm hocolim}_{\Cc}\, F)_{hW} \,,\qquad {\rm holim}_{\Cc_{hW}}(\tilde{F}) \,\cong\,({\rm holim}_\Cc\, F)^{hW}\ .$$
\end{lemma}
\begin{proof}
 This follows from Thomason's Theorem for homotopy colimits (see, e.g., \cite[Prop. 26.5, Theorem 26.8]{ChS02}).
\end{proof}
\subsection{Very small categories} \la{AA3}
We illustrate the above constructions on diagrams indexed by very small categories. Recall (see \cite[Section 10.13]{DS95})
\begin{defi} \la{vscat}
A category $\I$ is called {\it very small} if it has finitely many objects, finitely many morphisms, and
there exists an integer $N$ such that for any string of composable morphisms  in $\I$
$$ 
i_0 \stackrel{f_1}{\to} i_1 \to \ldots \stackrel{f_n}{\to} i_n 
$$
\noindent
of length $ n > N$ some $f_i$ is the identity morphism. 
\end{defi}
If $\I$ is a very small category, then the category of $\I$-diagrams $\Mc^{\I}$ 
with values in an {\it arbitrary} model category $ \Mc $ can be given two distinct model 
structures --- called projective and injective --- which share the same weak equivalences
(and hence have isomorphic homotopy categories) but differ in their fibrations and cofibrations.
In the projective model structure, the cofibrations in $ \Mc^{\I}$ have the following 
description ({\it cf.} \cite[Proposition~10.6]{DS95}). For every $i \in \I$, let $\partial i$ denote the full subcategory of $\I/i$ generated by all objects {\it except} the identity on $i$. Let $u_i: \partial i \to \I$ denote the functor $(j \to i) \mapsto j$. For any $X \in \Mc^\I$, let $\partial i(X):=\colim_{\partial i} u_i^{\ast}X$. There is a natural map $\partial_i(X) \to X(i)$ in $\Mc$ for all $i \in \I$. For any morphism $f: X \to Y$ in $\Mc^\I$ let $\partial_i(f):=\colim(\partial_i(Y) \leftarrow \partial_i(X) \to X)$. For all $i \in \I$, there is a natural map $\partial_i(f) \to Y(i)$. A morphism $f: X \to Y$ in $\Mc^\I$ is a cofibration if and only if for every $i \in \I$ the natural map $\partial_i(f) \to Y(i)$ is a cofibration in $\Mc$. In the injective model structure on $\Mc^{\I}$, the fibrations have a dual description ({\it cf.} \cite[Proposition~10.11]{DS95}). We remark that the projective model structure on $ \Mc^{\I} $
is adapted to constructing homotopy colimits, while the injective one to constructing homotopy limits.

\begin{lemma} \la{vsdiag} Let $\I$ be a very small diagram equipped with a $W$-action. Let $F: \I \to \Mc$ be a $W$-functor such that $F$ is cofibrant with respect to the projective model structure on $\Mc^\I$. Then the natural map in $\Mc^W$
$$\hocolim_{\I,W}(F) \to \colim_{\I,W}(F) $$
is a weak equivalence. 
\end{lemma}
\begin{proof}
There are weak equivalences in $\Mc$
\begin{equation} \la{comp1} \hocolim_{\I,W}(F) \stackrel{\sim}{\to} \hocolim_\I(F) \stackrel{\sim}{\to} \colim_\I(F) \ .\end{equation}
The first weak equivalence above is by Proposition \ref{Whob} and the second weak-equivalence is because $F$ is cofibrant as an $\I$-diagram. Since the composition of the weak equivalences in \eqref{comp1} coincides with the composition 
$$ \hocolim_{\I,W}(F) \to \colim_{\I,W}(F) \stackrel{\cong}{\to} \colim_{\I}(F)\,,$$
in which the second arrow is an isomorphism in $\Mc$, the map $\hocolim_{\I,W} F \to \colim_{\I,W} F$ is a morphism in $\Mc^W$ that is a weak equivalence in $\Mc$. Hence, it is a weak equivalence in $\Mc^W$.
\end{proof}
Dually, we have
\begin{lemma} \la{vsdiagdual}
Let $\I$ be a very small diagram equipped with a $W$-action. Let $F: \I \to \Mc$ be a $W$-functor such that $F$ is fibrant with respect to the injective model structure on $\Mc^\I$. Then the natural map in $\Mc^W$
$${\rm lim}_{\I,W}(F) \to \holim_{\I,W}(F) $$
is a weak equivalence. 
\end{lemma}
\begin{example} \la{star}
Let $S$ be a finite $W$-set. Let $\I$ be the generalized pushout category with a $W$-action whose set of objects is $S \sqcup \{0\}$ with the $W$-action on objects extending the action on $S$ by fixing the object $0$. The only non-identity morphisms in $\I$ are unique maps $f_s:0 \to s$ for each $s \in S$. We define the $W$-action on morphisms by $\hat{g}(f_s):=f_{gs}: 0 \to g \cdot s$ for all $g \in W, s \in S$. 

Let $F: \I \to \Mc$ be a $W$-functor. If $Q \stackrel{\sim}{\to} F$ is a morphism of $W$-diagrams $\I \to \Mc$ that is a cofibrant resolution with respect to the projective model structure on $\Mc^\I$, then Lemma \ref{vsdiag} implies that there are weak equivalences in $\Mc^W$
\begin{equation*}
\hocolim_{\I,W}(F) \xleftarrow{\sim} \hocolim_{\I,W}(Q) \simeq 
\colim_{\I,W}(Q) \simeq 
\colim\bigl(Q(0) \leftarrow \coprod_{s \in S} Q(0) \to \coprod_{s \in S} Q(s)\bigr) \,,
\end{equation*}
where the action of each $g \in W$ on $Q(0)$ is given by $\vartheta_g:Q(0) \to Q(0)$, the action of $g \in W$ on $\coprod_{s \in S} Q(0)$ is given by $\vartheta_g:Q(0)_s \to Q(0)_{g\cdot s}$, where $Q(0)_s$ denotes the copy of $Q(0)$ in $\coprod_{s \in S} Q(0)$ indexed by $s$, and the action of $g \in W$ on $\coprod_{s \in S} Q(s)$ is given by $\vartheta_g: Q(s) \to Q(g \cdot s)$. Here the natural transforms $\vartheta_g: Q \to Q \circ \hat{g}, g \in W$, are part of the structure of $Q$ as a $W$-functor.
 \end{example}
\begin{example} \la{starlim}
Let $S$ be a finite $W$-set. Let $\I$ be the generalized pullback category with a $W$-action whose set of objects is $S \sqcup \{0\}$ with the $W$-action on objects extending the action on $S$ by fixing the object $0$. The only non-identity morphisms in $\I$ are unique maps $f_s:s \to 0$ for each $s \in S$. We define the $W$-action on morphisms by $\hat{g}(f_s):=f_{gs}:  g \cdot s \to 0$ for all $g \in W, s \in S$. 

Let $F: \I \to \Mc$ be a $W$-functor. If $F \stackrel{\sim}{\to} R$ is a morphism of $W$-diagrams $\I \to \Mc$ that is a fibrant resolution with respect to the injective model structure on $\Mc^\I$, then Lemma \ref{vsdiagdual} implies that there is a weak equivalence in $\Mc^W$
$$ \holim_{\I,W}(F) \,\cong\, \lim\biggl(R(0) \rightarrow \prod_{s \in S} R(0) \leftarrow \prod_{s \in S} R(s)\biggr)  $$
where the $W$ actions on $R(0)$, $\prod_{s \in S} R(0)$ and $\prod_{s \in S} R(s)$ are dual to the analogous $W$-actions in example \ref{star} above.
\end{example}

\section{Coaffine Stacks and Rational Homotopy Theory}
\la{AB}
Rational homotopy theory describes homotopy types of topological spaces up to rational equivalence (i.e., 
disregarding torsion in cohomology and homotopy groups). There are two classical approaches to the subject differing in the choice of algebraic models: one -- due to D. Quillen \cite{Qui69} --- makes use of differential graded chain Lie algebras, and the other  -- due to D. Sullivan \cite{Su77} ---  commutative cochain algebras and polynomial differential forms. In recent years, both approaches have been  incorporated into the framework of Derived Algebraic Geometry (DAG) which led to some spectacular advances in homotopy theory (see, e.g., \cite{Heu21} and \cite{Yuan23}). 

In this Appendix, we briefly outline the DAG interpretation of Sullivan's construction in terms of 
{\it coaffine stacks}. This notion was originally introduced by B. To\"en \cite{To06}, and developed further --- in the framework of spectral algebraic geometry --- by J. Lurie (see \cite{DAGVIII}, \cite{DAGXIII}). 
The natural language for Lurie's approach is that of $\infty$-categories which we use freely in this Appendix. 
For all unexplained notation and terminology related to $\infty$-categories
we refer the reader to \cite{HTT}.


\subsection{$\bE_\infty$-algebras}
\la{A1}
Recall that a {\it spectrum} is a sequence of pointed topological spaces $ X = (X_0, X_1, \ldots)$ equipped with (weak) homotopy equivalences $ X_n \xrightarrow{\sim} \Omega X_{n+1}$, one for each $ n \ge 0 $, where $ \Omega $ stands for the based loop functor. The spectra form an $\infty$-category $ \Sp $, where the
1-morphisms are represented by maps of pointed spaces commuting with equivalences up to coherent higher homotopies. 
To each  $ X \in \Sp $ one can associate a sequence of abelian groups  $ \{\pi_i(X)\}_{i \in \Z}$ --- called the (stable) homotopy groups of $X$ --- by the formula $ \pi_i(X) := \lim_{n\to \infty}[\pi_{i+n}(X_n)]$. A spectrum $X$ is then called {\it connective} if $ \pi_i(X) = 0 $ for $ i < 0 $ and {\it discrete} if $ \pi_i(X) = 0 $ for all $ i \not=0 $. We write $ \Sp^{\rm cn}$ and $ \Sp^0$
for the full $\infty$-subcategories of connective and discrete spectra in $\Sp$, respectively.
The assignment $ X \mapsto \pi_0(X) $  yields an equivalence of $\infty$-categories $ \Sp^0 \xrightarrow{\sim} {\mathscr Ab} $, where $ {\mathscr Ab}$ is the ordinary category of abelian groups (viewed as an $\infty$-category via the nerve construction). Thus, $ \Sp $ can be thought of as a topological analogue (refinement) of $ {\mathscr Ab} $. Many natural algebraic structures on $ {\mathscr Ab} $ `lift' to $ \Sp $. In particular,
$ \Sp $ is a symmetric monoidal $\infty$-category relative to the smash product  $ \wedge: \Sp \times \Sp \to \Sp  $ that `lifts' the usual tensor product $\otimes$ on ${\mathscr Ab}$ (in the sense that  $\,\pi_0(X \wedge Y) \cong \pi_0(X) \otimes \pi_0(Y) \,$ for all  $\, X, Y \in \Sp^{\rm cn}\,$). The unit object in $ (\Sp, \wedge) $
 is the classical sphere spectrum $ \bS = (S^0, S^1, \ldots)\,$.

Now, an {\it $ \bE_\infty$-ring}  can be defined as a commutative monoid (i.e., commutative ring object) 
in  $ (\Sp, \wedge, \bS) $. We denote the $\infty$-category of such objects by $ \CAlg(\Sp)$. Thus, $ R \in \CAlg(\Sp) $ is a spectrum equipped with multiplication $ m: R \wedge R \to R $ which is unital, associative and commutative up to coherent homotopy. The functor $ \pi_0: \Sp \to {\mathscr Ab} $ restricts to commutative monoids $ \pi_0: \CAlg(\Sp) \to \CAlg({\mathscr Ab}) $, where $ \CAlg({\mathscr Ab}) \cong \Comm $ is (the nerve of) the ordinary category of commutative rings. Naturally, this becomes an equivalence $ \CAlg^0(\Sp) \simeq \Comm $, when we further restrict to the subcategory of discrete $ \bE_\infty$-rings in $ \CAlg(\Sp) $. 

Let $ \k$ be a field. Regarding $\k$ as a  discrete $ \bE_\infty$-ring, we define an {\it $ \bE_\infty$-algebra over $\k$} to be an $ \bE_\infty$-ring $A$ equipped with a morphism $ \k \to A $ in $ \CAlg(\Sp) $ called the unit map. The $ \bE_\infty$-algebras form the coslice $\infty$-category $ \CAlg(\Sp)_{\k/} $ under $\CAlg(\Sp)$, which we denote by $ \CAlg_{\k} $. Naturally, we say that an $ \bE_\infty$-algebra $ A \in \CAlg_\k $ is connective (resp., discrete) if its underlying spectrum is connective (resp., discrete). Furthermore, following \cite[4.1.1]{DAGVIII}, we make the following

\begin{defi}
\la{coconn}
An $ \bE_\infty$-algebra $ A \in \CAlg_\k $ is called {\it coconnective} if its unit map $ \k \to A $ induces an
isomorphism $ \k \xrightarrow{\sim} \pi_0(A) $ in $ {\mathscr Ab}$, and $ \pi_i(A) = 0 $ for all $ i > 0 $.
\end{defi}
As in \cite{DAGVIII}, we write $ \CAlg_\k^{\rm cn}$, $\,\CAlg^0_\k $ and  $ \CAlg_\k^{\rm cc}$ for the full $\infty$-subcategories of $ \CAlg_\k $ spanned by the connective, discrete and coconnective $ \bE_\infty$-algebras, respectively.

\subsection{Affine and coaffine stacks} \la{A2} In spectral algebraic geometry, the main objects of study are {\it derived stacks} represented by functors $\, X: \CAlg^{\rm cn}(\Sp) \to \sS \,$ on the {\it connective} $ \bE_\infty$-rings with values in the $\infty$-category of spaces. Such functors are required to satisfy descent conditions on flat \v{C}ech hypercovers of $ \bE_\infty$-rings (i.e. they behave as hypercomplete sheaves on a flat Grothendieck topology on $ \CAlg^{\rm cn}(\Sp)^{\rm op}$). We  restrict our consideration to
functors $\, X: \CAlg_\k^{\rm cn} \to \sS \,$ defined on the connective $ \bE_\infty$-algebras over a fixed field $\k$. 
We denote the $\infty$-category of such functors by $ \Fun(\CAlg_k^{\rm cn}, \,\sS)$. 

The basic examples of derived stacks are the {\it derived affine schemes} represented by
functors of the form
\begin{equation}
\la{RSpec}
\RSpec(A):\, \CAlg^{\rm cn}_\k \to \sS\ ,\quad B \mapsto \Map_{\CAlg_\k}(A,\,B)
\end{equation}
where $A$ is a connective $ \bE_\infty$-algebra.
By the ($\infty$-categorical version of) Yoneda Lemma (see \cite[5.1.3.1]{HTT}), the assignment
$ A \mapsto \RSpec(A) $ extends to a fully faithful functor 
$\,(\CAlg^{\rm cn}_\k)^{\rm op}\,\into\, \Fun(\CAlg_\k^{\rm cn}, \,\sS) $. This means that $A$ is determined
by $ \RSpec(A) $ uniquely, up to a contractible space of choices, and the $\infty$-category $\dAff_\k$
of derived affine $\k$-schemes can thus be identified with the opposite $\infty$-category of connective $ \bE_\infty$-algebras
\begin{equation}
\la{daff}
\dAff_\k \simeq (\CAlg^{\rm cn}_\k)^{\rm op}
\end{equation}

To define coaffine stacks we introduce the notation: for a functor $\, X: \CAlg^{\rm cn}_\k \to \sS \,$ 
we let $ X_0: \CAlg^0_\k \to \sS  $ denote its restriction to the subcategory of discrete $ \bE_\infty$-algebras,
i.e. $ X_0 := i^* X$, where $ i: \CAlg^0_\k \into \CAlg^{\rm cn}_\k $ is the natural inclusion. We note that
the restriction functor $\,i^*\,$ has a left adjoint
$i_!: \,\Fun(\CAlg^0_\k,\,\sS) \to \Fun(\CAlg^{\rm cn}_\k,\,\sS)\,$ given by the left Kan extension along $i$.
\begin{defi}
\la{coaff}
A functor $ X: \CAlg^{\rm cn}_\k \to \sS $ is a {\it coaffine stack} if its satisfies the properties
\begin{enumerate}
    \item[(a)] $X$ is the left Kan extension of $ X_0 $, i.e.  $\,X \cong i_! \,i^* X \,$.
    \item[(b)] $ X_0(B) $ is connected, i.e $ \pi_0[X_0(B)] \cong \ast \,$, for all $ B \in \CAlg_\k^0$.
    \item[(c)] For each $ n > 0 $, the functor $\,\pi_n(X_0(\mbox{--}),\,\ast): \CAlg_\k^0 \to {\rm Grp} \,$ is representable by a prounipotent affine group scheme over $\k$.
\end{enumerate}
\end{defi}
Property $(a)$ shows that the coaffine stacks are simpler objects than the derived affine 
schemes since as functors, they are completely determined by values on discrete 
$ \bE_\infty$-algebras. Because the discrete $ \bE_\infty$-algebras are equivalent to the ordinary commutative $\k$-algebras, this can be also viewed as a justification for omitting the adjective `derived' in Definition~\ref{coaff}. On the other hand, if $\k$ is a field of characteristic of zero, the category of coaffine stacks, $\cAff_\k$, admits a universal characterization (see Theorem~\ref{ToTh} below) which 
shows that it is, in a sense, dual (or complimentary) to $\dAff_\k$. This justifies the adjective "coaffine" in Definition~\ref{coaff}.

To state Theorem~\ref{ToTh} we observe that the functor $  A \mapsto \RSpec(A) $ defined by \eqref{RSpec} extends naturally to all (i.e., not necessarily connective) $ \bE_\infty$-algebras $A \in \CAlg_\k$.
%
%
This extended Yoneda functor is not fully faithful; however, it becomes so when restricted to
the {\it coconnective} $ \bE_\infty$-algebras. More precisely, we have the following basic observation,  originally due to B. To\"en \cite{To06} and --- in the  topological context --- due to J. Lurie 
({\it cf.} \cite[Theorem 4.4.1, Proposition 4.4.6, Proposition 4.4.8]{DAGVIII}). 
\begin{theorem}[To\"en, Lurie]
\la{ToTh}
Assume that $\k$ is a field of characteristic zero. Then the restricted Yoneda functor
\begin{equation}
\la{coSpec}
(\CAlg^{\rm cc}_\k)^{\rm op} \,\into \, (\CAlg_\k)^{\rm op} \xrightarrow{\RSpec} \Fun(\CAlg_\k^{\rm cn}, \,\sS)
\end{equation}
is fully faithful. The essential image of \eqref{coSpec} is precisely the $\infty$-category $\cAff_\k$ of coaffine stacks. 
\end{theorem}
Theorem~\ref{ToTh} says that, over a field $\k$ of characteristic zero,
a functor $\, X: \CAlg^{\rm cn}_\k \to \sS \,$ is a coaffine stack if and only if it is 
corepresentable by a coconnective $ \bE_\infty$-algebra. Hence, we have an equivalence
of $\infty$-categories
\begin{equation}
\la{caff}
\cAff_\k \simeq (\CAlg^{\rm cc}_\k)^{\rm op}
\end{equation}
This implies, in particular, that 
the functors $X$ satisfying the conditions of Definition~\ref{coaff} are indeed derived stacks as 
they are automatically hypercomplete as presheaves with respect to the flat topology on $\dAff_k $. 
Modifying slightly the notation of \cite{DAGVIII}, we denote the functor 
\eqref{coSpec} by $\,\Specc\,$: thus, every coconnective algebra $A$ determines and is determined 
by a coaffine stack $\Specc(A)$, which is given by the same formula as in \eqref{RSpec}. 

\subsection{Sullivan-To\"en-Lurie Theorem}
\la{Sullivan}
We say that a coaffine stack $ X = \Specc(A)$ is {\it simply connected} if $ \pi_{-1}(A) = 0 $.
The simply connected coaffine stacks can be characterized by the three 
conditions of Definition~\ref{coaff} with additional assumptions that
in (b), all spaces $X_0(B)$ are simply connected, and in (c), the functors 
$\,B \mapsto \pi_n[X_0(B)]\,$ for $ n \ge 2 $ are represented by {\it commutative} prounipotent
group $k$-schemes, i.e. are isomorphic to linear functors of the form $\, B \mapsto \Hom_{\k}(V_n, B) \,$, where $V_n$ are some vector spaces over $\k$ ({\it cf.} \cite[Remark 1.2.5]{DAGXIII}).
We say then that a simply connected coaffine stack $X$ is {\it of finite type} if the vector spaces 
$V_n$ are finite-dimensional for all $ n\ge 2 $.
We write $ \cAff_\k^{1, \,{\rm ft}} $ for the full $\infty$-subcategory of $ \cAff_\k$ consisting
of simply connected coaffine stacks of finite type over $\k$.

The next result is an updated version of a classical theorem of D. Sullivan~\cite{Su77} that
describes the category $ \sS^{1, {\rm ft}}_{\Q}  $ of  simply connected rational spaces
of finite type  in terms of cochain models.

\begin{theorem}[Sullivan, To\"en, Lurie]
\la{SuTh}
Let $ \k = \Q $. The natural evaluation functor
\begin{equation}\la{eval}
\Fun(\CAlg_\k^{\rm cn},\,\sS)\,\to\, \sS\ ,\quad X \mapsto X(\Q)\,,
\end{equation}
restricts to an equivalence of $\infty$-categories $\,\cAff_\k^{1, \,{\rm ft}} \xrightarrow{\sim} \sS^{1, {\rm ft}}_{\Q}$ given by
$$
\Specc(A) \,\mapsto \, \Map_{\CAlg_\k}(A, \Q)\,.
$$
The inverse functor $ \sS^{1, {\rm ft}}_{\Q} \to \cAff_\k^{1, \,{\rm ft}} $ can be written in the form
$$
Z \,\mapsto \,\Specc\,[C_{\Q}^*(Z)]\,,
$$
where $ C_{\k}^*(Z) $ is the classical cochain complex of a space $Z$ with coefficients in $\k$ viewed as an $ \bE_\infty$-algebra over $\k$.
\end{theorem}
\begin{remark}
In the above form, Sullivan's Theorem appears in \cite{DAGXIII} as Theorem 1.3.5.
In fact, this result is extended in \cite{DAGXIII} to all simply connected spaces. 
To this end, Lurie specifies a more general class of functors 
$X \in \Fun(\CAlg_\k^{\rm cn},\,\sS) $ than simply connected coaffine stacks that he refers to as $k$-rational homotopy 
types. \cite[Theorem 1.3.6]{DAGXIII} asserts that the evaluation functor $ X \mapsto X(\Q) $ 
provides an equivalence $ \rm RType(\Q) \simeq \sS^{1}_{\Q} $ between the $\infty$-category of
$\Q$-rational homotopy types and the $\infty$-category of all simply connected rational spaces,
without any finiteness assumptions.
\end{remark}

\subsection{DG models}
\la{DGModel}
We now explain the relation of the above construction to the classical approach in rational homotopy theory based on simplicial polynomial differential forms, due to Bousfield and Gugenheim \cite{BG76}. 
This will allow us to produce explicit models for coaffine stacks that we use in the main body of the paper.

If $\k$ is a commutative ring, we let $ \Com_k $ denote the category of (unbounded) cochain complexes over $\k$, and set $ \Mod_\k := \Com_\k[W^{-1}] $ to be the $\infty$-category obtained from $ \Com_\k $ by formally inverting the class $W$ of quasi-isomorphisms in $ \Com_k $. The last notation is motivated by the well-known fact  (see, e.g., \cite[Appendix B]{SS03})
that $\Mod_k$ can be identified with the $\infty$-category $\Mod_{k}(\Sp)$ of $k$-module spectra, via the classical Eilenberg-MacLane construction. Under this identification, the canonical functor $ \Com_\k \to \Mod_\k $ is symmetric monoidal, hence restricts to a functor on commutative monoids: 
$\, \CAlg(\Com_\k) \to \CAlg(\Mod_\k) \,$. By definition, $ \CAlg(\Com_k) $ is the category of commutative cochain dg algebras, while $\CAlg(\Mod_\k) \simeq \CAlg_k $ is the
$\infty$-category of $\bE_{\infty}$-algebras over $k$. Now, if $\k$ is a field 
of characteristic zero, it is well known that the above functor induces an equivalence of $\infty$-categories 
\begin{equation}
\la{DG}
\cdga_\k \simeq \CAlg_\k   
\end{equation}
where $ \cdga_\k = \CAlg(\Com_\k)[W^{-1}]$ stands for the $\infty$-category of commutative 
cochain dg algebras with quasi-isomorphisms inverted. Thus, when working over a field of characteristic zero, we can replace the $\bE_{\infty}$-algebras by the usual commutative cochain dg algebras. 

Under the equivalence \eqref{DG}, the connective $\bE_{\infty}$-algebras correspond to
the dg algebras $A \in  \cdga_\k^{\le 0} $ with cohomology $H^i(A)$ vanishing in all (strictly) positive degrees $i>0$, while the coconnective ones correspond to the dg algebras $ A \in \cdga_\k^{>0} $ with $H^i(A)=0 $ for all 
$ i < 0$ satisfying, in addition, the condition $ H^0(A) \cong \k\,$.
%
%
We write $ \cdga_\k^{\le 0} $ and $ \cdga_\k^{>0} $ for the full $\infty$-subcategories of $ \cdga_\k $
spanned by the dg algebras with above constraints on cohomology: thus, $\, \CAlg^{\rm cn}_\k \simeq \cdga^{\le 0}_\k $ and $ \CAlg^{\rm cc}_\k \simeq \cdga^{>0}_\k $.
Combining these equivalences with \eqref{daff} and \eqref{caff}, we conclude  
\begin{equation}
\la{DGcc}
\dAff_\k \simeq (\cdga^{\le 0}_\k)^{\rm op} \ ,\quad
\cAff_\k \simeq (\cdga^{>0}_\k)^{\rm op} 
\end{equation}

We will deal with various geometric constructions with coaffine stacks involving limits and colimits
in the $\infty$-category $ \cAff_\k $. To handle these effectively we `embed' $\cAff_\k $ into a convenient model category. Recall (see \cite[Appendix A.2]{HTT}) that every simplicial model 
category $ \M $ determines an $\infty$-category $ \underline{\M} := \mathcal{N}_{\Delta}(\M^{\rm cf}) $ called the {\it underlying $\infty$-category of $\M$}. Here $ \M^{\rm cf} $ is the full simplicial subcategory of $\M$ spanned by the fibrant-cofibrant objects (this simplicial category is always fibrant in the sense that its mapping spaces are Kan complexes); $ \mathcal{N}_{\Delta} $ denotes the homotopy coherent nerve functor. Not every $\infty$-category $\Cc$ arises 
this way; however, for every $\Cc$, there is a fully faithful embedding  $ \Cc \into \underline{\M} $ for an appropriately chosen simplicial category $ \M $.

The equivalence \eqref{DGcc} suggests that for the $\infty$-category $ \Cc = \cAff_\k $, we take $ \M $
to be the opposite of the category $ \ccdga_\k^{\ge 0} $ of non-negatively graded commutative cochain dg algebras.
This category carries a simplicial model structure (originally constructed in \cite{BG76}), where the weak equivalences are the quasi-isomorphisms and the fibrations are all degreewise surjective maps of dg algebras. This model structure is known to be cofibrantly generated (in fact, combinatorial), with all objects being fibrant and the cofibrant objects being the retracts of so-called KS-complexes, aka Sullivan algebras (see, e.g., \cite[Section 1.2]{Hess07}). Now, the mapping spaces in $ \ccdga_\k^{\ge 0} $ 
are defined in terms of simplicial polynomial differential forms as follows ({\it cf.} \cite[\S 5]{BG76}).

For each $n\ge 0$, let $ \dDel^n_{\k} \,$ denote the algebraic $n$-simplex over $\k$, which is the hyperplane $\, \sum_{i=0}^n t_i = 1 \,$ embedded in the affine $\k$-space $ \bA^{n+1}_\k$ with coordinates $(t_0,\ldots, t_n)$. Define $\,\Omega(\dDel^n_\k) := \Omega^{\ast}_{\rm dR}(\dDel^n_\k) $ to be the algebraic de Rham complex on $ \dDel^n_{\k} \,$: explicitly,
\begin{equation}
\la{Omegan} 
\Omega(\dDel^n_\k)\,\cong\, \frac{\Lambda_{\k}(t_0,\ldots,t_n,dt_0,\ldots,dt_n)}{(\sum_{i=0}^n t_i-1,\,\sum_{i=0}^n dt_i)} \,,
\end{equation}
where $\deg(t_i)=0$ and $\deg(dt_i)=1$, and the differential $d$ on $ \Omega(\dDel^n_\k) $ is defined on generators by $t_i \mapsto  dt_i$ and $dt_i \mapsto 0$ for all $0 \leq i \leq n$. The assignment $ [n] \mapsto
\Omega(\dDel^n_\k) $
makes $ \Omega_{\ast} := \{\Omega(\dDel^{n}_\k)\}_{n \ge 0} $  a simplicial cochain dg algebra,
i.e. a simplicial object in $\ccdga_\k^{\ge 0} $, with face and degeneracy maps determined  by
$$\partial_i\,:\,\Omega_n \to \Omega_{n-1}\,,\,\, t_{k} \mapsto 
\begin{cases}
 t_{k} & k<i\\
 0 & k=i\\
 t_{k-1} & k>i\\
\end{cases}\,\qquad s_i\,:\,\Omega_n \to \Omega_{n+1}\,,\,\, t_{k} \mapsto \begin{cases}
    t_{k} & k<i\\
    t_{k}+t_{k+1} & k=i\\
    t_{k+1} & k>i\\
\end{cases}\ .
$$
Now, for any objects $ A, B \in \ccdga_\k^{\ge 0} $, the simplicial mapping space from $A$ to $B$  in $\ccdga_\k^{\ge 0} $ is defined by
\begin{equation} 
\la{smap}
\Map(A,B)_{\ast} \,:=\, \Hom(A,\,B \otimes \Omega_*) 
\end{equation}
where $ \otimes = \otimes_\k $ and `$\Hom $' stand for the tensor product and the set of dg algebra homomorphisms between cochain algebras over $\k$. 

We write $ \cdga_\k^{\ge 0} := \underline{\ccdga}_\k^{\ge 0} $ for the underlying $\infty$-category of 
the simplicial model category $\ccdga_\k^{\ge 0} $ described above. By \eqref{DGcc}, we then identify
$\cAff_\k = (\cdga_\k^{>0})^{\rm op} $ as a full $\infty$-subcategory of $ (\cdga_\k^{\ge 0})^{\rm op} $ spanned by the cochain algebras $A$ with $ H^0(A) \cong \k $, and we will keep the notation $ \Specc(A) $ for the
coaffine stack associated to such an algebra $A$. Under this identification, the evaluation functor \eqref{eval}
corresponds to the {\it  Sullivan's realization functor}  
\begin{equation} 
\la{defsr} 
\langle\,\mbox{--}\,\rangle:\, (\ccdga_\k^{\ge 0})^{\rm op}\,\to\,\sS\ ,\quad A \mapsto  
|\Hom(A,\,\Omega_*)| 
\end{equation}
where  $\,|\,\mbox{--}\,|\,$ stands for the usual geometric realization of simplicial sets.
At the simplicial level, the functor \eqref{defsr} has a left adjoint $\, A_{\rm PL}: \,\sset \to  (\ccdga_\k^{\ge 0})^{\rm op}\,$ assigning to a simplicial set $X_* \in \sset$ the cochain 
algebra of so-called {\it piecewise linear} (PL) differential forms
\begin{equation} 
\la{defapl}
A_{\rm PL}(X)\,:=\,\Hom_{\sset}(X,\Omega_\ast) 
\end{equation}
with product and differential defined objectwise. The main theorem of \cite{BG76} asserts that 
\begin{equation} 
\la{SPL} 
A_{\rm PL}\,:\,\sset\, \rightleftarrows\, (\ccdga_\k^{\geqslant 0})^{\rm op}\,:\,\langle\, \mbox{--} \,\rangle \ .\end{equation}
is a Quillen pair of simplicial model structures. The functor $ A_{\rm PL} $ 
is known to be equivalent to the usual (simplicial) cochain functor over field 
$\k$: namely, for all $ X \in \sset $,  there is a quasi-isomorphism of DG algebras
(see, e.g., \cite[Theorem 10.9]{FHT}):
\begin{equation}
\la{CAPL}
C_{\k}^*(X) \,\simeq \,A_{\rm PL}(X)   
\end{equation}
Theorem~\ref{SuTh} can then be restated
by saying that \eqref{SPL} restricts to simply connected objects of finite type, inducing an equivalence of $\infty$-categories:
$\,\sset_{\Q}^{1, {\rm ft}}\,\simeq\, \cAff_{\Q}^{1, \,{\rm ft}}$.
\subsection{Limits and colimits}\la{B5}
Recall that if $ \cF: \Cc \to \Dc $ is a diagram (i.e., a functor) between two $\infty$-categories, with $ \Cc$ being small, the {\it limit} of $\cF$ is an object $\, \blim_{\Cc}(\cF) $ in $ \Dc $ that represents the vertex of a cone $ \, \blim_{\Cc}(\cF) \to \cF $ over $\cF$, which is universal (homotopically terminal) in the sense that, for all  objects $ D \in \Dc $, it induces a homotopy equivalence of mapping spaces
$$
\Map_{\Dc}(D,\,\blim_{\Cc} \,\cF) \xrightarrow{\sim} \Map_{\Fun(\Cc, \Dc)}(D, \cF)\,,
$$
where $ D $ on the right is identified with 
its constant diagram $ D: \Cc \to \Delta^0 \to \Dc $ in $ \Fun(\Cc, \Dc)$. The {\it colimit} of  $ \cF $ is 
defined dually by the formula: $\,\bcolim_{\Cc}(\cF) := \blim_{\Cc^{\rm op}}(\cF^{\rm op})\,$, and thus it represents the vertex of a universal cone $ \cF \to \bcolim_{\Cc}(\cF) $ under $ \cF $, inducing the equivalence
$$
\Map_{\Dc}(\bcolim_{\Cc}(\cF),\,D) \xrightarrow{\sim} \Map_{\Fun(\Cc, \Dc)}(\cF, D)\,.
$$

We are interested in limits and colimits of diagrams in $ \cAff_\k $. The next proposition expresses these in terms of classical  {\it homotopy}  (co)limits  of commutative dg algebras
which, in turn, can be computed in terms of polynomial differential forms \eqref{Omegan} (see Examples below).
\begin{prop}\la{PropA}
Let $\,\A: \I^{\rm op}\, \to \,\ccdga_\k^{\ge 0} $ be a $($contravariant$)$ diagram of commutative cochain  algebras indexed by a small category $\I $. Assume that $ H^0(A_i) \cong \k $ for each $ i\in \I $, and let $ \,\Specc(\A):\, \Cc \to \cAff_\k $  denote the associated $\infty$-diagram of coaffine stacks indexed by
$ \Cc = {\mathcal N}\I $.
 \begin{enumerate}
 \item[$(1)$]
If $\, H^0[\hocolim_{\I^{\rm op}}(\A)] \cong \k $, then
$$
\blim_{\Cc} \,[\,\Specc(\A)]\,\simeq\, \Specc\,[\hocolim_{\I^{\rm op}}(\A)]
$$

\item[$(2)$]
Dually, if $\, H^0[\holim_{\I^{\rm op}}(\A)] \cong \k $, then
$$
\bcolim_{\Cc} \,[\,\Specc(\A)]\,\simeq\, \Specc\,[\holim_{\I^{\rm op}}(\A)]\,,
$$    
\end{enumerate}
where $ \hocolim $ and $ \holim $ are taken in the  simplicial model structure on $ \ccdga_\k^{\ge 0} $. 
\end{prop}
\begin{proof}
By \cite[Corollary 4.2.4.8]{HTT}, the $\infty$-category $(\cdga_\k^{\ge 0})^{\rm op}$ associated to (the dual 
of) the  simplicial model category $ \ccdga_k^{\ge 0} $ admits limits and colimits over an arbitrary small $\infty$-category $\Cc$. By \cite[Theorem 4.2.4.1]{HTT}, these limits and colimits agree with the homotopy limits and colimits in $ (\ccdga_\k^{\ge 0})^{\rm op} $ indexed by a small category $\I$, 
provided $ \Cc  $ is equivalent to the nerve of $ \I $. 
Since we realized $ \cAff_\k $ as a full $\infty$-subcategory of $(\cdga_\k^{\ge 0})^{\rm op}$, by \cite[1.2.13.7]{HTT}, we can compute limits and colimits in $ \cAff_\k $ by working in  $(\cdga_\k^{\ge 0})^{\rm op}$ and then verifying that the resulting objects are in $ \cAff_\k $. The assumptions on the $0$th cohomology that we make
in part $(1)$ and $(2)$ are what is needed to verify this last condition.
\end{proof}
We apply Proposition~\ref{PropA} in two basic examples that provide explicit dg algebra models for
pullbacks and pushouts  in $\cAff_k$.
\begin{example}
\la{EA1}    
Let $ \Cc = (\Delta^1 \amalg_{\Delta^0} \Delta^1)^{\rm op} $, so that $ \Cc \cong \mathcal{N}\I $, where 
$ \I = \{1 \to 0 \leftarrow 2\} $ is the standard category indexing the pullback diagrams. The functor $\,\A: \I^{\rm op}\, \to \,\ccdga_\k^{\ge 0} $ is then represented by a pushout diagram
$\, \{A_1 \to A_0 \leftarrow A_2 \}\,$, and its
homotopy colimit is given by   the usual derived tensor products
of cochain dg algebras:
$$
\hocolim_{\I^{\rm op}}(\A)\,\cong\, A_1 \otimes_{A_0}^{\bL} A_2
$$
Since $H^0(A_1 \otimes_{A_0}^{\bL} A_2) \cong H^0(A_1) \otimes_{H^0(A_0)} H^0(A_2)$, condition $(1)$ of Proposition~\ref{PropA} obviously holds. Hence the pullbacks in the $\infty$-category $ \cAff_\k $ can be expressed by the formula
\begin{equation}\la{pullb}
\Specc(A_1) \times_{\Specc(A_0)} \Specc(A_2) 
\,\simeq\,\Specc(A_1 \otimes_{A_0}^{\bL} A_2)    
\end{equation}
\end{example}

\begin{example}
\la{EA2}
Let $ \Cc = \Delta^1 \amalg_{\Delta^0} \Delta^1 $, so that $ \Cc \cong \mathcal{N}\I $, where 
$ \I = \{1 \leftarrow 0 \to 2\} $ is the standard category indexing the pushout diagrams.
Then $\,\A: \I^{\rm op}\, \to \,\ccdga_\k^{\ge 0} $ is given by a pullback diagram of cochain dg algebras $\, \{A_1 \to A_0 \leftarrow A_2 \}\,$. To compute its homotopy limit, we recall that the model structure on the category $\ccdga_\k^{\ge 0}$ is fibrant and hence right proper (see \cite[Corollary 13.1.3 (2)]{Hir03}). 
Therefore, we can use standard results of abstract homotopy theory 
(see \cite{Hir03}, Section 13.3, in particular, Propositions 13.3.4 and 13.3.8) to express
$$
\holim(A_1 \to A_0 \leftarrow A_2)\,\simeq\, A_1 \times_{A_0} 
A_0^{\dDel^1} \times_{A_0} A_2 \,\cong\,
(A_1 \times A_2) \times_{(A_0 \times A_0)} A_0^{\dDel^1}
$$
where $\, A_0^{\dDel^1} $ is a path object on $A_0$. By definition, the diagonal map
$ A_0 \into A_0 \times A_0 $ factors  in 
$ \ccdga_k^{\ge 0} $ as $ 
A_0 \xrightarrow{\sim} A_0^{\dDel^1} \onto A_0 \times A_0
$, where the first arrow is a quasi-isomorphism. Hence
$H^0(A_0^{\dDel^1}) \cong H^0(A_0) \cong k $ and therefore
$$ 
H^0[A_1 \times A_2) \times_{(A_0 \times A_0)} A_0^{\dDel^1}] \cong
(H^0(A_1) \times H^0(A_2)) \times_{(H^0(A_0) \times H^0(A_0))} H^0(A_0^{\dDel^1}) \cong (\k\times \k) \times_{(\k \times \k)} \k \cong \k
$$
This shows that condition $(2)$ of Proposition~\ref{PropA} holds 
in this case. Thus the pushouts of coaffine stacks can be expressed by the formula
\begin{equation}\la{pusho}
\Specc(A_1)\, \amalg_{\Specc(A_0)} \Specc(A_2) 
\,\simeq\,\Specc(A_1 \times_{A_0} 
A_0^{\dDel^1} \times_{A_0} A_2)    
\end{equation}
We note that there is a canonical path object  in $ \ccdga_\k^{\ge 0} $ that we can use in \eqref{pusho}: 
$$ 
A_0^{\dDel^1} = A_0 \otimes \Omega(\dDel^1_\k)
$$
where $ \Omega(\dDel^1_\k)  $ is the de Rham complex of the algebraic $1$-simplex $ \dDel_\k^1 \subseteq \bA^2_\k $ defined explicitly in \eqref{Omegan}.
The corresponding path fibration $ (p_0, p_1): A_0 \otimes \Omega(\dDel^1_\k) \onto A_0 \times A_0 $ is then given by 
evaluation of the coordinate functions $ (x_0, x_1) $ on $ \bA^2_\k $ at points
$(1,0)$ and $(0,1)$, respectively.
\end{example}
%

\section{Homotopy decomposition of \texorpdfstring{$\infty$}{infinity}-categories over posets}
\la{AC}
In this Appendix, we prove a general theorem (Theorem~\ref{SlThmInf}) on 
homotopy decompositions of diagrams of spaces indexed by an $\infty$-category $\Cc$ equipped with a functor to a poset. 
This result is a key technical tool for computing homotopy colimits in the main body of the paper (see Section~\ref{S6}). In the special case when $\Cc$ is an ordinary 1-category, Theorem~\ref{SlThmInf} was first stated and proved by J. S{\l}omi{\'n}ska in \cite{Sl01} (see Theorem~\ref{SlThm} in Section~\ref{AC4}); however, various consequences of this result (e.g., its specialization to EI categories, see Corollary \ref{EIA}) were known and widely used in homotopy theory for quite some time (see \cite{DK85}, \cite{JM89}, \cite{Sl91}, \cite{JMO94}, \cite{Gr02}). The proof of Theorem~\ref{SlThmInf} that we offer in this Appendix --- and, in fact, the very approach --- are quite different from that of S{\l}omi{\'n}ska: our arguments do not specialize to \cite{Sl01} in the case when $\Cc$ is an ordinary category\footnote{The interested reader may skim the proof of the S{\l}omi{\'n}ska Theorem given in the first {\tt arXiv} version of this paper.}. The key new ingredient is provided by J. Lurie's Straightening/Unstraightening Equivalence for cartesian fibrations of $\infty$-categories and its variations studied in \cite{HTT} and the recent literature (see \cite{HM15, S17, HHR25}).

\vspace*{2ex}

\noindent
{\it Notation.}\ 
Throughout this Appendix, we use the notation and follow the conventions of \cite{HTT}: in particular, by an `$\infty$-category' we mean a quasi-category (in the sense of A. Joyal); we write $ \Cat_\infty $ for the $\infty$-category of $\infty$-categories and $\sS$ for the $\infty$-category of spaces. For the internal Hom of simplicial sets, we use the notation
$\Map(K,X) = \Fun(K, X)$ interchangeably (depending on a context). For a simplicial set $X \in \sset $, $|X| \in \sS$ denotes its homotopy type ($|N\Cc| = B\Cc$ for a category $\Cc$); the functor $|\,\cdot\,|: \Cat_\infty \to \sS$ is left adjoint to the inclusion, hence preserves colimits. A map of simplicial sets is {\it cofinal} if it is so in the sense of \cite[4.1.1.1]{HTT}: the key property of such a map is that the restriction along it preserves colimits of diagrams in any $\infty$-category; for nerves of ordinary categories, this is equivalent to (right) homotopy cofinality ({\it cf.} \cite[4.1.3.1]{HTT} and \cite[19.6.1]{Hir03}).

\subsection{The relative nerve over a poset subdivision}
\la{AC1}
Let $ (\Pc, \leq) $ be a poset. Recall that the {\it barycentric subdivision}
of $\Pc $ is the poset $ \sd(\Pc)$ consisting of all finite, strictly
increasing sequences $ \ux := (x_0, x_1, \ldots, x_n) $ of elements of $ \Pc $ ordered
by inclusion of the underlying sets: i.e.,
$$
(x_0, x_1, \ldots x_n)\,\leq\,(y_0, y_1,\ldots,y_m)\ \Longleftrightarrow\
\{x_0, x_1, \ldots x_n\} \subseteq \{y_0, y_1,\ldots,y_m\}\,.
$$
The poset $ \Pc $ is related to its barycentric subdivision
$ \sd(\Pc)$ through the canonical projection:
$$ \varepsilon:\ \sd(\Pc)^{\rm op} \to \Pc \ ,\quad
(x_0, x_1, \ldots, x_n) \mapsto x_0 \,.
$$
Now, the elements $ \ux := (x_0 < x_1 < \ldots < x_n) $ of $ \sd(\Pc) $ can be viewed as maps of posets $ [n] \to \Pc $, defining the functor
\begin{equation} \la{fungamma}
\gamma:\,\sd(\Pc) \,\to\, \Cat_{/\Pc}\,\xrightarrow{N}\, \sset_{/N\Pc}\ ,\quad
\ux \mapsto ([n] \to \Pc) \mapsto (\Delta[n] \to N\Pc)\,,
\end{equation}
where $ \sset_{/N\Pc} $ is the slice category of simplicial sets over the nerve of
$\Pc$. This slice category carries a natural --- the so-called {\it covariant} --- model structure, where the cofibrations are the monomorphisms of simplicial sets over $ N \Pc $ and the fibrant objects are the left fibrations $ X \to N\Pc$ (see \cite[Section 2.1.4]{HTT}). Then, by Dugger \cite{Du01}, the functor \eqref{fungamma} generates the universal Quillen adjunction of model categories
\begin{equation}
\la{Dugger}
{\rm Re}_{\gamma}:\,  {\rm sPre}(\sd\, \Pc)\,\leftrightarrows \sset_{/N\Pc}\, : {\rm Sing}_{\gamma} \,,
\end{equation}
where $  {\rm sPre}(\sd\,\Pc) = \sset^{\sd(\Pc)^{\rm op}} $ is the category of simplicial presheaves
on $ \sd(\Pc)$ equipped with standard projective model structure. On the other hand, a version of Lurie's
Straightening Theorem (see \cite[Proposition 3.2.5.18]{HTT}) asserts that there is a Quillen equivalence:
\begin{equation}
\la{Lurie}
\mathfrak F_{\bullet}(\sd \,\Pc):\,  \sset_{/N(\sd\, \Pc)^{\rm op}} \leftrightarrows  {\rm sPre}(\sd\, \Pc)\,: N_\bullet(\sd\, \Pc) \,,
\end{equation}
where the right adjoint is the {\it relative nerve functor} on $\sd(\Pc)^{\rm op} $ defined in \cite[3.2.5.2]{HTT}.
Combining \eqref{Dugger} and \eqref{Lurie} we thus get a Quillen adjunction relating the slice
categories of simplicial sets over  the poset $\Pc$ and its barycentric subdivision:
\begin{equation}
\la{Luriecomb}
{\rm Re}_{\gamma} \circ \mathfrak F_{\bullet}(\sd \,\Pc):\, \sset_{/N(\sd\, \Pc)^{\rm op}} \,\leftrightarrows\, \sset_{/N\Pc}\,:  N_\bullet(\sd\, \Pc) \circ {\rm Sing}_{\gamma}
\end{equation}
Our goal is to describe the right adjoint functor
of \eqref{Luriecomb} on objects (i.e., simplicial maps) of the form $F: \Cc \to N \Pc $, where $\Cc$ is a fixed $\infty$-category. Since $ N \Pc $ is an $\infty$-category, any such map is automatically a functor of $\infty$-categories (i.e, a morphism in $\Cat_{\infty}$).

First, observe that ${\rm Sing}_{\gamma}$ assigns to any object $ F: X \to N\Pc $ in $ \sset_{/N\Pc} $ the simplicial presheaf
\begin{equation}
\la{singg}
{\rm Sing}_{\gamma}(F):\ \sd(\Pc)^{\rm op} \to \sset\ ,\quad \ux \mapsto F_*^{-1}(\gamma(\ux)) \,,
\end{equation}
where $ F_*^{-1}(\gamma(\ux)) = \Map(\Delta[n], X) \times_{\Map(\Delta[n], N \Pc)} \Delta[0] $ is the fiber of the natural map
\begin{equation*}
F_*: \ \Map(\Delta[n], X) \to \Map(\Delta[n], N \Pc)
\end{equation*}
over the point $ \gamma(\ux): \Delta[n] \to N \Pc $. Now, if $ X = \Cc $ is an $\infty$-category, these fibers are well behaved. More precisely, 
\begin{lemma}
\la{Ffib}
Any functor $ F: \Cc \to N\Pc $ is a categorical fibration. Consequently, the map $F_*$ is a categorical fibration for every $n \ge 0$, and for every $ \ux = (x_0 < \ldots < x_n) \in \sd(\Pc) $ the strict fiber
\begin{equation}
\la{funfiberinf}
F_*^{-1}(\ux)\,:=\, \Fun(\Delta[n], \Cc) \times_{\Fun(\Delta[n], N\Pc)} \Delta[0]
\end{equation}
of $F_*$ over $ \gamma(\ux) $ is an $\infty$-category $($which is also the homotopy fiber of $F_*$ over $\gamma(\ux)$$)$.
\end{lemma}
\begin{proof}
By Joyal's criterion (see \cite[2.4.6.5]{HTT}), a functor between $\infty$-categories is a categorical fibration iff it is an inner fibration and lifts categorical equivalences. For $ n \geq 2 $, an inner horn $ \Lambda^n_k \subset \Delta[n] $ contains all the vertices of $\Delta[n]$, and a simplex of $ N\Pc $ is determined by its vertices; hence any filler in $\Cc$ of a horn $ \Lambda^n_k \to \Cc $ lies over the prescribed simplex $ \Delta[n] \to N\Pc $, and $F$ is an inner fibration. The only equivalences in $N\Pc$ are the identities, which lift trivially. The remaining claims follow, since $\Fun(\Delta[n],\,-\,) = \Map(\Delta[n],\,-\,)$ preserves categorical fibrations, 
and the strict fibers of categorical fibrations are homotopy fibers.
\end{proof}
Thus, for an $\infty$-category $\Cc$, the simplicial presheaf ${\rm Sing}_{\gamma}(F)$ of \eqref{singg} is the functor
\begin{equation}
\la{Fun*inf}
F_*^{-1}:\ \sd(\Pc)^{\rm op} \to \sset\ , \quad \ux \mapsto F_*^{-1}(\ux)\,,
\end{equation}
which takes values in $\infty$-categories: the objects of $ F_*^{-1}(\ux) $ are the $n$-simplices $ \underline{c}: \Delta[n] \to \Cc $ with $ F\underline{c} = \gamma(\ux) $, i.e. the (coherent) chains $ c_0 \to c_1 \to \ldots \to c_n $ in $\Cc$ with $ F(c_i) = x_i $, and the morphisms are the natural transformations of such chains lying over the identity of $ \gamma(\ux) $. By \cite[3.2.5.21]{HTT}, the composite functor $ N_\bullet(\sd\,\Pc) \circ {\rm Sing}_\gamma $ of \eqref{Luriecomb} evaluated at $F$ is the relative nerve
\begin{equation}
\la{relnerve}
N_\bullet(\sd\,\Pc)[{\rm Sing}_\gamma(F)]\,=\,N_{F_*^{-1}}(\sd(\Pc)^{\rm op})\ \longrightarrow\ N\sd(\Pc)^{\rm op}\,,
\end{equation}
which is a cocartesian fibration with fibers \eqref{funfiberinf}. 

\begin{example}
\la{SlEx}
Assume that $ F = NF_0 $ is the nerve of an ordinary functor $ F_0: \Cc_0 \to \Pc $ from a $1$-category $\Cc_0 $ to the poset $\Pc$.
Then, \eqref{Fun*inf} is the nerve of the $\Cat$-valued presheaf 
\begin{equation}
\la{funfiber}
(F_0)_*^{-1}: \ \sd(\Pc)^{\rm op} \to \Cat\ ,\quad \ux \mapsto  \Fun([n], \Cc_0)\times_{\Fun([n], \Pc)} [0]\,,
\end{equation}
and \eqref{relnerve} is the nerve of the Grothendieck construction on this presheaf (see \cite[3.2.5.6]{HTT}):
\begin{equation*} 
p:\,\Cc_{F_0}(\Pc) := \sd(\Pc)^{\rm op} \smallint\,(F_0)_*^{-1} \,\to\,\sd(\Pc)^{\rm op}
\end{equation*}
Explicitly, the category $ \Cc_{F_0}(\Pc) $ has objects $ \,([n],\,\ux, \,\underline{c})\, $, where
$ [n] \in \Delta $, $\,\ux = (x_0 < x_1 < \ldots <x_n ) \in \sd(\Pc)\,$
and $ \underline{c} = (c_0 \to c_1 \to \ldots \to c_n ) $ are the chains of composable maps in $ \Cc_0$ such that $ F(c_i) = x_i $ for all $i=0,\ldots, n$, and the morphisms are natural commutative diagrams in $\Cc$ indexed by maps in $ \sd(\Pc)^{\rm op}$. In this explicit form,  $ \Cc_{F_0}(\Pc) $ was originally defined by J. S{\l}omi{\'n}ska \cite{Sl01}.
She observed that --- besides the cocartesian fibration $p$ --- there is another natural functor
\begin{equation}
\la{p_F0}
p_{F_0}: \Cc_{F_0}(\Pc) \to \Cc \ ,\quad ([n],\,\ux,\, \underline{c})\, \mapsto \, c_0    \,,
\end{equation}
that defines a categorical correspondence (span) with $p$ over $\Pc$:
%
\begin{equation}
\la{Slomcor}
\begin{tikzcd}[scale cd= 0.9]
	& {\Cc_{F_0}(\Pc)} \\
	\Cc_0 && {\sd(\Pc)^{\rm op}} \\
	& \Pc
	\arrow["{p_{F_0}}"', from=1-2, to=2-1]
	\arrow["p", from=1-2, to=2-3]
	\arrow["F_0"', from=2-1, to=3-2]
	\arrow["\varepsilon", from=2-3, to=3-2]
\end{tikzcd}
\end{equation}
The main theorem of \cite{Sl01} (see Theorem~\ref{SlThm} below) states that \eqref{p_F0} is a (right) homotopy cofinal functor. Thus, the  correspondence \eqref{Slomcor} allows one to `transform' homotopy colimits of diagrams of spaces over $\Cc_0$ to those over ${\sd(\Pc)^{\rm op}}$. We return to this situation in Section~\ref{AC4}.
\end{example}

The S{\l}omi{\'n}ska functor \eqref{p_F0} does not admit an {\it obvious} extension from 1-categories to $\infty$-categories. To see such an extension we construct a different --- though equivalent --- model for 
the relative nerve \eqref{relnerve} that comes with a natural (strict) simplicial map generalizing \eqref{p_F0}.

\subsection{The \texorpdfstring{$\infty$}{infinity}-category \texorpdfstring{$\Cc_F(\Pc)$}{C\_F(P)}}
\la{AC2}
For a fixed $\Pc$, we define the poset
\begin{equation}
\la{Upos}
\U(\Pc)\,:=\,\{(\ux, x)\ :\ \ux \in \sd(\Pc)\,,\ x \in \ux\}\ ,\quad (\ux, x) \leq (\ux', x')\ \Longleftrightarrow\ \{\ux\} \supseteq \{\ux'\}\ \ \text{and}\ \ x \leq x' ,
\end{equation}
with three natural poset maps 
\begin{equation}
\la{qpis}
\begin{array}{lll}
q:\ \U(\Pc) \to \sd(\Pc)^{\rm op}\,, & (\ux, x) \mapsto \ux\ , & \\*[1ex]
\pi:\ \U(\Pc) \to \Pc\,, & (\ux, x) \mapsto x\ , & \\*[1ex]
s:\ \sd(\Pc)^{\rm op} \to \U(\Pc)\,, & \ux \mapsto (\ux, \,\min(\ux))\,, &
\end{array}
\end{equation}
%
The poset $\U(\Pc)$ together with the first map in \eqref{qpis} can be obtained by applying the {\it contravariant}\footnote{This should not be confused with the opposite of the  {\it covariant} Grothendieck construction on $\Phi$ over $\sd(\Pc)$, in which the partial order on $\U$  would be reversed in the second coordinate.}  Grothendieck construction to the tautological functor $ \Phi: (\sd(\Pc)^{\rm op})^{\rm op} = \sd(\Pc) \to \Cat $, $\,\ux \mapsto \ux$. Hence, $q$ is a Grothendieck fibration with fiber $\ux$ over $\ux$, whose cartesian morphisms are $ (\ux, x) \to (\ux', x) $; its nerve 
$ Nq $ is therefore a cartesian fibration of simplicial sets with fiber $ \Delta[n] $ over $ \ux $. 
The base-change functor $ q^*: \sset_{/N\sd(\Pc)^{\rm op}} \to \sset_{/N\U(\Pc)} $ along such a fibration admits a {\it right} adjoint 
$$
q_*: \sset_{/N\U(\Pc)} \,\to\,  \sset_{/N\sd(\Pc)^{\rm op}} \,,
$$
which can be explicitly described as follows (see \cite[3.2.2]{HTT}): 
for $ X \to N\U(\Pc) $, an $m$-simplex of $ q_*X $ over $ \sigma: \Delta[m] \to N\sd(\Pc)^{\rm op} $ is a map $ \Delta[m] \times_{N\sd(\Pc)^{\rm op}} N\U(\Pc) \to X $ over $ N\U(\Pc)$.

\begin{defi}
\la{CFPinfdef}
For an $\infty$-category $\Cc$ and a functor $ F: \Cc \to N\Pc $, we define
\begin{equation}
\la{CFPinf}
\Cc_F(\Pc)\,:=\, q_*\big(N\U(\Pc) \times_{N\Pc} \Cc\big)\ \in\ \sset_{/N\sd(\Pc)^{\rm op}}\,,
\end{equation}
the fiber product being taken along $ N\pi $ and $F$. To simplify the notation we set 
$ \Dc := N\U(\Pc) \times_{N\Pc} \Cc $, so that $ \Cc_F(\Pc) = q_* \Dc $, and denote by $ p: \Cc_F(\Pc) \to N\sd(\Pc)^{\rm op} $ the structure map. Next, using the last map $s$ in \eqref{qpis} (which is a canonical section of $q$), we define
\begin{equation}
\la{p_Finf}
p_F:\ \Cc_F(\Pc)\, =\, s^* q^* q_* \Dc \,\xrightarrow{\ s^*(\eta)\ }\, s^* \Dc \,=\, N\sd(\Pc)^{\rm op} \times_{N\Pc} \Cc \,\to\, \Cc\,,
\end{equation}
where $ \eta: q^*q_* \Dc \to \Dc $ is the counit of the adjunction $ (q^*, \,q_*)$, and the fiber product $ s^*\Dc $ is taken along the map $ \pi \circ s = \varepsilon $.
\end{defi}

Explicitly, an $m$-simplex of $ \Cc_F(\Pc) $ is a pair $ (\sigma, T) $, where $ \sigma = (\ux^{(0)} \supseteq \ux^{(1)} \supseteq \ldots \supseteq \ux^{(m)}) $ is an $m$-simplex of $ N\sd(\Pc)^{\rm op} $ and
\begin{equation}
\la{simpCFP}
T:\ N(V_\sigma) \to \Cc\ , \qquad V_\sigma := [m] \times_{\sd(\Pc)^{\rm op}} \U(\Pc) = \{(i, x)\ :\ 0 \le i \le m\,,\ x \in \ux^{(i)}\}\,,
\end{equation}
is a map of simplicial sets with $ F\,T(i,x) = x $, the poset $ V_\sigma \subseteq [m] \times \ux^{(0)} $ being given the product order; $ p $ sends $(\sigma, T)$ to $\sigma$, and $ p_F $ sends $(\sigma, T)$ to the $m$-simplex $ T \circ s_\sigma $ of $ \Cc$, where $ s_\sigma: [m] \to V_\sigma $, $ i \mapsto (i, \min \ux^{(i)}) $, is induced by $s$. In particular, the objects of $ \Cc_F(\Pc) $ are the pairs $ (\ux, \underline{c}) $ with $ \ux \in \sd(\Pc) $, $ \underline{c} \in F_*^{-1}(\ux) $, and $ p_F(\ux, \underline{c}) = c_0 $. Thus, $\Cc_F(\Pc)$ defines a natural correspondence (span) $\, \Cc \xleftarrow{\,p_F\,} \Cc_F(\Pc) \xrightarrow{\,p\,} \sd(\Pc)^{\rm op} \,$ over the poset $\Pc$, i.e. $ F \circ p_F = \varepsilon \circ p $; by \eqref{p_Finf}, it factors through the fiber product $ \sd(\Pc)^{\rm op} \times_{\Pc} \Cc $. 

We now verify that \eqref{CFPinf} has the expected properties.
\begin{prop}
\la{CFPprop}
Let $\Cc$ be an $\infty$-category and $ F: \Cc \to N\Pc $ a functor.
\begin{enumerate}[label=$(\roman*)$, topsep=2pt, itemsep=1pt]
\item The map $ p: \Cc_F(\Pc) \to N\sd(\Pc)^{\rm op} $ is a cocartesian fibration whose fiber over $\ux$ is the $\infty$-category $ F_*^{-1}(\ux) $ of \eqref{funfiberinf}; an edge $ (\sigma, T) $ with $ \sigma = (\ux \supseteq \ux') $ is $p$-cocartesian iff $ T $ carries the edges $ (0, x) \to (1, x) $, $ x \in \ux' $, to equivalences in $\Cc$. In particular, $ \Cc_F(\Pc) $ is an $\infty$-category, and $ p $ is classified by a functor $ \sd(\Pc)^{\rm op} \to \Cat_\infty $ with values \eqref{funfiberinf}.

\item There is a natural categorical equivalence of cocartesian fibrations over $ N\sd(\Pc)^{\rm op}$
\begin{equation}
\la{relnervecomp}
\Cc_F(\Pc)\ \xrightarrow{\,\sim\, }\ N_{F_*^{-1}}(\sd(\Pc)^{\rm op})\,=\, N_\bullet(\sd\,\Pc)\big[{\rm Sing}_\gamma(F)\big]\,.
\end{equation}

\item If $\Cc = N\Cc_0$ is the nerve of a small category $\Cc_0$ and $F = NF_0$, then \eqref{relnervecomp} is an isomorphism, so that $\Cc_F(\Pc)$ is the nerve of the classical category $ \sd(\Pc)^{\rm op}\smallint F_*^{-1} $ and $p_F$ is the nerve of the first vertex functor.

\end{enumerate}
\end{prop}
\begin{proof}
$(i)$ The fiber of $ p $ over $\ux$ is given by the constant simplices $\sigma$, for which $ V_\sigma = [m] \times \ux $. The $m$-simplices of $V_\sigma $ are the maps $ \Delta[m] \times \Delta[n] \to \Cc $ lying over $ \gamma(\ux) \circ {\rm pr}_2 $, i.e. the $m$-simplices of the simplicial set \eqref{funfiberinf}. To see that $p$ is a cocartesian fibration we recall that $ \Cc_F(\Pc) = q_* \Dc $ and
observe that $\, \Dc \cong  (N\U(\Pc) \times \Cc) \times_{N\U(\Pc) \times N\Pc} N\U(\Pc) \,$, the fiber product being taken along $ N\U(\Pc) \times F $ and the graph $ (\id, N\pi) $ of the map $\pi$ in \eqref{qpis}. Since $q_*$ is a right adjoint, it preserves fiber products and terminal objects: hence, the above isomorphism implies
$$ 
\Cc_F(\Pc)\,\cong\,q_*(N\U \times \Cc) \times_{q_*(N\U \times N\Pc)} N\sd(\Pc)^{\rm op} \,, 
$$
where the map $\, N\sd(\Pc)^{\rm op} \to q_*(N\U \times N\Pc) \,$ is the tautological section 
$ t := q_*(\id, N\pi) $, taking an $m$-simplex $\sigma$ to $ N\pi: N(V_\sigma) \to N\Pc$. By \cite[3.2.2.13]{HTT}, 
since $Nq$ is a cartesian fibration and the natural projections of $ N\U(\Pc) \times \Cc$ and $ N\U(\Pc) \times N\Pc $ onto  $ N\U(\Pc) $ are cocartesian fibrations, both $ q_*(N\U(\Pc) \times \Cc) $ and $ q_*(N\U(\Pc) \times N\Pc) $ are cocartesian fibrations over $ N\sd(\Pc)^{\rm op}$, an edge $(\sigma, T)$ being cocartesian iff $T$ carries the $q$-cartesian edges $ (0,x) \to (1,x) $ of $ N(V_\sigma) $ to edges whose second component is an equivalence. 
Next, we consider the functor  $ F_* := q_*(N\U(\Pc) \times F) $.
Since $F$ preserves equivalences, $ F_* $ preserves cocartesian edges, and so does $t$, as $ N\pi $ takes $ (0,x) \to (1,x) $ to $ \id_x$. We claim that $F_*$ is a categorical fibration. To see this we equip the categories  $ \sset_{/N\U(\Pc)} $ and $\sset_{/N\sd(\Pc)^{\rm op}}$ with model structures induced from the Joyal model structure on $\sset$.  In these model structures, the cofibrations are monomorphisms, and the weak equivalences are categorical equivalences of total simplicial sets. The functor $q^*$ preserves monomorphisms and categorical equivalences: indeed, if $ f: A \to B $ is a categorical equivalence over $ N\sd(\Pc)^{\rm op} $, then $ q^* A = A \times_B q^*B $ is a base change of the cartesian fibration $ q^*B \to B $ along $f$, and a base change of a cartesian fibration along a categorical equivalence is a categorical equivalence by \cite[3.3.1.3]{HTT}. 
Hence, $ (q^*,\, q_*) $ is a Quillen adjunction, and $q_*$ preserves fibrations. Since $ N\U(\Pc) \times F $ is a base change of the categorical fibration $F$ (see Lemma~\ref{Ffib}),  $F_*$ is a fibration between $\infty$-categories, i.e. a categorical fibration  (\cite[2.4.6.5]{HTT}) as claimed. Now, the map $ p: \Cc_F(\Pc) \to N\sd(\Pc)^{\rm op} $ is obtained as the base change of $F_*$ along $t$. Hence, $p$ is an inner fibration, and moreover, it admits cocartesian lifts. The latter can be described as follows: given an edge $e$ of $ N\sd(\Pc)^{\rm op} $ and a vertex $a$ of $\Cc_F(\Pc)$ over its source, choose a cocartesian lift $ \tilde{e}: a \to a' $ of $e$ in $ q_*(N\U(\Pc) \times \Cc)$; then $ F_*(\tilde{e}) $ and $ t(e) $ are cocartesian lifts of $e$ with the same source, so they differ by an equivalence $ \theta: F_*(a') \to t(e)_1 $ in the fiber, which, together with the $2$-simplex witnessing this \cite[2.4.1.7]{HTT}, lifts along the categorical fibration $F_*$; the resulting edge $ a \to a'' $ lies over $t(e)$ and is cocartesian for $ q_*(N\U(\Pc) \times \Cc) \to N\sd(\Pc)^{\rm op} $, hence $F_*$-cocartesian \cite[2.4.1.3(3)]{HTT}, hence $p$-cocartesian \cite[2.4.1.3(2)]{HTT}. 
We conclude that $p$
is a cocartesian fibration of simplicial sets.

$(ii)$ Unwinding \cite[3.2.5.2]{HTT}, an $m$-simplex of $ N_{F_*^{-1}}(\sd(\Pc)^{\rm op}) $ over $ \sigma $ is a map $ T': S_\sigma \to \Cc $ over $ N\Pc$, where $ S_\sigma := \bigcup_{j=0}^{m}\ \Delta[\{0,\ldots,j\}] \times \Delta[\ux^{(j)}] \subseteq N(V_\sigma) $, with $ \Delta[\ux^{(j)}] \subseteq \Delta[\ux^{(0)}] $ the face spanned by $\ux^{(j)}$. Restriction along $ S_\sigma \subseteq N(V_\sigma) $ defines \eqref{relnervecomp} over $ N\sd(\Pc)^{\rm op}$. Both sides are cocartesian fibrations ($(i)$ and \cite[3.2.5.21]{HTT}), the map is an isomorphism on fibers, and it preserves cocartesian edges: by \cite[3.2.5.21]{HTT}, $(\sigma, T')$ with $\sigma = (\ux \supseteq \ux')$ is cocartesian iff $T'|_{\Delta[1] \times \Delta[\ux']}$, viewed as an edge $ \underline{c}|_{\ux'} \to \underline{c}' $ of $ F_*^{-1}(\ux')$, is an equivalence, i.e. iff its components $ T'((0,x) \to (1,x)) $, $x \in \ux'$, are equivalences in $\Cc$. Hence \eqref{relnervecomp} is a categorical equivalence by \cite[3.3.1.5]{HTT}. 

$(iii)$ Maps $ S_\sigma \to N\Cc_0 $ correspond to functors $ \tau_1(S_\sigma) \to \Cc_0 $, and $ \tau_1(S_\sigma) = V_\sigma$, as every relation $ (i,x) \leq (i', x') $ in $V_\sigma$ factors as $ (i,x) \leq (i, x') \leq (i', x') $ through $ [i] \times \ux^{(i)} $ and $ [i'] \times \ux^{(i')}$. Hence \eqref{relnervecomp} is an isomorphism, and the claim follows from \cite[3.2.5.6]{HTT}; one checks directly that $p_F$ becomes the first vertex functor.
\end{proof}
\begin{remark}
By the same argument, a categorical equivalence $ \Cc \to \Cc' $ over $ N\Pc$ induces a categorical equivalence $ \Cc_F(\Pc) \to \Cc'_{F'}(\Pc) $ compatible with $p$ and $p_F$.
\end{remark}

\subsection{The main theorem}
\la{AC3}
We are now in position to state the main theorem of this Appendix.
\begin{theorem}
\la{SlThmInf}
Let $\Cc$ be an $\infty$-category, $\Pc$ a poset, and $ F: \Cc \to N\Pc $ a functor. Then the functor $ p_F: \Cc_F(\Pc) \to \Cc $ of \eqref{p_Finf} is cofinal.
\end{theorem}
The proof of this result relies on two technical Lemmas (\ref{sliceinf} 
and~\ref{LemmaABinf}), which we prove first. 

For an object $ c \in \Cc $, we let $ \Cc_{c/} $ denote 
the undercategory of $\Cc$, whose $m$-simplices are the maps $ \Delta[0] \star \Delta[m] \to \Cc $ restricting to $c$ on the cone point (see \cite[1.2.9.5]{HTT}). The projection $\Cc_{c/} \to \Cc $ is a left fibration \cite[2.1.2.2]{HTT}; in particular, $ \Cc_{c/} $ is an $\infty$-category, and we write 
$$ 
F_c: \Cc_{c/} \to \Cc \xrightarrow{F} N\Pc 
$$ 
Note that the undercategory$ (F_c)^{-1}_*(\ux) $ is empty unless $ F(c) \leq x_0 $; thus $ (\Cc_{c/})_{F_c}(\Pc) $ lies over the full subposet $ \sd(\Pc_c)^{\rm op} \subseteq \sd(\Pc)^{\rm op} $, where $ \Pc_c := \{x \in \Pc : F(c) \leq x\} $, and by Proposition~\ref{CFPprop}$(i)$ applied to $F_c$ it is a cocartesian fibration
\begin{equation}
\la{compfuncinf}
p_c:\ (\Cc_{c/})_{F_c}(\Pc) \to N\sd(\Pc_c)^{\rm op}
\end{equation}
whose fiber over $ \ux = (x_0 < \ldots < x_n) $ is the $\infty$-category $ (F_c)_*^{-1}(\ux) $ of coherent chains $\, c \to c_0 \to c_1 \to \ldots \to c_n \,$ in $\Cc$ (i.e. maps $ \Delta[0] \star \Delta[n] \to \Cc $ with cone point $c$) such that $ F(c_i) = x_i $ for all $i$.
\begin{lemma}
\la{sliceinf}
For every object $ c \in \Cc$, there is a trivial fibration of simplicial sets over $ N\sd(\Pc)^{\rm op}$
\begin{equation}
\la{sliceeq}
\rho_c:\ (\Cc_{c/})_{F_c}(\Pc)\ \xrightarrow{\ \sim\ }\ (p_F)_{c/}\,:=\, \Cc_F(\Pc) \times_{\Cc} \Cc_{c/}\,,
\end{equation}
where the fiber product is taken along $p_F$ and the projection $ \Cc_{c/} \to \Cc$.
\end{lemma}
\begin{proof}
Let $ \sigma = (\ux^{(0)} \supseteq \ldots \supseteq \ux^{(m)}) $ be an $m$-simplex of $N\sd(\Pc)^{\rm op}$. An $m$-simplex of $ (\Cc_{c/})_{F_c}(\Pc) $ over $\sigma$ is a map $ N(V_\sigma) \to \Cc_{c/} $ over $N\Pc$, i.e. a map $\, \Delta[0] \star N(V_\sigma) \to \Cc \,$ restricting to $c$ on the cone point and to a map over $N\Pc$ on $N(V_\sigma)$. By \eqref{simpCFP}, an $m$-simplex of $ (p_F)_{c/} $ over $\sigma$ is a map $ T: N(V_\sigma) \to \Cc $ over $N\Pc$ together with a map $ H: \Delta[0] \star \Delta[m] \to \Cc $ with $ H|_{\Delta[m]} = T \circ s_\sigma $ and cone point $c$; that is, a map 
$$ 
W_\sigma := N(V_\sigma) \cup_{\Delta[m]} (\Delta[0] \star \Delta[m]) \to \Cc \,,
$$
where $\Delta[m]$ is embedded in $N(V_\sigma)$ via $s_\sigma$. As $ \Delta[0] \star (\,\cdot\,) $ preserves monomorphisms and unions, $ W_\sigma = (\Delta[0] \star s_\sigma\Delta[m]) \cup (\emptyset \star N(V_\sigma)) $ is a subcomplex of $ \Delta[0] \star N(V_\sigma)$, and restriction along this inclusion defines \eqref{sliceeq}. 
 To see that $\rho_c$ is a trivial fibration we must solve, for every $\sigma$, the lifting problem
\begin{equation}
\la{liftrho}
\begin{tikzcd}
W_\sigma \cup \big(\Delta[0] \star N(V_{\partial\sigma})\big) \arrow[r] \arrow[d, hook] & \Cc \arrow[d, "F"] \\
\Delta[0] \star N(V_\sigma) \arrow[r] \arrow[ur, dashed] & N\Pc
\end{tikzcd}
\end{equation}
where $ N(V_{\partial \sigma}) := N(V_\sigma) \times_{\Delta[m]} \partial\Delta[m] $ (we use that $\Delta[0] \star (\,\cdot\,)$ preserves unions of subcomplexes). The bottom map in \eqref{liftrho} is well defined: it is the nerve of the map of posets $ \{\bot\} \sqcup V_\sigma \to \Pc $, $ \bot \mapsto F(c)$, $ (i,x) \mapsto x$, which is monotone because the given map already contains the edges $ c \to T s_\sigma(i) $, whence $ F(c) \leq \min \ux^{(i)} \leq x $. Since $F$ is an inner fibration (Lemma~\ref{Ffib}), it suffices to show that the left vertical map in \eqref{liftrho} is inner anodyne. Writing $ B := N(V_\sigma) $ and $ B_0 := N(V_{\partial\sigma}) \cup s_\sigma \Delta[m] $, this map is the inclusion $ (\Delta[0] \star B_0) \cup (\emptyset \star B) \into \Delta[0] \star B $, which is inner anodyne provided $ B_0 \into B $ is left anodyne, by \cite[2.1.2.3]{HTT}. To prove the latter, consider the map of posets
$$
h:\ [1] \times V_\sigma \to V_\sigma\ , \qquad h(0, (i,x)) = (i, \min \ux^{(i)})\,,\quad h(1, (i,x)) = (i,x)\,,
$$
which is monotone since $ \min\ux^{(i)} \leq \min\ux^{(i')} \leq x' $ whenever $ i \leq i' $ and $ x' \in \ux^{(i')}$. The map $ Nh: \Delta[1] \times B \to B $ commutes with the projections to $\Delta[m]$, hence preserves $ N(V_{\partial\sigma}) $, it restricts to the projection on $ \Delta[1] \times s_\sigma\Delta[m] $, and $ Nh|_{\{0\} \times B} = s_\sigma \circ {\rm pr} $ takes values in $ s_\sigma \Delta[m] \subseteq B_0 $. Therefore $B_0 \into B$ is a retract of the inclusion $ (\{0\} \times B) \cup (\Delta[1] \times B_0) \into \Delta[1] \times B $ (embed via $ b \mapsto (1,b) $ and retract via $Nh$), which is left anodyne by \cite[2.1.2.6]{HTT}. 
This completes the proof.
\end{proof}
Next, we record a formal lemma on colimits of diagrams indexed by $\infty$-categories in adjunction.
\begin{lemma}
\la{LemmaABinf}
Let $ X: {\mathscr B} \to \M $ be a diagram in an $\infty$-category $\M$ indexed by an $\infty$-category ${\mathscr B}$, and let $ i: {\mathscr A} \to {\mathscr B} $ be a functor admitting a right adjoint $ r: {\mathscr B} \to {\mathscr A} $ with counit $ \varepsilon: i r \to \id_{\mathscr B} $. If $ X(\varepsilon): (ir)^* X \to X $ is an equivalence in $ \Fun({\mathscr B}, \M) $ $($i.e. $ X(\varepsilon_b) $ is an equivalence for every $b$$)$, then $ \colim_{\mathscr B}(X) $ exists iff $ \colim_{\mathscr A}(i^*X) $ exists, in which case they are equivalent.
\end{lemma}
\begin{proof}
There is a `zig-zag' of natural maps of colimits
\begin{equation*}
    \begin{tikzcd}
        &   \colim_{{\mathscr B}}((ir)^*X) = \colim_{{\mathscr B}}(r^* i^* X) \arrow[dl, "\varepsilon_*"']\arrow[dr,"r^*"] & \\
       \colim_{{\mathscr B}}(X) && \colim_{{\mathscr A}}(i^*X) 
    \end{tikzcd}
\end{equation*}
where $\varepsilon_*$ is induced by $X(\varepsilon)$ and $r^*$ by restriction along $r$. The first is an equivalence because $X(\varepsilon)$ is, and the second because $r$, being a right adjoint, is cofinal: by \cite[4.1.3.1]{HTT} it suffices that $ {\mathscr B} \times_{\mathscr A} {\mathscr A}_{a/} $ be weakly contractible for every $a \in \mathscr{A}$, and the adjunction identifies this left fibration over ${\mathscr B}$ with $ {\mathscr B}_{i(a)/} $, which has an initial object.
\end{proof}
\begin{proof}[Proof of Theorem~\ref{SlThmInf}]
By Joyal's version of Quillen's Theorem A \cite[4.1.3.1]{HTT}, it suffices to show that, for every object $ c \in \Cc $, the simplicial set $ (p_F)_{c/} $ is weakly contractible; by Lemma~\ref{sliceinf}, it suffices to prove that $ |(\Cc_{c/})_{F_c}(\Pc)| \simeq {\rm pt} $. Let $ (F_c)^{-1}_*: \sd(\Pc_c)^{\rm op} \to \Cat_\infty $ be a functor classifying the cocartesian fibration \eqref{compfuncinf}. By \cite[3.3.4.3]{HTT}, 
the colimit of $ (F_c)^{-1}_* $ in $\Cat_\infty$ is obtained from $ (\Cc_{c/})_{F_c}(\Pc) $ by inverting the $p_c$-cocartesian edges; applying the colimit-preserving functor $ |\,\cdot\,|$, which does not see the inverted edges, we get
\begin{equation}
\la{Bdecompinf}
|(\Cc_{c/})_{F_c}(\Pc)|\,\simeq\, \colim_{\sd(\Pc_c)^{\rm op}} \big[\,|(F_c)^{-1}_*|\,\big]\,.
\end{equation}
We compute the colimit in \eqref{Bdecompinf} using Lemma~\ref{LemmaABinf}, with ${\mathscr A} := \{\,\ux \in {\rm sd}(\Pc_c)^{\rm op}\, :\, x_0 = F(c) \}$ the full subposet of $ {\mathscr B} = \sd(\Pc_c)^{\rm op} $ and $\,i: {\mathscr A} \into {\rm sd}(\Pc_c)^{\rm op} $ the inclusion. It has a right adjoint $ r $, given by $ r(\ux) = \ux $ if $ x_0 = F(c) $ and $ r(\ux) = (F(c) < \ux) := (F(c) < x_0 < \ldots < x_n) $ if $ x_0 \neq F(c) $; the unit is $ ri = \id$, and the counit $\varepsilon_{\ux}: ir(\ux) \to \ux $ is the identity for $ \ux \in {\mathscr A} $ and the map in $ \sd(\Pc_c)^{\rm op} $ represented by the inclusion $ \{\ux\} \subset \{F(c), \ux\} $ otherwise. Thus, for $ X := |(F_c)^{-1}_*| $, we need only check that $ X(\varepsilon_{\ux}) $ is an equivalence for $ \ux = (x_0 < \ldots < x_n) $ with $ x_0 > F(c)$. The functor
\begin{equation}
\la{FunF_cxinf}
(F_c)^{-1}_*(\varepsilon_{\ux}):\ (F_c)^{-1}_*(F(c) < \ux)\ \to\ (F_c)^{-1}_*(\ux)
\end{equation}
is the restriction $ (d^0)^*: \Fun(\Delta[n+1], \Cc_{c/}) \to \Fun(\Delta[n], \Cc_{c/}) $ along the face $ d^0: \Delta[n] \into \Delta[n+1] $ omitting the vertex $F(c)$, restricted to the fibers over $ N\Pc $ (for the model \eqref{CFPinf}, morphisms of the base act on fibers by restriction of sections); in terms of chains, 
$$ 
(d^0)^*:\ (c \to c_0 \to c_1 \to \ldots \to c_{n+1})\ \mapsto\ (c \to c_1 \to \ldots \to c_{n+1}) \,,
$$ 
where $ F(c_0) = F(c) $ and $ F(c_i) = x_{i-1} $ for $ i \geq 1$. Now, $\Cc_{c/}$ has an initial object $ \id_c $, and $ d^0 $ is fully faithful. Hence, $ (d^0)^* $ admits a fully faithful left adjoint --- the left Kan extension $ (d^0)_! $ (see \cite[4.3.2.16]{HTT} --- 
which sends 
$$
(d^0)_!:\  (c \to c_1 \to \ldots \to c_{n+1}) \ \mapsto\   (c \xrightarrow{\id} c \to c_1 \to \ldots \to c_{n+1}) \,,
$$
the value at the new vertex being the colimit of the empty diagram. 
The unit $ \id \to (d^0)^* (d^0)_! $ is an isomorphism, and the components of the counit $ (d^0)_! (d^0)^* \to \id $ are the natural transformations $ (c \xrightarrow{\id} c \to c_1 \to \ldots) \to (c \to c_0 \to c_1 \to \ldots) $ whose only non-identity component $ (\id_c) \to (c \to c_0) $ lies over $ \id_{F(c)}$. Since $ \Fun(\Delta[n+1], N\Pc) $ is the nerve of a poset, this counit lies over the identity of $ \gamma(F(c) < \ux)$, and $ (d^0)_! \dashv (d^0)^* $ restricts to an adjunction between the fibers $ (F_c)^{-1}_*(F(c) < \ux) $ and $ (F_c)^{-1}_*(\ux) $. As an adjunction between $\infty$-categories induces inverse homotopy equivalences on the underlying homotopy types (a natural transformation $ \Delta[1] \times {\mathscr D} \to {\mathscr D}' $ being a homotopy), \eqref{FunF_cxinf} induces an equivalence $ X(\varepsilon_{\ux}): |(F_c)^{-1}_*(F(c) < \ux)| \xrightarrow{\sim} |(F_c)^{-1}_*(\ux)| $, which verifies the hypothesis of Lemma~\ref{LemmaABinf}. Applying it, we conclude that $\, \colim_{\sd(\Pc_c)^{\rm op}}(X) \simeq \colim_{\mathscr A}(i^* X) \simeq X([0], F(c)) \,$, the last equivalence because $ ([0], F(c)) $ is a terminal object of ${\mathscr A}$. Finally, $ X([0], F(c)) $ is the homotopy type of $ (F_c)^{-1}_*([0], F(c)) = \Cc_{c/} \times_{N\Pc} \{F(c)\} $, the fiber of $F_c$ over $F(c)$: this is a full subcategory of $\Cc_{c/}$ (a simplex of $ \Cc_{c/} $ lies over $F(c)$ iff all its vertices do) containing the initial object $\id_c$ of $\Cc_{c/}$, hence weakly contractible. By \eqref{Bdecompinf}, $ (\Cc_{c/})_{F_c}(\Pc) $ is weakly contractible, which completes the proof.
\end{proof}
Theorem~\ref{SlThmInf} yields a decomposition of colimits over $\Cc$ into iterated colimits over $\sd(\Pc)^{\rm op}$ and the fibers \eqref{funfiberinf}.
\begin{cor}
\la{Sldecinf}
Let $ \M $ be an $\infty$-category admitting small colimits. Then, for every $ X: \Cc \to \M $, there are natural equivalences
\begin{equation}
\la{Sldecinfeq}
\colim_{\Cc}(X)\,\simeq\, \colim_{\Cc_F(\Pc)}(p_F^*X)\,\simeq\, \colim_{\ux \in \sd(\Pc)^{\rm op}}\ \colim_{F_*^{-1}(\ux)}(p_F^* X)\,,
\end{equation}
where the inner colimit is the value at $\ux$ of the left Kan extension $ p_!(p_F^*X) $ along the cocartesian fibration $p$.
\end{cor}
\begin{proof}
The first equivalence is Theorem~\ref{SlThmInf}. The second follows from \cite[4.2.3.10]{HTT}: 
for the cocartesian fibration $p$, the left Kan extension $ p_!(p_F^*X) $ exists and is computed fiberwise, and $ \colim_{\Cc_F(\Pc)}(p_F^* X) \simeq \colim_{\sd(\Pc)^{\rm op}}\, p_!(p_F^*X) $.
\end{proof}
We apply Corollary~\ref{Sldecinf} in a universal situation. For an $\infty$-category $\Cc$ with homotopy category $ h\Cc $, let $ \Pc(\Cc) := \Pc(h\Cc) $ be the universal poset under $h\Cc$ (see \eqref{erel} below) and $ F_{\Cc}: \Cc \to N(h\Cc) \to N\Pc(\Cc) $ the canonical functor. We say that $\Cc$ is an {\it EIA $\infty$-category} if $ h\Cc $ is an EIA category, i.e. every endomorphism in $\Cc$ is an equivalence; then every morphism $ f: c \to c' $ between equivalent objects is an equivalence (if $ g: c' \to c $ is an equivalence, the endomorphisms $gf$ and $fg$ are equivalences, hence so is $f$). Also, $ \Pc(\Cc) $ is the set of equivalence classes of objects, and $ \sd[\Pc(\Cc)] $ consists of the sequences $ \underline{c} = ([c_0], \ldots, [c_n]) $ of pairwise distinct classes with $ \Map_{\Cc}(c_{i-1}, c_i) \neq \varnothing$. For such $\underline{c}$, with chosen representatives $c_i$, we write $ {\rm Aut}_{\Cc}(c_i) \subseteq \Map_{\Cc}(c_i, c_i) $ for the space of self-equivalences of $c_i$ (a group-like $E_1$-space), $ {\rm Aut}_{\Cc}(\underline{c}) := \prod_i {\rm Aut}_{\Cc}(c_i)$, and $ \Map_{\Cc}(\underline{c}) := \prod_{i \geq 1} \Map_{\Cc}(c_{i-1}, c_i) $ (with $ \Map_{\Cc}(\underline{c}) = \ast $ if $n = 0$).
\begin{cor}
\la{EIAinf}
Let $\Cc$ be an EIA $\infty$-category and $ F = F_{\Cc}$. Then:
\begin{enumerate}[label=$(\roman*)$, topsep=2pt, itemsep=1pt]
\item For every $ \underline{c} \in \sd[\Pc(\Cc)] $, the fiber $ F_*^{-1}(\underline{c}) $ is an $\infty$-groupoid, and there is a natural equivalence of spaces
\begin{equation}
\la{fibEIAinf}
F_*^{-1}(\underline{c})\,\simeq\, \Map_{\Cc}(\underline{c})_{h{\rm Aut}_{\Cc}(\underline{c})}\,,
\end{equation}
where $ {\rm Aut}_{\Cc}(\underline{c}) $ acts on $ \Map_{\Cc}(\underline{c}) $ by conjugation, $ (g_i)\cdot(f_i) = (g_i f_i g_{i-1}^{-1})$.

\item For every $ X: \Cc \to \sS $, there is a natural equivalence
\begin{equation}
\la{Xnatinf}
\colim_{\Cc}(X)\,\simeq\, \colim_{\sd[\Pc(\Cc)]^{\rm op}}(X^\natural)\,, \qquad X^{\natural}(\underline{c})\,=\,\big[\Map_{\Cc}(\underline{c}) \times X(c_0)\big]_{h{\rm Aut}_{\Cc}(\underline{c})}\,,
\end{equation}
with $ {\rm Aut}_{\Cc}(\underline{c}) $ acting diagonally, through the projection $ {\rm Aut}_{\Cc}(\underline{c}) \to {\rm Aut}_{\Cc}(c_0) $ on the second factor.
\end{enumerate}
\end{cor}
\begin{proof}
$(i)$ Since $ \Fun(\Delta[n], N\Pc) $ is a poset nerve, $ F_*^{-1}(\underline{c}) $ is the full subcategory of $ \Fun(\Delta[n], \Cc) $ spanned by the chains $ (d_0 \to \ldots \to d_n) $ with $ [d_i] = [c_i] $. The components of a morphism in it are morphisms between equivalent objects, hence equivalences; since natural equivalences are detected pointwise, 
$ F_*^{-1}(\underline{c}) $ is a Kan complex \cite[1.2.5.1]{HTT}, namely the union of the components of $ \Fun(\Delta[n], \Cc)^{\simeq} $ containing chains with vertices in the classes $[c_i]$. Since the spine inclusion $ {\rm Sp}[n] \into \Delta[n] $ is inner anodyne, evaluation at the vertices $ \Fun(\Delta[n], \Cc)^{\simeq} \simeq \Fun({\rm Sp}[n], \Cc)^{\simeq} \to (\Cc^{\simeq})^{n+1} $ is a Kan fibration with fiber $ \Map_{\Cc}(\underline{c}) $ over $ (c_0, \ldots, c_n) $. Restricting to the component $ \prod_i B{\rm Aut}_{\Cc}(c_i) $ of $ (\Cc^\simeq)^{n+1} $, we see that $ F_*^{-1}(\underline{c}) $ is the total space of a Kan fibration over $ B{\rm Aut}_{\Cc}(\underline{c}) $ with fiber $ \Map_{\Cc}(\underline{c}) $, i.e. the homotopy quotient \eqref{fibEIAinf} of the induced (conjugation) action.

$(ii)$ By Corollary~\ref{Sldecinf}, it suffices to identify $ \colim_{F_*^{-1}(\underline{c})}(p_F^*X) $ with $ X^\natural(\underline{c})$. The restriction of $ p_F^* X $ to $ F_*^{-1}(\underline{c}) $ is $ X \circ {\rm ev}_0 $, where ${\rm ev}_0$ factors through $ B{\rm Aut}_{\Cc}(\underline{c}) \to B{\rm Aut}_{\Cc}(c_0) \subseteq \Cc^{\simeq}$; hence, by $(i)$ and \cite[4.2.3.10]{HTT} applied to $ F_*^{-1}(\underline{c}) \to B{\rm Aut}_{\Cc}(\underline{c})$, $ \colim_{F_*^{-1}(\underline{c})}(p_F^*X) \simeq \colim_{B{\rm Aut}_{\Cc}(\underline{c})} [\Map_{\Cc}(\underline{c}) \times X(c_0)] $, which is \eqref{Xnatinf}.
\end{proof}

\subsection{Homotopy decompositions of categories}
\la{AC4}
We now specialize the above results to the classical case, where $F = NF_0$ is the nerve of a functor $ F_0: \Cc_0 \to \Pc $ from an ordinary (small) category $ \Cc_0 $ to a poset $\Pc$ (see Example~\ref{SlEx}). By Proposition~\ref{CFPprop}$(iii)$, the $\infty$-category \eqref{CFPinf} is then the nerve of the Grothendieck construction $ \Cc_{F_0}(\Pc) $ on the $\Cat$-valued fiber functor \eqref{funfiber}:
\begin{equation}
\la{CFPcat}
\Cc_{F}(\Pc)\,\cong\, N[\Cc_{F_0}(\Pc)]\,,
\end{equation}
and the $\infty$-functor \eqref{p_Finf} is isomorphic to
the nerve of the S{\l}omi{\'n}ska functor $p_{F_0} $ defined  by \eqref{p_F0}. To simplify the notation, we will drop `$N$' and write $\Cc, \,F, \ldots $ for $\Cc_0,\, F_0, \ldots $ above.

Since a functor between ordinary categories is homotopy cofinal iff its nerve is cofinal in the sense of \cite[4.1.1.1]{HTT} (see \cite[4.1.3.1]{HTT}), our Theorem~\ref{SlThmInf} specializes to
\begin{theorem}[\cite{Sl01}, Theorem 2.1]
\la{SlThm}
For any functor $F: \Cc \to \Pc $ from a small category to a poset, the functor $ p_F: \Cc_{F}(\Pc) \to \Cc $ is homotopy cofinal.
\end{theorem}
The key property of a homotopy cofinal functor is that the associated restriction functor induces a weak equivalence on all homotopy colimits (see \cite[19.6.13]{Hir03}). Thus, for any $ X: \Cc \to \M $ (where $\M$ is a model category) and any functor $ F: \Cc \to \Pc$, Theorem~\ref{SlThm} provides a weak equivalence $\, p_F^*:\ \hocolim_{\Cc_F(\Pc)}(p_F^* X)\,\xrightarrow{\sim}\, \hocolim_{\Cc}(X) \,$, which, in combination with Thomason's Theorem \cite{ChS02} ({\it cf.} \eqref{Sldecinfeq}), yields the decomposition
\begin{equation}
\la{Sldec}
\hocolim_{\Cc}(X)\,\simeq\, \hocolim_{\ux \in \sd(\Pc)^{\rm op}}\
\hocolim_{F_{\ast}^{-1}(\ux)} (p^*_{F} X)
\end{equation}
We apply \eqref{Sldec} in the universal situation: namely, given a category $\Cc$, we take $ F $ to be the functor $ F_{\Cc}: \Cc \to \Pc({\Cc}) $ that is initial among all functors from $\Cc$ to poset categories (the unit of the adjunction between the inclusion $ {\tt Poset} \into \Cat $ and its left adjoint $ \Pc: \Cat \to {\tt Poset}$). Its target --- the universal poset under $\Cc$ --- is $\,\Pc(\Cc) := {\rm Ob}(\Cc)/\!\sim\,$, where
\begin{equation}
\la{erel}
c \sim c' \quad \Leftrightarrow \quad \Hom_{\Cc}(c, c') \not= \varnothing \quad \mbox{and}\quad \Hom_{\Cc}(c', c) \not= \varnothing
\end{equation}
with the partial order $\, [c] \leq [c'] \ \Leftrightarrow \ \Hom_{\Cc}(c, c') \not= \varnothing \,$ (well defined modulo \eqref{erel}), and $ F_{\Cc} $ is the natural projection $\,c \mapsto [c] \,$ onto equivalence classes.

Now, let us assume that $\Cc$ is an EIA-category, i.e. a category in which every endomorphism and every isomorphism is an automorphism. Then, it follows from \eqref{erel} that $\, c \sim c' \,$ iff $ c \cong c' $ in $\Cc$, and moreover, every isomorphism class in $ \Cc $ contains only one object. Hence $ \Pc(\Cc) = {\rm Ob}(\Cc) $, and the poset $ \sd[\Pc(\Cc)] $ consists of finite sequences $\underline{c} = (c_0, c_1, \ldots, c_n) $ of pairwise distinct objects in $\Cc$ such that $\Hom_{\Cc}(c_{i-1}, c_i) \not= \varnothing $ for all $ i =1,\ldots n$.
In this case, the fiber functor $\, F_*^{-1}: \sd[\Pc(\Cc)]^{\rm op} \to \Cat $ associated to $ F = F_{\Cc}$ takes values in groupoids: more precisely,
for $ \underline{c} = (c_0, c_1, \ldots, c_n) \in \sd[\Pc(\Cc)] $,
$F_*^{-1}(\underline{c})$ is given by the action groupoid
\begin{equation}
\la{fibEIA}
F_*^{-1}(\underline{c}) = {\rm Aut}_{\Cc}(\underline{c}) \ltimes \Hom_{\Cc}(\underline{c}) \,,
\end{equation}
where  $\, {\rm Aut}_{\Cc}(\underline{c}) :=  {\rm Aut}_{\Cc}(c_0) \times \ldots \times  {\rm Aut}_{\Cc}(c_n)\,$ and $\, \Hom_{\Cc}(\underline{c}) := \Hom_{\Cc}(c_0, c_1) \times \ldots \times \Hom_{\Cc}(c_{n-1}, c_n) \,$ (with $ \Hom_{\Cc}(\underline{c}) := \ast $ if $ n = 0$). The action of ${\rm Aut}_{\Cc}(\underline{c})$ on $  \Hom_{\Cc}(\underline{c}) $ is  the natural adjoint action of automorphisms on morphisms in $\Cc$: i.e., for $ \underline{g} = (g_0, g_1, \ldots, g_n) \in {\rm Aut}_{\Cc}(\underline{c})\,$,
\begin{equation}\la{Autact}
\underline{g} \cdot (f_1, f_2, \ldots, f_n)\,=\, (g_1 f_1 g_0^{-1},\,g_2 f_2 g_1^{-1},\,\ldots,\,g_n f_n g_{n-1}^{-1})\,.
\end{equation}
Since homotopy colimits over the action groupoids are represented by the Borel construction, we have $\, \hocolim_{F_*^{-1}(\underline{c})}(p_F^* X) \simeq E{\rm Aut}_{\Cc}(\underline{c}) \times_{{\rm Aut}_{\Cc}(\underline{c})}[\Hom_{\Cc}(\underline{c}) \times X(c_0)] \,$ for any $ \underline{c} \in \sd[\Pc(\Cc)]$. As a consequence of Theorem~\ref{SlThm} (or of Corollary~\ref{EIAinf}), we thus get
\begin{cor}
\la{EIA}
If $\Cc$ is an EIA-category, then for any $ X: \Cc \to \M $, there is a
natural equivalence
\begin{equation}
    \hocolim_{\Cc}(X)\,\simeq\,\hocolim_{\sd[\Pc(\Cc)]^{\rm op}}(X^{\natural})
\end{equation}
where the functor $\,X^{\natural}: \sd[\Pc(\Cc)]^{\rm op} \to \M $ is defined by
\begin{equation}
\la{Xnat}
X^{\natural}(\underline{c}) =  E{\rm Aut}_{\Cc}(\underline{c}) \times_{{\rm Aut}_{\Cc}(\underline{c})}\left[\Hom_{\Cc}(\underline{c}) \times X(c_0)\right]
\end{equation}
with ${\rm Aut}_{\Cc}(\underline{c})$ acting diagonally on $\, \Hom_{\Cc}(\underline{c}) \times X(c_0) $ by \eqref{Autact} on the first factor and via the canonical projection
${\rm Aut}_{\Cc}(\underline{c}) \onto {\rm Aut}_{\Cc}(c_0) $ on the second.
\end{cor}
%
\bibliography{secondbibtexfile}{}
\bibliographystyle{amsalpha}
\end{document}